\documentclass[opre]{informs3} %
\pdfoutput=1

\OneAndAHalfSpacedXI
\usepackage{booktabs,enumitem}
\usepackage{amsmath}
\usepackage{bm}
\usepackage{comment,soul}
\usepackage{algpseudocode,algorithm}
\usepackage{graphicx}
\usepackage{subcaption}
\graphicspath{{./Figures/}}
\usepackage{tikz}

\usepackage{ifdraft}

\usepackage{shortcuts}

\usepackage{natbib}
 \bibpunct[, ]{(}{)}{,}{a}{}{,}%
 \def\bibfont{\small}%
 \def\bibsep{\smallskipamount}%
\TheoremsNumberedBySection  %
\ECRepeatTheorems

\EquationsNumberedBySection %

\MANUSCRIPTNO{}

\let\oldendproof\endproof
\renewcommand{\endproof}{\halmos\oldendproof}

\usepackage[capitalise]{cleveref}

\Crefname{assumption}{Assumption}{Assumptions}
\Crefname{appsec}{Appendix}{Appendices}

\begin{document}
 \RUNAUTHOR{Gupta and Kallus}
\RUNTITLE{Data-Pooling in Stochastic Optimization}

\TITLE{Data-Pooling in Stochastic Optimization}

\ARTICLEAUTHORS{%
\AUTHOR{Vishal Gupta}
\AFF{Data Science and Operations, USC Marshall School of Business, Los Angles, CA 90089,\\ \EMAIL{guptavis@usc.edu} } 
\AUTHOR{Nathan Kallus}
\AFF{School of Operations Research and Information Engineering and Cornell Tech, Cornell University, New York, NY 10044,\\ \EMAIL{kallus@cornell.edu} } %
} %

\ABSTRACT{%
Managing large-scale systems often involves simultaneously solving thousands of unrelated stochastic optimization problems, each with limited data. Intuition suggests one can decouple these unrelated problems and solve them separately without loss of generality. We propose a novel data-pooling algorithm called Shrunken-SAA that disproves this intuition.  In particular, we prove that combining data across problems can outperform decoupling, even when there is no a priori structure linking the problems and data are drawn independently.  Our approach does not require strong distributional assumptions and applies to constrained, possibly non-convex, non-smooth optimization problems such as vehicle-routing, economic lot-sizing or facility location.  We compare and contrast our results to a similar phenomenon in statistics (Stein's Phenomenon), highlighting unique features that arise in the optimization setting that are not present in estimation.  We further prove that as the number of problems grows large, Shrunken-SAA learns \emph{if} pooling can improve upon decoupling \emph{and} the optimal amount to pool, even if the average amount of data per problem is fixed and bounded. Importantly, we highlight a simple intuition based on stability that highlights \emph{when} and \emph{why} data-pooling offers a benefit, elucidating this perhaps surprising phenomenon.  This intuition further suggests that data-pooling offers the most benefits when there are many problems, each of which has a small amount of relevant data.  Finally, we demonstrate the practical benefits of data-pooling using real data from a chain of retail drug stores in the context of inventory management. 
}


\maketitle
%


%
%
%

\section{Introduction}
\label{sec:Intro}
The stochastic optimization problem 
\begin{equation} \label{eq:StochOpt}
\min_{\bx \in \X} \quad \E^{\P}[ c(\bx, \bxi) ]
\end{equation}
is a fundamental model with applications ranging from inventory management to personalized medicine.  
In typical data-driven settings, the measure $\P$ governing the random variable $\bxi$ is unknown.  Instead, we have access to a dataset $\S = \{ \bxihat_1, \ldots, \bxihat_N \}$ drawn i.i.d. from $\P$ and seek a decision $\bx \in \X$ depending on these data.  This model and its data-driven variant have been extensively studied in the literature (see \citealp{shapiro2009lectures} for an overview).  

Managing real-world, large-scale systems, however, frequently involves solving thousands of potentially unrelated stochastic optimization problems like Problem~\eqref{eq:StochOpt} simultaneously. For example, inventory management often requires optimizing stocking levels for many distinct products across categories, not just a single product.  Firms typically determine staffing and capacity for many warehouses and fulfillment centers across the supply-chain, not just at a single location.  Logistics companies often divide large territories into many small regions and solve separate vehicle routing problems, one for each region, rather than solving a single monolithic problem.  In such applications, a more natural model than Problem~\eqref{eq:StochOpt} might be 
\begin{equation} \label{eq:MultiStochOpt}
\frac{1}{K}  \sum_{k=1}^K \frac{\lambda_k}{\lambdabar} 
	\ \min_{ \bx_k \in \X_k} \
	\E^{\P_k}[ c_k( \bx_k, \bxi_k) ],
\end{equation}
where we solve a separate subproblem of the form \eqref{eq:StochOpt} for each $k$, e.g., setting a stocking level for each product.  Here, $\lambda_k > 0$ represents the frequency with which the decision-maker incurs costs from problems of type $k$, and $\lambdabar=\frac1K\sum_{k=1}^K\lambda_k$.  Thus, this formulation captures the fact that our total costs in such systems are driven by the frequency-weighted average of the costs of many distinct optimization problems. 

Of course, intuition strongly suggests that since there are no coupling constraints across the feasible regions $\X_k$ in Problem~\eqref{eq:MultiStochOpt}, one can and should decouple the problem into $K$ unrelated subproblems and solve them separately.  Indeed, when the measures $\P_k$ are known, this procedure is optimal.  When the $\P_k$ are unknown and unrelated, but one has access to a dataset $\S_k = \{ \bxihat_{k,1}, \ldots, \bxihat_{k, \Nhat_k} \}$ drawn i.i.d. from $\P_k$ independently across $k$, intuition \emph{still} suggests decoupling is without loss of generality and that \mbox{data-driven procedures can be applied separately by subproblem.}
\vskip 5pt
\begin{center}
	{\bf \textit{A key message of this paper is that this intuition is \textul{false}}.}
\end{center}
\vskip 5pt

\noindent In the data-driven setting, when solving many stochastic optimization problems, we show there exist algorithms which pool data across sub-problems that outperform decoupling \emph{even} when the underlying problems are unrelated, and data are \emph{independent.} This phenomenon holds despite the fact that  the $k^\text{th}$ dataset $\S_k$ tells us nothing about $\P_l$ for $l \neq k$, and there is no a priori relationship between the $\P_k$.   
We term this phenomenon the \emph{data-pooling phenomenon in stochastic optimization.}

\begin{figure}[t!]%
\centering%
\begin{minipage}[m]{0.5\textwidth}
	\input{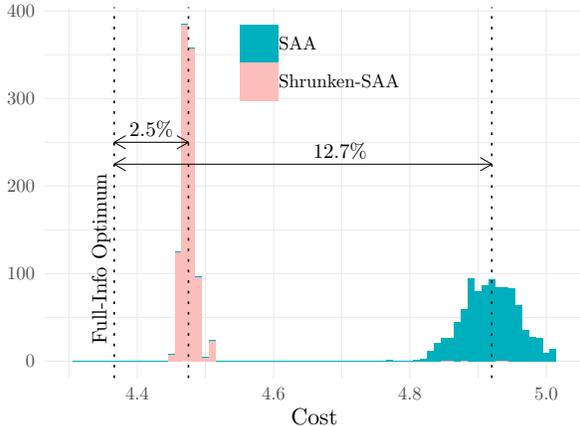}
\end{minipage}%
\begin{minipage}[m]{0.5\textwidth}
    \captionof{figure}{\textbf{The Data-Pooling Phenomenon} 
    Consider $K = 10,000$ data-driven newsvendor problems each with critical fractile $90\%$ and $20$ data points drawn independently across problems. SAA decouples the problems and orders the $90^{\text{th}}$-sample quantile in each. Shrunken-SAA (cf. \cref{alg:ssaa} in \cref{sec:ShrunkenSAA}), leverages data-pooling.  Indicated percentages are losses to the full-information optimum.     
    Additional details in Appendix~\ref{sec:SimSetUp}.  
	}
	\label{fig:IllustratingDataPooling}
\end{minipage}%
\vspace{-20pt}
\end{figure}

\Cref{fig:IllustratingDataPooling} illustrates the data-pooling phenomenon with a simulated example for emphasis.  Here $K = 10,000$, and the $k^\text{th}$ subproblem is a newsvendor problem with critical quantile $90\%$, i.e.,  $c_k(x;\xi)=\max\braces{9 (\xi-x), (x-\xi)}$.  The measures $\P_k$ are fixed and in each run we simulate $\Nhat_k=20$ data points per subproblem.  For the decoupled benchmark, we use a standard method, Sample Average Approximation (SAA; \cref{def:SAA}) which is particularly well-suited to the data-driven newsvendor problem \citep{levi2015data}.  
For comparison, we use our novel Shrunken-SAA algorithm which exploits the data-pooling phenomenon.  We motivate and formally define Shrunken-SAA in \cref{sec:ShrunkenSAA}, but, loosely speaking Shrunken-SAA proceeds by replacing the $k^\text{th}$ dataset $\S_k$ with a ``pooled" dataset which is a weighted average of the original $k^\text{th}$ dataset and all of the remaining $l \neq k$ datasets.  It then applies SAA to these  each of these new pooled datasets.  Perhaps surprisingly, by pooling data across the unrelated subproblems, Shrunken-SAA reduces the loss to full-information optimum by over 80\% compared to SAA in this example.   

\vskip 5pt
\noindent{\textbf{Our Contributions: }}  
We describe and study the data-pooling phenomenon in stochastic optimization in context of Problem~\eqref{eq:MultiStochOpt}.  Our analysis applies to constrained, potentially non-convex, non-smooth optimization problems under fairly mild assumptions on the data-generating process.  In particular, \edit{we assume only that each $\P_k$ has finite support (potentially differing across $k$); in some cases, we can even relax this assumption.  We contrast the data-pooling phenomenon to a similar phenomenon in statistics (Stein's phenomenon), highlighting unique features that arise in the optimization setting (cf. \cref{thm:DataPoolingDoesntHelp,ex:DifferentNewsvendor}).  In particular, and in contrast to traditional statistical settings, we show that the potential benefits of data-pooling depend strongly on the structure of the underlying optimization problems, and, in some cases, data-pooling may offer no benefit over decoupling.}

This observation raises important questions:  Given a particular data-driven instance of Problem~\eqref{eq:MultiStochOpt}, should we data-pool, and, if so, how?  More generally, does data-pooling \emph{typically} offer a significant benefit over decoupling, or are instances like \cref{fig:IllustratingDataPooling} somehow the exception to the rule?

To help resolve these questions, we propose a simple, novel algorithm we call Shrunken Sample Average Approximation (Shrunken-SAA).  Shrunken-SAA generalizes the classical SAA algorithm and, consequently, inherits many of its excellent large-sample asymptotic properties (cf. \cref{rem:LargeSample}).  Moreover, Shrunken-SAA is incredibly versatile and can be tractably applied to a wide variety of optimization problems with computational requirements similar to traditional SAA (cf. \cref{rem:Computational}).  Unlike traditional SAA, however, Shrunken-SAA exploits the data-pooling phenomenon to improve performance over SAA, as seen in \cref{fig:IllustratingDataPooling}.
Moreover, Shrunken-SAA exploits the structure of the optimization problems and strictly improves upon an estimate-then-optimize approach using traditional statistical shrinkage estimators (cf. \cref{ex:DifferentNewsvendor,sec:Numerics}).  

Shrunken-SAA data-pools by combining data across subproblems in a particular fashion motivated by an empirical Bayesian argument.  We prove that (under frequentist assumptions)
 for many classes of optimization problems, as the number of subproblems $K$ grows large, Shrunken-SAA determines \emph{if} pooling in this way can improve upon decoupling and, if so, also determines the optimal amount to pool (cf. \cref{thm:FixedPointShrinkage,thm:SmoothThmGeneralAnchor,thm:FixedPointShrinkageDiscrete,thm:ShrinkageDiscreteGeneralAnchor}).  
These theoretical results study Problem~\cref{eq:MultiStochOpt} when the random variables $\bxi_k$ have finite, discrete support and the amount of data available for the $k^\text{th}$ subproblem is, itself, random (see \cref{ass:RandomData}). \edit{Some of our results do extend to the case of continuous $\bxi_k$ (cf. \cref{remark:continuousanalysis} and \cref{thm:ctsdist_fixed_anchor,thm:ctsdist_general_anchor,thm:ctsdist_loo_anchor} in \cref{sec:ContinuousExtension}), and numerical experiments suggest our results are generally robust to the assumption of a random amount of data.}    
 
 More interestingly, our theoretical performance guarantees for Shrunken-SAA hold even when the expected amount of data per subproblem is small and fixed, and the number of problems $K$ is large, as in \cref{fig:IllustratingDataPooling}, i.e., they hold in the so-called small-data, large-scale regime \citep{gupta2017small}.  Indeed, since many traditional data-driven methods (including SAA) converge to the full-information optimum in the large-sample regime, the small-data, large-scale regime is arguably the more interesting regime in which to study the benefits of data-pooling.  

In light of the above results, Shrunken-SAA provides an algorithmic approach to deciding if, and, by how much to pool.  To develop an intuitive understanding of \emph{when} and \emph{why} data-pooling might improve upon decoupling, we also introduce the \emph{Sub-Optimality-Instability Tradeoff}, a decomposition of the benefits of data-pooling.  We show that the performance  of a data-driven solution to Problem~\eqref{eq:MultiStochOpt} (usually called its out-of-sample performance in machine learning settings) can be decomposed into a sum of two terms:  a term that roughly depends on its in-sample sub-optimality, and a term that depends on its instability, i.e., how much does in-sample performance change when training with one fewer data points?
As we increase the amount of data-pooling, we increase the in-sample sub-optimality because we ``pollute" the $k^\text{th}$ subproblem with data from other, unrelated subproblems.  At the same time, however, we decrease  the instability of the $k^\text{th}$ subproblem, because the solution no longer relies on its own data so strongly.  Shrunken-SAA works by navigating this tradeoff, seeking a ``sweet spot" to improve performance.  (See \cref{sec:StabilityIntuition} for discussion.)  

In many ways, the Sub-Optimality-Instability Tradeoff resembles the classical bias-variance tradeoff from statistics.  However, they differ in that the Sub-Optimality-Instability tradeoff applies to general optimization problems, while the bias-variance tradeoff applies specifically to the case of mean-squared error.  Moreover, even in the special case when Problem~\eqref{eq:MultiStochOpt} models mean-squared error, we prove that these two tradeoffs are distinct (cf. \cref{sec:BiasVariance}).  In this sense, the Sub-Optimality-Instability Tradeoff may be of independent interest outside data-pooling.  

Stepping back, this simple intuition suggests that Shrunken-SAA, and data-pooling more generally, offer significant benefits whenever the decoupled solutions to the subproblems are sufficiently unstable, which typically happens when there is only a small amount of relevant data per subproblem.  It is in this sense that the behavior in \cref{fig:IllustratingDataPooling} is typical and not pathological.  
Moreover, this intuition also naturally extends beyond Shrunken-SAA, paving the way to developing and analyzing new algorithms which also exploit the, hitherto underutilized, data-pooling phenomenon. 

Finally, we present numerical evidence in an inventory management context using real-data from a chain of European Drug Stores showing that Shrunken-SAA can offer significant benefits over decoupling when the amount of data per subproblem is small to moderate.  These experiments also suggest that Shrunken-SAA's ability to identify an optimal amount of pooling and improve upon decoupling are relatively robust to violations of our assumptions on the data-generating process.  

\vskip 5pt
\noindent \textbf{Connections to Prior Work:}  
As shown in \cref{sec:ShrunkenSAA}, our proposed algorithm Shrunken-SAA generalizes SAA.  
In many ways, SAA is \emph{the} most fundamental approach to solving Problem~\eqref{eq:StochOpt} in a data-driven setting.  SAA proxies $\P$ in \eqref{eq:StochOpt} by the empirical distribution $\hat \P$ on the data and optimizes against $\hat \P$.  It enjoys strong theoretical and practical performance in the large-sample limit, i.e., when $N$ is large \citep{kleywegt2002sample,shapiro2009lectures}.  \edit{For data-driven newsvendor problems, specifically -- an example we use throughout our work -- SAA is the maximum likelihood estimate of the optimal solution and also is the distributionally robust optimal solution when using a Wasserstein ambiguity set \citep[pg. 151]{esfahani2018data}.\label{line:NewsvendorWassersteinDRO}} SAA is incredibly versatile and applicable to a wide-variety of classes of optimization problems.   This combination of strong performance and versatility has fueled SAA's use in practice.  

When applied to Problem~\eqref{eq:MultiStochOpt}, SAA by construction decouples the problem into its $K$ subproblems.  
Because of this strong theoretical and practical performance, we use SAA throughout as the natural, ``apples-to-apples" decoupled benchmark to which we compare our data-pooling procedure Shrunken-SAA.  

More generally, the data-pooling phenomenon for stochastic optimization is closely related to Stein's phenomenon in statistics (\citealp{stein1956inadmissibility}; see also \citealp{efron2016computer} for a modern overview).  \citet{stein1956inadmissibility} considered estimating the mean of $K$ normal distributions, each with known variance $\sigma^2$, from $K$ datasets.  The $k^\text{th}$ dataset is drawn i.i.d. from the $k^\text{th}$ normal distribution and draws are independent across $k$. The natural decoupled solution to the problem (and the maximum likelihood estimate) is to use the $k^\text{th}$ sample mean as an estimate for the $k^\text{th}$ distribution. 
Surprisingly, while this estimate is optimal for each problem separately in a very strong sense (uniformly minimum variance unbiased and admissible), \citet{stein1956inadmissibility} describes a pooled procedure that \emph{always} outperforms this decoupled procedure with respect to total mean-squared error whenever $K \geq 3$.  

The proof of Stein's landmark result is remarkably short, but arguably opaque.  Indeed, many textbooks refer to it as ``Stein's Paradox," perhaps because it is not immediately clear what drives the result.  Why does it always improve upon decoupling, and what is special about $K=3$?  Is this a feature of normal distributions?  The known variance assumption?  The structure of mean-squared error loss?  All of the above?  

Many authors have tried to develop simple intuition for Stein's result (e.g., \citealp{efron1977stein,stigler19901988,brown2012geometrical,brown1971admissible,beran1996stein}) with mixed success.  As a consequence, although Stein's phenomenon has had tremendous impact in statistics, it has, in our humble opinion, had fairly limited impact on data-driven optimization.  
It is simply not clear how to generalize Stein's original algorithm to optimization problems different from minimizing mean-squared error.  
Indeed, the few data-driven optimization methods that attempt to leverage shrinkage apply either to quadratic optimization (e.g., \citealp{davarnia2017estimation,jorion1986bayes,demiguel2013size}) or else under Gaussian or near-Gaussian assumptions \citep{gupta2017small,mukherjee2015efficient}, both of which are very close to Stein's original setting.  

By contrast, our analysis of the data-pooling phenomenon \label{ReviewerWeakerPhrasingDistribution}\editt{requires very mild distributional assumptions} and applies to constrained, potentially non-convex, non-smooth optimization problems.  Numerical experiments in \cref{sec:Numerics} further suggest that even our few assumptions are not crucial to the data-pooling phenomenon.  Moreover, our proposed algorithm, Shrunken-SAA, is extremely versatile, and can be applied in any setting in which SAA can be applied.  

Finally, we note that (in)stability has been well-studied in the machine-learning community (see, e.g., \citealp{bousquet2002stability,shalev2010learnability,yu2013stability} and references therein).   \citet{shalev2010learnability}, in particular, argues that stability is the fundamental feature of data-driven algorithms that enables learning.  Our Sub-Optimality-Instability Tradeoff connects the data-pooling phenomenon in stochastic optimization to this larger statistical concept.  To the best of our knowledge, however, existing theoretical analyses of stability focus on the large-sample regime.  Ours is the first work to leverage stability concepts in the small-data, large-scale regime.  From a technical perspective, this analysis requires somewhat different tools.

\vskip 5pt
\noindent \textbf{Notation:}
Throughout the document, we use boldfaced letters $(\bp, \bfm, \ldots)$  to denote vectors and matrices, and ordinary type to denote scalars.  We use ``hat" notation ($\bphat, \bfmhat, \ldots)$ to denote observed data, i.e., an observed realization of a random variable.  We reserve the index $k$ to denote parameters for the $k^\text{th}$ subproblem.  For any random variable $X$ and $p \geq 1$, let $\| X \|_p \equiv \sqrt[p]{\E[\abs{X}^p]}$ denote the $p^\text{th}$ norm of $X$.  {\blockedit \label{NotationForReviewer} Finally, $\be_i$ refers to the $i^\text{th}$ unit vector and $\rightarrow_p$ denotes convergence in probability.}
%

\section{Model Setup and the Data-Pooling Phenomenon}
\label{sec:Model}
As discussed in the introduction, we assume throughout that $\P_k$ has finite, discrete support, i.e.,  $\bxi_k \in \{\ba_{k1},\dots, \ba_{kd}\}$ with $d \geq 2$.  Notice that while the support may in general be distinct across subproblems, without loss of generality $d$ is common.\footnote{\edit{\Cref{remark:continuousanalysis} below discusses relaxing this discrete support assumption.}} To streamline the notation, we write 
\[
p_{ki} \equiv \P_k( \bxi_k = \ba_{ki})
\ \ \text{ and } \ \ 
 c_{ki}(\bx) \equiv c_k(\bx, \ba_{ki}),
\quad  \quad   i = 1\, \ldots, d. \]
For each $k$, we let $\S_k = \{ \bxihat_{kj} : j = 1, \ldots, \Nhat_k \}$ be the $k^\text{th}$ dataset with $\bxihat_{kj} \sim\P_k$ drawn i.i.d.  
Since $\P_k$ is discrete, we can equivalently represent the $k^\text{th}$ dataset $\S_k$ via counts, $\bfmhat_k = (\mhat_{k1}, \ldots, \mhat_{kd})$, where $\mhat_{ki}$ denotes the number of times that $\ba_{ki}$ occurs in $\S_k$, and $\be^\top\bfmhat_k = \Nhat_k$.  \edit{In what follows, we will use $\bfmhat_k$ and $\S_k$ interchangeably to refer to the $k^\text{th}$ dataset.}  

Note that because $\bxihat_{kj}$ are i.i.d., 
\begin{equation} \label{eq:MultinomialCounts}
\bfmhat_k \mid \Nhat_k \sim \text{Multinomial}(\Nhat_k, \bp_k), \quad k = 1, \ldots K.
\end{equation}
Let $\S = (\S_1,\dots,\S_K)$, or equivalently, $\bfmhat=( \bfmhat_1,\dots, \bfmhat_K)$, denote all the data across all $K$ subproblems, and let $\bfNhat=(\Nhat_1,\dots,\Nhat_K)$ denote the total observation counts.  For convenience, we let $\Nhat_{\max} = \max_k \Nhat_k$ and \edit{$\Nhatbar \equiv \frac{1}{K} \sum_{k=1}^K \Nhat_k$}.  
Finally, let $\bphat_k \equiv \bfmhat_k / \Nhat_k$ denote the empirical distribution for the $k^\text{th}$ subproblem. 

Notice we have used $\hat{\cdot}$ notation when denoting $\Nhat_k$ and conditioned on its value in specifying the distribution of $\bfmhat_k$. This is because in our subsequent analysis, we will sometimes view the amount of data available for each problem as random (see Sec.~\ref{sec:UnbiasedRiskEstimate} below).  When the amount of data is fixed and \emph{non-random}, we condition on $\Nhat_k$ explicitly to emphasize this fact.  

With this notation, we can rewrite our target optimization problem:
\begin{align} \label{eq:TargetProblem}
Z^* \equiv \min_{\bx_1\in\X_1,\, \ldots,\, \bx_K\in\X_K } \quad \frac{1}{K}\sum_{k=1}^K \frac{\lambda_k}{\lambdabar} \;{\bp_k}^\top \bc_k( \bx_k)  
\end{align}
Our goal is to identify a data-driven policy, i.e., a function $\bx(\bfmhat)=( \bx_1(\bfmhat), \ldots,  \bx_K(\bfmhat) )$ mapping $\bfmhat$ to $\X_1\times\cdots\times \X_K$ for which 
\(
\frac{1}{K}\sum_{k=1}^K \frac{\lambda_k}{\lambdabar} {\bp_k}^\top \bc_k( \bx_k(\bfmhat))
\)
 is small.  We stress that the performance of a data-driven policy is random because it depends on the data.  

As mentioned with full information of $\bp_k$, Problem~\eqref{eq:TargetProblem} decouples across $k$, and, after decoupling, no longer depends on the frequency weights $\frac{\lambda_k}{ K\lambdabar}$.  
Our proposed algorithms will also \emph{not} require knowledge of the weights $\lambda_k$.
For convenience we let $\lambdamin = \min_k \lambda_k$, and $\lambdamax=\max_k \lambda_k$.  

A canonical policy to which we will compare  is the \emph{Sample Average Approximation} (SAA) policy which proxies the solution of these de-coupled problems by replacing $\bp_k$ with $\bphat_k$:

\begin{definition}[\textbf{Sample Average Approximation}] \label{def:SAA}
Let $ \bx^{\rm SAA}_k (\bfmhat_k) \in \arg\min_{\blockedit \bx_k \in \X_k} {\bphat_k}^\top \bc_k( \bx_k)$ denote the SAA policy for the $k^{th}$ problem and let   $\bx^{\rm SAA}(\bfmhat)  = (\bx_1^{\rm SAA}(\bfmhat_1), \ldots, \bx_K^{\rm SAA}(\bfmhat_K))$.
\end{definition}
As we will see, SAA is closely related to our proposed algorithm Shrunken-SAA, and hence provides a natural (decoupled) benchmark when assessing the value of data-pooling. 

{ \label{Reviewer:DefineNewsvendor}
\blockedit Finally, we use the newsvendor problem as a running example in what follows.  We say the $k^\text{th}$ subproblem is a \emph{newsvendor problem} with critical fractile $0 < s < 1$ if $c_k(x;\xi)=\max\braces{\frac{s}{1-s} (\xi-x), (x-\xi)}$.  Its full-information solution is the $s^\text{th}$ quantile of the $k^\text{th}$ distribution.
}

\subsection{A Bayesian Perspective of Data-Pooling}
\label{sec:Bayes}
To motivate data-pooling, we first consider a Bayesian approximation to our problem.  Specifically, suppose that each $\bp_k$ were independently drawn from a common Dirichlet prior, i.e., 
\[
\bp_k \sim \text{Dir}(\bp_0, \alpha_0), \quad k = 1, \ldots, K, 
\]
with $\alpha_0 > 0$ and $\bp_0 \in \Delta_d$, the $d$-dimensional simplex.  
The Bayes-optimal decision minimizes the posterior risk, which is
\(
\Eb{ \frac{1}{K}\sum_{k=1}^K \frac{\lambda_k}{\lambdabar} {\bp_k}^\top \bc_k( \bx_k) \mid \bfmhat} = 
\frac{1}{K}\sum_{k=1}^K \frac{\lambda_k}{\lambdabar} \Eb{{\bp_k} \mid \bfmhat}^\top \bc_k( \bx_k)
\),
by linearity.  Furthermore, by independence and conjugacy, respectively, 
\[
\E\left[ \bp_k  \mid \bfmhat \right]  
\ = \  
\E\left[ \bp_k  \mid \bfmhat_k \right] 
 \ =  \ \frac{\alpha_0}{\Nhat_k + \alpha_0} \bp_0 + \frac{\Nhat_k}{\Nhat_k + \alpha_0} \bphat_k.
\]
Hence, a Bayes-optimal solution is
\(
\bx(\alpha_0, \bp_0, \bfmhat_k) = \left( \bx_1(\alpha_0, \bp_0, \bfmhat_1), \ldots, \bx_K(\alpha_0, \bp_0, \bfmhat_K) \right)
\),
where
\begin{align}
\bphat_{k}(\alpha) &= \left(\frac{\alpha}{\Nhat_k + \alpha} \bp_{0} + \frac{\Nhat_k}{\Nhat_k + \alpha} \bphat_{k} \right), \quad  k = 1, \ldots, K\label{eq:ShrunkPHat}
\\ \label{eq:ShrunkXSol}
 \bx_k(\alpha, \bp_0, \bfmhat_k) &\in \arg\min_{ \bx_k \in \X_k} \quad \bphat_k(\alpha)^\top \bc_k( \bx_k), \quad k = 1, \ldots, K.
\end{align}
For any fixed (non-data-driven) $\alpha$ and $\bp_0$, $ \bx_k(\alpha, \bp_0, \bfmhat_k)$ only depends on the data through $\bfmhat_k$, but not on $\bfmhat_l$ for $l\neq k$.   

This policy has an appealing, intuitive structure. Notice $\bphat_k(\alpha)$ \editt{overloads notation slightly}\label{ReviewerOverloadingNotation}
and is a convex combination between $\bphat_k\editt{=\bphat_k(0)}$, a data-based estimated of $\bp_k$, and $\bp_0$, an a priori estimate of $\bp_k$. In traditional statistical parlance, we say $\bphat_k(\alpha)$ \emph{shrinks} the empirical distribution $\bphat_k$ toward the anchor $\bp_0$.  The Bayes-optimal solution is the plug-in solution when using this shrunken empirical measure, i.e., it optimizes $\bx_k$ as though that were the known true measure.  Note in particular, this differs from the SAA solution, which is the plug-in solution when using the ``unshrunken" $\bphat_k$.  

The parameter $\alpha$ controls the degree of shrinkage.  As $\alpha \rightarrow 0$, $ \bx_k(\alpha, \bp_0, \bfmhat)$ converges to an SAA solution, and as $\alpha \rightarrow \infty$, $ \bx_k(\alpha, \bp_0, \bfmhat)$ converges to the (non-random) solution to the fully-shrunken $k^{\text{th}}$ subproblem.  
In this sense the Bayes-optimal solution ``interpolates" between the SAA solution and the fully-shrunken solution.  
The amount of data $\Nhat_k$ attenuates the amount of shrinkage, i.e., subproblems with more data are shrunk less aggressively for the same $\alpha$.  

Alternatively, we can give a data-pooling interpretation of $\bx_k(\alpha, \bp_0, \bfmhat_k)$ via the Bayesian notion of pseudocounts.  Observe
\(
\bx_k(\alpha, \bp_0, \bfmhat_k) \in \arg\min_{ \bx_k \in \X_k} \  \left(\frac{\alpha\bp_{0}  + \bfmhat_k}{\Nhat_k + \alpha}  \right)^\top \bc_k( \bx_k)
\)
and that $\frac{\alpha\bp_{0}  + \bfmhat_k}{\Nhat_k + \alpha}$ is a distribution on $\{\ba_{k1},\dots, \ba_{kd}\}$.  In other words, we can interpret $\bx_k(\alpha, \bp_0, \bfmhat_k)$ as the solution obtained when we augment each of our original $K$ datasets with $\alpha$ additional ``synthetic" data points with counts $\alpha \bp_0$.  %
As we increase $\alpha$, we add more synthetic data.

For $\alpha > 0$, $ \bx_k(\alpha, \bp_0, \bm 0)$ is the solution to the fully shrunken $k^{\text{th}}$ subproblem.  For emphasis, let
\[
\ts \bx_k(\infty, \bp_0) \in \arg \min_{ \bx_k \in \X_k } \sum_{i=1}^d  p_{0i} c_{ki}( \bx_k),
\]
so that $ \bx_k(\alpha, \bp_0, \bm 0) =  \bx_k(\infty, \bp_0)$ for all $\alpha > 0$.  For completeness,  we also define 
\(
 \bx_k(0, \bp_0, \bm 0 ) =  \bx_k(\infty, \bp_0), 
\)
so that $ \bx_k(\alpha, \bp_0, \cdot)$ is continuous in $\alpha$.

In summary, $\bx_k(\alpha, \bp_0, \bfmhat_k)$ has an intuitive structure that is well-defined \emph{regardless of the precise structure of the cost functions $\bc_k(\cdot)$ or feasible region $\mathcal X$}.  
Importantly, this analysis shows that when the $\bp_k$ follow a Dirichlet prior, data-pooling by $\alpha$ is never worse than decoupling, and will be strictly better whenever $\bx^{\rm SAA}_k(\bfmhat_k)$ is not an optimal solution to the problem defining $\bx_k(\alpha, \bp_0, \bfmhat_k)$.

\subsection{Data-Pooling in a Frequentist Setting}
It is perhaps not surprising that data-pooling (or shrinkage) improves upon the decoupled SAA solution in the Bayesian setting because 
problems $l \neq k$ contain information about $\alpha$ and $\bp_0$ which in turn contain information about $\bp_k$.  
What may be surprising is that even in frequentist settings, i.e., when the $\bp_k$ are fixed constants that may have no relationship to one another and there is no ``ground-truth" values for $\alpha$ or $\bp_0$, policies like $\bx(\alpha, \bp_0, \bfmhat)$ can still improve upon the decoupled SAA solution through a careful choice of $\alpha$ and $\bp_0$ that depend on \emph{all} the data.  Indeed, this is the heart of Stein's result for Gaussian random variables and mean-squared error.  

To build intuition, we first study the specific case of minimizing mean-squared error
and show that data-pooling can improve upon the decoupled SAA solution in the frequentist framework of \cref{eq:MultinomialCounts}.  This result is thus reminiscent of Stein's classical result, but does not require the Gaussian assumptions.  
Consider the following example:
\begin{example}[A Priori-Pooling for Mean-Squared Error]  \label{ex:MSE}
Consider a special case of Problem~\eqref{eq:TargetProblem} such that for all $k$ that $\lambda_k = \lambdabar$, $\Nhat_k = \Nhat \geq 2$, $\bp_k$ is supported on $\{a_{k1}, \ldots, a_{kd}\} \subseteq \mathbb R$, $\X_k = \mathbb R$ and  $c_{ki}(x) = (x - a_{ki})^2$.  In words, the $k^\text{th}$ subproblem estimates the unknown mean $\mu_k=\bp_k^\top\bm a_k$ by minimizing the mean-squared error.  Let $\sigma_k^2=\bp_k^\top(\bm a_k-\mu_k \edit{\be})^2$.

Fix any $\bp_0 \in \Delta_d$ and $\alpha \geq 0$ (not depending on the data).  A direct computation shows that 
\[\ts
x_k(\alpha, \bp_0, \bfmhat_k) 
\ \ \equiv  \ \
\muhat_k(\alpha)
\  \ \equiv \ \ 
 \frac{\Nhat}{\Nhat + \alpha} \muhat_k + \frac{\alpha}{\Nhat + \alpha} \mu_{k0},
\]
where $\muhat_k = \frac{1}{\Nhat} \sum_{i=1}^{\Nhat} \xihat_{ki}$ is the usual sample mean, and $\mu_{k0} = \bp_0^\top \ba_k$.  Notice in particular that the 
decoupled SAA solution is $\bx^{\sf SAA} = (\muhat_1, \ldots, \muhat_K)$, corresponding to $\alpha = 0$.  

For any $\bp_0$ and $\alpha$, the objective value of $\bx(\alpha, \bp_0, \bfmhat)$ is 
\begin{align*}
\frac{1}{K} \sum_{k=1}^K \bp_k^\top\bc_k(x_k(\alpha, \bp_0, \bfmhat_k)) 
\ =  \ 
	\frac{1}{K} \sum_{k=1}^K \Eb{( \muhat_k(\alpha) - \xi_k)^2  \mid \bfmhat }
\ = \
\frac{1}{K} \sum_{k=1}^K \prns{\sigma^2_k + (\mu_k - \muhat_k(\alpha))^2},
\end{align*} 
by the usual bias-variance decomposition of mean-squared error (MSE).  This objective is the average of $K$ independent random variables.  Hence, we might intuit that under appropriate regularity conditions (see Theorem~\ref{thm:DataPoolingJS} below) that, \edit{conditional on $\Nhat$}, as $K \rightarrow \infty$, 
\label{line:FixingConvNotation}
\begin{equation} \label{eq:MSEInLimit}{\blockedit\ts
\frac{1}{K} \sum_{k=1}^K \prns{\sigma^2_k + (\mu_k - \muhat_k(\alpha))^2} 
-
\frac{1}{K} \prns{\sum_{k=1}^K \sigma^2_k + \Eb{(\mu_k - \muhat_k(\alpha))^2 \mid \Nhat}    } \  \rightarrow_p  \ 0.}
\end{equation}
Moreover, 
\(
\frac{1}{K} \prns{\sum_{k=1}^K \sigma^2_k + \Eb{(\mu_k - \muhat_k(\alpha))^2 \mid \Nhat}    }
 \ =  \ 
\frac{1}{K} \sum_{k=1}^K \prns{\sigma^2_k + \left(\frac{\alpha}{\Nhat + \alpha}\right)^2 (\mu_k - \mu_{k0})^2 + \left(\frac{\Nhat}{\Nhat + \alpha}\right)^2 \frac{\sigma^2_k}{\Nhat}},
\)
again using the bias-variance decomposition of MSE.  We can minimize the righthand side over $\alpha$ explicitly, yielding the value
\[
\alphaAP_{\bp_0} \ = \ \frac{ \sum_{k=1}^K \sigma_k^2 }{\sum_{k=1}^K (\mu_k - \mu_{k0})^2}  \ > 0,
\]
where $\sf AP$ stands for \textit{a priori}, meaning $\alphaAP_{\bp_0}$ is the on-average-best a priori choice of shrinkage before observing any data.
In particular, \edit{substituting $\alpha=0$ and $\alpha = \alphaAP_{\bp_0}$ into the second term of \cref{eq:MSEInLimit}}
\edit{shows that, up to a term that is vanishing as $K\to\infty$, shrinking by $\alphaAP_{\bp_0}$ decreases the MSE by}
\begin{equation}\label{eq:mseimproveap}
\left( \frac{1}{K}  \sum_{k=1}^K \frac{\sigma_k^2}{\Nhat}\right)    \frac{ \alphaAP_{\bp_0} } { \Nhat+ \alphaAP_{\bp_0}} 
\ \ = \ \ 
\frac{  \left( \frac{1}{K \Nhat } \sum_{k=1}^K \sigma_k^2 \right)^2 } { \frac{1}{K\Nhat} \sum_{k=1}^K \sigma_k^2 + \frac{1}{K} \sum_{k=1}^K (\mu_k - \mu_{k0})^2 }
\ \ >  \ \
 0.
 \end{equation}
This benefit is strictly positive for any values of $\bp_k$ and $\bp_0$, and increasing in $\alphaAP_{\bp_0}$.
\end{example}

Unfortunately, we cannot implement $x(\alphaAP_{\bp_0}, \bp_0, \bfmhat)$ in practice because $\alphaAP_{\bp_0}$ is not computable from the data; it depends on the unknown $\mu_k$ and $\sigma^2_k$.  The next theorem shows that we can, however, estimate $\alphaAP_{\bp_0}$ from the data in a way that achieves the same benefit as $K \rightarrow \infty$, even if $\Nhat$ is fixed and small.   See Appendix~\ref{sec:ProofofThmDataPoolingJS} for proof.
\begin{theorem}[Data-Pooling for MSE] \label{thm:DataPoolingJS}
Consider a sequence of subproblems, indexed by $k = 1, 2, \ldots.$  Suppose for each $k$, the $k^\text{th}$ subproblem minimizes mean-squared error, i.e., $\bp_k$ is supported on $\{a_{k1}, \ldots, a_{kd}\} \subseteq \mathbb R$, $\X_k = \mathbb R$ and  $c_{ki}(x) = (x - a_{ki})^2$.  Suppose further that there exists $\lambdabar$, $\Nhat\geq 2$ and $a_{\rm max} < \infty$ such that $\lambda_k = \lambdabar$, $\Nhat_k = \Nhat$, and $\| \ba_k \|_\infty \leq a_{\rm max}$ for all $k$.  
Fix any $\bp_0 \in \Delta_d$, and let 
\[
\alphaJS_{\bp_0} = \frac{ \frac{1}{K} \sum_{k=1}^K  \frac{1}{\Nhat - 1} \sum_{i=1}^{\Nhat} (\xihat_{ki} - \muhat_k)^2 }
				  {   \frac{1}{K} \sum_{k=1}^K (\mu_{k0} - \muhat_k)^2  - \frac{1}{K\Nhat} \sum_{k=1}^K  \frac{1}{\Nhat - 1} \sum_{i=1}^{\Nhat} (\xihat_{ki} - \muhat_k)^2 }.
\]  
Then, conditional on $\Nhat$, as $K \rightarrow \infty$, 
\begin{align*} 
\underbrace{\frac{1}{K} \sum_{k=1}^K \bp_k^\top \bc_k(\bx_k^{\sf SAA}) 
	- \frac{1}{K} \sum_{k=1}^K \bp_k^\top \bc_k(\bx_k(\alphaJS_{\bp_0}, \bp_0, \bfmhat_k ) )}_{\text{Benefit over decoupling of $\alpha=\alphaJS_{\bp_0}$}}
\ - \ 
\underbrace{\frac{  \left( \frac{1}{K} \sum_{k=1}^K \sigma_k^2 / \Nhat \right)^2 } { \frac{1}{K} \sum_{k=1}^K \sigma_k^2 / \Nhat + \frac{1}{K} \sum_{k=1}^K (\mu_k - \mu_{k0})^2 }}_{\text{Expected benefit over decoupling of $\alpha=\alphaAP_{\bp_0}$}}
\ \rightarrow_p \ 0.
\end{align*}
\end{theorem}

Note that $x_k(\alphaJS_{\bp_0}, \bp_0, \bfmhat) = (1-\theta) \muhat_k + \theta \muhat_{k0}$ where 
\(
\theta  \ =  \  \frac1{\Nhat}\frac{\frac{1}{K} \sum_{k=1}^K  \frac{1}{\Nhat - 1} \sum_{i=1}^{\Nhat} (\xihat_{ki} - \muhat_k)^2 } {\frac{1}{K} \sum_{k=1}^K (\mu_{k0} - \muhat_k)^2 }.
\)
In this form, we can see that the resulting estimator with pooling $\alphaJS_{\bp_0}$ strongly resembles the classical James-Stein mean estimator (cf. \citealp[Eq. (7.51)]{efron2016computer}), with the exception that we have replaced the 
variance $\sigma_k^2$, which is assumed to be $1$ in Stein's setting,  with the usual, unbiased estimator of that variance.
This resemblance  motivates our ``$\sf JS$'' notation.  Theorem~\ref{thm:DataPoolingJS} is neither stronger nor weaker than the James-Stein theorem.   \edit{Our} result applies to non-Gaussian random variables and holds in probability, but is asymptotic;  the James-Stein theorem requires Gaussian distributions and holds in expectation, but applies to any fixed $K \geq 3$.

Theorem~\ref{thm:DataPoolingJS} shows that data-pooling for mean-squared error always offers a benefit over decoupling for sufficiently large $K$, no matter what the $\bp_k$ may be.  
Data-pooling for general optimization problems, however, exhibits more subtle behavior.  In particular, as shown in the following example and theorem, there exist instances where data-pooling offers no benefit over decoupling, and instances where data-pooling may be worse than decoupling.  

\begin{example}[Data-Pooling for Simple Newsvendor] 
\label{ex:BabyNewsvendor}
Consider a special case of Problem~\eqref{eq:TargetProblem} such that for all $k$, $\lambda_k = \lambdabar$,
$\edit{\xi_k}$ is supported on $\{ 1, 0 \}$, $\X_k = [0,1]$ and $c_k(x, \xi_k) = \abs{x - \xi_k}$ so that $\bp_k^\top\bc_k(x) = p_{k1} + x(1 - 2p_{k1})$.  In words, the $k^\text{th}$ subproblem estimates the median of a Bernoulli random variable by minimizing mean absolute deviation, or, equivalently, is a newsvendor problem with critical fractile $0.5$ for Bernoulli demand.  
We order the support so that $p_{k1} = \P( \xi_k = 1)$, as is typical for a Bernoulli random variable.  
Suppose further for each $k$, $p_{k1} > \frac{1}{2}$, and fix any $p_{01} < \frac{1}{2}$. 

Note $x_k(\alpha, \bp_0, \bfmhat_k) = \I{  \phat_{k1} \geq \frac{1}{2} + \frac{\alpha}{\Nhat_k}(\frac{1}{2} - p_{01})}$.\footnote{This solution is non-unique, and the solution $ \I{  \phat_{k1} > \frac{1}{2} + \frac{\alpha}{\Nhat_k}(\frac{1}{2} - p_{01})}$ is also valid.  We adopt the former solution in what follows, but our comments apply to either solution.}  
Further, for any $\alpha$ (possibly depending on $\bfmhat$),
\begin{align*}
\bp_k^\top \left( \bc_k(\bx_k(\alpha, \bp_0, \bfmhat_k)) - \bc_k(\bx_k(0, \bp_0, \bfmhat_k)) \right)
&=
(2 p_{k1} - 1) \left(\I{\phat_{k1} \geq 1/2} - \I{\phat_{k1} \geq \frac{1}{2} + \frac{\alpha}{\Nhat_k}\left(\frac{1}{2} - p_{01}\right)} \right)
\\ &=
(2 p_{k1} - 1) \I{1/2  \ \leq  \  \phat_{k1}  \ <  \ \frac{1}{2} + \frac{\alpha}{\Nhat_k}\left(\frac{1}{2} - p_{01}\right)} , 
\end{align*}
where the last equality follows since $\phat_{k1} < 1/2 \implies \phat_{k1} < \frac{1}{2} + \frac{\alpha}{2}(\frac{1}{2} - p_{01}) $.  Notice $p_{k1} > \frac{1}{2} \implies (2 p_{k1} - 1) > 0$, so this last expression is nonnegative.  
It follows that path by path, shrinkage by any $\alpha > 0$ cannot improve upon the decoupled solution ($\alpha = 0$).  Moreover, if $\bx_k(\alpha, \bp_0, \bfmhat_k) \neq \bx_k(0, \bp_0, \bfmhat_k)$, the performance is strictly worse.

One can check directly that if we had instead chosen $p_{01} \geq \frac{1}{2}$ and $p_{k1} < \frac{1}{2}$, a similar result holds.
\end{example}
We summarize this example in the following theorem:
\begin{theorem}[Data-Pooling Does Not Always Offer Benefit] 
\label{thm:DataPoolingDoesntHelp}
Given any $\bp_0$, there exist instances of Problem~\eqref{eq:TargetProblem} such that \edit{shrinkage does not outperform the (decoupled) SAA solution.}
\edit{Moreover, if $\bx(\alpha, \bp_0, \bfmhat)$ performs the same as SAA, then $\bx(\alpha, \bp_0, \bfmhat)$ is, itself, an SAA solution.}  
\end{theorem}

On the other hand, there exist examples where the James-Stein estimator and {traditional statistical reasoning} might suggest the benefits of pooling are marginal, but, by data-pooling in way that exploits the optimization structure, we can achieve significant benefits.  \edit{Specifically, our Bayesian motivation in \cref{sec:Bayes} suggests pooling offers little benefit when the $\bp_k$ are very dispersed, i.e., the Dirichlet prior has high variance and $\alpha_0$ is small.  Similarly, \cref{thm:DataPoolingJS} and \citet{efron1977stein} both suggest that the benefits of pooling over decoupling for MSE are marginal if the subproblem means are quite dispersed (cf. \cref{eq:mseimproveap}). \label{line:DispersedMeans1} }  Nonetheless, for general optimization problems, we observe pooling might still offer substantive benefits in these situations:  
\begin{example}[Pooling Can Offer Benefit Even When $\bp_k$ are Dispersed]  \label{ex:DifferentNewsvendor}
Let $d > 3$ and fix some $0 < s < 1$. Suppose the $k^\text{th}$ subproblem is a newsvendor problem with critical fractile $f_k > s$ and demand distribution supported on the integers $1, \dots, d$.    For each $k$,  let $p_{k1} = 0$ , $p_{kd} = 1-s$, and $p_{kj_k} = s$ for some $ 1 < j_k < d$.  Consider the fixed anchor $p_{01} = s$, $p_{0d} = 1-s$, and $p_{0j} = 0$ for $ 1 < j < d$.  Notice typical $\bp_k$'s are very far from $\bp_0$ since $\| \bp_k - \bp_0 \|_2 = \sqrt{2} s$.  For $s$ sufficiently close to $1$, this value is close to $\sqrt 2$, which is the maximal distance between two points on the simplex.  In other words, the $\bp_k$ are not very similar.   \edit{Moreover, the means are also dispersed for $s$ close to $1$ since $\frac{1}{K} \sum_{k=1}^K (\mu_k - \mu_0)^2 = s^2 \frac{1}{K} \sum_{k=1}^K (j_k-1)^2 \approx s^2 d/2$ if the $j_k$ are chosen uniformly.\label{line:DispersedMeans2}}  

Consequently, the James-Stein estimator does not shrink very much in this example.  A straightforward computation shows that for $K$ sufficiently large, $\alpha^{\sf JS}_{\bp_0} \leq \frac{(1-s)d^2}{s}$ with high probability, which is close to $0$ for $s$ close to $1$.  However, the full-information solution for the $k^{\text{th}}$ problem is $\bx^*_k=d$, which \emph{also} equals the fully-pooled ($\alpha=\infty$) solution, $\bx_k(\infty, \bp_0)$.  Hence, pooling in an optimization-aware way can achieve full-information performance, while both decoupling and an ``estimate-then-optimize'' approach using James-Stein shrinkage \emph{necessarily} perform worse. 
In other words, pooling offers significant benefits despite the $\bp_k$ being as dispersed as possible, because of the optimization structure, and leveraging this structure is necessary to obtain the best shrinkage.  \Halmos
\end{example}

\Cref{thm:DataPoolingJS,thm:DataPoolingDoesntHelp,ex:BabyNewsvendor,ex:DifferentNewsvendor} highlight the fact that data-pooling for general optimization is more complex
than Stein's phenomenon. In particular, in Stein's classical result for mean-squared error and Gaussian data, data-pooling \emph{always} offers a benefit for $K\geq3$.  For other optimization problems and data distributions, data-pooling may \emph{not} offer a benefit, or may offer a benefit but requires a new way of choosing the pooling amount. An interplay between $\bp_0$, $\bp_k$ and $\bc_k$ determines if data-pooling can improve upon decoupling and how much pooling is best.  

This raises two important questions:  First, how do we identify if an instance of Problem~\eqref{eq:TargetProblem} would benefit from data-pooling?  Second, if it does, how do we compute the ``optimal" amount of pooling? 
In the next sections, we show how our Shrunken-SAA algorithm can be used to address both questions in the relevant regime, where $K$ is large but the average amount of data per subproblem remains small.  Indeed, we show that Shrunken-SAA \edit{achieves} the 
best-possible shrinkage in an optimization-aware fashion for many types of problems \edit{and choices of anchor}.

%
\section{The Shrunken SAA Algorithm}
\label{sec:ShrunkenSAA}
{\blockedit
\begin{algorithm}[t!]  \small\blockedit 
\caption{\textbf{The Shrunken-SAA Algorithm.} 
}\label{alg:ssaa}
\begin{algorithmic}
\State \textbf{Input:} Data $\S_k = \{\bxihat_{k1}, \ldots, \bxihat_{k\Nhat_k}\}$, $k=1, \ldots, K$, and an anchor distribution $h(\S)$
\State Fix a finite grid $\mathcal A \subseteq [0, \infty)$
\State \textbf{for} $\alpha\in\mathcal A$, \ \  $k=1,\dots,K$,  \ \ $j=1,\dots,\Nhat_k$ \textbf{ define:}
    \State \quad 
    \(
    \bx_{k, -j}(\alpha, h(\S))\ \gets\ \arg\min_{ \bx_k \in \X_k} \  
	\sum_{\ell \neq j} c_k(\bx_k, \bxihat_{k\ell} ) + \alpha \E_{\bxi_k \sim h(\S)}\left[ c_k(\bx_k, \bxi_k) \right]
   \)
   \hfill   // Compute LOO solutions
\State \textbf{end for}

 \State $\alphaLOO_h\ \gets\ \arg\min_{\alpha\in\mathcal A}\sum_{k=1}^K \sum_{j=1}^{\Nhat_k}c_{k}(\bx_{k,-j}(\alpha, h(\S)), \bxihat_{kj})$  \hfill   // Modified LOO-Cross-Validation
  \ForAll{$k=1,\dots,K$}
    \State 
    \(
    \bx_k^\SSAA\ \gets\ \arg\min_{ \bx_k \in \X_k} \  
	\sum_{j=1}^{\Nhat_k} c_k(\bx_k, \bxihat_{kj} ) + \alphaLOO_h \E_{\bxi_k \sim h(\S)}\left[ c_k(\bx_k, \bxi_k) \right]
    \)\hfill // Compute Pooled solution
 \EndFor
    \State \Return $\left( \bx_1^\SSAA, \ldots, \bx_K^\SSAA\right)$
\end{algorithmic}
\end{algorithm}

}

{\blockedit
\Cref{alg:ssaa} formally defines Shrunken-SAA. 
\editt{The crucial step is the ``Modified LOO-Cross-Validation," which we discuss in detail in \cref{sec:UnbiasedRiskEstimate,sec:MotivatingViaLOO} below.}
To highlight similarities to SAA, we have stated the algorithm in terms of the datasets $\S_k$ and $\S=(\S_1,\dots,\S_K)$. 
Here $h(\S)$ represents an arbitrary, possibly data-driven anchor distribution (see below for examples).
Recall that we can equivalently express $\S_k$ in terms of the counts $\bfmhat_k$.  In that notation, we recognize that if the $j^\text{th}$ data point of $\S_k$ is $\ba_{ki}$, then $\bx_{k, -j}(\alpha, h(\bfmhat)) = \bx_k(\alpha, h(\bfmhat), \bfmhat_k - \be_i)$  and $\bx_k^\SSAA = \bx_k(\alphaLOO, h(\bfmhat), \bfmhat_k)$.  In other words, Shrunken-SAA retains the particular pooling structure of \cref{eq:ShrunkXSol} suggested by our Bayesian argument, but allows for a data-dependent anchor $h(\S)$ (equiv. $h(\bfmhat)$) and chooses the amount of pooling via a particular cross-validation scheme.  We present \Cref{alg:ssaa} using a finite grid of $\alpha\in\mathcal A$, but our theory below will study the algorithm with $\mathcal A=[0,\infty)$.
}

\begin{remark}[Computational Complexity]  \label{rem:Computational}
{\blockedit
Computationally, \cref{alg:ssaa} does not depend on $d$, the size of the support of $\bxi_k$.
Its bottleneck is computing $\bx_{k, -j}$ which is similar to solving the $k^\text{th}$ subproblem by SAA with an augmented data set described by $h(\S)$.
More specifically, \cref{alg:ssaa} depends on the data only through $h(\S)$ and averages of functions over subsets of $\S$, neither of which explicitly depend upon $d$.  
Consequently, although our setup and analysis assumes $\bxi_k$ has finite discrete support, from an implementation perspective, we can apply Shrunken-SAA when $\bxi_k$ has continuous support \emph{without} discretization so long as we can efficiently solve these augmented SAA problems
(cf. our empirical study in \cref{sec:Infinited}).  From a theoretical perspective, some of our analysis extends to this continuous setting (see \cref{remark:continuousanalysis} below).  
In the remainder, we follow \cref{sec:Model} and  treat the data as discrete, referring to the data by $\bfmhat_k$ and $\bfmhat$.
}

We consider Shrunken-SAA to be \emph{roughly} as tractable as SAA.  \edit{We say ``roughly" because, in the worst-case, one must solve at most $\abs{\mathcal A}\sum_{k=1}^K\min(d, \Nhat_k)$ 
problems in the LOO-cross-validation step, which, if we sample from $h(\bfmhat)$, have a similar structure to SAA.} 
Fortunately, we can parallelize these problems in distributed computing environments and  use previous iterations to ``warm-start" solvers. \edit{Moreover, in \cref{sec:KFoldCrossVal} we observe empirically that less computationally expensive $\kappa$-fold cross-validation procedures can be used in place of LOO with similar performance.}  \Halmos
\end{remark}

{\blockedit For clarity, the $\alphaLOO_{h}$ parameter (with $\mathcal A=[0,\infty)$) computed by \cref{alg:ssaa} is 
\begin{align}\label{eq:alphaloo}
\alphaLOO_{h}&
\textstyle
\in \arg \min_{\alpha \geq 0 } \  \sum_{k=1}^K  \bfmhat_{k}^\top \bc_{k}( \bx_k(\alpha,h(\bfmhat), \bfmhat_k - \be_i ) ).
\end{align}

\vskip 6pt
\noindent \textbf{The Anchor Distribution $h(\bfmhat)$}

As stated, the anchor in \cref{alg:ssaa}, $h(\bfmhat)$, is an input.  We think of $h(\bfmhat)$ as a function that selects an anchor distribution from a candidate set of distributions $\mathcal P$.  In what follows, we will focus on two types of anchors and corresponding candidate sets $\mathcal P$:  
\begin{itemize}[leftmargin=*]
\item \textbf{Fixed Anchors}: In this case, $h(\bfmhat)=\bp_0$, $\mathcal P = \{ \bp_0 \}$ for some fixed $\bp_0$, e.g., 
the uniform distribution  $\bp_0=\be/d$.  In general, fixed-anchors might be used 
for computational/statistical simplicity or when there is strong a priori knowledge of a good anchor.  In this special case, we abuse notation slightly, replacing the map $h: \bfmhat \mapsto \bp_0$ with the constant $\bp_0$ when it is clear from context, e.g., we write $\alphaLOO_{\bp_0}$ for $\alphaLOO_h$.
\item \textbf{Data-Driven Anchors}: In this case $h(\bfmhat)$ is any procedure that uses the data $\bfmhat$ to select a distribution, and $\mathcal P$ is the image of $h(\cdot)$.  One example might be to use all the data to fit a parametric distribution, e.g., a lognormal distribution, via maximum likelihood and use this fitted distribution as the anchor.  Then, $\mathcal P$ would be the set of lognormal distributions.
\end{itemize}

\vskip 8pt 
We also pay particular focus to two special cases of data-driven anchors in what follows:
\begin{itemize} [leftmargin=*]
\item \textbf{LOO-Optimized Anchor}: For a given $\mathcal P\subseteq\Delta_d$, let
\begin{equation}\label{eq:hloo} \textstyle
\hloo(\bfmhat) \in \argmin_{\bq\in\mathcal P}\min_{\alpha\in\mathcal A} \ 
	\sum_{k=1}^K\bfmhat_k^\top\bc_k(\bx_k(\alpha, \bq, \bfmhat_k - \be_i).
\end{equation}
We will see below that $\hloo$ satisfies stronger optimality properties than general data-driven anchors and, hence, we treat it separately.  From an implementation point of view, when applying \cref{alg:ssaa}, we only ever require the value of $\hloo(\bfmhat)$, not the full-function $\hloo(\cdot)$.  Thus, \cref{alg:ssaa} with $\hloo(\cdot)$ amounts to replacing the ``Modified LOO-Cross-Validation" step by a joint optimization over anchor and pooling amount:
\begin{equation} \label{eq:ComputingHlooAlphaLOO}  \textstyle
( \alphaLOO_{\hloo}, \hloo(\bfmhat) )\ \gets\ \argmin_{\alpha\in\mathcal A, \bq \in \mathcal P} \ 
 \sum_{k=1}^K \sum_{j=1}^{\Nhat_k}c_{k}(\bx_{k,-j}(\alpha, \bq), \bxihat_{kj}).
\end{equation}
We note that the multivariate optimization problem in \cref{eq:ComputingHlooAlphaLOO} may be challenging depending on the structure of $\mathcal P$, motivating our second special case below.  

\item \textbf{GM-Anchor}  
\label{Reviewer:GMAnchorChoice}
We also consider a computationally simpler ``grand-mean" anchor $h(\bfmhat) = \bphat^{\sf GM}$ where $\bphat^{\sf GM} \equiv \sum_{k=1}^K  \bphat_k \I{\Nhat_k > 0}  / \sum_{k=1}^K \I{\Nhat_k > 0}$ if $\Nhat_{\max} > 0$ and $\be/d$ otherwise.  (For this data-driven anchor, $\mathcal P = \Delta_d$.)
This choice is motivated by our Bayesian perspective on data-pooling from \cref{sec:Bayes}.  In the Bayesian setting $\bphat^{\sf GM}$ is an unbiased estimator of the prior mean.  We observe empirically in \cref{sec:Numerics} that $\bphat^{\sf GM}$ is a strong and computationally-efficient heuristic.  
\end{itemize}
}

\subsection{\edit{Oracle Benchmarks}}
\label{sec:OracleBenchmarks}
From \cref{thm:DataPoolingDoesntHelp}, data-pooling need not improve upon decoupling for a given $h(\cdot)$.  To establish appropriate benchmarks, we first define the \emph{oracle} pooling for given $h(\cdot)$, i.e., 
\begin{align} \label{eq:alphaOR}
\alphaOR_{h} \in \argmin_{\alpha \geq 0 } 
\Zperf(\alpha,h(\bfmhat)),
\quad\text{where}
\quad 
\Zperf(\alpha,\bq)
	&=  \ts
	\frac{1}{K} \sum_{k=1}^K Z_{k}(\alpha,\bq),
\\  \notag
Z_{k}(\alpha,\bq)
	&= \ts
	\frac{\lambda_k}{\lambdabar}{\bp_k}^\top \bc_k( \bx_k(\alpha, \bq,  \bfmhat_k) ).
\vspace{-5pt}
\end{align}

Notice $\alphaOR_{h}$ is random, depending on the entire data-sequence.  By construction, $\Zperf(\alphaOR_{h}, h(\bfmhat))$ lower bounds the performance of \emph{any} other data-driven pooling policy with anchor $h(\bfmhat)$ path-by-path.  Hence, it serves as a strong performance benchmark.  However, $\alphaOR_{h}$ also depends on the unknown $\bp_k$ and $\lambda_k$, and hence, is not implementable in practice.  In this sense, it is an oracle. 

Given any $\alpha$ (possibly depending on the data), we measure the sub-optimality of pooling by $\alpha$ relative to the oracle \edit{pooling for $h(\cdot)$} on a particular data-realization by
\[
\edit{{\sf Sub Opt}_{h,K}(\alpha)}=
\Zperf(\alpha,h(\bfmhat))
- 
\Zperf(\alphaOR_h,h(\bfmhat)).
\]
Good pooling procedures will have small sub-optimality with high-probability with respect to the data.  
Note we allow for the possibility that $\alphaOR_{h} = 0$, as is the case in \cref{ex:BabyNewsvendor}.  \edit{Thus, procedures that have small sub-optimality will still have good performance in instances where data-pooling is not beneficial.\label{reviewerSubOptNonneg}}  Moreover, studying when $\alphaOR_{h} > 0$ gives intuition into when and why data-pooling is helpful, a task we take up in \cref{sec:StabilityIntuition}.  

{\blockedit
The above oracle is defined with respect to a given anchor.  One might also seek to benchmark performance relative to the best-possible anchor.  
Given any $\mathcal P\subseteq\Delta_d$,
we define the oracle choice of anchor and pooling amount for anchors in $\mathcal P$ and for a particular data realization by 
\begin{equation} \label{eq:DefinesORValuesmP}
\ts 
(\alphaOR_\mP,\,\bq^\OR_\mP) \in \argmin_{\alpha \geq 0,\,\bq\in\mP} \Zperf(\alpha,\bq).
\end{equation}
Then, given any anchor $\bq \in \mathcal P$ and pooling amount $\alpha$ (both possibly depending the data), we measure the sub-optimality of shrinking by $\alpha$ towards $\bq$ by
\[
\edit{{\sf Sub Opt}_{\mP,K}(\alpha,\bq)} \ = \ 
\Zperf(\alpha,\bq)
- 
\Zperf(\alphaOR_\mP,\bq^\OR_\mP).
\]
For clarity, we observe that by construction  $\alphaOR_\mP=\alphaOR_{\bq^\OR_\mP}$.

}

\subsection{Motivating $\alphaLOO$ through Unbiased Estimation}
\label{sec:UnbiasedRiskEstimate}

\edit{ 
We first consider a \edit{fixed} anchor $h(\bfmhat) = \bp_0$.  Recall in this case, we abuse notation slightly, writing 
\begin{equation}\ts
\label{eq:alphaORFixedAnchor}
\alphaOR_{{\bp_0}} \in \argmin_{\alpha \geq 0 } \Zperf(\alpha, \bp_0)
\end{equation}}%
One approach to choosing $\alpha_{\bp_0}$ might be to construct a suitable proxy for $\Zperf(\alpha, \bp_0)$ in \cref{eq:alphaORFixedAnchor} based only on the data, and then choose the $\alpha_{\bp_0}$ that optimizes this proxy.  

If we knew the values of $\lambda_k$, a natural proxy might be to replace the unknown $\bp_k$ with $\bphat_k$, i.e., optimize 
\(
\frac{1}{K} \sum_{k=1}^K 
\frac{\lambda_k}{\lambdabar}{\bphat_k}^\top \bc_k( \bx_k(\alpha, \bq,  \bfmhat_k) ).
\)
Unfortunately, even for a fixed, non-data-driven $\alpha$, this proxy is \emph{biased}, i.e. $\Eb{ \frac{1}{K} \sum_{k=1}^K 
\frac{\lambda_k}{\lambdabar}{\bphat_k}^\top \bc_k( \bx_k(\alpha, \edit{\bp_0},  \bfmhat_k) )
 }  \neq \Eb{\Zperf(\alpha, \bp_0)}$, since both $\bphat_k$ and $\bx_k(\alpha, \bp_0, \bfmhat_k)$ depend on the data $\bfmhat_k$.  Worse, this bias wrongly suggests $\alpha = 0$, i.e. decoupling, is always a good policy, because $\bx_k(0, \bp_0, \bfmhat_k)$ always optimizes this proxy, by construction.  By contrast, \cref{thm:DataPoolingJS} shows data-pooling can offer significant benefits.  This type of bias and its consequences are well-known in other contexts and are often termed 
 the ``optimizer's curse" --  in-sample costs are optimistically biased and may not generalize well.  

These features motivate us to seek an unbiased estimate of $\Zperf(\alpha, \bp_0)$.  At first glance, however, $Z_K(\alpha, \bp_0)$, which depends on both the unknown $\bp_k$ and unknown $\lambda_k$, seems particularly intractable unless $\bx_k(\alpha, \bp_0, \bfmhat_k)$ admits a closed-form solution as in \cref{ex:MSE}.  A key observation is that, in fact, $\Zperf(\alpha, \bp_0)$ does more generally admit an unbiased estimator, \emph{if} we also introduce an additional assumption on our data-generating mechanism, i.e., that the amount of data is random.  
\begin{assumption}[Randomizing Amount of Data]  \label{ass:RandomData}
There exists an $N$  such that $\Nhat_k \sim \op{Poisson}(N\lambda_k)$ for each $k = 1, \ldots, K$.  
\end{assumption}
Under \cref{ass:RandomData}, (unconditional) expectations and probabilities should be interpreted as over both the random draw of $\Nhat_k$ and the counts $\bfmhat_k$.  

{\blockedit \label{WhyRandomData} 
Analytically, the benefit of \cref{ass:RandomData} is that it breaks the dependence across $i$ in $\bfmhat_k$.  Namely, by the Poisson-splitting property, under \cref{ass:RandomData},
\[
\mhat_{ki} \sim \op{Poisson}(m_{ki}) \ \  \text{where} \ \ m_{ki} \equiv N \lambda_k p_{ki}, \quad  i = 1, \ldots, d,  \quad k = 1, \ldots, K,
\]
and, furthermore, the $\mhat_{ki}$ are independent across $i$ and $k$. Notice if $\Nhat_k$ were non-random, these $\mhat_{ki}$ would be dependent.}

\label{ApplicationsRandomData}
\edit{Beyond its analytical convenience, we consider \cref{ass:RandomData} to be reasonable in many applications.  Consider 
for instance a retailer optimizing the price of $k$ distinct products, i.e., $x_k$ represents the price of product $k$, $\xi_k$, represents the (random) valuation of a typical customer, and $c_k(x_k, \xi_k)$ is the (negative) profit earned. 
In such settings, one frequently ties data collection to time, i.e., one might collect $N=6$ months worth of data.  To the extent that customers arrive seeking product $k$ in a random fashion, the number of arrivals $\Nhat_k$ that one might observe in $N$ months is, itself, random, and reasonably modeled as Poisson with rate proportional to $N$. Similar statements apply whenever data for problem $k$ is generated by an event which occurs randomly, e.g., when observing response time of emergency responders (disasters occur intermittently), effectiveness of a new medical treatment (patients with the relevant disease arrive sequentially), or any aspect of a customer service interaction (customers arrive randomly to service). }

In some ways, this perspective tacitly underlies the formulation of Problem~\eqref{eq:TargetProblem}, itself.  Indeed, one way to interpret the subproblem weights $\frac{\lambda_k}{K \lambdabar} = \frac{\lambda_k}{\sum_{j=1}^K \lambda_j}$ is that the decision-maker incurs costs $c_k(x_k, \xi_k)$ at rate $\lambda_k$, so that problems of type $k$ contribute a $\frac{\lambda_k}{\sum_{j=1}^K \lambda_j}$ fraction of the total long-run costs.  However, if problems of type $k$ occur at rate $\lambda_k$, it should be that observations of type $k$, i.e. realizations of $\bxi_k$, also occur at rate $\lambda_k$, supporting \cref{ass:RandomData}.

In settings where data-collection is not tied to randomly occurring events, modeling $\Nhat_k$ as Poisson may still be a reasonable approximation if $d$ is large relative to $\Nhat_k$ and each of the individual $p_{ki}$ are small.  Indeed, under such assumptions, a $\op{Multinomial}(\Nhat_k, \bp_k)$ is well-approximated by independent Poisson random variables with rates $\Nhat_kp_{ki}$, $i = 1, \ldots d$ (see \citealp{mcdonald1980poisson,deheuvels1988poisson} for a formal statement).  In this sense, we can view the consequence of \cref{ass:RandomData} as a useful approximation to the setting where $\Nhat_k$ are fixed, even if it is not strictly true.

In any case, under \cref{ass:RandomData}, we develop an unbiased estimate for $\Zperf(\alpha, \bp_0, \bfmhat)$.  We use the following identity \citep{chen1975poisson}.  For any $f:\mathbb Z_+\to\Rl$, for which the expectations exist, 
\begin{equation}\label{eq:steinchen}
W\sim\op{Poisson}(\lambda) \implies \lambda\E[f(W+1)]=\E[Wf(W)].
\end{equation}
The proof of the identity is immediate from the Poisson probability mass function.\footnote{In particular,
$\E[Wf(W)] = \sum_{w=0}^\infty w f(w) e^{-\lambda} \frac{\lambda^w}{w!} 
= 
\lambda \sum_{w=0}^\infty f(w) e^{-\lambda} \frac{\lambda^{w-1}}{(w-1)!}
= \lambda \E[f(W+1)]$.}   

\enlargethispage*{20pt}
Now, for any $\alpha \geq 0$ and $\bq \in \Delta_d$, define 
\begin{align}  \label{eq:DefZLOO}
Z^{\sf LOO}_{k}(\alpha,\bq) \equiv \frac{1}{N\lambdabar} \sum_{i=1}^d  \mhat_{ki} c_{ki}( \bx_k(\alpha, \bq,  \bfmhat_k - \be_i ) ),
\ \ \text{ and } \ \ 
\Zloo(\alpha, \bq) \equiv    \frac{1}{K}  \sum_{k=1}^K   Z^{\sf LOO}_{k}(\alpha, \bp_0).
\end{align}
\begin{lemma}[An Unbiased Estimator for $\Zperf(\alpha, \bp_0)$]  %
\label{cor:steinchen}
Under \cref{ass:RandomData}, we have  for any $\alpha\geq0,$ and $\bq\in\Delta_d$ that 
\(
\Eb{Z^{\sf LOO}_{k}(\alpha,\bq)}  =   \Eb{Z_{k}(\alpha,\bq)}.
\)
In particular, 
\(
\Eb{ \Zloo(\alpha, \bq) } = \Eb{ \Zperf(\alpha, \bq)}.
\)  
\end{lemma}
\proof{Proof.} 
Recall that
\(
	Z_{k}(\alpha,\bq)=\frac1{N\lambdabar}\sum_{i=1}^d m_{ki}c_{ki}( \bx_k(\alpha, \bq,  \bfmhat_k ) )
\)
and that under \cref{ass:RandomData} $\mhat_{ki}\sim\op{Poisson}(m_{ki})$ independently over $i=1,\dots,d$. Let $\mhat_{k,-i}$ denote $\prns{\mhat_{k,j}}_{j\neq i}$.  Then, 
by \cref{eq:steinchen}, 
\[
\Eb{m_{ki}
c_{ki}( \bx_k(\alpha, \bq,  \bfmhat_k ) )
\mid \mhat_{k,-i}
}
 \  =   \ 
 \Eb{
\mhat_{ki}c_{ki}( \bx_k(\alpha, \bq,  \bfmhat_k - \be_i) )
\mid \mhat_{k,-i}
}.
\]
Taking expectations of both sides, summing over $i = 1, \ldots, d$ and scaling by $N \lambdabar$ proves $\Eb{Z^{\sf LOO}_{k}(\alpha,\bq)}  =   \Eb{Z_{k}(\alpha,\bq)}$.  Finally, averaging this last equality over $k$ completes the lemma. 
\endproof

We therefore propose selecting $\alpha$ by minimizing the estimate $\Zloo(\alpha, \bp_0)$.  As written, $\Zloo(\alpha, \bp_0)$ still depends on the unknown $N$ and $\lambdabar$, however, these values occur multiplicatively and are positive, and so do not affect the optimizer.  \edit{Hence, the optimizer is exactly $\alphaLOO_h$ as in \cref{eq:alphaloo}.}
\subsection{Motivating $\alphaLOO$ via {Modified} Leave-One-Out Cross-Validation}
\label{sec:MotivatingViaLOO}
Although we motivated \cref{eq:alphaloo} via an unbiased estimator, we can alternatively motivate it through leave-one-out cross-validation.  This latter perspective informs our ``{\sf LOO}'' notation above.  Indeed, consider again our decision-maker, and assume in line with \cref{ass:RandomData} that subproblems of type $k$ arrive randomly according to a Poisson process with rate $\lambda_k$, independently across $k$.  When a problem of type $k$ arrives, she incurs  a cost $c_k(\bx_k, \bxi)$.  Again, the objective of Problem~\eqref{eq:TargetProblem} thus represents her expected, long-run costs.

We can alternatively represent her costs via the modified cost function 
$C\left(\bx_1, \ldots, \bx_K,\kappa,\bxi \right)=c_\kappa( \bx_\kappa,\bxi)$, where $\kappa$ is a random variable indicating which of the $k$ subproblems she is currently facing.  In particular, letting $\P( \kappa = k) = \frac{\lambda_k}{K \lambdabar}$ and $\P(\bxi=a_{ki}\mid \kappa=k)=p_{ki}$,
the objective of Problem~\eqref{eq:TargetProblem} can be more compactly written
\(
\Eb{ C\left( \bx_1, \ldots, \bx_K, \kappa, \bxi \right)}.
\)

Now consider pooling all the data into a single ``grand" data set of size $\Nhat_1+\cdots+\Nhat_K$:
\[
\braces{(k,\,\bxi_{kj}):j=1,\dots,\Nhat_k,\,k=1,\dots,K}.
\]
The grand dataset can be seen as i.i.d. draws of $(\kappa, \bxi)$.  

For a fixed $\alpha$ and $\bp_0$, the leave-one-out estimate of $\Eb{ C\left( \bx_1(\alpha, \bp_0, \bfmhat), \ldots, \bx_K(\alpha, \bp_0, \bfmhat), \kappa, \bxi \right)}$ is given by removing one data point from the grand data set, training $\bx_1(\alpha, \bp_0, \cdot), \ldots, \bx_K(\alpha, \bp_0, \cdot)$ on the remaining data, and evaluating $C(\cdot)$ on the left-out point using these policies.  
We repeat this procedure for each point in the grand data set and average.  After some bookkeeping, we can write this leave-one-out estimate as
\[
\frac{1}{\sum_{k=1}^K \Nhat_k} \sum_{k=1}^K \sum_{i=1}^{d} \mhat_{ki} c_{ki}( \bx_k(\alpha,\bp_0, \bfmhat_k - \be_i ) ),
\]
which agrees with the objective of \cref{eq:alphaloo} up to a positive multiplicative constant.  Although this multiplicative constant does not affect the choice of $\alphaLOO$, it \emph{does} cause the traditional leave-one-out estimator to be \emph{biased}.  This bias agrees with folklore results in machine learning that assert that leave-one-out does generally exhibit a small bias \citep{friedman2001elements}.

{
For data-driven anchors, we stress that, unlike traditional leave-one-out validation, we do \emph{not} use one fewer points when computing the anchor in \cref{alg:ssaa}; we use $h(\bfmhat)$ for all iterations.  Hence, Shrunken-SAA is \emph{not} strictly a leave-one-out procedure, motivating our qualifier ``Modified."  
}

\section{Performance Guarantees for Shrunken-SAA}
\label{sec:PerformanceGuarantees}
In this section, we show that in the limit where the number of subproblems $K$ grows, shrinking by $\alphaLOO_{h}$ is essentially best possible.
More precisely, for any $K \geq 2$ and any $0 < \delta < 1/2$, with probability at least $1-\delta$, we prove that
\begin{equation} \label{eq:StrongConvSubOpt}
\edit{{\sf Sub Opt}_{h,K}(\alphaLOO_{h})} \leq \tilde{\mathcal O}\left( \frac{ \log^\beta(1/\delta) }{ \sqrt K }  \right),
\end{equation}
where the $\tilde{\mathcal O}(\cdot)$ notation suppresses logarithmic factors in $K$, and $1 < \beta < 2$ is a constant that depends on the particular class of optimization problems under consideration. 
Imporantly, by Borel-Cantelli lemma, \cref{eq:StrongConvSubOpt}  implies 
\(
\edit{{\sf Sub Opt}_{h,K}(\alphaLOO_{h})} \rightarrow 0,
\)
almost surely as $K \rightarrow \infty$, 
\emph{even} if the expected amount of data per subproblem remains fixed.  

\edit{\Cref{eq:StrongConvSubOpt} asserts that for a \emph{given} anchor $h(\cdot)$, Shrunken-SAA achieves the best possible shrinkage amount as $K\rightarrow \infty$.  We will also prove similar bounds on ${{\sf Sub Opt}_{\mP,K}(\alphaLOO_{h},\hloo(\bfmhat))}$.  Such bounds assert that for a given class $\mP$, Shrunken-SAA with $\hloo(\cdot)$ achieves the best possible anchor and shrinkage amount \emph{simultaneously}.  }

\subsection{Overview of Proof Technique}\label{sec:ProofTechniqueOverview}
To prove performance guarantees like \cref{eq:StrongConvSubOpt}, we first bound the sub-optimality of Shrunken-SAA in terms of the maximal stochastic deviations of $\Zperf(\alpha, h)$ and $\Zloo(\alpha, h)$ from their means.  
\begin{lemma}[Bounding Sub-Optimality] \label{lem:ConditionsForOptimality}
\edit{Suppose \cref{ass:RandomData} holds.\label{Reviewer:StateAssumptions}}   
\\For a non-data-driven anchor $h(\bfmhat) = \bp_0$, 
\begin{align}\notag
\edit{{\sf Sub Opt}_{\bp_0,K}(\alphaLOO_{\bp_0})}
 &\leq  
2\underbrace{\sup_{\alpha \geq 0} \abs{ \Zperf(\alpha, \bp_0) - \Eb{ \Zperf(\alpha, \bp_0)} } }_{\substack{\text{Maximal Stochastic Deviation in $\Zperf(\cdot, \bp_0)$} } }
 +
2\underbrace{\sup_{\alpha \geq 0} \abs{ \Zloo(\alpha, \bp_0,) - \Eb{\Zloo(\alpha, \bp_0)} }}_{\substack{\text{Maximal Stochastic Deviation in $\Zloo(\cdot, \bp_0)$} } }.
\end{align}
Similarly, for a general data-driven anchor with $h(\bfmhat)\in\mP$,
\begin{align}\label{eq:suboptbounddatadriven}
\edit{{\sf Sub Opt}_{h,K}(\alphaLOO_{h})}
 &\leq  \edit{2}
\underbrace{\sup_{\substack{\alpha \geq 0 \\ \ \bq \in \edit{\mP}}} \abs{ \Zperf(\alpha, \bq) - \Eb{ \Zperf(\alpha, \bq)} } }_{\substack{\text{Maximal Stochastic Deviation in $\Zperf(\cdot,\cdot)$} } }
 +
2\underbrace{\sup_{\substack{\alpha \geq 0 \\ \ \bq \in \edit{\mP}}} \abs{ \Zloo(\alpha, \bq) - \Eb{\Zloo(\alpha, \bq)} }}_{\substack{\text{Maximal Stochastic Deviation in $\Zloo(\cdot,\cdot)$} } }.
\end{align}
\edit{Finally, for $h=\hloo$, ${\sf Sub Opt}_{\mP,K}(\alphaLOO_{\hloo},\hloo(\bfmhat))$ is also bounded by the right-hand side of \cref{eq:suboptbounddatadriven}.}
\end{lemma}
\proof{Proof.}
By definition of $\alphaLOO_{\bp_0}$,
$\Zloo(\alphaOR_{\bp_0}, \bp_0) - \Zloo(\alphaLOO_{\bp_0}, \bp_0) \geq 0$.
Therefore,
\begin{align*}
\ts\edit{{\sf Sub Opt}_{\bp_0,K}(\alphaLOO_{\bp_0})}
&\ts\leq
\Zperf(\alphaLOO_{\bp_0}, \bp_0)
- \Zperf(\alphaOR_{\bp_0}, \bp_0)
+
\Zloo(\alphaOR_{\bp_0}, \bp_0)-\Zloo(\alphaLOO_{\bp_0}, \bp_0)
\\ 
&\ts\leq 
    2\sup_{\alpha\geq0}\abs{ \Zperf(\alpha, \bp_0) -\Zloo(\alpha, \bp_0) }
\\ 
&\ts\leq 
2\sup_{\alpha\geq0}\abs{ \Zperf(\alpha, \bp_0) - \E\Zperf(\alpha, \bp_0) }
+2\sup_{\alpha\geq0}\abs{ \Zloo(\alpha, \bp_0) -\E\Zloo(\alpha, \bp_0) }
\\&\ts\phantom{\leq}+2\sup_{\alpha\geq0}\abs{ \E\Zperf(\alpha, \bp_0) -\E\Zloo(\alpha, \bp_0) }.
\end{align*}
By \Cref{cor:steinchen}, the last term is zero\edit{, which establishes the first statement.  The proof of the second statement is similar, but in the second inequality, we take an additional supremum over $\bq\in\edit{\mP}$ to replace $h(\bfmhat)$. The proof of the third statement is similar, using $\Zloo(\alphaLOO_{\hloo}, \hloo(\bfmhat)) \leq \Zloo(\alphaOR_\mP, \bq^\OR_\mP)$, and taking a supremum over  $\alpha\geq0$, $\bq\in\edit{\mP}$ in the second inequality.}
\endproof

Proving a performance guarantee for $\alphaLOO_{h}$ thus reduces to bounding the maximal deviations in the lemma.  
Recall $\Zperf(\alpha, \bq) = \frac{1}{K} \sum_{k=1}^K Z_k(\alpha, \bq)$ and $\Zloo(\alpha, \bq) = \frac{1}{K} \sum_{k=1}^K Z^{\sf LOO}_{k}(\alpha,\bq)$.  Both processes have a special form: they are the empirical average of $K$ independent stochastic processes (indexed by $k$).
Fortunately, there exist standard tools to bound the maximal deviations of such empirical processes that rely on bounding their metric entropy.  

To keep our paper self-contained, we summarize one such approach presented in \cite{pollard1990empirical}, specifically in Eq.~(7.5) of that work.  
Recall, for any set $S\subseteq\R d$, the $\epsilon$-packing number of $S$, denoted by \edit{$D(\epsilon,S)$}, is the largest number of elements of $S$ that can be chosen so that the Euclidean distance between any two is at least $\epsilon$.  Intuitively, packing numbers describe the size of $S$ at scale $\epsilon$.  

\begin{theorem}[A Maximal Inequality; \citealp{pollard1990empirical}] \label{thm:pollard}
Let $\mathbf W(t)=(W_1(t),\dots,W_K(t))\in\R K$ be a stochastic process indexed by $t\in\mathcal T$ and let $\overline W_K(t)=\frac1K\sum_{k=1}^KW_k(t)$.  Let $\mathbf F\in\R K_+$ be a random variable such that $\abs{W_k(t)}\leq F_k$ for all $t\in\mathcal T,\,k=1,\dots,K$. 
Finally, define the random variable 
\begin{equation}\label{eq:JDudley}
J
 \ \equiv \ 
 J\left(\{\mathbf W(t):t\in\mathcal T\} , \bm F \right)  
 \ \equiv \ 
9 \| \mathbf F\|_2 \int_0^{1} \sqrt{ \log D  \left( \| \mathbf F \|_2u, \  \big\{\mathbf W(t):t\in\mathcal T \big\} \right) } du.
\end{equation}
Then, for any $p \geq1$ and any $0 < \delta < 1$, with probability at least $1-\delta$,\footnote{
Strictly speaking, eq. (7.5) of \cite{pollard1990empirical} shows that 
\(
\E\left[\left| \sup_{t\in\mathcal T } \abs{ \overline W_K(t) - \E[\overline W_K(t)] } \right|^p\right]
 \leq 2^p C^p_p \E\left[J^p\right]K^{-p},
\)
for some constant $C_p$ that relates the $\ell_p$ norm of a random variable and a particular Orlicz norm. In \Cref{lem:RelatingOrlicz4}, we prove that it suffices to take $C_p = 5^{1/p} \sqrt{ \frac{p}{2e}}$.  The result then follows from Markov's Inequality.
} 
\[
\ts\sup_{t\in\mathcal T} \abs{\overline W_K(t) - \E[\overline W_K(t)] } 
\ \leq  \ 
5^{1/p} \sqrt{p}     
{\| J \|_p}{K^{-1} \delta^{-1/p}}.
\]
\end{theorem}
If $\mathcal T$ is finite, one can bound the maximal deviation with a union bound.  
\Cref{thm:pollard} extends beyond this simple case to cases where $\abs{\mathcal T} = \infty$.
The random variable $\bm F$ in the theorem is called an \emph{envelope} for the process $\mathbf W(t)$.  
The random variable $J$ is often called the \emph{Dudley integral}.
While packing numbers describe the size of a set at scale $\epsilon$, the Dudley integral roughly describes the size of the set at varying scales. 
We again refer the reader to \citet{pollard1990empirical} for discussion.

{\blockedit \label{proofStrategy}
Our overall proof strategy is to use \cref{thm:pollard} to bound the two suprema in \cref{lem:ConditionsForOptimality}, and thus obtain a bound on the sub-optimality. 
Specifically, define the following stochastic processes:
\[
\bZperf(\alpha,\bq)=(Z_1(\alpha,\bq),\dots,Z_K(\alpha,\bq)),
\quad  \quad 
\bZloo(\alpha,\bq)=(Z^{\sf LOO}_1(\alpha,\bq),\dots,Z^{\sf LOO}_K(\alpha,\bq)).
\]
Our proof strategy will be to 1) Compute envelopes for both processes 2) Compute the packing numbers and Dudley integrals for the relevant sets above 3) Apply \cref{thm:pollard} to bound the relevant maximal deviations and 4) Use these bounds in \cref{lem:ConditionsForOptimality} to bound the sub-optimality.  We execute this strategy for several special cases in the remainder of the section.  
}

As a first step, we identify envelopes for each process.  We restrict attention to the case where the optimal value of each subproblem is bounded for any choice of anchor and shrinkage.
\begin{assumption}[Bounded Optimal Values]
\label{asn:bounded}
 There exists $\Cmax$ such that for all $i= 1, \ldots, d$, and $k=1\, \ldots, K$,
$\sup_{\bq\in\Delta_d} \abs{c_{ki}(\bx_k(\infty,\bq))}\leq\Cmax$.
\end{assumption}

Notice that $\sup_{\alpha \geq 0, \ \bq \in \Delta_d} \abs{ c_{ki}(\bx_k(\alpha, \bq ) )  } = \sup_{\bq \in \Delta_d} \abs{c_{ki}(\bx_k(\infty, \bq)}$, so that the assumption bounds the optimal value associated to every policy.  Assumption 4.1 is a mild assumption, and follows for example if $c_{ki}( \cdot)$ is continuous and $\mathcal X_k$ is compact. However, the assumption also holds, e.g, if $c_{ki}(\cdot)$ is unbounded but coercive.
 With it, we can easily compute envelopes.  Recall, $\Nhat_{\max} \equiv \max_k \Nhat_k$.  
\begin{lemma}[Envelopes for $\bZperf,\bZloo$]
\label{lem:Envelopes}
Under \cref{asn:bounded}, 
\begin{enumerate}
    \item The vector $\mathbf F^{\sf Perf} \equiv \Cmax\lambdavec/\lambdabar$ is an envelope for $\mathbf Z(\alpha, \bq)$ with 
    \(
    \| \mathbf F^{\sf Perf} \|_2  =   \frac{\Cmax }{\lambdabar} \| \lambdavec \|_2 
    .\)
    \item The random vector $\mathbf F^{\sf LOO}  = \Cmax \frac{\bm{\Nhat}}{N \lambdabar}$ is an envelope for $\mathbf Z^{\sf LOO}(\alpha, \bq)$ with 
    \(
    \| \mathbf F^{\sf LOO} \|_2 =   \ \frac{\Cmax}{N \lambdabar} \| \bm{\Nhat} \|_2 
    .\)
\end{enumerate}
\end{lemma}
\edit{The proof is immediate from the definitions and omitted.\label{line:OmittedProof}}

\edit{Our next step is to }bound the packing numbers (and Dudley integrals) for the sets
\(
\left\{ \bZperf(\alpha,\bp_0) : \alpha \geq 0 \right\} \subseteq \mathbb R^K, 
\)
and 
\(
\left\{ \bZloo(\alpha,\bp_0) : \alpha \geq 0 \right\} \subseteq \mathbb R^K, 
\)
for the case of fixed anchors and the sets
\(
\left\{ \bZperf(\alpha,\bq) : \alpha \geq 0, \ \bq \in \edit{\mP}  \right\} \subseteq \mathbb R^K,
\) 
and 
 \(
\left\{ \bZloo(\alpha,\bq) : \alpha \geq 0, \ \bq \in \edit{\mP} \right\} \subseteq \mathbb R^K, 
\)
for the case of data-driven anchors.  
Bounding these packing numbers is subtle and requires exploiting the specific structure of the optimization problem \eqref{eq:TargetProblem}.   
We \emph{separately} consider two general classes of optimization problems -- strongly convex optimization problems and discrete optimization problems -- in the remainder.  Although we focus on these classes, we expect a similar proof strategy and technique might be employed to attack other classes of optimization problems.

{ \blockedit 
\begin{remark}[Performance of $\alphaLOO$ in the Large-Sample Regime]  \label{rem:LargeSample}
Although we focus on performance guarantees for $\alphaLOO$ in settings where $K$ is large and the expected amount of data per problem is fixed, one could also ask how $\alphaLOO$ performs in the large-sample regime, i.e., where $K$ is fixed and $\Nhat_k \rightarrow \infty$ for all $k$.  Using similar techniques, i.e., reducing the problem to bounding a certain maximal stochastic deviation, one can show that $\bx_k(\alphaLOO, \bp_0, \bfmhat)$ performs comparably to the full-information solution in Problem~\eqref{eq:TargetProblem} in this limit.  The proof uses somewhat standard arguments for empirical processes.  Moreover, the result is perhaps unsurprising; many data-driven methods converge to full-information performance in the large-sample regime (see, e.g., \cite{kleywegt2002sample} for the case of SAA) since $\bphat_k$ is consistent for $\bp_k$ for all $k$ in this regime. 
Consequently, we focus on the small-data, large-scale regime, where Shrunken SAA enjoys strong suboptimality guarantees not enjoyed by SAA. 
\editt{This small-data, large-scale focus, however, causes the $N$ dependence in our bounds to be looser than that obtained from a direct large-sample analysis.  Developing a unified analysis of data-pooling for \emph{any} sequence of $N,K$ remains an open question.}
\label{Reviewer:NDependence}
\Halmos
\end{remark}
}

\subsection{Fixed Anchors and Strongly-Convex Optimization Problems}
\label{sec:SmoothCostsFixedAnchor}

In this section, we treat the case where the $K$ subproblems are smooth enough so that $\bx_k(\alpha, \bq,  \bfmhat_k)$ is smooth in $\alpha$ and $\bq$ for each $k$.  
Specifically, in this section we assume:
\begin{assumption}[\textbf{Lipschitz, Strongly-Convex Optimization}]\label{asn:smooth}
There exists $L,\gamma$ such that 
$c_{ki}(\bx)$ are $\gamma$-strongly convex and $L$-Lipschitz over $\X_k$, and, moreover, $\X_k$ is non-empty and convex, for all $k=1,\dots,K$, and $i=1,\dots,d$.
\end{assumption}
{\blockedit  \begin{theorem}[Shrunken-SAA with Fixed Anchors for Strongly Convex Problems]  \label{thm:FixedPointShrinkage} 
Fix any $\bp_0$.  \edit{\label{Reviewer:StateAssumptions2} Suppose \cref{ass:RandomData,asn:bounded,asn:smooth} hold}, $K \geq 2$ and $N\lambdamin \geq 1$.
Then, there exists a universal constant $\const$ such that for any $0 < \delta < 1/2$, with probability at least $1-\delta$, we have that
 \[
	\edit{{\sf Sub Opt}_{\bp_0,K}(\alphaLOO_{\bp_0})}   \ \leq  \ 
	\const \cdot
	\max\left(\Cmax, L \sqrt{\frac{\Cmax}{\gamma}} \right) \cdot 
	\left( \frac{\lambdamax}{\lambdamin}\right)^{5/4} \cdot
	\frac{ \log^{2} (1/\delta)\cdot  \log^{3/2}(K)}{ \sqrt K }.
\]
\end{theorem} }
 \enlargethispage*{0pt}
The proof follows our strategy from \cref{sec:ProofTechniqueOverview}. (See \cref{sec:FixedPointShrinkageAppendix}.) We sketch the main ideas:

We first bound the packing numbers of $\Fperf$ and $\Floo$. 
The key observation is that since the subproblems are strongly-convex, the optimal solutions $\bx_k(\alpha, \bp_0, \bfmhat_k)$ are continuous as functions of $\alpha$.  We utilize this continuity to construct a packing.  

\begin{figure}[t!]%
\centering%
\begin{minipage}[m]{0.45\textwidth}\includegraphics[width=\textwidth]{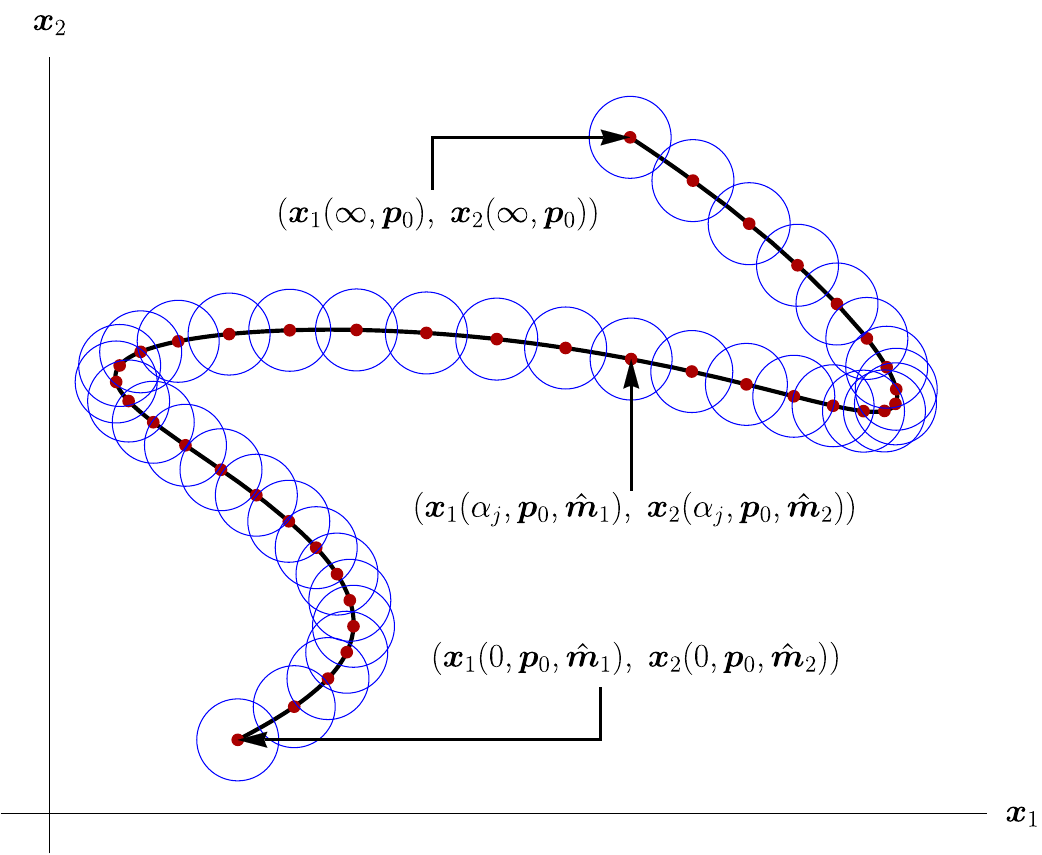}\end{minipage}\hfill%
\begin{minipage}[m]{0.525\textwidth}\captionof{figure}{\textbf{Covering a continuous process.} 
The set $\{(\bx_k(\alpha, \bp_0, \bfmhat_k))_{k=1}^K:\alpha\geq0\}$ can be thought of as a parametric curve indexed by $\alpha$ in the space $\prod_{k=1}^K\X_k$. 
Because of the continuity in $\alpha$ (cf. \cref{lem:ContinuityAlpha}, part~\ref{part:contAlpha}), to cover this curve for any compact set $\alpha\in[0,\alpha_{\max}]$ requires \edit{$\mathcal O(1/\epsilon)$} balls of size $\epsilon$. Because of the continuity at $\alpha=\infty$ (cf. \cref{lem:ContinuityAlpha}, part~\ref{part:contInfinity}), it suffices to take \edit{$\alpha_{\max}=\mathcal O(1/\epsilon)$}. This yields a packing number bound of  \edit{$\mathcal O(1/\epsilon^2)$} (cf. \cref{lem:PackSmooth}).
}\label{fig:ContinuousPackingFixed}\end{minipage}%
\end{figure}

Specifically, consider $\Fperf$.  
Continuity in $\alpha$ implies that by evaluating $\bx(\alpha, \bp_0, \bfmhat)$ on a sufficiently dense grid of $\alpha$'s, we can construct a covering of $ \left\{ \left( \bx_k(\alpha, \bp_0, \bfmhat_k) \right)_{k=1}^K : \alpha \geq 0 \right\}$, which in turn yields a covering of $\Fperf$.  
By carefully choosing the initial grid of $\alpha$'s, we can ensure that this last covering is a valid $(\epsilon/2)$-covering.  
By \cite[pg. 10]{pollard1990empirical}, the size of this covering bounds the $\epsilon$-packing number as desired.  
\Cref{fig:ContinuousPackingFixed} illustrates this intuition and further argues the initial grid of $\alpha$'s should be of size \edit{$\mathcal O(1/\epsilon^2)$}.  
A similar argument holds for $D(\epsilon, \Floo)$, using a grid of $\alpha$'s to cover $ \left\{ \big( x_k(\alpha, \bp_0, \bfmhat_k - \be_i)  : i = 1, \ldots, d, \ k = 1, \ldots, K \big) : \alpha \geq 0 \right\}$. \edit{The packing is also of size $\mathcal O(1/\epsilon^2)$.}

To complete the proof, we use these packing numbers in \cref{thm:pollard} to bound the maximal deviations of $\Zperf(\cdot,\bp_0),\Zloo(\cdot,\bp_0)$. Substituting into \cref{lem:ConditionsForOptimality} proves \cref{thm:FixedPointShrinkage} above.  Again, please see \cref{sec:FixedPointShrinkageAppendix} for details.
%
%

%

%

%
%
%
%

%
%
%
%
%
%
%
%
%
%
%
%
%
%
%
%
%
%
%
%
%
%
%
%
%
%

%
%
%
%
%
%
%
%
%
%
%
%
%
%
%
%
%
%
%
%
%
%
%
%
%
%
%
%
%
%
%
%
%
%
%
%
%
%
%
%
%
%
%
%

%
%
%
%
%
%
%
%
%
%
%
%
%

%
%
%
%
%
%
%
%
%
%
%
%
%
%
%
%

%
%
%
%
%
%

%
%

%
%
%
%
%
%
%
%
%

%
%
%
%
%
%
%
%
%
%
%
%
%
%
%
%

%
%
%
%
%

\subsection{Data-Driven Anchors and Strongly Convex Problems}  \label{sec:SmoothCostsGeneral}

{\blockedit 
We next consider the case of a data-driven anchor $h(\bfmhat) \in \mathcal P$.  Our performance guarantees will depend on the complexity of $\mathcal P$ as measured by the size of its $\ell_1$-packing numbers.  Namely, we let $D_1(\epsilon, \mathcal P)$ be the largest number of elements of $\mathcal P$ that can be chosen so that the $\ell_1$-distance between any two is at least $\epsilon$.\footnote{Recall $D(\epsilon, S)$ is defined with respect to $\ell_2$-distance.}  Then,
\begin{theorem} \label{thm:SmoothThmGeneralAnchor} {\rm \textbf{(Shrunken-SAA with Data-Driven Anchors for Strongly Convex Problems)}}
    \edit{\label{Reviewer:StateAssumptions3} Suppose \cref{ass:RandomData,asn:bounded,asn:smooth} hold}, $K \geq 2$.   \editt{ Let $d_0 \geq 1$ be such that for any $0 < \epsilon < 1/2$, 
	$\log D_1(\epsilon,  \mathcal P) \leq d_0 \log(1/\epsilon)$.}
    Then, there exists a universal constant $\const$ such that for any $0<\delta<1/2$, with probability at least $1-\delta$, we have that 
    \[
    \edit{{\sf Sub Opt}_{h,K}(\alphaLOO_{h})}
      \ \leq \
	\edit{    \const  \cdot 
	    \max\left( \Cmax, \ \frac{L^2}{\gamma}  \ + \  L\sqrt{\frac{\Cmax}{\gamma}} \right) \left( \frac{\lambdamax}{\lambdamin} \right)^{5/4}
	    \frac{d_0^2 \log^{7/2}(K)  \log^2(1/\delta) }{\sqrt K}
	    }.
    \]
\end{theorem}

In the special case of $\hloo(\cdot)$, we can prove an even stronger result, i.e., that Shrunken-SAA with $\hloo$ performs comparably to pooling in an optimal way to the best anchor within the class $\mP$.
\begin{theorem}[Shrunken-SAA with $\hloo$ for Strongly Convex Problems]
    \label{thm:SmoothThmLooAnchor}
    \editt{Under the assumptions of \cref{thm:SmoothThmGeneralAnchor},}
    there exists a universal constant $\const$ such that for any $0<\delta<1/2$, with probability at least $1-\delta$, we have that 
    \[
    {\sf Sub Opt}_{\mP,K}(\alphaLOO_{\hloo},\hloo(\bfmhat))
      \ \leq \
    {    \const  \cdot 
        \max\left( \Cmax, \ \frac{L^2}{\gamma}  \ + \  L\sqrt{\frac{\Cmax}{\gamma}} \right) \left( \frac{\lambdamax}{\lambdamin} \right)^{5/4}
        \frac{d_0^2 \log^{7/2}(K)  \log^2(1/\delta) }{\sqrt K}
        }.
    \]
\end{theorem}

In both theorems, the constant $d_0$ measures the {complexity} of $\mP$. 
Without loss of generality, $d_0 \leq 3d$ since $\mP\subseteq\Delta_d$ and $\log D_1(\epsilon,  \Delta_d) \leq 3d \log(1/\epsilon)$ \citep[Lemma 4.1]{pollard1990empirical}.
In practice, we might choose flexible, parametric families for $\mathcal P$ with small $d_0$ that do not scale with $d$.  An example might be when $\mP$ consists of all (truncated) Poisson distributions with mean at most $\Lambda$, in which case one can take $d_0 = 2\max(1,\log (\Lambda))$, \emph{independently} of $d$ (and the truncation). Another example is given in \Cref{sec:Numerics} using Beta-distributions.  
In general, we expect that our performance bounds must depend on the complexity of $\mP$ in some way, because we impose no assumptions on the function $h(\bfmhat)$ that selects the anchor, and, hence, must control behavior across all of $\mP$.  
}

\edit{Both proofs follow the strategy of \cref{sec:ProofTechniqueOverview} (see \cref{sec:SmoothCostsGeneralAppendix}).  The key idea to bounding the packing numbers is again to leverage continuity and cover the set $\{(\alpha, \bq) : \alpha \geq 0, \bq \in \mP \}$.}   \editt{Since both proofs leverage \cref{lem:ConditionsForOptimality}, the right hand sides of the bounds are the same.} \label{Reviewer:DifferentRHS}

\edit{\editt{By contrast,} the left-hand sides of \cref{thm:SmoothThmGeneralAnchor,thm:SmoothThmLooAnchor} are different: the first measures suboptimality relative to an oracle with a pre-specified anchor, while the second is relative to an oracle that can optimize the choice of anchor.  This distinction mirrors the difference between ``estimate-then-optimize" procedures and those which choose parameters in an optimization-aware fashion.   Continuing our example where $\mP$ is a set of Poisson distributions, \cref{thm:SmoothThmGeneralAnchor} bounds the suboptimality of Shrunken-SAA when  using (all) the data to fit a Poisson distribution without regard to the downstream optimization, e.g., by maximum likelihood, and then choosing $\alpha$ and $\bx_k(\cdot)$ to optimize.  By contrast, \cref{thm:SmoothThmLooAnchor} bounds the performance of Shrunken-SAA when choosing the anchor, $\alpha$ and $\bx_k(\cdot)$ simultaneously to optimize the downstream optimization.}

\subsection{Fixed Anchors and Discrete Optimization Problems}
\label{sec:DisceteFixedAnchor} 

In this section we consider the case where the $K$ subproblems are discrete optimization problems. Specifically, we require $\abs{ \X_k} < \infty$ for each $k=1,\dots,K$.
\edit{This encompasses, e.g., binary linear or non-linear optimization and linear optimization over a polytope, since we may restrict to its vertices.\label{line:linearprograms}}

Unlike the case of strongly convex problems, the optimization defining $\bx_k(\alpha, \bp_0, \bfmhat_k)$ (cf. \cref{eq:ShrunkXSol}) may admit multiple optima, and hence, $\bx_k(\alpha, \bp_0, \bfmhat_k)$ requires a tie-breaking rule. For our results below, we assume this tie-breaking rule is consistent in the sense that 
if the set of minimizers to \cref{eq:ShrunkXSol} is the same for two distinct values of $(\alpha, \bp_0)$, then the tie-breaking minimizer is also the same for both.  We express this requirement by representing the tie-breaking rule as a function from a set of minimizers to a chosen minimizer: 
\begin{assumption}[\textbf{Consistent Tie-Breaking}]\label{asn:tiebreak}
For each $k$, there exists $\sigma_k:2^{\X_k}\to\X_k$ such that
\[\ts
\bx_k(\alpha, \bp_0, \bfmhat_k) = \sigma_k\prns{\arg\min_{ \bx_k \in \X_k} ~ \bphat_k(\alpha)^\top \bc_k( \bx_k)}.
\]
\end{assumption}
Then,
\begin{theorem}[Shrunken-SAA with Fixed Anchors for Discrete Problems]
\label{thm:FixedPointShrinkageDiscrete}
Suppose that $\abs{ \X_k } < \infty$ for each $k$\edit{, $K\geq2$,} and \edit{\label{Reviewer:StateAssumptions4} that \cref{ass:RandomData,asn:bounded,asn:tiebreak} hold}.
Then, there exists a universal constant $\const$ such that \edit{for any $0<\delta<1/2$ we have that,} with probability at least $1-\delta$, 
\[ \edit{
\edit{{\sf Sub Opt}_{\bp_0,K}(\alphaLOO_{\bp_0})}
\ \leq \ \const \cdot \Cmax  \frac{\lambdamax}{\lambdamin} \cdot 
\sqrt{ \log\left( 2 N_{\max} \sum_{k=1}^K \abs{\X_k} \right)} \cdot 
	 \frac{ \log^{3/2}(K) \cdot \log^{3/2}(1/\delta)}{\sqrt K}.
}
\]
\end{theorem}
We stress that $\abs{\X_k}$ occurs logarithmically in the bound, so that the bound is reasonably tight even when the number of feasible solutions per subproblem may be large.  For example, consider binary optimization.  Then, $\abs{\X_k}$ often scales exponentially in the number of binary variables, so that $\log(\abs{\X_k})$ scales like the number of binary variables.  Thus, as long as the number of binary variables per subproblem is much smaller than $K$, the sub-optimality will be small with high probability. 

We also note that, unlike \cref{thm:FixedPointShrinkage}, \edit{the above bound depends on $\log(N_{\max})$}.  This mild dependence stems from the fact that we have made \emph{no assumptions of continuity} on the functions $\bc_{k}(\bx, \bxi)$ in $\bx$ or $\bxi$.  Since these functions could be arbitrarily non-smooth, \edit{we need to control their behavior separately across all of the LOO iterations, which introduces the $N_{\max}$ dependence.  With stronger assumptions, it might be possible to remove this dependence.  
However, since we are mostly interested in the setting where $N_k$ is moderate to small for all $k$, we do not pursue this idea. \label{NDependence}}

To prove \cref{thm:FixedPointShrinkageDiscrete}, we again follow the approach outlined in \cref{sec:ProofTechniqueOverview}. 
Since the policy $\bx(\alpha, \bp_0, \bfmhat)$ need not be smooth in $\alpha$, however, we adopt a different strategy than in \cref{sec:SmoothCostsFixedAnchor}.  Specifically, we bound the cardinality of $\Fperf$, $\Floo$, directly.  (Recall that the cardinality of a set bounds its $\epsilon$-packing number for any $\epsilon$.)  

\begin{figure}[t!]%
\centering%
\begin{minipage}[m]{0.55\textwidth}\includegraphics[width=\textwidth]{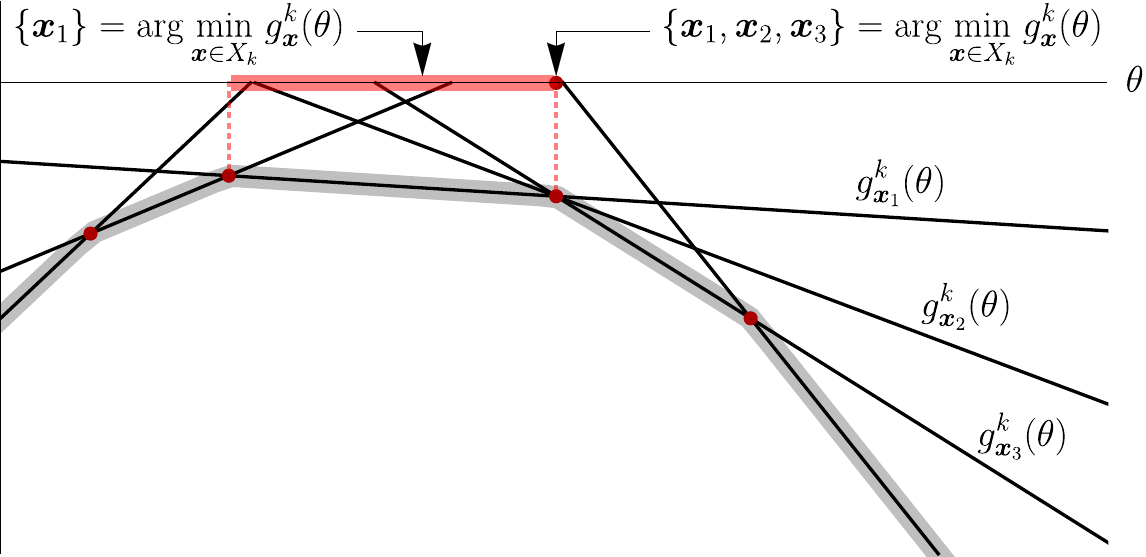}\end{minipage}\hfill%
\begin{minipage}[m]{0.425\textwidth}\captionof{figure}{\textbf{Counting Discrete Solutions.}  A \edit{concave} piecewise-linear function consisting of $\abs{\X_k}$ lines has at most $\abs{\X_k}-1$ breakpoints, between which the set of active supporting lines is constant. Any function of this set of active supporting lines is piecewise constant with at most $\abs{\X_k}-1$ discontinuities.}\label{fig:DiscretePackingFixed}\end{minipage}%
\end{figure}
First note  the cardinality of $\Fperf$ is at most that of $\braces{\prns{\bx_k(\alpha,\bp_0,\bfmhat_k)}_{k=1}^K:\alpha\geq0}$. A trivial bound on this latter set's cardinality is $\prod_{k=1}^K \abs{\X_k}$. This bound is too crude for our purposes; it grows exponentially in $K$ even if $\abs{\X_k}$ is bounded for all $k$. Intuitively, this bound is crude because it supposes we can vary each solution $\bx_k(\alpha, \bp_0, \bfmhat_k)$ independently of the others to achieve all $\prod_{k=1}^K \abs{\X_k}$ possible combinations.  In reality, we can only vary a single parameter, $\alpha$, that simultaneously controls all $K$ solutions, rather than varying them separately.  We use this intuition to show that a much smaller bound, i.e., $2\sum_{k=1}^K \abs{\X_k}$, is valid.

To this end, we fix $k$ and study the dependence of $\bx_k(\alpha,\bp_0,\bfmhat_k)$ on $\alpha$. In the trivial case $\Nhat_k=0$, $\bx_k(\alpha,\bp_0,\bfmhat_k)$ takes only one value: $\bx_k(\infty,\bp_0)$.  Hence we focus on the case 
$\Nhat_k\geq1$. 

Consider reparameterizing the solution in terms of $\theta = \frac{\alpha}{\alpha + \Nhat_k}\in[0,1)$ and let $\alpha(\theta)=\frac{\theta}{1-\theta}\Nhat_k$.  Then for any $\bx \in \X_k$, define the linear function 
\begin{align*}\ts
g_{k \bx}(\theta) &= ( (1-\theta) \bphat_k + \theta \bp^0 )^\top \bc_k(\bx), \quad \ \theta \in [0, 1).
\end{align*}
Since $g_{k \bx}(\cdot)$ is linear, the function $\theta \mapsto \min_{\bx\in\X_k}g_{k \bx}(\theta)$ is \edit{concave},\label{line:ConvexityTypo} piecewise-linear with at most $\abs{\X_k}-1$ breakpoints.  By construction, $\bx_k(\alpha(\theta),\bp_0,\bfmhat_k)\in\argmin_{\bx_k\in\X_k}g_{k \bx}(\theta)$.  
More precisely, for any $\theta$, the set of active supporting hyperplanes of $\min_{\bx\in\X_k}g_{k \bx}(\cdot)$ at $\theta$ is $\left\{ (\bp^0 - \bphat_k)^\top \bc_k(\bx)  \ : \  \bx \in \argmin_{\bx_k\in\X_k}g_{k \bx}(\theta) \right\}$.  

Since the set of active supporting hyperplanes is constant between breakpoints, the set of minimizers $\argmin_{\bx_k\in\X_k}g_{k \bx}(\theta)$ is also constant between breakpoints.  By \cref{asn:tiebreak}, this implies $\theta\mapsto \bx_k(\alpha(\theta),\bp_0,\bfmhat_k)$ is piecewise constant with at most $\abs{\X_k}-1$ points of discontinuity.  (See also \cref{fig:DiscretePackingFixed}.)  Viewed in the original parameterization in terms of $\alpha$, it follows that $\alpha \mapsto \bx_k(\alpha, \bp_0, \bfmhat_k)$ is also piecewise constant with at most $\abs{\X_k}-1$ points of discontinuity.  
Thus,
\begin{lemma} \label{thm:XPiecewise}
Suppose \cref{asn:tiebreak} holds.
Fix any $\bp_0$ and $\bfmhat_k$.
Then, the function \break$\alpha \mapsto \bx_k(\alpha, \bp_0, \bfmhat_k)$ is piecewise constant with at most $\abs{\X_k}-1$ points of discontinuity.  
\end{lemma}
Taking the union of all these points of discontinuity over $k$ proves that $\prns{\bx_k(\alpha, \bp_0, \bfmhat_k)}_{k=1}^K$ is also piecewise constant with at most $\sum_{k=1}^K(\abs{\X_k}-1)$ points of discontinuity. {Therefore, it takes at most $2 \sum_{k=1}^K \abs{\X_k} - 2K + 1$ different values --  a distinct value for each of the $\sum_{k=1}^K(\abs{\X_k}-1)$ breakpoints plus a distinct value for the $\sum_{k=1}^K(\abs{\X_k}-1) + 1$ regions between breakpoints.} This gives the desired cardinality bound on $\abs{ \Fperf  }$. A similar argument considering the larger \edit{$\prns{\bx_k(\alpha, \bp_0, \bfmhat_k-\be_i)}_{i\in\mathcal I_k, k=1,\dots,K}$, where $\mathcal I_k=\{i=1,\dots,d:\mhat_{ki}>0\}$,} gives a corresponding cardinality bound on $\abs{ \Floo }$. \edit{Noting $\abs{\mathcal I_k}\leq\min(d,\Nhat_k)$ 
gives the following (proof omitted): \enlargethispage*{20pt}}
\begin{corollary} [Size of Discrete Solutions Sets]
\label{thm:SizeOfDiscreteSets}
Suppose \cref{asn:tiebreak} holds.
 Then, 
\[\ts
\abs{ \Fperf  } \ \leq \ 2 \sum_{k=1}^K \abs{\X_k},
\quad
\edit{\abs{ \Floo }  \   \leq \  1 + 2 \sum_{k=1}^K \min(d,\Nhat_k)\abs{\X_k}}.
\]
\end{corollary}
\edit{The additional ``1" in the case of $\abs{\Floo}$ covers the case where $\Nhat_{\max} = 0$ and $\Floo = \{ \bm 0 \}$.}  
Although these bounds may appear large, an important feature is that they are only linear in $K$ as long as $\abs{\X_k}$ are bounded over $k$.

We use these cardinality bounds to bound the packing numbers and then apply our usual strategy via \cref{thm:pollard} and \cref{lem:ConditionsForOptimality} to prove \cref{thm:FixedPointShrinkageDiscrete}.   The details are in \cref{sec:FixedPointShrinkageDiscreteAppendix}.


\subsection{Data-Driven Anchors and Discrete Optimization Problems}  \label{sec:DiscreteGeneral}
We next extend the results of \cref{sec:DisceteFixedAnchor} to the case of a data-driven anchor, $h(\bfmhat)$.  \edit{As in \cref{sec:SmoothCostsGeneral}, our bounds will depend on a measure of complexity of $\mP$, namely, the dimension of 
$\spn(\mP) \equiv \{ \sum_{\ell = 1}^d \theta_\ell \bq_\ell \ : \theta_\ell \in \Rl, \ \bq_\ell \in \mP,\ \ell=1,\dots,d \}$ when viewed as a linear subspace.  Denote this dimension by $d_0$ and note $1\leq d_0 \leq d$.  A canonical example might be when $\mP$ consists of mixture distributions with $d_0$ (specified) components.}
We prove that:
\begin{theorem}[Shrunken-SAA with Data-Driven Anchors for Discrete Problems] 
\label{thm:ShrinkageDiscreteGeneralAnchor}
Suppose that $\abs{ \X_k } < \infty$ for each $k$, \editt{that $\spn(\mP)$ has dimension $d_0$,} and that \cref{asn:bounded,asn:tiebreak} hold.
Then, there exists a universal constant $\const$ such that \edit{for all $0<\delta<1/2$, we have that,} with probability at least $1-\delta$, 
\[
\edit{
\edit{{\sf Sub Opt}_{h,K}(\alphaLOO_{h})}
 \ \leq   \ 
 \const \cdot 
\Cmax  \frac{\lambdamax}{\lambdamin}  \sqrt{ \edit{d_0} \log\left( N_{\max} \sum_{k=1}^K \abs{\X_k}  \right) }  
	    \cdot  \frac{\log^{3/2}(K)\log^{2}(1/\delta) }{\sqrt K}.
}	    
\]
\end{theorem}
\edit{
\begin{theorem}[Shrunken-SAA with $\hloo$ for Discrete Problems]
    \label{thm:ShrinkageDiscreteLooAnchor}
    \editt{Under the assumptions of \cref{thm:ShrinkageDiscreteGeneralAnchor},}
    there exists a universal constant $\const$ such that for any $0<\delta<1/2$, with probability at least $1-\delta$, we have that 
    \[
    {\sf Sub Opt}_{\mP,K}(\alphaLOO_{\hloo},\hloo(\bfmhat))
      \ \leq \
    \const \cdot 
\Cmax  \frac{\lambdamax}{\lambdamin}  \sqrt{ \edit{d_0} \log\left( N_{\max} \sum_{k=1}^K \abs{\X_k}  \right) }  
        \cdot  \frac{\log^{3/2}(K)\log^{2}(1/\delta) }{\sqrt K}
        .
    \]
\end{theorem}}
\noindent Both proofs follow the strategy from \cref{sec:ProofTechniqueOverview} (see \cref{sec:DiscreteGeneralAppendix}) and, \editt{hence, lead to the same right hand sides.  However, the left hand sides are distinct.    We sketch the main ideas of the proof:} \label{Reviewer:SecondTwoTheorems} 

\begin{figure}[t!]%
\centering%
\begin{minipage}[m]{0.5\textwidth}\includegraphics[width=\textwidth]{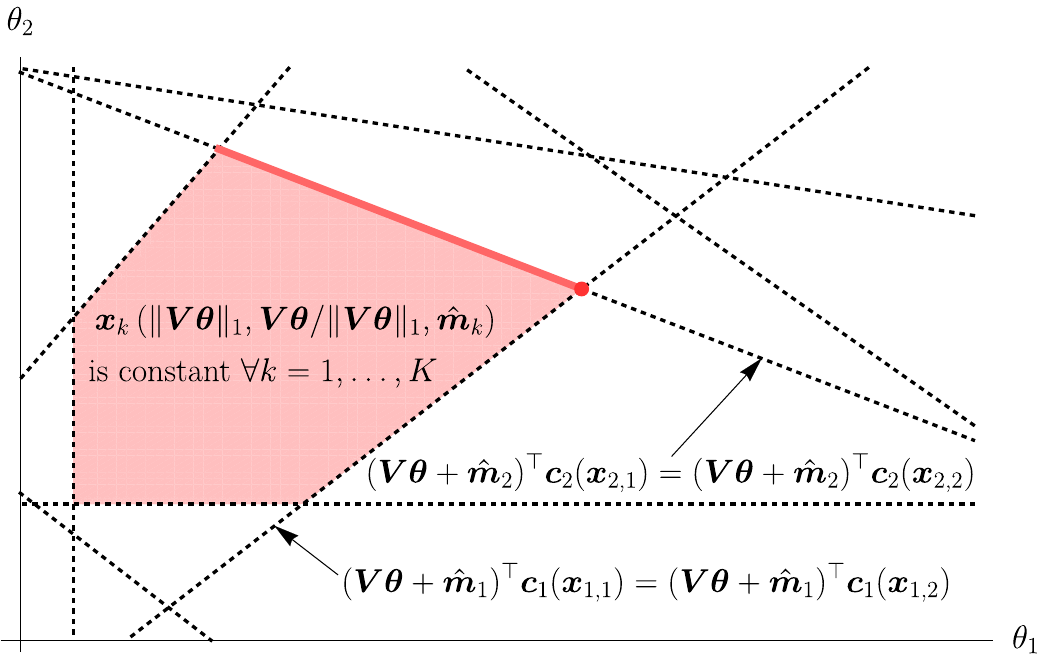}\end{minipage}\hfill%
\begin{minipage}[m]{0.45\textwidth}\captionof{figure}{\textbf{Solution Induced Hyperplane Arrangement.} 
The hyperplanes $H_{kij}$ (cf. \cref{eq:totalorderinghyperplanes}) in $\R d$ are indifference curves between solutions $\bx_{ki}$ and $\bx_{kj}$ in \cref{eq:ThetaParam}.  The total ordering on each set $\X_k$ induced by the objective of \cref{eq:ThetaParam}
is thus constant on the interior of the fully-specified polyhedra defined by the hyperplanes.}
\label{fig:DiscretePackingGeneral}
\end{minipage}%
\end{figure}

We first bound the cardinality of \edit{$\FperfqmP$}, \edit{$\FlooqmP$}.  
The key is to generalize the argument of \cref{sec:DisceteFixedAnchor} from counting breakpoints in a univariate piecewise affine function to counting the pieces in a multivariate piecewise affine function.
\edit{First, we reparameterize our policies.  
Let the columns of $\bm V \in \mathbb R^{d \times d_0}$ be a basis of $\spn(\mP)$.
Then, intrepreting $\bm 0/0$ as an arbitrary point in $\Delta_d$ (e.g., $\be/d$),} 
{\blockedit
\begin{align}  \notag
 \abs{ \FperfqmP }  
&\ \leq \ 
 \abs{ \left\{ \left(\bx_k(\alpha, \bq, \bfmhat_k) \right)_{k=1}^K  \ : \ \bq \in \mP, \alpha \geq 0 \right\}  } 
\\ \notag \
& \ \leq  \ 
 \abs{ \left\{ \left(\bx_k(\| \bm w \|_1,  \bm w / \| \bm w \|_1 , \bfmhat_k) \right)_{k=1}^K  \ : \ \bm w \in \spn(\mP)\cap\R d_+ \right\}  } 
\\ \label{eq:ySolPerf} \
 & \ = \ 
\abs{ \left\{ \left( \bx_k(\| \bm V\btheta \|_1,  \bm V\btheta / \| \bm V\btheta \|_1, \bfmhat_k) \right)_{k=1}^K \ : \ \btheta \in \R{d_0},\,\bm V\btheta\in\R d_+  \right\} }.
\end{align}}%
Hence, it suffices to bound the right most side of 
\cref{eq:ySolPerf}.
An advantage of this $\btheta$-parameterization over the original $(\alpha, \bq)$-parameterization is that, for $\Nhat_k > 0$,
\begin{equation} \label{eq:ThetaParam}
\ts\edit{\bx_k(\| \bm V\btheta \|_1,  \bm V\btheta / \| \bm V\btheta \|_1, \bfmhat_k) \in \arg\min_{\bx \in \X_k}\  ( \bm V\btheta + \bfmhat_k)^\top \bc_k(\bx)},
\end{equation}
and $\btheta$ occurs linearly in this representation. 
{\blockedit The set of $\btheta$ where we are indifferent between $\bx_{ki}, \bx_{kj} \in \X_k$ in \cref{eq:ThetaParam} is the hyperplane
\begin{equation} \label{eq:totalorderinghyperplanes}
\edit{H_{kij}=\braces{\btheta\in\R{d_0}:
\left( \bm V \btheta +\bfmhat_{k} \right)^\top \left( \bc_k( \bx_{ki}) - \bc_k( \bx_{kj}) \right)
= \bm 0
}.
}
\end{equation}
Consider drawing all $\sum_{k=1}^K \binom{\abs{ \X_k}}{2}$ such hyperplanes, as in \cref{fig:DiscretePackingGeneral}. Then, for any $\btheta \in \R{d_0}$, consider the polyhedron given by the equality constraints of those hyperplanes containing $\btheta$, and the inequality constraints defined by the side on which $\btheta$ lies for the remaining hyperplanes. 
The relative ordering of $\{\left( \bm V\btheta +\bfmhat_{k} \right)^\top \bc_k( \bx_k):  \bx_k \in \X_k \}$ is constant for all $\btheta$ in this polyhedron's interior.  Hence, 
$(\bx_k(\| \bm V\btheta \|_1,  \bm V\btheta / \| \bm V\btheta \|_1, \bfmhat_k))_{k=1}^K$ is also constant. Thus, to bound $\FperfqmP$, it suffices to count the number of such polyhedra.  We do this counting in \cref{sec:DiscreteGeneralAppendix}. 
A similar argument (with a different hyperplane arrangement) can be used to bound the cardinality of \edit{$\FlooqmP$}.} \edit{We summarize the results as:}
{\blockedit
\begin{lemma}[Size of Discrete Solutions Sets]  \label{thm:SizeOfDiscreteSetsGeneral}
Under the assumptions of \cref{thm:ShrinkageDiscreteGeneralAnchor},
\[\ts
\edit{\abs{ \FperfqmP  } \leq \prns{\sum_{k=1}^K \abs {\X_k}^2}^{d_0},~
\abs{ \FlooqmP }   \leq 1 + \Nhat_{\max}^{d_0}  \prns{\sum_{k=1}^K \abs{\X_k}^2}^{d_0}.}
\]
\end{lemma}
Importantly, both bounds 
are polynomial in $K$ if $\abs{\X_k}$ are bounded over $k$.
We then apply \cref{thm:pollard} to bound the maximal deviations in \cref{lem:ConditionsForOptimality}, proving the theorems.   Again, see \cref{sec:DiscreteGeneralAppendix} for details.  }

{\blockedit
\subsection{Performance Guarantees for Continuous Distributions}\label{remark:continuousanalysis}
Notice that none of our previous theorems (cf. \cref{thm:FixedPointShrinkage,thm:SmoothThmGeneralAnchor,thm:SmoothThmLooAnchor,thm:FixedPointShrinkageDiscrete,thm:ShrinkageDiscreteGeneralAnchor,thm:ShrinkageDiscreteLooAnchor}) depend explicitly on 
$d$, the size of the support of $\bp_k$.   Recall also that \cref{alg:ssaa} does not depend on $d$.  These observations beg the question of whether similar performance guarantees hold for Shrunken-SAA when  $\bxi_k$ are not discrete with finite support.  

For the case of strongly-convex optimization problems, the short answer is ``yes." One simply applies \cref{alg:ssaa} as written to the potentially continous $\bxi_k$, but \emph{analyzes} a discretized system where the discretization is chosen sufficiently fine that the two systems behave similarly.  The details are somewhat tedious.  See \cref{sec:ContinuousExtension} in the appendix for a formal statement and proof.

Unfortunately, for the case of discrete optimization problems, the answer is more subtle, and it is not clear that similar performance guarantees hold without additional assumptions.  Again, see \cref{sec:ContinuousExtension} for a discussion of the key issues.  
}

%
\section{The Sub-Optimality-Stability Tradeoff:  An Intuition for Data-Pooling}
\label{sec:StabilityIntuition}

\enlargethispage{10pt}

In the previous section, we established that for various classes of optimization problems, Shrunken SAA pools the data in the best possible way for a given anchor, 
\edit{or, when used with $\hloo$, pools the data in the best possible way to the best-in-class anchor}, asymptotically as $K \rightarrow \infty$.
In this section, we show how Shrunken SAA can also be used to build a strong intuition into \emph{when} and \emph{why} data-pooling improves upon decoupling.  

We focus first on the case of a non-data-driven anchor $\bp_0$ for simplicity.  \Cref{cor:steinchen} shows that (under \cref{ass:RandomData}) $\Eb{\Zperf(\alpha, \bp_0)} = \Eb{\Zloo(\alpha, \bp_0)}$.  
\Cref{thm:FixedPointShrinkage,thm:FixedPointShrinkageDiscrete} establish that under mild conditions, we often have the stronger statement
\[
\underbrace{\Zperf(\alpha, \bp_0)}_{\text{True Performance of } \alpha} 
\textstyle 
\ \ = \ \ 
\underbrace{\Zloo(\alpha, \bp_0)}_{\text{LOO Performance of } \alpha} 
\ + \ 
\underbrace{\tilde{\mathcal O}_p(1/ \sqrt K)}_{\text{Stochastic Error}}, 
\]
where the error term is uniformly small in $\alpha$.  In these two senses, optimizing $\Zperf(\alpha, \bp_0)$ over $\alpha$ is roughly equivalent to optimizing $\Zloo(\alpha, \bp_0)$ over $\alpha$, especially for large $K$.  

A simple algebraic manipulation then shows that 
\[ \textstyle
\Zloo(\alpha, \bp_0)  \ \ = \ \ \frac{1}{N \lambdabar} \big( \text{SAA-SubOpt}(\alpha) \ + \ \text{Instability}(\alpha) \ +  \ \text{SAA}(0)\big) , 
\vspace{-10pt}
\]
\begin{align*}
\text{where} \hspace{20pt}
\text{SAA-SubOpt}(\alpha)  &\textstyle\ \equiv \ 
	\frac{1}{K} \sum_{k=1}^K  \sum_{i=1}^d \mhat_{ki} \Big( c_{ki}\big(x_k(\alpha, \bp_0, \bfmhat_k) \big) - c_{ki}\big(x_k(0, \bp_0, \bfmhat_k)  \big) \Big)
\\
\text{Instability}(\alpha) &\ \textstyle \equiv \ 
\frac{1}{K} \sum_{k=1}^K  \sum_{i=1}^d \mhat_{ki}  \Big( c_{ki}\big(x_k(\alpha, \bp_0, \bfmhat_k - \be_i) \big)
- c_{ki}\big(x_k(\alpha, \bp_0, \bfmhat_k)  \big) \Big),
\\
\text{SAA}(0) &  \textstyle \ \equiv \ 
	\frac{1}{K} \sum_{k=1}^K  \sum_{i=1}^d \mhat_{ki} c_{ki}\big(x_k(0, \bp_0, \bfmhat_k)  \big).
\end{align*}
Note $\text{SAA}(0)$ does not depend on $\alpha$.  In other words, optimizing $\Zperf(\alpha, \bp_0)$ over $\alpha$ is roughly equivalent to optimizing $\Zloo(\alpha, \bp_0)$, which in turn is equivalent to optimizing 
\begin{equation}\ts\tag{\textbf{Sub-Optimality-Instability Tradeoff}}
\min_{\alpha \geq 0 } \quad \text{SAA-SubOpt}(\alpha) + \text{Instability}(\alpha)  .
\end{equation}
We term this last optimization the ``Sub-Optimality-Instability Tradeoff.''

To develop some intuition, notice $\text{SAA-SubOpt}(\alpha)$ is nonnegative, and measures the average degree to which each $\bx_k(\alpha, \bp_0, \bfmhat_k)$ is sub-optimal with respect to a (scaled) SAA objective.  In particular, $\text{SAA-SubOpt}(\alpha)$ is minimized at $\alpha = 0$, and we generally expect it is increasing in $\alpha$.  
By contrast, $\text{Instability}(\alpha)$ measures the average degree to which the (scaled) performance of $\bx_k(\alpha, \bp_0, \bfmhat_k)$ changes on the training sample if we were to use one fewer data points.  It is minimized at $\alpha = \infty$, since the fully-shrunken solution $\bx_k(\infty, \bp_0, \bfmhat_k)$ does not depend on the data and is, hence, completely stable.  Intuitively, we might expect $\text{Instability}(\alpha)$ to be decreasing since as $\alpha$ increases, the shrunken measure $\bphat_k(\alpha)$ depends less and less on the data. In reality, $\text{Instability}(\alpha)$ is often decreasing for large enough $\alpha$, but for smaller $\alpha$ can have subtle behavior depending on the optimization structure.  (See below for examples.)  

This tradeoff is intuitive in light of our data-pooling interpretation of $\bx_k(\alpha, \bp_0, \bfmhat_k)$ from \cref{sec:Bayes}.  Recall, we interpret $\bx_k(\alpha, \bp_0,  \bfmhat_k)$ as the solution when we augment our original dataset with a synthetic dataset of size $\alpha$ drawn from $\bp_0$.  As we increase $\alpha$, we introduce more SAA-sub-optimality into $\bx_k(\alpha, \bp_0, \bfmhat_k)$ because we ``pollute" the $k^\text{th}$ dataset with draws from a distinct distribution.  However, we also increase the stability of $\bx_k(\alpha, \bp_0, \bfmhat_k)$ because we reduce its dependence on $\bfmhat_k$.  Shrunken-SAA seeks an $\alpha$ in the ``sweet spot" that balances these two effects.  

Importantly, this tradeoff also illuminates \emph{when} data-pooling offers an improvement, i.e., when $\alphaLOO > 0$.  Intuitively, $\alphaLOO > 0$ only if $\text{Instability}(0)$ is fairly large and decreasing.  Indeed, in this setting, the SAA-sub-optimality incurred by choosing a small positive $\alpha$ is likely outweighed by the increased stability.  However, if $\text{Instability}(0)$ is already small, the marginal benefit of additional stability likely won't outweigh the cost of sub-optimality.  

More precisely, we intuit that data-pooling offers a benefit whenever i) the SAA solution is unstable, ii) the fully-shrunken solution $\bx_k(\infty, \bp_0, \bfmhat)$ is not too sub-optimal, and iii) $K$ is sufficiently large.  In particular, when $\Nhat_k$ is relatively small for most $k$, the SAA solution is likely to be \emph{very} unstable.  Hence, intuition suggests data-pooling likely provides a benefit whenever $\Nhat_k$ is small but $K$ is large, i.e., the small-data, large-scale regime.  

The intuition for a data-driven anchor $h(\bfmhat)$ is essentially the same. The proofs of \Cref{thm:SmoothThmGeneralAnchor,thm:ShrinkageDiscreteGeneralAnchor} show that the approximation 
$\Zperf(\alpha, \bp_0) \approx \Zloo(\alpha, \bp_0)$ holds uniformly in $\alpha$ \emph{and} $\bp_0$.  Consequently, the Sub-Optimality-Instability Tradeoff also holds for all $\bp_0$.  Hence, it holds for the specific realization of $h(\bfmhat)$, and changing $\alpha$ balances these two sources of error for this anchor.  
We recall in contrast to traditional leave-one-out validation, however, Shrunken-SAA does not remove a data point and retrain the anchor.  This detail is important because it ensure\edit{s} the fully-shrunken solution $\bx_k(\infty, h(\bfmhat), \bfmhat)$ is still completely stable per our definition, i.e., has instability equal to zero, despite depending on the data.  

The Sub-Optimality-Instability Tradeoff resembles the classical bias-variance tradeoff for MSE.  Both tradeoffs decompose performance into a systematic loss (bias or SAA-sub-optimality) and a measure of dispersion (variance or instability).  An important distinction, however, is that the Sub-Optimality-Instability tradeoff applies to general optimization problems, not just mean-squared error.  Even if we restrict to the case of MSE (cf. \cref{ex:MSE}), however, the two tradeoffs still differ and are two different ways to split the ``whole" into ``pieces."  See \cref{sec:BiasVariance}.

\subsection{Sub-Optimality-Instability Tradeoff as a Diagnostic Tool}
\begin{figure}
    \centering
    \begin{subfigure}[b]{0.3\textwidth}
	\ifdraft{Figure Goes Here}
	{
\begin{tikzpicture}[x=1pt,y=1pt]
\definecolor{fillColor}{RGB}{255,255,255}
\path[use as bounding box,fill=fillColor,fill opacity=0.00] (0,0) rectangle (144.54,144.54);
\begin{scope}
\path[clip] ( 25.23, 23.41) rectangle (140.54,140.54);
\definecolor{drawColor}{gray}{0.92}

\path[draw=drawColor,line width= 0.2pt,line join=round] ( 25.23, 44.29) --
	(140.54, 44.29);

\path[draw=drawColor,line width= 0.2pt,line join=round] ( 25.23, 75.42) --
	(140.54, 75.42);

\path[draw=drawColor,line width= 0.2pt,line join=round] ( 25.23,106.54) --
	(140.54,106.54);

\path[draw=drawColor,line width= 0.2pt,line join=round] ( 25.23,137.67) --
	(140.54,137.67);

\path[draw=drawColor,line width= 0.2pt,line join=round] ( 43.57, 23.41) --
	( 43.57,140.54);

\path[draw=drawColor,line width= 0.2pt,line join=round] ( 69.78, 23.41) --
	( 69.78,140.54);

\path[draw=drawColor,line width= 0.2pt,line join=round] ( 95.99, 23.41) --
	( 95.99,140.54);

\path[draw=drawColor,line width= 0.2pt,line join=round] (122.20, 23.41) --
	(122.20,140.54);

\path[draw=drawColor,line width= 0.4pt,line join=round] ( 25.23, 28.73) --
	(140.54, 28.73);

\path[draw=drawColor,line width= 0.4pt,line join=round] ( 25.23, 59.85) --
	(140.54, 59.85);

\path[draw=drawColor,line width= 0.4pt,line join=round] ( 25.23, 90.98) --
	(140.54, 90.98);

\path[draw=drawColor,line width= 0.4pt,line join=round] ( 25.23,122.10) --
	(140.54,122.10);

\path[draw=drawColor,line width= 0.4pt,line join=round] ( 30.47, 23.41) --
	( 30.47,140.54);

\path[draw=drawColor,line width= 0.4pt,line join=round] ( 56.68, 23.41) --
	( 56.68,140.54);

\path[draw=drawColor,line width= 0.4pt,line join=round] ( 82.88, 23.41) --
	( 82.88,140.54);

\path[draw=drawColor,line width= 0.4pt,line join=round] (109.09, 23.41) --
	(109.09,140.54);

\path[draw=drawColor,line width= 0.4pt,line join=round] (135.30, 23.41) --
	(135.30,140.54);
\definecolor{drawColor}{RGB}{248,118,109}

\path[draw=drawColor,line width= 0.6pt,line join=round] ( 30.47, 44.23) --
	( 31.53, 44.23) --
	( 32.59, 44.23) --
	( 33.65, 44.23) --
	( 34.71, 44.23) --
	( 35.77, 44.23) --
	( 36.82, 44.23) --
	( 37.88, 44.23) --
	( 38.94, 44.23) --
	( 40.00, 44.23) --
	( 41.06, 44.23) --
	( 42.12, 44.23) --
	( 43.18, 44.23) --
	( 44.24, 49.43) --
	( 45.30, 49.43) --
	( 46.35, 49.43) --
	( 47.41, 49.43) --
	( 48.47, 49.43) --
	( 49.53, 49.43) --
	( 50.59, 49.43) --
	( 51.65, 49.43) --
	( 52.71, 49.43) --
	( 53.77, 49.43) --
	( 54.82, 49.43) --
	( 55.88, 49.43) --
	( 56.94, 55.26) --
	( 58.00, 55.26) --
	( 59.06, 55.26) --
	( 60.12, 55.26) --
	( 61.18, 55.26) --
	( 62.24, 55.26) --
	( 63.30, 55.26) --
	( 64.35, 55.26) --
	( 65.41, 55.26) --
	( 66.47, 55.26) --
	( 67.53, 55.26) --
	( 68.59, 55.26) --
	( 69.65, 55.26) --
	( 70.71, 59.06) --
	( 71.77, 59.06) --
	( 72.83, 59.06) --
	( 73.88, 59.06) --
	( 74.94, 59.06) --
	( 76.00, 59.06) --
	( 77.06, 59.06) --
	( 78.12, 59.06) --
	( 79.18, 59.06) --
	( 80.24, 59.06) --
	( 81.30, 59.06) --
	( 82.36, 59.06) --
	( 83.41, 61.86) --
	( 84.47, 61.86) --
	( 85.53, 61.86) --
	( 86.59, 61.86) --
	( 87.65, 61.86) --
	( 88.71, 61.86) --
	( 89.77, 61.86) --
	( 90.83, 61.86) --
	( 91.89, 61.86) --
	( 92.94, 61.86) --
	( 94.00, 61.86) --
	( 95.06, 61.86) --
	( 96.12, 63.37) --
	( 97.18, 63.37) --
	( 98.24, 63.37) --
	( 99.30, 63.37) --
	(100.36, 63.37) --
	(101.41, 63.37) --
	(102.47, 63.37) --
	(103.53, 63.37) --
	(104.59, 63.37) --
	(105.65, 63.37) --
	(106.71, 63.37) --
	(107.77, 63.37) --
	(108.83, 63.37) --
	(109.89, 61.68) --
	(110.94, 61.68) --
	(112.00, 61.68) --
	(113.06, 61.68) --
	(114.12, 61.68) --
	(115.18, 61.68) --
	(116.24, 61.68) --
	(117.30, 61.68) --
	(118.36, 61.68) --
	(119.42, 61.68) --
	(120.47, 61.68) --
	(121.53, 61.68) --
	(122.59, 56.34) --
	(123.65, 56.34) --
	(124.71, 56.34) --
	(125.77, 56.34) --
	(126.83, 56.34) --
	(127.89, 56.34) --
	(128.95, 56.34) --
	(130.00, 56.34) --
	(131.06, 56.34) --
	(132.12, 56.34) --
	(133.18, 56.34) --
	(134.24, 56.34) --
	(135.30, 50.99);
\definecolor{drawColor}{RGB}{0,191,196}

\path[draw=drawColor,line width= 0.6pt,dash pattern=on 2pt off 2pt ,line join=round] ( 30.47, 28.73) --
	( 31.53, 28.73) --
	( 32.59, 28.73) --
	( 33.65, 28.73) --
	( 34.71, 28.73) --
	( 35.77, 28.73) --
	( 36.82, 28.73) --
	( 37.88, 28.73) --
	( 38.94, 28.73) --
	( 40.00, 28.73) --
	( 41.06, 28.73) --
	( 42.12, 28.73) --
	( 43.18, 28.73) --
	( 44.24, 30.74) --
	( 45.30, 30.74) --
	( 46.35, 30.74) --
	( 47.41, 30.74) --
	( 48.47, 30.74) --
	( 49.53, 30.74) --
	( 50.59, 30.74) --
	( 51.65, 30.74) --
	( 52.71, 30.74) --
	( 53.77, 30.74) --
	( 54.82, 30.74) --
	( 55.88, 30.74) --
	( 56.94, 35.71) --
	( 58.00, 35.71) --
	( 59.06, 35.71) --
	( 60.12, 35.71) --
	( 61.18, 35.71) --
	( 62.24, 35.71) --
	( 63.30, 35.71) --
	( 64.35, 35.71) --
	( 65.41, 35.71) --
	( 66.47, 35.71) --
	( 67.53, 35.71) --
	( 68.59, 35.71) --
	( 69.65, 35.71) --
	( 70.71, 44.96) --
	( 71.77, 44.96) --
	( 72.83, 44.96) --
	( 73.88, 44.96) --
	( 74.94, 44.96) --
	( 76.00, 44.96) --
	( 77.06, 44.96) --
	( 78.12, 44.96) --
	( 79.18, 44.96) --
	( 80.24, 44.96) --
	( 81.30, 44.96) --
	( 82.36, 44.96) --
	( 83.41, 58.12) --
	( 84.47, 58.12) --
	( 85.53, 58.12) --
	( 86.59, 58.12) --
	( 87.65, 58.12) --
	( 88.71, 58.12) --
	( 89.77, 58.12) --
	( 90.83, 58.12) --
	( 91.89, 58.12) --
	( 92.94, 58.12) --
	( 94.00, 58.12) --
	( 95.06, 58.12) --
	( 96.12, 75.29) --
	( 97.18, 75.29) --
	( 98.24, 75.29) --
	( 99.30, 75.29) --
	(100.36, 75.29) --
	(101.41, 75.29) --
	(102.47, 75.29) --
	(103.53, 75.29) --
	(104.59, 75.29) --
	(105.65, 75.29) --
	(106.71, 75.29) --
	(107.77, 75.29) --
	(108.83, 75.29) --
	(109.89, 95.40) --
	(110.94, 95.40) --
	(112.00, 95.40) --
	(113.06, 95.40) --
	(114.12, 95.40) --
	(115.18, 95.40) --
	(116.24, 95.40) --
	(117.30, 95.40) --
	(118.36, 95.40) --
	(119.42, 95.40) --
	(120.47, 95.40) --
	(121.53, 95.40) --
	(122.59,116.47) --
	(123.65,116.47) --
	(124.71,116.47) --
	(125.77,116.47) --
	(126.83,116.47) --
	(127.89,116.47) --
	(128.95,116.47) --
	(130.00,116.47) --
	(131.06,116.47) --
	(132.12,116.47) --
	(133.18,116.47) --
	(134.24,116.47) --
	(135.30,135.22);
\end{scope}
\begin{scope}
\path[clip] (  0.00,  0.00) rectangle (144.54,144.54);
\definecolor{drawColor}{gray}{0.30}

\node[text=drawColor,anchor=base east,inner sep=0pt, outer sep=0pt, scale=  0.64] at ( 21.63, 26.53) {0.0};

\node[text=drawColor,anchor=base east,inner sep=0pt, outer sep=0pt, scale=  0.64] at ( 21.63, 57.65) {0.1};

\node[text=drawColor,anchor=base east,inner sep=0pt, outer sep=0pt, scale=  0.64] at ( 21.63, 88.77) {0.2};

\node[text=drawColor,anchor=base east,inner sep=0pt, outer sep=0pt, scale=  0.64] at ( 21.63,119.90) {0.3};
\end{scope}
\begin{scope}
\path[clip] (  0.00,  0.00) rectangle (144.54,144.54);
\definecolor{drawColor}{gray}{0.30}

\node[text=drawColor,anchor=base,inner sep=0pt, outer sep=0pt, scale=  0.64] at ( 30.47, 15.40) {0};

\node[text=drawColor,anchor=base,inner sep=0pt, outer sep=0pt, scale=  0.64] at ( 56.68, 15.40) {5};

\node[text=drawColor,anchor=base,inner sep=0pt, outer sep=0pt, scale=  0.64] at ( 82.88, 15.40) {10};

\node[text=drawColor,anchor=base,inner sep=0pt, outer sep=0pt, scale=  0.64] at (109.09, 15.40) {15};

\node[text=drawColor,anchor=base,inner sep=0pt, outer sep=0pt, scale=  0.64] at (135.30, 15.40) {20};
\end{scope}
\begin{scope}
\path[clip] (  0.00,  0.00) rectangle (144.54,144.54);
\definecolor{drawColor}{RGB}{0,0,0}

\node[text=drawColor,anchor=base,inner sep=0pt, outer sep=0pt, scale=  0.80] at ( 82.88,  5.94) {$\alpha$};
\end{scope}
\begin{scope}
\path[clip] (  0.00,  0.00) rectangle (144.54,144.54);
\definecolor{drawColor}{RGB}{248,118,109}

\path[draw=drawColor,line width= 0.6pt,line join=round] ( 47.09,124.34) -- ( 58.65,124.34);
\end{scope}
\begin{scope}
\path[clip] (  0.00,  0.00) rectangle (144.54,144.54);
\definecolor{drawColor}{RGB}{0,191,196}

\path[draw=drawColor,line width= 0.6pt,dash pattern=on 2pt off 2pt ,line join=round] ( 47.09,109.89) -- ( 58.65,109.89);
\end{scope}
\begin{scope}
\path[clip] (  0.00,  0.00) rectangle (144.54,144.54);
\definecolor{drawColor}{RGB}{0,0,0}

\node[text=drawColor,anchor=base west,inner sep=0pt, outer sep=0pt, scale=  0.64] at ( 60.10,122.14) {Instability};
\end{scope}
\begin{scope}
\path[clip] (  0.00,  0.00) rectangle (144.54,144.54);
\definecolor{drawColor}{RGB}{0,0,0}

\node[text=drawColor,anchor=base west,inner sep=0pt, outer sep=0pt, scale=  0.64] at ( 60.10,107.68) {SAA-SubOpt};
\end{scope}
\end{tikzpicture} \vspace{-20pt}}
        \caption{$p_{01} = .3$, $s = .5$}
    \end{subfigure}
    ~ %
    \begin{subfigure}[b]{0.3\textwidth}
    	\ifdraft{Figures Goes Here}
        {
\begin{tikzpicture}[x=1pt,y=1pt]
\definecolor{fillColor}{RGB}{255,255,255}
\path[use as bounding box,fill=fillColor,fill opacity=0.00] (0,0) rectangle (144.54,144.54);
\begin{scope}
\path[clip] ( 28.43, 23.41) rectangle (140.54,140.54);
\definecolor{drawColor}{gray}{0.92}

\path[draw=drawColor,line width= 0.2pt,line join=round] ( 28.43, 39.43) --
	(140.54, 39.43);

\path[draw=drawColor,line width= 0.2pt,line join=round] ( 28.43, 60.84) --
	(140.54, 60.84);

\path[draw=drawColor,line width= 0.2pt,line join=round] ( 28.43, 82.24) --
	(140.54, 82.24);

\path[draw=drawColor,line width= 0.2pt,line join=round] ( 28.43,103.64) --
	(140.54,103.64);

\path[draw=drawColor,line width= 0.2pt,line join=round] ( 28.43,125.05) --
	(140.54,125.05);

\path[draw=drawColor,line width= 0.2pt,line join=round] ( 46.26, 23.41) --
	( 46.26,140.54);

\path[draw=drawColor,line width= 0.2pt,line join=round] ( 71.74, 23.41) --
	( 71.74,140.54);

\path[draw=drawColor,line width= 0.2pt,line join=round] ( 97.22, 23.41) --
	( 97.22,140.54);

\path[draw=drawColor,line width= 0.2pt,line join=round] (122.70, 23.41) --
	(122.70,140.54);

\path[draw=drawColor,line width= 0.4pt,line join=round] ( 28.43, 28.73) --
	(140.54, 28.73);

\path[draw=drawColor,line width= 0.4pt,line join=round] ( 28.43, 50.13) --
	(140.54, 50.13);

\path[draw=drawColor,line width= 0.4pt,line join=round] ( 28.43, 71.54) --
	(140.54, 71.54);

\path[draw=drawColor,line width= 0.4pt,line join=round] ( 28.43, 92.94) --
	(140.54, 92.94);

\path[draw=drawColor,line width= 0.4pt,line join=round] ( 28.43,114.35) --
	(140.54,114.35);

\path[draw=drawColor,line width= 0.4pt,line join=round] ( 28.43,135.75) --
	(140.54,135.75);

\path[draw=drawColor,line width= 0.4pt,line join=round] ( 33.52, 23.41) --
	( 33.52,140.54);

\path[draw=drawColor,line width= 0.4pt,line join=round] ( 59.00, 23.41) --
	( 59.00,140.54);

\path[draw=drawColor,line width= 0.4pt,line join=round] ( 84.48, 23.41) --
	( 84.48,140.54);

\path[draw=drawColor,line width= 0.4pt,line join=round] (109.96, 23.41) --
	(109.96,140.54);

\path[draw=drawColor,line width= 0.4pt,line join=round] (135.44, 23.41) --
	(135.44,140.54);
\definecolor{drawColor}{RGB}{248,118,109}

\path[draw=drawColor,line width= 0.6pt,line join=round] ( 33.52,135.22) --
	( 34.55, 97.09) --
	( 35.58, 97.09) --
	( 36.61, 97.09) --
	( 37.64, 97.09) --
	( 38.67, 97.09) --
	( 39.70, 97.09) --
	( 40.73, 97.09) --
	( 41.76, 97.09) --
	( 42.79, 97.09) --
	( 43.82, 69.35) --
	( 44.85, 69.35) --
	( 45.88, 69.35) --
	( 46.91, 69.35) --
	( 47.94, 69.35) --
	( 48.97, 69.35) --
	( 50.00, 69.35) --
	( 51.03, 69.35) --
	( 52.06, 69.35) --
	( 53.08, 69.35) --
	( 54.11, 51.18) --
	( 55.14, 51.18) --
	( 56.17, 51.18) --
	( 57.20, 51.18) --
	( 58.23, 51.18) --
	( 59.26, 51.18) --
	( 60.29, 51.18) --
	( 61.32, 51.18) --
	( 62.35, 51.18) --
	( 63.38, 51.18) --
	( 64.41, 40.22) --
	( 65.44, 40.22) --
	( 66.47, 40.22) --
	( 67.50, 40.22) --
	( 68.53, 40.22) --
	( 69.56, 40.22) --
	( 70.59, 40.22) --
	( 71.62, 40.22) --
	( 72.65, 40.22) --
	( 73.67, 40.22) --
	( 74.70, 33.87) --
	( 75.73, 33.87) --
	( 76.76, 33.87) --
	( 77.79, 33.87) --
	( 78.82, 33.87) --
	( 79.85, 33.87) --
	( 80.88, 33.87) --
	( 81.91, 33.87) --
	( 82.94, 33.87) --
	( 83.97, 33.87) --
	( 85.00, 31.06) --
	( 86.03, 31.06) --
	( 87.06, 31.06) --
	( 88.09, 31.06) --
	( 89.12, 31.06) --
	( 90.15, 31.06) --
	( 91.18, 31.06) --
	( 92.21, 31.06) --
	( 93.24, 31.06) --
	( 94.26, 31.06) --
	( 95.29, 29.86) --
	( 96.32, 29.86) --
	( 97.35, 29.86) --
	( 98.38, 29.86) --
	( 99.41, 29.86) --
	(100.44, 29.86) --
	(101.47, 29.86) --
	(102.50, 29.86) --
	(103.53, 29.86) --
	(104.56, 29.86) --
	(105.59, 28.92) --
	(106.62, 28.92) --
	(107.65, 28.92) --
	(108.68, 28.92) --
	(109.71, 28.92) --
	(110.74, 28.92) --
	(111.77, 28.92) --
	(112.80, 28.92) --
	(113.82, 28.92) --
	(114.85, 28.92) --
	(115.88, 28.73) --
	(116.91, 28.73) --
	(117.94, 28.73) --
	(118.97, 28.73) --
	(120.00, 28.73) --
	(121.03, 28.73) --
	(122.06, 28.73) --
	(123.09, 28.73) --
	(124.12, 28.73) --
	(125.15, 28.73) --
	(126.18, 28.73) --
	(127.21, 28.73) --
	(128.24, 28.73) --
	(129.27, 28.73) --
	(130.30, 28.73) --
	(131.33, 28.73) --
	(132.36, 28.73) --
	(133.39, 28.73) --
	(134.41, 28.73) --
	(135.44, 28.73);
\definecolor{drawColor}{RGB}{0,191,196}

\path[draw=drawColor,line width= 0.6pt,dash pattern=on 2pt off 2pt ,line join=round] ( 33.52, 28.73) --
	( 34.55, 28.73) --
	( 35.58, 28.73) --
	( 36.61, 28.73) --
	( 37.64, 28.73) --
	( 38.67, 28.73) --
	( 39.70, 28.73) --
	( 40.73, 28.73) --
	( 41.76, 28.73) --
	( 42.79, 28.73) --
	( 43.82, 34.21) --
	( 44.85, 34.21) --
	( 45.88, 34.21) --
	( 46.91, 34.21) --
	( 47.94, 34.21) --
	( 48.97, 34.21) --
	( 50.00, 34.21) --
	( 51.03, 34.21) --
	( 52.06, 34.21) --
	( 53.08, 34.21) --
	( 54.11, 40.93) --
	( 55.14, 40.93) --
	( 56.17, 40.93) --
	( 57.20, 40.93) --
	( 58.23, 40.93) --
	( 59.26, 40.93) --
	( 60.29, 40.93) --
	( 61.32, 40.93) --
	( 62.35, 40.93) --
	( 63.38, 40.93) --
	( 64.41, 45.87) --
	( 65.44, 45.87) --
	( 66.47, 45.87) --
	( 67.50, 45.87) --
	( 68.53, 45.87) --
	( 69.56, 45.87) --
	( 70.59, 45.87) --
	( 71.62, 45.87) --
	( 72.65, 45.87) --
	( 73.67, 45.87) --
	( 74.70, 49.64) --
	( 75.73, 49.64) --
	( 76.76, 49.64) --
	( 77.79, 49.64) --
	( 78.82, 49.64) --
	( 79.85, 49.64) --
	( 80.88, 49.64) --
	( 81.91, 49.64) --
	( 82.94, 49.64) --
	( 83.97, 49.64) --
	( 85.00, 51.14) --
	( 86.03, 51.14) --
	( 87.06, 51.14) --
	( 88.09, 51.14) --
	( 89.12, 51.14) --
	( 90.15, 51.14) --
	( 91.18, 51.14) --
	( 92.21, 51.14) --
	( 93.24, 51.14) --
	( 94.26, 51.14) --
	( 95.29, 52.04) --
	( 96.32, 52.04) --
	( 97.35, 52.04) --
	( 98.38, 52.04) --
	( 99.41, 52.04) --
	(100.44, 52.04) --
	(101.47, 52.04) --
	(102.50, 52.04) --
	(103.53, 52.04) --
	(104.56, 52.04) --
	(105.59, 52.64) --
	(106.62, 52.64) --
	(107.65, 52.64) --
	(108.68, 52.64) --
	(109.71, 52.64) --
	(110.74, 52.64) --
	(111.77, 52.64) --
	(112.80, 52.64) --
	(113.82, 52.64) --
	(114.85, 52.64) --
	(115.88, 52.64) --
	(116.91, 52.64) --
	(117.94, 52.64) --
	(118.97, 52.64) --
	(120.00, 52.64) --
	(121.03, 52.64) --
	(122.06, 52.64) --
	(123.09, 52.64) --
	(124.12, 52.64) --
	(125.15, 52.64) --
	(126.18, 52.64) --
	(127.21, 52.64) --
	(128.24, 52.64) --
	(129.27, 52.64) --
	(130.30, 52.64) --
	(131.33, 52.64) --
	(132.36, 52.64) --
	(133.39, 52.64) --
	(134.41, 52.64) --
	(135.44, 52.64);
\end{scope}
\begin{scope}
\path[clip] (  0.00,  0.00) rectangle (144.54,144.54);
\definecolor{drawColor}{gray}{0.30}

\node[text=drawColor,anchor=base east,inner sep=0pt, outer sep=0pt, scale=  0.64] at ( 24.83, 26.53) {0.00};

\node[text=drawColor,anchor=base east,inner sep=0pt, outer sep=0pt, scale=  0.64] at ( 24.83, 47.93) {0.01};

\node[text=drawColor,anchor=base east,inner sep=0pt, outer sep=0pt, scale=  0.64] at ( 24.83, 69.33) {0.02};

\node[text=drawColor,anchor=base east,inner sep=0pt, outer sep=0pt, scale=  0.64] at ( 24.83, 90.74) {0.03};

\node[text=drawColor,anchor=base east,inner sep=0pt, outer sep=0pt, scale=  0.64] at ( 24.83,112.14) {0.04};

\node[text=drawColor,anchor=base east,inner sep=0pt, outer sep=0pt, scale=  0.64] at ( 24.83,133.55) {0.05};
\end{scope}
\begin{scope}
\path[clip] (  0.00,  0.00) rectangle (144.54,144.54);
\definecolor{drawColor}{gray}{0.30}

\node[text=drawColor,anchor=base,inner sep=0pt, outer sep=0pt, scale=  0.64] at ( 33.52, 15.40) {0};

\node[text=drawColor,anchor=base,inner sep=0pt, outer sep=0pt, scale=  0.64] at ( 59.00, 15.40) {5};

\node[text=drawColor,anchor=base,inner sep=0pt, outer sep=0pt, scale=  0.64] at ( 84.48, 15.40) {10};

\node[text=drawColor,anchor=base,inner sep=0pt, outer sep=0pt, scale=  0.64] at (109.96, 15.40) {15};

\node[text=drawColor,anchor=base,inner sep=0pt, outer sep=0pt, scale=  0.64] at (135.44, 15.40) {20};
\end{scope}
\begin{scope}
\path[clip] (  0.00,  0.00) rectangle (144.54,144.54);
\definecolor{drawColor}{RGB}{0,0,0}

\node[text=drawColor,anchor=base,inner sep=0pt, outer sep=0pt, scale=  0.80] at ( 84.48,  5.94) {$\alpha$};
\end{scope}
\begin{scope}
\path[clip] (  0.00,  0.00) rectangle (144.54,144.54);
\definecolor{drawColor}{RGB}{248,118,109}

\path[draw=drawColor,line width= 0.6pt,line join=round] ( 49.01,124.34) -- ( 60.57,124.34);
\end{scope}
\begin{scope}
\path[clip] (  0.00,  0.00) rectangle (144.54,144.54);
\definecolor{drawColor}{RGB}{0,191,196}

\path[draw=drawColor,line width= 0.6pt,dash pattern=on 2pt off 2pt ,line join=round] ( 49.01,109.89) -- ( 60.57,109.89);
\end{scope}
\begin{scope}
\path[clip] (  0.00,  0.00) rectangle (144.54,144.54);
\definecolor{drawColor}{RGB}{0,0,0}

\node[text=drawColor,anchor=base west,inner sep=0pt, outer sep=0pt, scale=  0.64] at ( 62.02,122.14) {Instability};
\end{scope}
\begin{scope}
\path[clip] (  0.00,  0.00) rectangle (144.54,144.54);
\definecolor{drawColor}{RGB}{0,0,0}

\node[text=drawColor,anchor=base west,inner sep=0pt, outer sep=0pt, scale=  0.64] at ( 62.02,107.68) {SAA-SubOpt};
\end{scope}
\end{tikzpicture}\vspace{-20pt}}
        \caption{$p_{01} = .75$, $s = .5$}
    \end{subfigure}
    ~ %
    \begin{subfigure}[b]{0.3\textwidth}
	\ifdraft{Figure Goes Here}
	{
\begin{tikzpicture}[x=1pt,y=1pt]
\definecolor{fillColor}{RGB}{255,255,255}
\path[use as bounding box,fill=fillColor,fill opacity=0.00] (0,0) rectangle (144.54,144.54);
\begin{scope}
\path[clip] ( 28.43, 23.41) rectangle (140.54,140.54);
\definecolor{drawColor}{gray}{0.92}

\path[draw=drawColor,line width= 0.2pt,line join=round] ( 28.43, 38.89) --
	(140.54, 38.89);

\path[draw=drawColor,line width= 0.2pt,line join=round] ( 28.43, 59.22) --
	(140.54, 59.22);

\path[draw=drawColor,line width= 0.2pt,line join=round] ( 28.43, 79.54) --
	(140.54, 79.54);

\path[draw=drawColor,line width= 0.2pt,line join=round] ( 28.43, 99.86) --
	(140.54, 99.86);

\path[draw=drawColor,line width= 0.2pt,line join=round] ( 28.43,120.19) --
	(140.54,120.19);

\path[draw=drawColor,line width= 0.2pt,line join=round] ( 28.43,140.51) --
	(140.54,140.51);

\path[draw=drawColor,line width= 0.2pt,line join=round] ( 46.26, 23.41) --
	( 46.26,140.54);

\path[draw=drawColor,line width= 0.2pt,line join=round] ( 71.74, 23.41) --
	( 71.74,140.54);

\path[draw=drawColor,line width= 0.2pt,line join=round] ( 97.22, 23.41) --
	( 97.22,140.54);

\path[draw=drawColor,line width= 0.2pt,line join=round] (122.70, 23.41) --
	(122.70,140.54);

\path[draw=drawColor,line width= 0.4pt,line join=round] ( 28.43, 28.73) --
	(140.54, 28.73);

\path[draw=drawColor,line width= 0.4pt,line join=round] ( 28.43, 49.05) --
	(140.54, 49.05);

\path[draw=drawColor,line width= 0.4pt,line join=round] ( 28.43, 69.38) --
	(140.54, 69.38);

\path[draw=drawColor,line width= 0.4pt,line join=round] ( 28.43, 89.70) --
	(140.54, 89.70);

\path[draw=drawColor,line width= 0.4pt,line join=round] ( 28.43,110.02) --
	(140.54,110.02);

\path[draw=drawColor,line width= 0.4pt,line join=round] ( 28.43,130.35) --
	(140.54,130.35);

\path[draw=drawColor,line width= 0.4pt,line join=round] ( 33.52, 23.41) --
	( 33.52,140.54);

\path[draw=drawColor,line width= 0.4pt,line join=round] ( 59.00, 23.41) --
	( 59.00,140.54);

\path[draw=drawColor,line width= 0.4pt,line join=round] ( 84.48, 23.41) --
	( 84.48,140.54);

\path[draw=drawColor,line width= 0.4pt,line join=round] (109.96, 23.41) --
	(109.96,140.54);

\path[draw=drawColor,line width= 0.4pt,line join=round] (135.44, 23.41) --
	(135.44,140.54);
\definecolor{drawColor}{RGB}{248,118,109}

\path[draw=drawColor,line width= 0.6pt,line join=round] ( 33.52,135.22) --
	( 34.55,135.22) --
	( 35.58,128.30) --
	( 36.61,128.30) --
	( 37.64,119.48) --
	( 38.67,119.48) --
	( 39.70,111.98) --
	( 40.73,111.98) --
	( 41.76,100.89) --
	( 42.79,100.89) --
	( 43.82, 90.60) --
	( 44.85, 90.60) --
	( 45.88, 78.64) --
	( 46.91, 78.64) --
	( 47.94, 67.97) --
	( 48.97, 67.97) --
	( 50.00, 59.22) --
	( 51.03, 59.22) --
	( 52.06, 49.25) --
	( 53.08, 49.25) --
	( 54.11, 44.20) --
	( 55.14, 44.20) --
	( 56.17, 38.38) --
	( 57.20, 38.38) --
	( 58.23, 35.01) --
	( 59.26, 35.01) --
	( 60.29, 32.25) --
	( 61.32, 32.25) --
	( 62.35, 30.34) --
	( 63.38, 30.34) --
	( 64.41, 29.87) --
	( 65.44, 29.87) --
	( 66.47, 28.83) --
	( 67.50, 28.83) --
	( 68.53, 28.77) --
	( 69.56, 28.77) --
	( 70.59, 29.12) --
	( 71.62, 29.12) --
	( 72.65, 28.73) --
	( 73.67, 28.73) --
	( 74.70, 28.73) --
	( 75.73, 28.73) --
	( 76.76, 28.73) --
	( 77.79, 28.73) --
	( 78.82, 28.73) --
	( 79.85, 28.73) --
	( 80.88, 28.73) --
	( 81.91, 28.73) --
	( 82.94, 28.73) --
	( 83.97, 28.73) --
	( 85.00, 28.73) --
	( 86.03, 28.73) --
	( 87.06, 28.73) --
	( 88.09, 28.73) --
	( 89.12, 28.73) --
	( 90.15, 28.73) --
	( 91.18, 28.73) --
	( 92.21, 28.73) --
	( 93.24, 28.73) --
	( 94.26, 28.73) --
	( 95.29, 28.73) --
	( 96.32, 28.73) --
	( 97.35, 28.73) --
	( 98.38, 28.73) --
	( 99.41, 28.73) --
	(100.44, 28.73) --
	(101.47, 28.73) --
	(102.50, 28.73) --
	(103.53, 28.73) --
	(104.56, 28.73) --
	(105.59, 28.73) --
	(106.62, 28.73) --
	(107.65, 28.73) --
	(108.68, 28.73) --
	(109.71, 28.73) --
	(110.74, 28.73) --
	(111.77, 28.73) --
	(112.80, 28.73) --
	(113.82, 28.73) --
	(114.85, 28.73) --
	(115.88, 28.73) --
	(116.91, 28.73) --
	(117.94, 28.73) --
	(118.97, 28.73) --
	(120.00, 28.73) --
	(121.03, 28.73) --
	(122.06, 28.73) --
	(123.09, 28.73) --
	(124.12, 28.73) --
	(125.15, 28.73) --
	(126.18, 28.73) --
	(127.21, 28.73) --
	(128.24, 28.73) --
	(129.27, 28.73) --
	(130.30, 28.73) --
	(131.33, 28.73) --
	(132.36, 28.73) --
	(133.39, 28.73) --
	(134.41, 28.73) --
	(135.44, 28.73);
\definecolor{drawColor}{RGB}{0,191,196}

\path[draw=drawColor,line width= 0.6pt,dash pattern=on 2pt off 2pt ,line join=round] ( 33.52, 28.73) --
	( 34.55, 28.73) --
	( 35.58, 31.38) --
	( 36.61, 31.38) --
	( 37.64, 36.75) --
	( 38.67, 36.75) --
	( 39.70, 43.45) --
	( 40.73, 43.45) --
	( 41.76, 53.05) --
	( 42.79, 53.05) --
	( 43.82, 62.95) --
	( 44.85, 62.95) --
	( 45.88, 74.23) --
	( 46.91, 74.23) --
	( 47.94, 85.79) --
	( 48.97, 85.79) --
	( 50.00, 95.91) --
	( 51.03, 95.91) --
	( 52.06,105.84) --
	( 53.08,105.84) --
	( 54.11,112.90) --
	( 55.14,112.90) --
	( 56.17,119.55) --
	( 57.20,119.55) --
	( 58.23,123.70) --
	( 59.26,123.70) --
	( 60.29,126.67) --
	( 61.32,126.67) --
	( 62.35,128.38) --
	( 63.38,128.38) --
	( 64.41,128.91) --
	( 65.44,128.91) --
	( 66.47,129.64) --
	( 67.50,129.64) --
	( 68.53,129.64) --
	( 69.56,129.64) --
	( 70.59,129.64) --
	( 71.62,129.64) --
	( 72.65,130.03) --
	( 73.67,130.03) --
	( 74.70,130.03) --
	( 75.73,130.03) --
	( 76.76,130.03) --
	( 77.79,130.03) --
	( 78.82,130.03) --
	( 79.85,130.03) --
	( 80.88,130.03) --
	( 81.91,130.03) --
	( 82.94,130.03) --
	( 83.97,130.03) --
	( 85.00,130.03) --
	( 86.03,130.03) --
	( 87.06,130.03) --
	( 88.09,130.03) --
	( 89.12,130.03) --
	( 90.15,130.03) --
	( 91.18,130.03) --
	( 92.21,130.03) --
	( 93.24,130.03) --
	( 94.26,130.03) --
	( 95.29,130.03) --
	( 96.32,130.03) --
	( 97.35,130.03) --
	( 98.38,130.03) --
	( 99.41,130.03) --
	(100.44,130.03) --
	(101.47,130.03) --
	(102.50,130.03) --
	(103.53,130.03) --
	(104.56,130.03) --
	(105.59,130.03) --
	(106.62,130.03) --
	(107.65,130.03) --
	(108.68,130.03) --
	(109.71,130.03) --
	(110.74,130.03) --
	(111.77,130.03) --
	(112.80,130.03) --
	(113.82,130.03) --
	(114.85,130.03) --
	(115.88,130.03) --
	(116.91,130.03) --
	(117.94,130.03) --
	(118.97,130.03) --
	(120.00,130.03) --
	(121.03,130.03) --
	(122.06,130.03) --
	(123.09,130.03) --
	(124.12,130.03) --
	(125.15,130.03) --
	(126.18,130.03) --
	(127.21,130.03) --
	(128.24,130.03) --
	(129.27,130.03) --
	(130.30,130.03) --
	(131.33,130.03) --
	(132.36,130.03) --
	(133.39,130.03) --
	(134.41,130.03) --
	(135.44,130.03);
\end{scope}
\begin{scope}
\path[clip] (  0.00,  0.00) rectangle (144.54,144.54);
\definecolor{drawColor}{gray}{0.30}

\node[text=drawColor,anchor=base east,inner sep=0pt, outer sep=0pt, scale=  0.64] at ( 24.83, 26.53) {0.00};

\node[text=drawColor,anchor=base east,inner sep=0pt, outer sep=0pt, scale=  0.64] at ( 24.83, 46.85) {0.01};

\node[text=drawColor,anchor=base east,inner sep=0pt, outer sep=0pt, scale=  0.64] at ( 24.83, 67.17) {0.02};

\node[text=drawColor,anchor=base east,inner sep=0pt, outer sep=0pt, scale=  0.64] at ( 24.83, 87.50) {0.03};

\node[text=drawColor,anchor=base east,inner sep=0pt, outer sep=0pt, scale=  0.64] at ( 24.83,107.82) {0.04};

\node[text=drawColor,anchor=base east,inner sep=0pt, outer sep=0pt, scale=  0.64] at ( 24.83,128.14) {0.05};
\end{scope}
\begin{scope}
\path[clip] (  0.00,  0.00) rectangle (144.54,144.54);
\definecolor{drawColor}{gray}{0.30}

\node[text=drawColor,anchor=base,inner sep=0pt, outer sep=0pt, scale=  0.64] at ( 33.52, 15.40) {0};

\node[text=drawColor,anchor=base,inner sep=0pt, outer sep=0pt, scale=  0.64] at ( 59.00, 15.40) {5};

\node[text=drawColor,anchor=base,inner sep=0pt, outer sep=0pt, scale=  0.64] at ( 84.48, 15.40) {10};

\node[text=drawColor,anchor=base,inner sep=0pt, outer sep=0pt, scale=  0.64] at (109.96, 15.40) {15};

\node[text=drawColor,anchor=base,inner sep=0pt, outer sep=0pt, scale=  0.64] at (135.44, 15.40) {20};
\end{scope}
\begin{scope}
\path[clip] (  0.00,  0.00) rectangle (144.54,144.54);
\definecolor{drawColor}{RGB}{0,0,0}

\node[text=drawColor,anchor=base,inner sep=0pt, outer sep=0pt, scale=  0.80] at ( 84.48,  5.94) {$\alpha$};
\end{scope}
\begin{scope}
\path[clip] (  0.00,  0.00) rectangle (144.54,144.54);
\definecolor{drawColor}{RGB}{248,118,109}

\path[draw=drawColor,line width= 0.6pt,line join=round] ( 71.43, 89.20) -- ( 82.99, 89.20);
\end{scope}
\begin{scope}
\path[clip] (  0.00,  0.00) rectangle (144.54,144.54);
\definecolor{drawColor}{RGB}{0,191,196}

\path[draw=drawColor,line width= 0.6pt,dash pattern=on 2pt off 2pt ,line join=round] ( 71.43, 74.75) -- ( 82.99, 74.75);
\end{scope}
\begin{scope}
\path[clip] (  0.00,  0.00) rectangle (144.54,144.54);
\definecolor{drawColor}{RGB}{0,0,0}

\node[text=drawColor,anchor=base west,inner sep=0pt, outer sep=0pt, scale=  0.64] at ( 84.44, 87.00) {Instability};
\end{scope}
\begin{scope}
\path[clip] (  0.00,  0.00) rectangle (144.54,144.54);
\definecolor{drawColor}{RGB}{0,0,0}

\node[text=drawColor,anchor=base west,inner sep=0pt, outer sep=0pt, scale=  0.64] at ( 84.44, 72.54) {SAA-SubOpt};
\end{scope}
\end{tikzpicture} \vspace{-20pt}}
        \caption{$p_{01} = .3$, $s = .2$}
    \end{subfigure}
    \caption{\textbf{Sub-Optimality-Instability Curves}.  We consider $K = 10,000$ newsvendors where $p_{k1} \sim \text{Uniform}[.6, .9]$, $\Nhat_k \sim \text{Poisson}(10)$, and a single data draw.  The values of $p_{01}$ and the critical fractile $s$ is given in each panel.  In the first panel, instability initially increases, and there is no benefit to pooling.  In the second and third,  instability is decreasing and there is a benefit to pooling.}
    \label{fig:TradeoffPlot}
	\vspace{-10pt}
\end{figure}
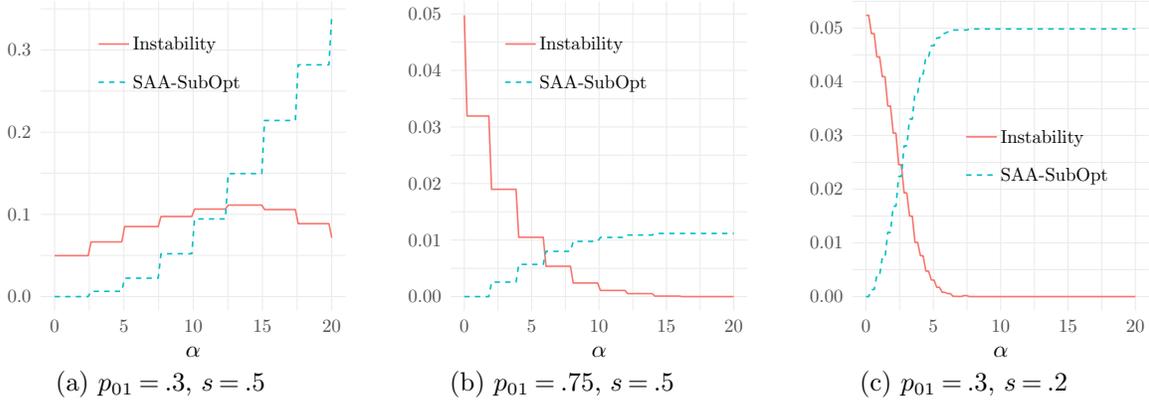

Our comments above are qualitative, focusing on developing intuition.  However, the Sub-Optimality-Instability Tradeoff also 
provides a quantitative diagnostic tool for studying data-pooling. Indeed, for simple optimization problems such as minimizing MSE, it may be possible to analytically study the effects of pooling (cf. \cref{thm:DataPoolingJS}), but for more complex optimization problems where $\bx_k(\alpha, h(\bfmhat), \bfmhat_k)$ is not known analytically, such a study is not generally possible.  Fortunately, both $\text{SAA-SubOpt}(\alpha)$ and $\text{Instability}(\alpha)$ can be evaluated \emph{directly from the data}.  Studying their dependence on $\alpha$ for a particular instance provides insight into how data-pooling improves (or does not improve) solution quality.
We illustrate with  \cref{ex:BabyNewsvendor}:
\begin{example}[Simple Newsvendor Revisited]
\label{ex:BabyNewsvendorTake2}
We revisit \cref{ex:BabyNewsvendor} and simulate an instance with $K = 1000$, $p_{k1}$ distributed uniformly on $[.6, .9]$ and $p_{01} = .3$.  One can confirm that as in \cref{ex:BabyNewsvendor}, data-pooling offers no benefit over decoupling (regardless of the choice of $\Nhat_k$) for these parameters.  We take $\Nhat_k \sim \text{Poisson}(10)$ for all $k$, and simulate a single data realization $\bfmhat$.  

Using the data, we can evaluate $\text{SAA-SubOpt}(\alpha)$ and $\text{Instability}(\alpha)$ explicitly.  We plot them in the first panel of Fig.~\ref{fig:TradeoffPlot}.  Notice that as expected, $\text{SAA-SubOpt}(\alpha)$ increases steadily in $\alpha$, however, perhaps surprisingly, $\text{Instability}(\alpha)$ \emph{increases} at first, before ultimately decreasing.  The reason is that as in \cref{ex:BabyNewsvendor}, $\bx_k(\alpha, \bp_0, \bfmhat_k) = \I{ \phat_{k1}(\alpha) \geq 1/2}$.  For small positive $\alpha$, $\phat_{k1}(\alpha)$ is generally closer to $\frac{1}{2}$ than $\phat_{k1}$, and since $ \frac{1}{2}$ is the critical threshold where $\bx_k(\alpha, \bp_0, \bfmhat)$ changes values, the solution is less stable.  Hence, $\text{Instability}(\alpha)$ \emph{increases} for small $\alpha$.  Because of this initial increasing behavior, the ``gains" in stability never outweigh the costs of sub-optimality, and hence decoupling is best.  Indeed, the first panel of \cref{fig:TradeOffLOO} in the appendix shows $\alphaLOO_{\bp_0} = \alphaOR_{\bp_0} = 0.0$.  

We earlier observed that the benefits of pooling depend on the anchor.  We next consider the same parameters and data as above but let $p_{01} = .75$.  The second panel of \cref{fig:TradeoffPlot} shows the Sub-Optimality-Instability tradeoff.  We see here that again $\text{Sub-Optimality}(\alpha)$ is increasing, and, perhaps more intuitively, $\text{Instability}(\alpha)$ is decreasing.  Hence, there is a positive $\alpha$ that minimizes their sum, and  the second panel \cref{fig:TradeOffLOO} shows $\alphaLOO_{\bp_0} \approx \alphaOR_{\bp_0} \approx 16.16$.  

Finally, as mentioned previously, the potential benefits of data-pooling also depends on the problem structure.  The Sub-Optimality-Instability tradeoff allows us to study this dependence.  Consider again letting $p_{01} = .3$, but now consider newsvendor problems with critical fractile $s = .2$.  We again see a benefit to pooling.  The Sub-Optimality-Instability tradeoff is in the last panel of \cref{fig:TradeoffPlot}.  The last panel of \cref{fig:TradeOffLOO} shows $\alphaLOO_{\bp_0} \approx 2.42$ and $\alphaOR_{\bp_0} \approx 2.22$.
\end{example}

In summary, while $\alphaLOO_{h}$ identifies a good choice of shrinkage in many settings,  Sub-Optimality and Instability graphs as above often illuminate \emph{why} this is a good choice of shrinkage, providing insight.  This is particularly helpful for complex optimization problems for which it may be hard to reason about $\bx_k(\alpha, h(\bfmhat), \bfmhat_k)$.

\section{Computational Experiments}
\label{sec:Numerics}
In this section we study the empirical performance of Shrunken-SAA on synthetic and real data.  All code for reproducing these experiments and plots is available at 
\blinded{\url{https://github.com/vgupta1/JS_SAA}}
We focus on assessing the degree to which Shrunken-SAA is robust to violations of the assumptions underlying 
\cref{thm:FixedPointShrinkage,thm:SmoothThmGeneralAnchor,thm:FixedPointShrinkageDiscrete,thm:ShrinkageDiscreteGeneralAnchor,thm:SmoothThmLooAnchor,thm:ShrinkageDiscreteLooAnchor}.  
Specifically, we ask how Shrunken-SAA performs when 
i) $K$ is small to moderate, and not growing to infinity;
ii) \cref{ass:RandomData} is violated, i.e., each $\Nhat_k$ is fixed and non-random;
iii) the true  $\mathbb P_k$ do not have finite, discrete support; or
iv) $N$ grows large.

For simplicity, we take each subproblem to be a newsvendor problem with critical fractile $s = 95\%$.  Since the performance of Shrunken-SAA depends on the true distributions $\bp_k$, we use real sales data from a chain of European pharmacies.  (See \cref{sec:DataBackground} for more details.)  

{\blockedit  \label{newDROBenchmarks}
We compare several policies:   
The first two, \textbf{SAA} and \textbf{KS}, are decoupled-benchmarks.  Recall that for the newsvendor problem, {SAA}, i.e., $\bx(0, \bp_0, \bfmhat)$, is also the optimal solution to a distributionally robust formulation using a Wasserstein ambiguity set (\cite{esfahani2018data}).  We define {KS} to be  an optimal solution to a distributionally robust formulation of the newsvendor problem using the Kolmogorov-Smirnov ambiguity set (see \cref{sec:DataCleaning} for formal definition).  This set enjoys strong large-sample statistical guarantees \citep{bertsimas2018robust}.  

\label{line:ClarifyingJS}
The next three policies, \textbf{JS-Fixed}, \textbf{S-SAA-Fixed} and \textbf{Oracle-Fixed}, each shrink towards the uniform distribution, i.e., a fixed anchor.  They differ in the amount of shrinkage.  {JS-Fixed}, i.e., $\bx(\alphaJS_{\bp_0}, \bp_0, \bfmhat)$, pools according to \cref{thm:DataPoolingJS}; {S-SAA-Fixed}, i.e.,  $\bx(\alphaLOO_{\bp_0}, \bp_0, \bfmhat)$, is our Shrunken-SAA algorithm; and {Oracle-Fixed}, i.e., $\bx(\alphaOR_{\bp_0}, \bp_0, \bfmhat)$ is the oracle shrinkage.  

The next two policies, \textbf{S-SAA-Beta} and \textbf{Oracle-Beta}, each shrink towards a data-driven choice of anchor in $\mP$, where $\mP$ consists of scaled beta-distributions (cf. \cref{sec:DataCleaning}).  S-SAA, i.e., $\bx(\alphaLOO_{\mP}, \hloo(\bfmhat), \bfmhat)$,
uses $\hloo$, while Oracle-Beta, i.e., $\bx(\alphaOR_\mP, \bq^{\sf OR}_{\mP}, \bfmhat)$, uses the oracle anchor.  

Finally, the last set of policies, \textbf{JS-GM}, \textbf{S-SAA-GM} and \textbf{Oracle-GM} each shrink towards the grand-mean distribution, $\bphat^{\sf GM}$.  They differ in the amount of shrinkage. JS-GM, pools according to  \cref{thm:DataPoolingJS}, S-SAA-GM is our Shrunken-SAA Algorithm, and Oracle-GM is the oracle pooling.  

Intuitively, the difference between the JS policies and the decoupled policies illustrates the value of data-pooling in a ``generic" fashion that does not account for the shape of the cost functions.  By contrast, the difference between the Shrunken-SAA policies and the JS policies quantifies the additional benefit of tailoring the amount of pooling to the specific newsvendor cost function.  \edit{Similarly, the difference between the ``Beta" anchor versions and the Fixed versions help quantify the value of a good choice of anchor, and, as we will see, the GM variants highlight that the grand-mean is often a good heuristic choice of anchor.}

}

Before presenting the details, we summarize our main findings.  When $N$ is moderate to large, all methods (including Shrunken-SAA) perform comparably to the full-information solution.  When $N$ is small to moderate, however, our Shrunken-SAA policies provide a significant benefit over SAA and a substantial benefit over JS variants that do not leverage the optimization structure.  This is true even for moderate $K$ ($K \leq 100$) and even when $\Nhat_k$ are fixed (violating \cref{ass:RandomData}).  The value of $d$ has little effect on the performance of Shrunken-SAA; it strongly outperforms decoupling even as $d \rightarrow \infty$.  \edit{Finally, our GM heuristic has very strong performance, comparable to the Beta variants which optimize the choice of anchor, at a much smaller computational cost.}

For ease of comparison in what follows, we present all results as ``\% Benefit over SAA," i.e., bigger values are better.  \edit{In many cases, to aid readability, we only present a subset of benchmark policies on a graph.  In these cases, larger tables with all benchmarks are available in \cref{sec:AdditionalFigures}.}

\subsection{Data Description}
\label{sec:DataBackground}
Our dataset consists of daily sales at the store level for a European pharmacy chain with locations across $7$ countries.  We treat these aggregated store sales as if they were the realized daily demand of a single product.  Although this is clearly a simplification of the underlying inventory management problem, we do not believe it significantly impacts the study of our key questions outlined above.   Additionally, aggregating over products makes demand censoring insignificant.

The original dataset contains $942$ days of data across $1115$ stores.  After some preliminary data-cleaning (see \cref{sec:DataCleaning}), we are left with $629$ days.  Due to local holidays, individual stores may still be closed on these $629$ days.  Almost all ($1105$) stores have at least one missing day, and 16\% of stores have 20\% of days missing.  

Stores vary in size, available assortment of products, promotional activities and prices, creating significant heterogeneity in demand.  The average daily demand ranges from 3,183 to 23,400.  The first panel of \cref{fig:AvgDailyDemands} in Appendix~\ref{sec:AdditionalFigures} plots the average daily demand by store.    The second panel provides a more fine-grained perspective, showing the distribution of daily demand for a few representative stores.  The distributions are quite distinct, at least partially because the overall scale of daily sales differs wildly between stores.  

Finally, with the exception of \cref{sec:Infinited}, we discretize demand by dividing the range of observations into $d$ equally-spaced bins to form the true distributions $\bp_k$.  \Cref{fig:SampleDemands} plots $\bp_k$ for some representative stores when $d = 20$.  We consider these distributions to be quite diverse and far from the uniform distribution (our fixed anchor).  We also plot the distribution of the 95\% quantile with respect to this discretization in the second panel of \cref{fig:SampleDemands}.  Note that it is not the case that 95\% quantile occurs in the same (discretized) bin for each $\bp_k$, i.e., the quantile itself displays some heterogeneity, unlike \cref{ex:DifferentNewsvendor}.

\enlargethispage{10pt}
\subsection{An Idealized Synthetic Dataset}
\label{sec:SyntheticData}
\begin{figure}
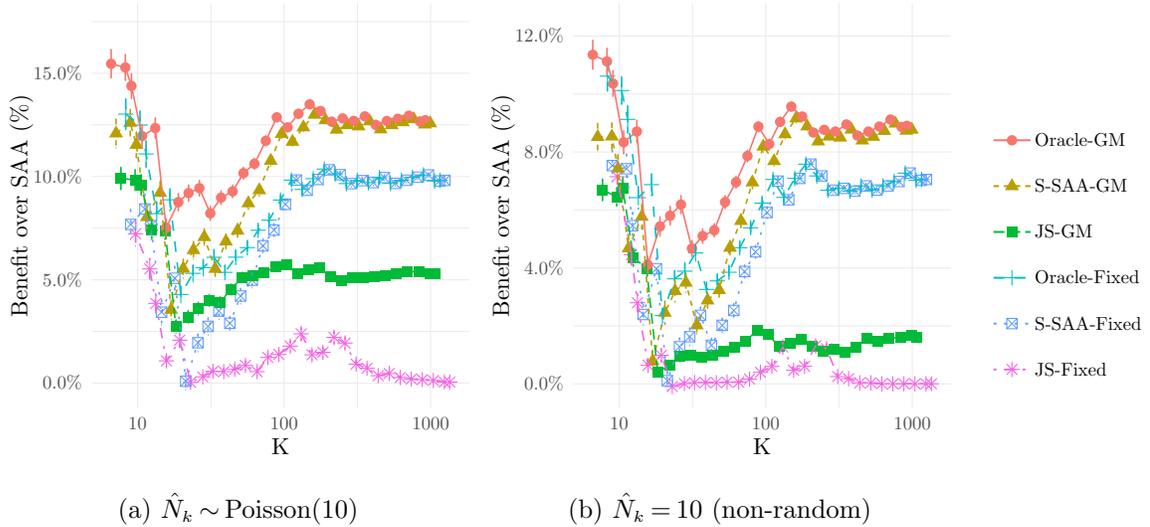

    \centering
    \begin{subfigure}[b]{0.38\textwidth}
	\ifdraft{FiguresGoesHere}{ 
	\input{./Figures/synthDataTrue.tex}\vspace{-10pt} 
	}
        \caption{$\Nhat_k \sim \op{Poisson}(10)$}
    \end{subfigure}
    \begin{subfigure}[b]{0.38\textwidth}
	\ifdraft{FiguresGoesHere}{ 
	\input{./Figures/synthDataFalse.tex} \vspace{-10pt}
	}
        \caption{$\Nhat_k = 10$ (non-random)}
    \end{subfigure}
    \begin{subfigure}[b]{.1\textwidth}
	\ifdraft{FiguresGoesHere}{ 
\begin{tikzpicture}[x=1pt,y=1pt]
\definecolor{fillColor}{RGB}{255,255,255}
\path[use as bounding box,fill=fillColor,fill opacity=0.00] (0,0) rectangle ( 72.27,180.67);
\begin{scope}
\path[clip] (  0.00,  0.00) rectangle ( 72.27,180.67);
\definecolor{fillColor}{RGB}{248,118,109}

\path[fill=fillColor] ( 17.57,133.70) circle (  1.96);
\end{scope}
\begin{scope}
\path[clip] (  0.00,  0.00) rectangle ( 72.27,180.67);
\definecolor{drawColor}{RGB}{248,118,109}

\path[draw=drawColor,line width= 0.6pt,line join=round] ( 10.63,133.70) -- ( 24.51,133.70);
\end{scope}
\begin{scope}
\path[clip] (  0.00,  0.00) rectangle ( 72.27,180.67);
\definecolor{drawColor}{RGB}{248,118,109}

\path[draw=drawColor,line width= 0.6pt,line join=round] ( 10.63,133.70) -- ( 24.51,133.70);
\end{scope}
\begin{scope}
\path[clip] (  0.00,  0.00) rectangle ( 72.27,180.67);
\definecolor{fillColor}{RGB}{183,159,0}

\path[fill=fillColor] ( 17.57,119.41) --
	( 20.21,114.83) --
	( 14.93,114.83) --
	cycle;
\end{scope}
\begin{scope}
\path[clip] (  0.00,  0.00) rectangle ( 72.27,180.67);
\definecolor{drawColor}{RGB}{183,159,0}

\path[draw=drawColor,line width= 0.6pt,dash pattern=on 2pt off 2pt ,line join=round] ( 10.63,116.35) -- ( 24.51,116.35);
\end{scope}
\begin{scope}
\path[clip] (  0.00,  0.00) rectangle ( 72.27,180.67);
\definecolor{drawColor}{RGB}{183,159,0}

\path[draw=drawColor,line width= 0.6pt,dash pattern=on 2pt off 2pt ,line join=round] ( 10.63,116.35) -- ( 24.51,116.35);
\end{scope}
\begin{scope}
\path[clip] (  0.00,  0.00) rectangle ( 72.27,180.67);
\definecolor{fillColor}{RGB}{0,186,56}

\path[fill=fillColor] ( 15.61, 97.05) --
	( 19.53, 97.05) --
	( 19.53,100.97) --
	( 15.61,100.97) --
	cycle;
\end{scope}
\begin{scope}
\path[clip] (  0.00,  0.00) rectangle ( 72.27,180.67);
\definecolor{drawColor}{RGB}{0,186,56}

\path[draw=drawColor,line width= 0.6pt,dash pattern=on 4pt off 2pt ,line join=round] ( 10.63, 99.01) -- ( 24.51, 99.01);
\end{scope}
\begin{scope}
\path[clip] (  0.00,  0.00) rectangle ( 72.27,180.67);
\definecolor{drawColor}{RGB}{0,186,56}

\path[draw=drawColor,line width= 0.6pt,dash pattern=on 4pt off 2pt ,line join=round] ( 10.63, 99.01) -- ( 24.51, 99.01);
\end{scope}
\begin{scope}
\path[clip] (  0.00,  0.00) rectangle ( 72.27,180.67);
\definecolor{drawColor}{RGB}{0,191,196}

\path[draw=drawColor,line width= 0.4pt,line join=round,line cap=round] ( 14.80, 81.67) -- ( 20.35, 81.67);

\path[draw=drawColor,line width= 0.4pt,line join=round,line cap=round] ( 17.57, 78.89) -- ( 17.57, 84.44);
\end{scope}
\begin{scope}
\path[clip] (  0.00,  0.00) rectangle ( 72.27,180.67);
\definecolor{drawColor}{RGB}{0,191,196}

\path[draw=drawColor,line width= 0.6pt,dash pattern=on 4pt off 4pt ,line join=round] ( 10.63, 81.67) -- ( 24.51, 81.67);
\end{scope}
\begin{scope}
\path[clip] (  0.00,  0.00) rectangle ( 72.27,180.67);
\definecolor{drawColor}{RGB}{0,191,196}

\path[draw=drawColor,line width= 0.6pt,dash pattern=on 4pt off 4pt ,line join=round] ( 10.63, 81.67) -- ( 24.51, 81.67);
\end{scope}
\begin{scope}
\path[clip] (  0.00,  0.00) rectangle ( 72.27,180.67);
\definecolor{drawColor}{RGB}{97,156,255}

\path[draw=drawColor,line width= 0.4pt,line join=round,line cap=round] ( 15.61, 62.36) rectangle ( 19.53, 66.28);

\path[draw=drawColor,line width= 0.4pt,line join=round,line cap=round] ( 15.61, 62.36) -- ( 19.53, 66.28);

\path[draw=drawColor,line width= 0.4pt,line join=round,line cap=round] ( 15.61, 66.28) -- ( 19.53, 62.36);
\end{scope}
\begin{scope}
\path[clip] (  0.00,  0.00) rectangle ( 72.27,180.67);
\definecolor{drawColor}{RGB}{97,156,255}

\path[draw=drawColor,line width= 0.6pt,dash pattern=on 1pt off 3pt ,line join=round] ( 10.63, 64.32) -- ( 24.51, 64.32);
\end{scope}
\begin{scope}
\path[clip] (  0.00,  0.00) rectangle ( 72.27,180.67);
\definecolor{drawColor}{RGB}{97,156,255}

\path[draw=drawColor,line width= 0.6pt,dash pattern=on 1pt off 3pt ,line join=round] ( 10.63, 64.32) -- ( 24.51, 64.32);
\end{scope}
\begin{scope}
\path[clip] (  0.00,  0.00) rectangle ( 72.27,180.67);
\definecolor{drawColor}{RGB}{245,100,227}

\path[draw=drawColor,line width= 0.4pt,line join=round,line cap=round] ( 15.61, 45.01) -- ( 19.53, 48.94);

\path[draw=drawColor,line width= 0.4pt,line join=round,line cap=round] ( 15.61, 48.94) -- ( 19.53, 45.01);

\path[draw=drawColor,line width= 0.4pt,line join=round,line cap=round] ( 14.80, 46.98) -- ( 20.35, 46.98);

\path[draw=drawColor,line width= 0.4pt,line join=round,line cap=round] ( 17.57, 44.20) -- ( 17.57, 49.75);
\end{scope}
\begin{scope}
\path[clip] (  0.00,  0.00) rectangle ( 72.27,180.67);
\definecolor{drawColor}{RGB}{245,100,227}

\path[draw=drawColor,line width= 0.6pt,dash pattern=on 1pt off 3pt on 4pt off 3pt ,line join=round] ( 10.63, 46.98) -- ( 24.51, 46.98);
\end{scope}
\begin{scope}
\path[clip] (  0.00,  0.00) rectangle ( 72.27,180.67);
\definecolor{drawColor}{RGB}{245,100,227}

\path[draw=drawColor,line width= 0.6pt,dash pattern=on 1pt off 3pt on 4pt off 3pt ,line join=round] ( 10.63, 46.98) -- ( 24.51, 46.98);
\end{scope}
\begin{scope}
\path[clip] (  0.00,  0.00) rectangle ( 72.27,180.67);
\definecolor{drawColor}{RGB}{0,0,0}

\node[text=drawColor,anchor=base west,inner sep=0pt, outer sep=0pt, scale=  0.64] at ( 26.24,131.50) {Oracle-GM};
\end{scope}
\begin{scope}
\path[clip] (  0.00,  0.00) rectangle ( 72.27,180.67);
\definecolor{drawColor}{RGB}{0,0,0}

\node[text=drawColor,anchor=base west,inner sep=0pt, outer sep=0pt, scale=  0.64] at ( 26.24,114.15) {S-SAA-GM};
\end{scope}
\begin{scope}
\path[clip] (  0.00,  0.00) rectangle ( 72.27,180.67);
\definecolor{drawColor}{RGB}{0,0,0}

\node[text=drawColor,anchor=base west,inner sep=0pt, outer sep=0pt, scale=  0.64] at ( 26.24, 96.81) {JS-GM};
\end{scope}
\begin{scope}
\path[clip] (  0.00,  0.00) rectangle ( 72.27,180.67);
\definecolor{drawColor}{RGB}{0,0,0}

\node[text=drawColor,anchor=base west,inner sep=0pt, outer sep=0pt, scale=  0.64] at ( 26.24, 79.46) {Oracle-Fixed};
\end{scope}
\begin{scope}
\path[clip] (  0.00,  0.00) rectangle ( 72.27,180.67);
\definecolor{drawColor}{RGB}{0,0,0}

\node[text=drawColor,anchor=base west,inner sep=0pt, outer sep=0pt, scale=  0.64] at ( 26.24, 62.12) {S-SAA-Fixed};
\end{scope}
\begin{scope}
\path[clip] (  0.00,  0.00) rectangle ( 72.27,180.67);
\definecolor{drawColor}{RGB}{0,0,0}

\node[text=drawColor,anchor=base west,inner sep=0pt, outer sep=0pt, scale=  0.64] at ( 26.24, 44.77) {JS-Fixed};
\end{scope}
\end{tikzpicture}
	}
    \end{subfigure}
\caption{\textbf{Robustness to \cref{ass:RandomData}.}
Performance of policies on simulated data.  In the first panel, the amount of data per store follows \cref{ass:RandomData} with $N_k = 10$.  In the second panel, the amount of data is fixed at $\Nhat_k = 10$ for all runs.  Error bars show $\pm 1$ standard error.  
\vspace{-10pt}
}\label{fig:SynthData}
\end{figure}

We first consider an ideal setting for Shrunken-SAA.  Specifically, after discretizing demand for each store into $d = 20$ buckets, we set $\bp_k$ to be the empirical distribution of demand over the \emph{entire} dataset with respect to these buckets.  We then simulate synthetic data according to \cref{eq:MultinomialCounts} under \cref{ass:RandomData}.  We train each of our methods using this data, and then evaluate their true performance using the $\bp_k$.  We repeat this process $200$ times.  The left panel of \cref{fig:SynthData} shows the average results \edit{for a subset of the policies.  \Cref{tab:SynthDataTrueAllPolicies} in the appendix includes all policies.}  

As suggested by \cref{thm:FixedPointShrinkageDiscrete,thm:ShrinkageDiscreteGeneralAnchor}, Shrunken-SAA significantly outperforms decoupling even for $K$ as small as $10$.  For large $K$, the benefit is as large as $10-15\%$.  Both of our Shrunken-SAA policies converge quickly to their oracle benchmarks.  We note the JS policies also outperform the decoupled solutions, but by a smaller amount (5-10\%).    
For both sets of policies, shrinking to the grand mean outperforms shrinking to the uniform distribution, since, as observed earlier, the true distributions are far from uniform and have quantiles far from the uniform quantile.  \edit{Indeed, the grand-mean policies perform comparably to our Beta policies (cf. \cref{tab:SynthDataTrueAllPolicies}).}

We also illustrate the standard deviation of the performance for each of these methods in \cref{fig:StdDevPlot} in Appendix~\ref{sec:AdditionalFigures}.  For all approaches, the standard deviation tends to zero as $K \rightarrow \infty$, because the true performance concentrates at its expectation for each method.  For small $K$, our Shrunken-SAA approaches exhibit significantly smaller standard deviation than SAA, and, for larger $K$, the standard deviation is comparable to the oracle values, and much less than JS variants.  The reduction in variability compared to SAA follows intuitively since pooling increases stability.  

Finally, we plot the average amount of shrinkage across runs as a function of $K$ for each method in \cref{fig:AlphaPlot} in Appendix~\ref{sec:AdditionalFigures}.  We observe that the shrinkage amount converges quickly as $K \rightarrow \infty$, and that our Shrunken-SAA methods pool much more than the JS variants.  In particular, when shrinking to the grand-mean or to an optimized Beta distribution, our Shrunken-SAA methods use a value of $\alpha \geq 30$ for large $K$, i.e., placing  $3$ times more weight on the anchor than the data, itself. By contrast, JS variants eventually engage in almost no pooling.

\subsection{Relaxing \cref{ass:RandomData}}
\label{sec:NumericallyTestingAssRandomData}
We next consider robustness to \cref{ass:RandomData}.  Specifically, we repeat the experiment of the previous section but now simulate data with $\Nhat_k = 10$ for all $k$ and all runs.  Results are shown in the second panel of \cref{fig:SynthData}, and \edit{\cref{fig:StdDevPlot,fig:AlphaPlot,tab:SynthDataFalseAllPolicies} in \cref{sec:AdditionalFigures}}.  We see the same qualitative features.  Specifically, our Shrunken-SAA methods converge to oracle performance, and, even for moderate $K$, they significantly outperform decoupling.  The JS methods offer a much smaller improvement over SAA.  Many of the other features with respect to convergence in $\alpha$ and standard deviation of the performance are also qualitatively similar.  

\subsection{Historical Backtest}
\label{historical}
For our remaining tests we consider a more realistic setting for Shrunken-SAA.  Specifically, we employ {repeated random subsampling validation} with our data to assess each method: for each store we select $\Nhat_k = 10$ days randomly from the dataset, then train each method with these points, and finally evaluate their out-of-sample performance on $N_{\rm test} = 10$ data points, again chosen randomly from the dataset.  Note that unlike the previous experiment, it is possible that some of sampled training days have missing data for store $k$.  In this cases, we will have fewer than $\Nhat_k$ points when training store $k$.  Similar missing data occur for the $N_{\rm test}$ testing points. 
{We prefer repeated, random subsampling validation to more traditional $5$-fold cross-validation \edit{when evaluating our methods}, in order to finely control the number of data points $\Nhat_k$ used in each subproblem.\label{MCValidation}}
  
\begin{figure}[t!]%
\centering%
\begin{minipage}[m]{0.69\textwidth}
	\ifdraft{Figures Goes Here}{
	\input{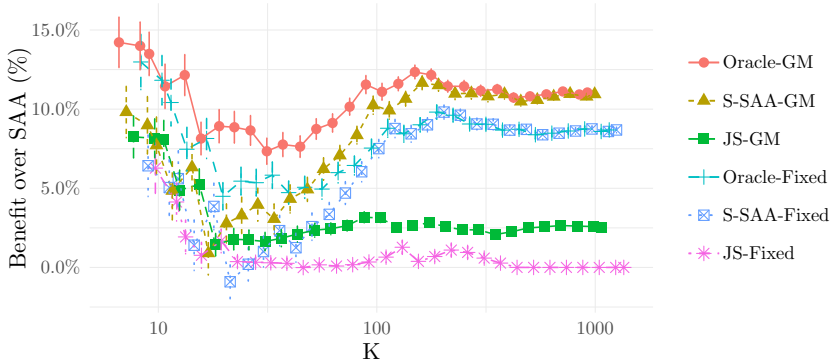}
	}
\end{minipage}%
\begin{minipage}[m]{0.29\textwidth}
    \captionof{figure}{\textbf{Historical Backtest.}
    We evaluate our policies on historical data using $d=20$.  Error bars show $\pm 1$ standard error.  
	}
	\label{fig:BacktestBaseCase}
\end{minipage}%
\vspace{-20pt}
\end{figure}

We evaluate each of our policies using our historical backtest set-up with $d = 20$ in \cref{fig:BacktestBaseCase}.  \edit{For readability, the figure shows a subset of policies.  \Cref{tab:HistDataTrueAllPolicies20} in the appendix shows all policies.} Importantly, we see the same features as in our synthetic data experiment:  our Shrunken-SAA methods converge to oracle optimality and offer a substantive improvement over SAA for large enough $K$. They also outperform JS variants that do not leverage the optimization structure.  

{\blockedit
\subsection{Other Experiments with Synthetic and Real Data}
\label{sec:OtherExperiments}
\Cref{sec:Infinited,sec:LargeSample,sec:KFoldCrossVal} in the appendix study the robustness of Shrunken-SAA to the number of support points $d$, its performance as $N \rightarrow \infty$, and compares computationally cheaper variants of the algorithm that substitute $2\textendash$fold or $5\textendash$fold cross-validation for the LOO validation step.  We omit details for space.  Generally, we find that: i) Shrunken-SAA is quite robust to $d$.  ii) As $N$ increases
Shrunken-SAA retains many of SAA's strong large-sample properties.  Namely, both methods approach full-information optimum, so there is less ``room" to improve upon decoupling, but Shrunken-SAA offers some marginal benefit for large $K$.  iii) Other forms of cross-validation perform quite well and are viable alternatives in computationally limited settings.
}

\section{Conclusion and Future Directions}
\label{sec:Conclusion}
In this paper, we introduce and study the data-pooling phenomenon for stochastic optimization problems, i.e., that when solving many separate data-driven stochastic optimization subproblems, there exist algorithms which pool data across subproblems that outperform decoupling, even when 1) the underlying subproblems are distinct and unrelated, and 2) data for each subproblem are independent.  We propose a simple algorithm Shrunken-SAA that exploits this phenomenon by pooling data in a particular fashion motivated by a Bayes model.  We prove that under frequentist assumptions, in the limit as the number of subproblems grows large, Shrunken-SAA identifies whether pooling in this way can improve upon decoupling, and, if so, the ideal amount to pool,  even if the amount of data per subproblem is fixed and small.  In other words, Shrunken-SAA identifies an optimal level of pooling in the so-called small-data, large-scale regime.  In particular, we prove explicit high-probability bounds on the performance of Shrunken-SAA relative to an oracle benchmark that decay like $\tilde{\mathcal O}(1/\sqrt K)$ where $K$ is the number of subproblems.  

Shrunken-SAA need not offer a strict benefit over decoupling in all instances.  Hence, we also introduce the Sub-Optimality-Instability tradeoff, a decomposition of the benefits of data-pooling that provides strong intuition into the kinds of problems for which data-pooling offers a benefit.  Overall, this intuition and empirical evidence with real data suggest Shrunken-SAA offers significant benefits in the small-data, large-scale regime for  a variety of problems.

We hope our work inspires fellow researchers to think of data-pooling as an ``additional knob" that might be leveraged to improve  performance when designing algorithms for data-driven decision-making under uncertainty.  
\label{ReviewerLastPage}

\enlargethispage{20pt}

\ACKNOWLEDGMENT{The authors would like to thank the editorial team including 3 anonymous reviewers for the constructive comments on an earlier draft.  
Grant Funding: V.G. is partially supported by the National Science Foundation under Grant No. 1661732. N.K. is partially supported by the National Science Foundation under Grant No. 1656996.}

{
\bibliographystyle{ormsv080} %
\setlength{\bibsep}{-3pt plus .3ex}
\bibliography{reference.bib} %
}

\newpage
\ECSwitch
\ECHead{
\begin{center}
$ $\\
Online Appendix:
\vspace{8pt} Data-Pooling for Stochastic Optimization
\end{center}
}

\begin{APPENDICES}
\crefalias{section}{appsec}
\crefalias{subsection}{appsec}
\crefalias{subsubsection}{appsec}

\section{Proof of \cref{thm:DataPoolingJS}: Data-Pooling for MSE}
\label{sec:ProofofThmDataPoolingJS}
\begin{proof}{Proof of Theorem~\ref{thm:DataPoolingJS}.}  
First note that
\begin{align}\notag
\frac{1}{K} \sum_{k=1}^K& \bp_k^\top \bc_k(\bx_k^{\sf SAA}) 
    - \frac{1}{K} \sum_{k=1}^K \bp_k^\top \bc_k(\bx_k(\alphaJS_{\bp_0}, \bp_0, \bfmhat_k ) )
\ - \ 
{\frac{  \left( \frac{1}{K} \sum_{k=1}^K \sigma_k^2 / \Nhat \right)^2 } { \frac{1}{K} \sum_{k=1}^K \sigma_k^2 / \Nhat + \frac{1}{K} \sum_{k=1}^K (\mu_k - \mu_{k0})^2 }}\\\notag
=&
\prns{\frac1K\sum_{k=1}^K\prns{\sigma_k^2+(\mu_k - \muhat_k(0))^2}-
\frac1K\sum_{k=1}^K\prns{\sigma_k^2+(\mu_k - \muhat_k(\alphaJS))^2}}
\\\notag&-\Eb{ \left. \frac1K\sum_{k=1}^K\prns{\sigma_k^2+(\mu_k - \muhat_k(0))^2}-
\frac1K\sum_{k=1}^K\prns{\sigma_k^2+(\mu_k - \muhat_k(\alphaAP))^2} \right| \Nhat}
\\\notag\leq
&
\abs{ \frac{1}{K} \sum_{k=1}^K\prns{ (\mu_k - \muhat_k(0))^2  - \Eb{(\mu_k - \muhat_k(0))^2 \mid \Nhat}   }}
+
\abs{ \frac{1}{K} \sum_{k=1}^K\prns{ (\mu_k - \muhat_k(\alphaJS))^2  - \Eb{(\mu_k - \muhat_k(\alphaAP))^2 \mid \Nhat}   }}
\\\leq&\label{eq:UCofMSE}2\sup_{\alpha\geq0}\abs{ \frac{1}{K} \sum_{k=1}^K\prns{ (\mu_k - \muhat_k(\alpha))^2  - \Eb{(\mu_k - \muhat_k(\alpha))^2 \mid \Nhat}   }}
\\&+\label{eq:JSproofBias}
\abs{ \frac{1}{K} \sum_{k=1}^K \prns{\Eb{(\mu_k - \muhat_k(\alphaJS))^2 \mid \Nhat}  - \Eb{(\mu_k - \muhat_k(\alphaAP))^2 \mid \Nhat}   }}.
\end{align}

We begin by showing \cref{eq:UCofMSE} converges to zero in probability.
Notice Eq.~\eqref{eq:UCofMSE} is the maximal deviation of a stochastic process (indexed by $\alpha$) composed of averages of independent, but not identically distributed, random variables.  Such processes are discussed in \cref{thm:pollard}, and we follow that approach to establish convergence here.  

We first claim that the constants $F_k = 4 a_{\max}^2$ yield an envelope.  Specifically, 
\begin{align*}
\abs{ \mu_k - \muhat_k(\alpha) }  
&\ \leq   \ 
\abs{\bp^\top\ba_k}+\abs{\bphat(\alpha)^\top\ba_k}
\ \leq \ 
2\|  \ba_k \|_\infty.
\end{align*}
which is at most $2 a_{\max}$.
Hence $(\mu_k - \muhat_k(\alpha))^2 \leq F_k$.  %

We next show that the set $\left\{ \Big( (\mu_k - \muhat_k(\alpha))^2 \Big)_{k=1}^K  : \alpha \geq 0 \right\} \subseteq \mathbb R^K$ has pseudo-dimension at most $3$.  Indeed, this set is contained within the set 
\begin{align*}
&\left\{ \Big( \left(\theta(\mu_k - \mu_{k0})  + (1-\theta) (\mu_k - \muhat_k) \right)^2 \Big)_{k=1}^K  : \theta \in \mathbb R \right\}   \ \subseteq \  \R K
\end{align*} 
This set is the range of a quadratic function of $\theta$, and is hence contained within a linear subspace of dimension at most $3$. Thus, it has pseudo-dimension at most $3$.

Since this set has pseudo-dimension at most $3$, there exists a constant $\const_1$ (not depending on $K$ or other problem parameters) such that the corresponding Dudley integral can be bounded as $J \leq \const_1 \| \bm F \|_2$ \citep[pg. 37]{pollard1990empirical}.  \Cref{thm:pollard} with $p = 1$ thus implies there exists a constant $\const_2$ (not depending on $K$ or other problem parameters) such that
\[
\Eb{\sup_{\alpha \geq 0} \abs{ \frac{1}{K} \sum_{k=1}^K \prns{(\mu_k - \muhat_k(\alpha))^2  - \Eb{(\mu_k - \muhat_k(\alpha))^2}}}}
\  \leq \ 
\const_2 \cdot a^2_{\max} /\sqrt{K}.
\]
Markov's inequality then yields the convergence of \cref{eq:UCofMSE} to 0.

We will next show that \cref{eq:JSproofBias} converges to 0.
Let $\theta^{\sf JS} = \frac{\alphaJS}{\alphaJS + \Nhat}$ and $\theta^{\sf AP} = \frac{\alphaAP}{\alphaAP + \Nhat}$ and note $\theta^{\sf JS}, \theta^{\sf AP} \in [0, 1]$ almost surely.  Write, 
\begin{align*}
&\abs{
\frac{1}{K} \sum_{k=1}^K \left( \Eb{ (\mu_k - \muhat_k(\alphaJS) )^2  - (\mu_k - \muhat_k(\alphaAP))^2 \mid \Nhat} \right) 
}
\\
& \qquad \ \leq \ 
\frac{1}{K} \sum_{k=1}^K \Eb{ \abs{ (\mu_k - \muhat_k(\alphaJS) )^2  - (\mu_k - \muhat_k(\alphaAP))^2 }  \mid \Nhat}  
\\& \qquad   = \ 
\frac{1}{K} \sum_{k=1}^K \Eb{ \abs{ (\mu_k - \muhat_k + \theta^{\sf JS} (\muhat_k - \mu_{k0}) )^2  - (\mu_k - \muhat_k + \theta^{\sf AP} (\muhat_k - \mu_{k0}) )^2 }  \mid \Nhat} . 
\end{align*}

Consider the function $\theta \mapsto (\mu_k - \muhat_k + \theta (\muhat_k - \mu_{k0}) )^2$.  For $\abs{ \theta} \leq 1$, its derivative is bounded in magnitude by 
\begin{align*} \hspace{-10pt}
2  \abs{\mu_k - \muhat_k + \theta (\muhat_k - \mu_{k0}) }  \abs{ \muhat_k - \mu_{k0} }
&\ \leq \ 
2  \Big( \abs{\mu_k - \muhat_k} + \abs{\muhat_k - \mu_{k0}} \Big)    \abs{ \muhat_k - \mu_{k0} }
\\&\ \leq \ 
2  \Big( 2 a_{\max} + 2 a_{\max} \Big)    2 a_{\max}
\ = \ 
16 a^2_{\max}.
\end{align*}
Hence, by the mean-value theorem, 
\[
\Eb{\abs{ (\mu_k - \muhat_k + \theta^{\sf JS} (\muhat_k - \mu_{k0}) )^2  - (\mu_k - \muhat_k + \theta^{\sf AP} (\muhat_k - \mu_{k0}) )^2 } \mid \Nhat}
\ \leq \ 
16 a^2_{\max} \Eb{\abs{ \theta^{\sf JS} -  \theta^{\sf AP} } \mid \Nhat}.
\]

We will next show that, conditional on $\Nhat$, $\theta^{JS} - \theta^{\sf AP} \rightarrow_p 0$ as $K\rightarrow \infty$ .  Since $\abs{\theta^{\sf JS} - \theta^{\sf AP}} \leq 2$ almost surely, this will imply that  
\(
\Eb{\abs{ \theta^{\sf JS} -  \theta^{\sf AP} } \mid \Nhat} \rightarrow 0
\)
as $K \rightarrow\infty$, completing the proof.  

Since $\Nhat\geq1,\alphaJS\geq0,\alphaAP\geq0$, we have $\abs{\theta^{\sf JS} - \theta^{\sf AP}}= \frac{\Nhat}{(\alphaJS+\Nhat)(\alphaAP+\Nhat)}\abs{\alphaJS - \alphaAP}\leq \abs{\alphaJS - \alphaAP}$. We proceed to show $\alphaJS \rightarrow_p \alphaAP$.  
We show this second convergence by showing that both the numerator and denominator converge in probability.  For the numerator, 

\[ 0\leq{ \frac{1}{\Nhat - 1} \sum_{i=1}^{\Nhat} ( \xihat_{ki} - \muhat_k)^2 }
\ \leq \
\frac{\Nhat}{\Nhat - 1} 4a_{\max}^2
\ \leq \ 
8 a_{\max}^2,
\]
since $\Nhat \geq 2 \implies \frac{\Nhat}{\Nhat -1 } \leq 2$.  By Hoeffding's inequality, for any $ t > 0$, 
\[
\P\left( \abs{ \frac{1}{K} \sum_{k=1}^K  \frac{1}{\Nhat - 1} \sum_{i=1}^{\Nhat} (\xihat_{ki} - \muhat_k)^2 
		- \frac{1}{K} \sum_{k=1}^K \sigma_k^2 } > t  \ \mid \ \Nhat \right) 
\ \leq\  
		2 \exp\left( - \frac{K t^2 }{ 32 a_{\max}^4} \right) 
\ \rightarrow \ 
0,
\]
as $K \rightarrow \infty$.  
Thus, \(
\frac{1}{K} \sum_{k=1}^K  \frac{1}{\Nhat - 1} \sum_{i=1}^{\Nhat} (\xihat_{ki} - \muhat_k)^2 \rightarrow_p \frac{1}{K} \sum_{k=1}^K \sigma_k^2 .
\)  

Entirely analogously, $0\leq\prns{ \muhat_k - \mu_{k0}}^2 = \prns{(\bphat_k - \bp_0)^\top \ba_k}^2 \leq 4 a^2_{\max}$.  Hence, by Hoeffding's inequality, 
\[
\P \left( \abs{ \frac{1}{K} \sum_{k=1}^K \prns{(\mu_{k0}-\muhat_k)^2  - \Eb{(\mu_{k0}-\muhat_k)^2 }}} > t  \ \mid \ \Nhat \right) \ \leq \ 
2 \exp\left( -\frac{K t^2}{8a_{\max}^2} \right) \ \rightarrow \ 0,
\]
as $K \rightarrow \infty$.  Recall $\Eb{(\mu_{k0}-\muhat_k)^2 \mid \Nhat } = \sigma_k^2/\Nhat + (\mu_{k0} - \mu_k)^2$ by the bias-variance decomposition.  
Combining the numerator and denominator, we have by Slutsky's Theorem that $\alphaJS \rightarrow \alphaAP$.  
\end{proof}

\section{Auxiliary Lemmas}  
\label{sec:AuxLemmas}
\edit{In this section,, we first prove some auxiliary lemmas that we will need when proving our performance guarantees.  These results are largely elementary or well-known facts about tails of random variables.} 
{\blockedit
\begin{lemma}[Bounding a Gaussian Integral]
\label{lemma:integral}
Suppose $t\geq1$. Then
$$
\int_0^1 \sqrt{ \log ( t/ \epsilon) } d\epsilon\leq \sqrt{\log t}+\sqrt{\pi}/2
.$$
\end{lemma}
\proof{Proof.}
Make the substitution $u=\sqrt{2 \log(t/\epsilon)}$.  Then, 
\begin{align*}
\int_0^1 \sqrt{ \log ( t/ \epsilon) } d\epsilon 
& \ = \ 
 \frac{t}{\sqrt 2} \int_{\sqrt{2 \log t}}^\infty u^2  e^{-u^2/2} du 
\\ 
&\ = \ 
  \left.  \frac{tu}{\sqrt 2} e^{-u^2/2} \right|_\infty^{\sqrt{2\log t}}  + \frac{t}{\sqrt 2} \int_{\sqrt{2\log t}}^\infty e^{-u^2/2} du  
&& (\text{integration by parts})
\\
& = 
\sqrt{\log t}  +  \frac{t}{\sqrt 2} \int_{\sqrt{2 \log t}}^\infty   e^{-u^2/2} du .
\end{align*} 
Consider $t \mapsto \frac{t}{\sqrt 2} \int_{\sqrt{2 \log t}}^\infty   e^{-u^2/2} du$. 
Its derivative with respect to $t$ is 
\[
\frac{1}{\sqrt 2} \int_{\sqrt{2 \log t}}^\infty e^{-u^2/2} du -\frac{1}{2 t \sqrt{\log t }} 
\ \leq \ 
\frac{1}{ t \sqrt 2 \cdot \sqrt{2 \log t}}   -\frac{1}{2 t \sqrt{\log t }} 
\ = \ 
0,
\]
where we have a standard inequality for the tail CDF of the normal distribution: $\int_x^{\infty} e^{-u^2/2} du  \leq x^{-1} \cdot e^{-\frac{x^2}{2}}$.  Since the derivative is always non-positive, the integral is non-increasing in $t$.  Thus, 
\[
t \int_{\sqrt{\log t}}^\infty e^{-u^2} du  \ \leq \ 1 \int_{\sqrt{\log 1}}^\infty e^{-u^2} du   =  \frac{\sqrt {\pi} }{2}.
\]
Substituting above completes the proof. 
\endproof  
}

{\blockedit 
\begin{lemma}[$L_p$-norms of Products] \label{lem:NormProduct}
For any $p \geq 1$ and random variables $X$, $Y$.  Then,  
\( 
\| XY \|_p \ \leq \ \| {X} \|_{2p} \| {Y} \|_{2p}.
\)
\end{lemma}
\proof{Proof.}
By H\"older's inequality, 
\(
\E[ \abs{X Y}^p ] \ \leq \ \sqrt{ \E[ X^{2p}] } \cdot \sqrt{ \E[Y^{2p} ] }
\).
Taking the $p^\text{th}$ root of both side yields the result.
\endproof
}
{\blockedit
The following lemma is a specific case of Lemma 2.2.2 of \citet{van1996weak} with explicit constants:
\begin{lemma}[Tails of the Maximum]\label{lem:psi1norm}
Suppose the random variables $Y_1,\dots,Y_K$ satisfy $\E\exp(\beta_0Y_k)\leq 2$ for all $k=1,\dots,K$, $K\geq2$. Let $Y_{\max}=\max_{k=1,\dots,K}Y_k$, and define $\beta = \frac{\beta_0}{1+\log K}$.  Then, $\E\exp(\beta Y_{\max}) \leq 6$.
\end{lemma}
\proof{Proof.}
By definition of $\beta$,
\begin{equation} \label{eq:ImplicationBetaDef}
t \leq \beta Y_{\max}  \iff 1 \leq e^{ \beta_0 Y_{\max}  - t(1+ \log K)  }.
\end{equation}
Then, writing $\exp(\cdot)$ as an integral,
\begin{align*}
\exp( \beta Y_{\max} ) & = e + \int_{1}^{\beta Y_{\max} } e^t dt 
\\
& \leq 
e + \int_{1}^{\beta Y_{\max} } e^{\beta_0 Y_{\max} } \cdot e^{-t ( 1 + \log K )} \cdot e^{t} dt  &&(\text{\cref{eq:ImplicationBetaDef}})
\\
& \leq 
e + \int_{1}^{\beta Y_{\max} } e^{\beta_0 Y_{\max} } \cdot e^{-t \log K } dt 
\\
& \leq
e + \sum_{k=1}^K \int_{1}^{\infty } e^{ \beta_0 Y_{k} } \cdot e^{-t \log K } dt,
\end{align*}
where in the last step we have bounded the maximum by a sum and extended the limits of integration because the integrand is positive.  Now take expectations of both sides and evaluate the integral, yielding
\begin{align*}
\E\left[ \exp( \beta Y_{\max} ) \right]  
\  \leq \ 
e + 2 K \int_1^\infty e^{-t \log K} dt
\ = \ 
 e + \frac{2}{ \log K} \leq 6, 
\end{align*}
since $K \geq 2$.  
\endproof
}

\edit{Recall, for any random variable $Y$ and function $\Psi(\cdot)$, $\|Y \|_\Psi \equiv \inf \left\{ \beta > 0: \E\left[ \Psi\left( \abs{Y}\beta^{-1}  \right) \right] \leq 1 \right\}$ is the Orlicz norm of $Y$ with respect to $\Psi(\cdot)$.}
\begin{lemma}[Relating $\Psi$-norm and $L_p$-norm] \label{lem:RelatingOrlicz4}
Fix $p \geq 1$.  Let $\Psi(t) = \frac{1}{5} \exp(t^2)$, and $\| \cdot \|_\Psi$ be the corresponding Orlicz norm.  Then, 
\begin{enumerate}[label=\roman*)]
    \item For any $t\geq 0$, $t^p \leq \left(\frac{p}{e}\right)^p e^{t}$.  
    \item For any $t \geq 0$, $t^p \leq \left(\frac{p}{2}\right)^{\frac{p}{2}} e^{-\frac{p}{2}}  e^{t^2}$.  
    \item Let $C_p = 5^{1/p} \left( \frac{p}{2} \right)^{1/2} e^{-1/2}$.  For any random variable $Y$, $\| Y \|_p \leq C_p \|Y \|_\Psi$.  
    	\item \label{PsiNormSqrtLogY} 
	{\blockedit 
	    	For any random variable $Y \geq 1$, 
    		$ \| \sqrt{\log Y} \|_p \leq 5^{1/p} \left( \frac{p}{2e} \right)^{1/2}  \max(1, \sqrt{\E[Y]} / 2)$.
	}	
\end{enumerate}
\end{lemma}
\proof{Proof.}
Consider the optimization $\max_{t \geq 0} t^p e^{-t}$.  Taking derivatives shows the optimal solution is $t^* = p$, and the optimal value is $p^p e^{-p}$. Hence, $ t^p e^{-t} \leq p^p e^{-p}$ for all $t$.  Rearranging proves the first statement.  The second follows from the first since, $t^p = \left( t^2 \right)^{\left(\frac{p}{2}\right)} \leq \left(\frac{p}{2}\right)^{p/2} e^{-\frac{p}{2}}  e^{t^2}$.  

For the third, statement, let 
$\beta = \| Y \|_\Psi$, i.e., $\E\left[ \exp\left( \frac{Y^2}{\beta^2} \right) \right] \leq 5$.  Then, 
\[
\E \left[ \left(\frac{ \abs{Y} }{C_p \beta} \right)^p \right] 
\ = \ 
\frac{1}{C_p^p} \E \left[ \left(\frac{ \abs{Y} }{\beta} \right)^p  \right]
\ \leq  \ 
\frac{1}{C_p^p } \left(\frac{p}{2}\right)^{p/2} e^{-\frac{p}{2}}   \E\left[ e^{\frac{Y^2}{\beta^2}} \right]
\ \leq  \  1.
\]
Rearranging and taking the $p^\text{th}$ root of both sides proves the third statement.  

{\blockedit 
Finally, for the last statement, we will first bound
$\| \sqrt{\log Y} \|_\Psi$ where $\Psi(t) = \frac{1}{5} \exp(t^2)$.  
To this end, it suffices to find a $B > 0$ such that 
\begin{align*}
\frac{1}{5} \E\left[ \exp( \log(  Y ) / B^2 )  \right] \leq 1
~\text{ or, equivalently, }~
\E\left[  Y^{1/B^2}  \right] \leq 5.
\end{align*}

We have two possibilities:  Suppose first \(
 \E[Y] \leq 5.
\) Then $B=1$ is feasible above, and so $\| \sqrt{\log( Y) } \|_\Psi \leq 1$. 

On the other hand, suppose 
\(
 \E[Y] > 5.
\)
Consider $\theta = \frac{4}{\E[ Y ] - 1} \in (0, 1)$.  
Then, from convexity of the function $t \mapsto \E[Y^t]$,
\[
\E[Y^\theta] 
\ \leq \ \theta \E[Y^1] + (1-\theta) \E[Y^0] 
\ = \ \theta \E[Y] + (1-\theta)
\ = \ 5.
\]
Thus, if we let $B = \sqrt{\E[Y]}/2$, we have 
\[
\E\left[Y^{1/B^2}\right] \ = \ \E\left[ Y^{4 / \E[Y] } \right] \ \leq  \ \E\left[ Y^{4 / (\E[Y]  - 1) } \right] 
\ = \ 
\E[ Y^{\theta} ] \leq 5.
\]
Hence, $\| \sqrt{\log( Y) } \|_\Psi \leq \sqrt{\E[Y]}/2$.  Combining both cases  proves 
$\| \sqrt{\log( Y) } \|_\Psi 
 \leq \max(1, \sqrt{\E[Y]} / 2)$.  
 
 Apply Part iii) to complete the proof.  
}
\endproof

{\blockedit

\begin{lemma}[Properties of Poisson Random Variables] \label{lem:PropertiesOfPoisson} 
Suppose $\Nhat_k\sim\op{Poisson}(N_k)$, for $k=1, \ldots, K$, where $N_k \geq 1$ for all $k$, and $K \geq 2$. 
Let $\Nhat_{\max} \equiv \max_k \Nhat_k$, $N_{\max} \equiv \max_k N_k$, $\Nhat_{\min} \equiv \min_k \Nhat_k$ and $N_{\min} \equiv \min_k N_k$.  
Then for any $p \geq 1$:
\begin{enumerate}[label=\roman*)]
	\item \label{Psi1Poisson} 
	\( 
	\E\left[ \exp\prns{ \frac{\Nhat_{k}}{2N_k} }\right]  \ \leq \ 2,
	\)  
	\vspace{3pt}
	\item \label{Psi1InversePoisson}
	\(
	\E\left[\exp\prns{ \frac{N_k}{2(\Nhat_{k}+1)} } \right]  \ \leq  \ 2,
	\)
	\vspace{3pt}
	\item \label{Psi1MaxPoisson} 
	\( 
	\E \left[ \exp\left(\frac{ \Nhat_{\max}}{2(1+\log K) N_{\max} }  \right)\right] 
	\ \leq \ 6
	\), 
	\vspace{3pt}
	\item \label{Psi1InverseMinPoisson} 
	\(
	\E \left[ \exp\left( \frac{N_{\min}} {2(1+\log (K)) (\Nhat_{\min} + 1)}  \right) \right] 
	\  \leq \ 6
	\),
	\vspace{3pt}
	\item \label{PNormNhatmax}
	\(  \| \Nhat_{\max}\|_p \ \leq \ 
	6^{1/p} \left( \frac{2p}{e} \right) N_{\max} (1 + \log(K)) 
	\ \leq \
	6^{1/p} \left( \frac{6p}{e} \right) N_{\max} \log(K), 
	\)
	\vspace{3pt}
	\item \label{PNormSqrtRatio}
\(
	        	\magd{\sqrt{\frac{ \Nhat_{\max}}{\Nhat_{\min}+1} }}_p \ \leq  \ 
		6^{1/p} \left(\frac{6p}{e}\right)  \sqrt{\frac{\lambdamax}{\lambdamin} } \cdot \log(K). 
\)
\end{enumerate}
\end{lemma}
\proof{Proof.}\hfill \\
\noindent Part~\ref{Psi1Poisson}
Let $\beta_0 \equiv \log\left(1 + \frac{ \log 2}{N_k} \right)$.  From the Poisson moment generating function,
\[
	\E\left[ \exp(\beta_0 \Nhat_k) \right] 
	\ = \ 
	\exp( N_k ( e^\beta_0 - 1) ) = 2.
\]
Thus, to prove \ref{Psi1Poisson}, it suffices to show that $\beta_0 =  \log\left(1 + \frac{ \log 2}{N_k} \right) \geq \frac{1}{2N_k}$.  
The function $N \mapsto  \log\left(1 + \frac{ \log 2}{N} \right) - \frac{1}{2N}$ is positive at $N =1$ and tends to zero as $N \rightarrow \infty$.  By differentiating, we see it has one critical point at $N = \frac{ \log 2}{2 \log 2 - 1}$ which by inspection is a maximum.  Hence, it is always non-negative, proving the claim and the first statement.  

\vskip 8pt
\noindent Part~\ref{Psi1InversePoisson}
Use the Poisson probability mass function to write
\begin{align*}
\E\left[ \exp\left( \frac{N_k}{2(\Nhat_k +1)} \right) \right]
	& \ = \ e^{-N_k} \ \sum_{n=0}^\infty   \frac{N_k^n}{n!}  \cdot  \exp\left( \frac{N_k}{2(n+1)}\right)
\\
	&\ = \ e^{-N_k} \ \sum_{n=0}^\infty \frac{N_k^n}{n!} \cdot \sum_{j=0}^\infty \prns{\frac{N_k}{2(n+1)}}^j\frac1{j!}
\\
&\ = \ e^{-N_k} \ 
\sum_{j=0}^\infty\frac1{j!}\prns{\frac{N_k}{2}}^j
 \cdot \sum_{n=0}^\infty \frac{N_k^n}{n!} \prns{\frac{1}{n+1}}^j, 
\end{align*}
where the first equality uses the Taylor expansion of $\exp(\cdot)$ and the second from reversing the summations.  
Since 
\(
\frac{1}{n+1} \leq \frac{i}{n+i} \quad \text{for all } n, i \geq 1, 
\) 
we obtain that 
\[
\prns{\frac{1}{n+1}}^j\leq \frac{1}{n+1} \cdot \frac{ 2}{n+2}  \cdots  \frac{j}{n+j} 
\ = \ 
\frac{n!j!}{(n+j)!}.
\]
Substituting above yields 
\begin{align*}
\E\left[ \exp\left( \frac{N_k}{2(\Nhat_k +1)} \right) \right]
& \ \leq  \ e^{-N_k} \ 
\sum_{j=0}^\infty\frac1{j!}\prns{\frac{N_k}{2}}^j
\cdot  \sum_{n=0}^\infty \frac{N_k^n}{n!} \frac{n!j!}{(n+j)!}
\\
&\ = \ e^{-N_k} \ 
\sum_{j=0}^\infty\frac1{2^j}
\cdot \sum_{n=0}^\infty\frac{N_k^{n+j}}{(n+j)!}
\\
& \ = \ e^{-N_k} \ 
\sum_{j=0}^\infty\frac1{2^j} \cdot \sum_{n=j}^\infty\frac{N_k^n}{n!}
\\
& \ \leq  \ e^{-N_k}
\sum_{j=0}^\infty\frac1{2^j} \cdot \sum_{n=0}^\infty\frac{N_k^n}{n!}
\\& \ = \ 2.
\end{align*}

\vskip 8pt
\noindent Parts~\ref{Psi1MaxPoisson} and \ref{Psi1InverseMinPoisson}
These results follow by combining \cref{lem:psi1norm} with parts \ref{Psi1Poisson} and \ref{Psi1InversePoisson} respectively. 

\vskip 8pt
\noindent Part~\ref{PNormNhatmax}
Let $\beta = \frac{1}{2(1+ \log K ) N_{\max} }$.  Then, from \cref{lem:RelatingOrlicz4} Part i),
\begin{align*}
\E[ \Nhat_{\max}^p] 
\ = \
\beta^{-p}  \E[ ( \beta \Nhat_{\max} )^p ] 
\ \leq  \ 
\beta^{-p} \left( \frac{p}{e} \right)^p \E[ \exp(\beta\Nhat_{\max})] 
\ \leq \ 
6 \left( \frac{2p}{e} \right)^p N_{\max}^p (1 + \log K)^p,
\end{align*}
where the second inequality uses Part~\ref{Psi1MaxPoisson}.  Taking the $p^\text{th}$ root of both sides proves the first statement. The second follows because $K \geq 2$ implies that $1 + \log K \leq 3 \log K$.  

\vskip 8pt
\noindent Part~\ref{PNormSqrtRatio}
Applying an identical argument to the previous part but with Part~\ref{Psi1InverseMinPoisson} , we have 
\begin{align*}
	\E\left[ (\Nhat_{\min} + 1 )^{-p} \right]  \ \leq \ 6  \left(\frac{2p}{eN_{\min}}\right)^p   (1 + \log K )^p.
\end{align*}
Therefore, we have
\begin{align*}
	\magd{\sqrt{  \frac{ \Nhat_{\max} }{ \Nhat_{\min}+1} } }_p^p
	&=
	\E\left[ \prns{ \frac{ \Nhat_{\max} } { \Nhat_{\min}+1 } }^{p/2}\right] 
	\\
	&\leq
	\sqrt{\E\Nhat^p_{\max}} \cdot \sqrt{\E\left[ (\Nhat_{\min}+1)^{-p} \right]}  &&(\text{Cauchy-Schwarz Inequality})
	\\
	&\leq
	6\cdot2^{p}\left( \log(K) + 1 \right))^p e^{-p}p^p\prns{\frac{N_{\max}}{N_{\min}}}^{p/2}
	\\
	& = 
	6\cdot  \left(\frac{2p }{e}\right)^{p} \left( \log(K) + 1 \right)^p \prns{\frac{\lambdamax}{\lambdamin} }^{p/2}
	.
\end{align*}
Finally, since $K \geq 2$, we have $1 + \log K \leq 3 \log K$.  Making this substitution and simplifying completes the proof. 
\endproof

}
%

\section{Deferred Proofs for Sub-Optimality Guarantees from \cref{sec:PerformanceGuarantees}}

In this section, we provide the complete proofs for the high-probability sub-optimality bounds presented in \cref{sec:PerformanceGuarantees}.

\subsection{Proof of \cref{thm:FixedPointShrinkage}: Shrunken-SAA with Fixed Anchors for Strongly Convex Problems}\label{sec:FixedPointShrinkageAppendix}

We first prove the results summarized in \cref{sec:SmoothCostsFixedAnchor}.

\subsubsection{Proof of continuity lemma and packing number bounds}
As mentioned in the main text, the key idea is to establish continuity of the solutions $\bx_k(\alpha, \bp_0, \bfmhat_k)$ in the parameters.  
{\blockedit
\begin{lemma}[Continuity properties of $\bx_k(\alpha, \bp_0, \bfmhat_k)$] \label{lem:ContinuityAlpha}
Under the assumptions of \cref{thm:FixedPointShrinkage}, 
\begin{enumerate}[label=\roman*)]
\item \label{part:ContinuityinP0} (Continuity in anchor) For any $\alpha \geq 0$, and any $\bp, \overline \bp \in \Delta_d$,
\[ 
\|  \bx_k(\alpha, \bp,  \bfmhat_k) -  \bx_k(\alpha, \overline \bp,  \bfmhat_k) \|_2  \ \leq  \  \frac{L }{\gamma}  \cdot{ \| \bp - \overline \bp \|_1 }.
\]
\item \label{part:contLOO} (Continuity in $ \bfmhat_k$)  For any $ \bfmhat_k$ such that $\Nhat_k  \geq 1$ we have 
\[
\|  \bx_k (\alpha, \bp_0, \bfmhat_k) -  \bx_k(\alpha, \bp_0, \bfmhat_k - \be_i) \|_2 
 \ \leq \ 
\frac{ 4 L }{ \gamma \Nhat_k} .
 \]
\item \label{part:contAlpha} (Continuity in $\alpha$) 
For any $\alpha,\, \overline\alpha \geq 0$,  and $\bp_0\in\Delta_d$,
\[
\|  \bx_k(\alpha, \bp_0, \bfmhat_k) -  \bx_k(\overline\alpha, \bp_0, \bfmhat_k) \|_2 
\ \leq \ 
\frac{4L}{\gamma}  \cdot \frac{\abs{\alpha - \overline \alpha}}{\Nhat_k + 1}.
\]

\item \label{part:contInfinity} (Continuity at $\alpha = \infty$) 
For any $\alpha \geq 0$ and $\bp_0 \in \Delta_d$  such that $\max(\alpha, \Nhat_k) > 0$, 
\[
	\| \bx_k(\alpha, \bp_0, \bfmhat_k)  - \bx_k(\infty, \bp_0, \bfmhat_k) \|_2 
	\ \leq  \ 
	\frac{2L}{\gamma} \frac{\Nhat_k}{\Nhat_k + \alpha}.
\]
\end{enumerate}
\end{lemma}
}
\proof{Proof.}%
Fix $k$.  
For any $\bq \in \Delta_d$, define 
\[
f_{\bq}(\bx) \equiv \bq^\top \bc_k(\bx), \quad \quad \bx(\bq) \in \arg\min_{\bx \in \X_k} f_{\bq} (\bx).
\]
We first prove the general inequality for any $\bq, \overline \bq \in \Delta_d$, 
\begin{equation} \label{eq:GeneralContinuityIneq}
\blockedit
\| \bx( \bq) - \bx(\overline \bq ) \|_2 \leq { \frac{L}{\gamma} }  \cdot { \| \bq - \overline \bq \|_1 }.
\end{equation}
We will then use this general purpose inequality to prove the various parts of the lemma by choosing particular values for $\bq$ and $\overline \bq$.  

Note that since each $c_{ki}(\bx)$ is $\gamma$-strongly convex for each $i$, $f_{\bq}(\bx)$ is also $\gamma$-strongly convex.  From the first-order optimality conditions, 
\(
\nabla f_{\bq}( \bx (\bq) )^\top \left( \bx(\overline \bq ) - \bx(\bq ) \right) \geq 0.
\)
Then, from strong-convexity, 
\begin{align*}
f_{\bq}(\bx(\overline \bq) ) - f_{\bq}( \bx (\bq ) ) 
& \ \geq \ 
\nabla f_{\bq}( \bx (\bq) )^\top \left( \bx(\overline \bq ) - \bx(\bq ) \right)  + 
\frac{\gamma}{2} \| \bx(\bq) - \bx (\overline \bq)  \|_2^2
\\
&\  \geq  \ 
\frac{\gamma}{2} \| \bx(\bq) - \bx (\overline \bq)  \|_2^2.
\end{align*}
A symmetric argument holds switching $\bq$ and $\overline \bq$ yielding 
\[
f_{\overline \bq}(\bx(\bq) ) - f_{\overline \bq}( \bx (\overline \bq ) ) \geq \frac{\gamma}{2} \| \bx(\bq) - \bx (\overline \bq)  \|_2^2.
\]
Adding yields, 
\begin{align*}
\gamma \| \bx(\bq) - \bx (\overline \bq)  \|_2^2 
& \ \leq \ 
\Big( f_{\overline \bq}(\bx(\bq) ) - f_{\bq}( \bx (\bq ) )  \Big) + 
\Big( f_{\bq}(\bx(\overline \bq) ) - f_{\overline \bq}( \bx (\overline \bq ) ) \Big)
\\
& \ = \ 
\left( \overline \bq - \bq \right)^\top \big( \bc_k(\bx (\bq) ) - \bc_k(\bx (\overline \bq ) ) \big)
\\
&\ \leq \ 
\big\| \bc_k(\bx (\bq) ) - \bc_k(\bx (\overline \bq ) ) \big\|_\infty
\| \bq - \overline \bq \|_1
\\
&\ \leq \ 
L\| \bx(\bq) - \bx (\overline \bq)  \|_2
\| \bq - \overline \bq \|_1,
\end{align*}
\edit{by the H\"older inequality and assumed Lipschitz constant.}  Rearranging proves \cref{eq:GeneralContinuityIneq}.  

We can now prove each part of the lemma.  

\noindent Part~\ref{part:ContinuityinP0}  
First suppose $\alpha + \Nhat_k > 0$.  Take 
\begin{align*}
\bq  = \frac{\alpha}{\alpha + \Nhat_k} \bp + \frac{\Nhat_k}{\Nhat_k + \alpha} \bphat_k, 
\quad \text{ and } \quad
\overline \bq  = \frac{\alpha}{\alpha + \Nhat_k} {\overline \bp} + \frac{\Nhat_k}{\Nhat_k + \alpha} \bphat_k .
\end{align*}
Then, 
\(
\| \bq - \overline \bq \|_1 
 \ = \ 
\frac{\alpha}{\Nhat_k + \alpha} \| \bp - \overline \bp \|_1 
\ \leq \ \| \bp - \overline \bp \|_1.
\)
Substituting into \cref{eq:GeneralContinuityIneq} proves the result in this case.  \edit{Next, suppose $\alpha + \Nhat_k = 0$.  Then, applying \cref{eq:GeneralContinuityIneq}  with $\bq = \bp$ and $\overline \bq= \overline \bp$ yields the result.}  
\vskip 12pt
{\blockedit 
\noindent Part~\ref{part:contLOO}
First suppose $\Nhat_k \geq 2$.  Take
\begin{align*}
\bq = \frac{\alpha}{\Nhat_k + \alpha} \bp_0  + \frac{1}{\Nhat_k + \alpha} \bfmhat_k
\quad \text{ and } \quad 
\overline \bq = \frac{\alpha}{\Nhat_k + \alpha - 1} \bp_0 + \frac{1}{\Nhat_k + \alpha - 1} (\bfmhat_k - \be_i).
\end{align*}
Then, 
\begin{align*}
\| \bq - \overline \bq \|_1 
& \ \leq \ 
\abs{ \frac{\alpha}{\Nhat_k + \alpha} - \frac{\alpha}{\Nhat_k + \alpha - 1} }  \| \bp_0 \|_1  \ + \ \abs{\frac{1}{\Nhat_k + \alpha} - \frac{1}{\Nhat_k + \alpha - 1} } \| \bfmhat_k \|_1 + \frac{1}{\Nhat_k + \alpha - 1}
\\
& \ = \ 
\frac{2}{\Nhat_k -1 + \alpha}
\\ 
& \ \leq \ 
\frac{4}{\Nhat_k}.
\end{align*}
Substituting into \cref{eq:GeneralContinuityIneq} proves the result when $\Nhat_k \geq 2$.  

Next, when $\Nhat_k = 1$, let $\bq$ be as above and $\overline \bq = \bp_0$.  Then $\| \bq - \overline \bq \|_1 \leq 2 \leq \frac{4}{\Nhat_k}$.  Again, substituting into \cref{eq:GeneralContinuityIneq} proves the result.  
}
\vskip 12 pt
\noindent Part~\ref{part:contAlpha} 
Notice if $\Nhat_k = 0$, then $\|  \bx_k(\alpha, \bp_0, \bfmhat_k) -  \bx_k(\overline\alpha, \bp_0, \bfmhat_k) \|_2  = 0$ and the bounds holds trivially.  Hence, suppose $\Nhat_k \geq 1$.  Consider taking $\bq = \bphat_k(\alpha)$ and $\overline \bq = \bphat_k(\overline\alpha)$.  Then 
\begin{align*}
\| \bq - \overline \bq\|_1 
& \ = \ \left \| \left(\left(\frac{\alpha}{\Nhat_k + \alpha} - \frac{\overline\alpha}{\Nhat_k + \overline\alpha} \right) \bp_0  
+ \left(\frac{\Nhat_k}{\Nhat_k + \alpha} - \frac{\Nhat_k}{\Nhat_k + \overline\alpha} \right)  \bphat_k \right) \right\|_1
\\ 
& \ \leq \ \left( \abs{ \frac{\alpha}{\Nhat_k + \alpha} - \frac{\overline\alpha}{\Nhat_k + \overline\alpha} } + \abs{ \frac{\Nhat_k}{\Nhat_k + \alpha} - \frac{\Nhat_k}{\Nhat_k + \overline\alpha} } \right)
\\ \label{eq:WeightDiffs}
& = 2 \abs{ \frac{\Nhat_k}{\Nhat_k + \alpha} - \frac{\Nhat_k}{\Nhat_k + \overline\alpha} }, 
\\
& = \frac{ 2\Nhat_k \abs{\alpha - \overline \alpha}}{ (\Nhat_k + \alpha) (\Nhat_k + \overline \alpha) } ,
\end{align*}
where second equality follows because $\abs{ \frac{\alpha}{\Nhat_k + \alpha} - \frac{\overline \alpha}{\Nhat_k + \overline \alpha} } = \abs{ \frac{\Nhat_k}{\Nhat_k + \alpha} - \frac{\Nhat_k}{\Nhat_k + \overline \alpha} }$.  Next write, 
\[
\frac{ 2\Nhat_k \abs{\alpha - \overline \alpha}}{ (\Nhat_k + \alpha) (\Nhat_k + \overline \alpha) }
\ \leq \ 
\frac{ 2 \abs{\alpha - \overline \alpha}}{(\Nhat_k + \overline \alpha) }
\ \leq \ 
\frac{ 2 \abs{\alpha - \overline \alpha}}{\Nhat_k }
\ \leq \ 
\frac{ 4 \abs{\alpha - \overline \alpha}}{(\Nhat_k + 1)},
\]
where the last inequality follows because $\frac{1}{N} \leq \frac{2}{N+1}$ for $N \geq 1$.  

Substituting into \cref{eq:GeneralContinuityIneq} completes the proof of part~\ref{part:contAlpha}.

{\blockedit 
 \vskip 12pt
 \noindent Part~\ref{part:contInfinity} Take $\bq = \bp_0$ and $\overline \bq = \bphat_k(\alpha)$.  
 Then,
 \begin{align*}
 \| \bq - \overline \bq \|_1 
 & \ = \ 
 \left\| \left(1 - \frac{\alpha}{\Nhat_k + \alpha} \right) \bp_0  
 + \left(0 - \frac{\Nhat_k}{\Nhat_k + \alpha} \right)  \bphat_k \right\|_1
 \\ & \ \leq \
 \abs{ 1 - \frac{\alpha}{\Nhat_k + \alpha} } + \abs{ 0 - \frac{\Nhat_k}{\Nhat_k + \alpha} }
 \\
 &\  = \ 
 2\frac{\Nhat_k}{\Nhat_k + \alpha}.
 \end{align*}
 Again, substituting into \cref{eq:GeneralContinuityIneq} proves the inequality. 
\endproof 
}

{
{ \blockedit\begin{lemma}[Packing Numbers for Strongly-Convex Problems] \label{lem:PackSmooth}
Under the assumptions of \cref{thm:FixedPointShrinkage}, we have
for any $0 < \epsilon \leq 1$,
\begin{align} 
\label{eq:PackSmoothTrue}
D\left( \epsilon \| \mathbf F^{\sf Perf} \|_2, \ \Fperf \right) 
& \ \leq \ 
2 + 
    \frac{\Nhat_{\max}}{\Nhat_{\min} + 1}  \frac{32 L^4}{\Cmax^2 \gamma^2 \epsilon^2},
\\ 
\label{eq:PackSmoothLOO}
D\left(\epsilon \| \mathbf F^{\sf LOO} \|_2, \ \Floo \right) 
& \ \leq \ 
2 + 
    \frac{\Nhat_{\max}}{\Nhat_{\min} + 1}  \frac{32 L^4}{\Cmax^2 \gamma^2 \epsilon^2}.
\end{align}
\end{lemma}}
\proof{Proof.}%
We first prove \cref{eq:PackSmoothTrue}.
We proceed by constructing an $\frac{\epsilon}{2} \| \mathbf F^{\sf Perf} \|_2 $-covering. 
The desired packing number is at most the size of this covering.  Recall, by \cref{lem:Envelopes}, $\| \mathbf F^{\sf Perf}\|_2^2 = \frac{\Cmax^2}{\lambdabar^2} \| \bm{\lambda} \|_2^2$, and let $Z_k(\infty,\bp_0) = \frac{1}{\lambdabar} \sum_{i=1}^d \lambda_kp_{ki} c_{ki}( \bx_k(\infty,\bp_0) )$.

First, suppose \edit{$\Nhat_{\max} = 0$}, which implies $\Nhat_k = 0$ for all $k =1, \ldots, K$.  In this case, $ \bx_k(\alpha,\bp_0,  \bfmhat_k) =  \bx_k(\infty,\bp_0)$ for all $k$, whereby $\{ \mathbf Z(\alpha,\bp_0) : \alpha \geq 0\} = \{ \mathbf Z(\infty,\bp_0) \}$, and the covering number is $1$, so the above bound is valid.  

Now suppose $\Nhat_{\max} > 0$.  
Let  \edit{$\alpha_{\max} = \frac{4L^2\Nhat_{\max}}{\Cmax \gamma \epsilon}$}. \edit{For any $\alpha \geq \alpha_{\max} > 0$,}
\begin{align*}
\abs{Z_k(\alpha,\bp_0)  - Z_k(\infty,\bp_0) } 
&\ \leq \ \frac{\lambda_k}{\lambdabar}\sum_{i=1}^dp_{ki} \abs{c_{ki}(\bx_k(\alpha, \bfmhat)) - c_{ki}(\bx_k(\infty) ) }
\\
&\ \leq \ \frac{\lambda_k}{\lambdabar}
\sum_{i=1}^dp_{ki} L \| \bx_k(\alpha, \bfmhat)) - \bx_k(\infty) \|_2 && (\text{Lipschitz continuity})
\\
&\  \leq \ \edit{\frac{\lambda_k}{\lambdabar}
 {\frac{2L^2}{\gamma}} \cdot { \frac{ \Nhat_k}{\Nhat_k + \alpha} }} &&(\text{\cref{lem:ContinuityAlpha}, part~\ref{part:contInfinity} since } \alpha > 0).
\end{align*}
\edit{It follows that for all $\alpha\geq\alpha_{\max}$ we have}
{\blockedit\begin{align*}
\| \mathbf Z(\alpha,\bp_0) - \mathbf Z(\infty,\bp_0) \|_2
 \ &\leq  \ 
    \prns{\frac{4L^4}{\lambdabar^2\gamma^2} \sum_{k=1}^K \lambda_k^2\prns{\frac{\Nhat_k}{\Nhat_k + \alpha}}^2}^{1/2}
\\ \ &\leq \ 
	\frac{2L^2\|\lambdavec\|_2}{\lambdabar\gamma} \prns{\frac{\Nhat_{\max}}{\Nhat_{\max} + \alpha}}
\\ \ &\leq \ 
	\frac{2L^2\|\lambdavec\|_2}{\lambdabar\gamma} \prns{\frac{\Nhat_{\max}}{\Nhat_{\max} + \alpha_{\max}}}
\\ \ &\leq \ 
	\frac{2L^2\|\lambdavec\|_2}{\lambdabar\gamma} \prns{\frac{1}{1 + \frac{4L^2}{C\gamma\epsilon}}}
\\ \ &\leq \ 
	\frac{2L^2\|\lambdavec\|_2}{\lambdabar\gamma} \cdot{\frac{C\gamma\epsilon}{4L^2}}
\\ \ &= \ 
	\frac{\epsilon}{2} \| \mathbf F^{\sf Perf} \|.
\end{align*}}
Thus, in our covering, we place one point at $\mathbf Z(\infty,\bp_0)$ to cover all points $\mathbf Z(\alpha,\bp_0)$ with $\alpha \geq \alpha_{\max}$.  

Next let $\{ \alpha_1, \ldots, \alpha_M\}$ be a \edit{$\frac{\gamma(\Nhat_{\min} + 1) \Cmax \epsilon}{8 L^2}$} covering of $[0, \alpha_{\max}]$.  Note, \edit{$M \leq 1+
\frac{8 L^2 \alpha_{\max}}{\gamma(\Nhat_{\min} + 1) \Cmax \epsilon}$}.  We claim $\{ \mathbf Z(\alpha,\bp_0), \ldots, \mathbf Z(\alpha_M,\bp_0) \}$ is an \edit{$\frac{\epsilon}{2} \| \mathbf F^{\sf Perf} \|$}-covering of $\{ \mathbf Z(\alpha,\bp_0) : \alpha \in [0, \alpha_{\max}] \}$.  Indeed, for any $\alpha \in [0, \alpha_{\max}]$, let $\alpha_j$ be the nearest element of the $\alpha$-covering.  Then,  
\begin{align*}
\abs{ Z_k(\alpha,\bp_0) - Z_k(\alpha_j,\bp_0) } 
& \ \leq  \ 
    \frac{\lambda_k}{\lambdabar}\sum_{i=1}^d p_{ki} \abs{ c_{ki}( \bx_k(\alpha,\bp_0,  \bfmhat_j) ) - c_{ki} (  \bx_k(\alpha_j,\bp_0,  \bfmhat_k))} 
\\
& \ \leq \
    \frac{\lambda_k}{\lambdabar}\sum_{i=1}^d p_{ki} L \|  \bx_k(\alpha,\bp_0,  \bfmhat_j) -  \bx_k(\alpha_j,\bp_0,  \bfmhat_k) \|_2
\\
& \ \leq \
    \edit{\frac{\lambda_k}{\lambdabar}  { \frac{ 4L^2}{\gamma(\Nhat_{\min} + 1)} } { \abs{\alpha - \alpha_j }}}
    && (\text{\cref{lem:ContinuityAlpha}, part~\ref{part:contAlpha}})
\\
& \ \leq \ 
        \edit{\frac{\lambda_k}{\lambdabar}  { \frac{ 4L^2}{\gamma(\Nhat_{\min} + 1)} } \cdot \frac{\gamma(\Nhat_{\min} + 1) \Cmax \epsilon}{8 L^2}}
\\ &
\ = \  \frac{ \Cmax \epsilon \lambda_k}{2 \lambdabar}
\end{align*}
Thus, $\edit{\| \mathbf Z(\alpha,\bp_0) -\mathbf Z(\alpha_j,\bp_0) \|_2 \leq \frac{\Cmax  \epsilon\|\lambdavec\|_2 }{2\lambdabar}=\frac{\epsilon}{2} \| \mathbf F^{\sf Perf} \|}$ as was to be shown.   

The total size of the covering is thus 
\[
1 + M 
\ \leq \ \edit{
    2 + \frac{8 L^2 \alpha_{\max}}{\gamma(1+\Nhat_{\min}) \Cmax \epsilon}
\ =\ 
    2 + 
    \frac{\Nhat_{\max}}{1+\Nhat_{\min}}  \frac{32 L^4}{\Cmax^2 \gamma^2 \epsilon^2}.
    }
\]

We next prove \cref{eq:PackSmoothLOO}.
We again proceed by constructing an \edit{$\frac{\epsilon   }{2} \| \mathbf F^{\sf LOO} \|$}-covering, since the desired packing is at most the size of this covering.  Recall by \cref{lem:Envelopes}, $\| \mathbf F^{\sf LOO}\|_2^2 = \frac{\Cmax^2}{N^2\lambdabar^2} \| \bm{\Nhat} \|^2_2$. 

If $\Nhat_{\max} = 0$, then $\Nhat_k = 0$ for all $k$, and $\{ \mathbf Z^{\sf LOO}(\alpha,\bp_0) : \alpha \geq 0 \} = \{ \bm 0 \}$, so this covering number is $1$. 

Otherwise, $\Nhat_{\max} > 0$.  
Let \edit{$\alpha_{\max} =  \frac{4\Nhat_{\max}L^2}{\Cmax \gamma \epsilon}$}.
\edit{Then, for any $\alpha \geq \alpha_{\max} > 0$,}
\begin{align*}
\abs{ Z_k^{\sf LOO}(\alpha,\bp_0) - Z_k^{\sf LOO}(\infty,\bp_0) } 
& \leq 
\frac{1}{N\lambdabar} \sum_{i=1}^d \mhat_{ki}\abs{ c_{ki}(\bx_k(\alpha,  \bfmhat_k - \be_i) ) - c_{ki}(\bx_k(\infty) ) }
\\
&\leq 
\frac{L}{N\lambdabar} \sum_{i=1}^d \mhat_{ki}  \| \bx_k(\alpha,  \bfmhat_k - \be_i) - \bx_k(\infty) \|_2 && \text{(Lipschitz-Continuity)}
\\
& \leq 
\edit{\frac{L}{N\lambdabar} \sum_{i=1}^d \mhat_{ki}  {\frac{ 2L }{\gamma} } { \frac{\Nhat_k -1 }{\Nhat_k -1 + \alpha}}}  && \text{(\cref{lem:ContinuityAlpha}, part~\ref{part:contInfinity}} )
\\ & \leq 
\edit{\frac{2L^2\Nhat_k}{\gamma N\lambdabar} { \frac{\Nhat_k }{\Nhat_k + \alpha}}},  
\end{align*}
because $x \mapsto \frac{x}{x+ \alpha}$ is an increasing function.
Thus, \edit{for any $\alpha \geq \alpha_{\max}$}, 
\begin{align*}
\| \mathbf Z_k^{\sf LOO}(\alpha,\bp_0) - \mathbf Z_k^{\sf LOO}(\infty,\bp_0) \|_2
& \  \leq \ 
\frac{ 2 L^2 }{\gamma} \prns{\sum_{k=1}^K \frac{\Nhat_k^2}{N^2\lambdabar^2} \cdot \prns{\frac{\Nhat_k}{\Nhat_k + \alpha}}^2}^{1/2}
\\
&  \ \leq \ 
\frac{ 2 L^2 }{\gamma C} \prns{\sum_{k=1}^K \frac{C^2\Nhat_k^2}{N^2\lambdabar^2}}^{1/2} \cdot \frac{\Nhat_{\max}}{\Nhat_{\max} + \alpha}
\\
& \ =\ 
\frac{ 2 L^2 }{\gamma C} \| \mathbf F^{\sf LOO} \|_2 \  \frac{\Nhat_{\max}}{\Nhat_{\max} + \alpha}
\\
& \ \leq\ 
\frac{ 2 L^2 }{\gamma C} \| \mathbf F^{\sf LOO} \|_2 \  \frac{\Nhat_{\max}}{\Nhat_{\max} + \alpha_{\max}}
\\
& \ =\ 
\frac{ 2 L^2 }{\gamma C} \| \mathbf F^{\sf LOO} \|_2  \ \frac{1}{1 + \frac{4L^2}{\Cmax \gamma \epsilon}}
\\
& \ \leq\ 
\frac{ 2 L^2 }{\gamma C} \| \mathbf F^{\sf LOO} \|_2  \ \frac{\Cmax \gamma \epsilon}{4L^2}
\\
& \ =\ \frac{\epsilon}{2}\| \mathbf F^{\sf LOO} \|_2
\end{align*}
Thus, in our covering, we place one point at $\mathbf Z^{\sf LOO}(\infty,\bp_0)$ to cover all points $\mathbf Z^{\sf LOO}(\alpha,\bp_0)$ for $\alpha \geq \alpha_{\max}$.  

Next let $\{ \alpha_1, \ldots, \alpha_M \}$ be a $\edit{\frac{\gamma(\Nhat_{\min} + 1)\Cmax\epsilon}{8 L^2}}$-covering of $[0, \alpha_{\max}]$.  Note, \edit{$M \leq 
1+\frac{8 L^2 \alpha_{\max}}{\gamma(\Nhat_{\min} + 1) \Cmax \epsilon}$}.  We claim this covering induces an $\frac{\epsilon}{2} \| \mathbf F^{\sf LOO}\|_2 $-covering of $\{ \mathbf Z^{\sf LOO}(\alpha,\bp_0) : \alpha \in [0, \alpha_{\max}] \}$.  Indeed, for any $\alpha \in [0, \alpha_{\max}]$, let $\alpha_j$ be the nearest element of the $\alpha$-covering.  Then, for any $k$ such that $\Nhat_k \geq 1$, 
\begin{align*}
\Big|Z^{\sf LOO}_k(\alpha,\bp_0) &- Z^{\sf LOO}_k(\alpha_j,\bp_0) \Big| 
\\
&\leq 
\frac{1}{N\lambdabar} \sum_{i=1}^d \mhat_{ki} \abs{ c_{ki}(\bx_k(\alpha,   \bfmhat_{ki} - \be_i) ) - c_{ki}(\bx_k(\alpha_j,   \bfmhat_{ki} - \be_i) ) } 
\\ &\leq
\frac{L}{N\lambdabar} \sum_{i=1}^d \mhat_{ki} \| \bx_k(\alpha,   \bfmhat_{ki} - \be_i) ) - \bx_k(\alpha_j,   \bfmhat_{ki} - \be_i) \|_2 &&\text{(Lipschitz Continuity)} 
\\
& \leq
\edit{\frac{\Nhat_k }{N\lambdabar} \cdot {\frac{4L^2}{\gamma(\Nhat_{\min} + 1)} } \cdot { \abs{\alpha - \alpha_j } }}&&\text{(\cref{lem:ContinuityAlpha}, part~\ref{part:contAlpha}} )
\\
& \leq
\edit{\frac{\Nhat_k }{N\lambdabar} \cdot {\frac{4L^2}{\gamma(\Nhat_{\min} + 1)} } \cdot \frac{\gamma(\Nhat_{\min} + 1) \Cmax\epsilon}{8L^2}}
\\
&= 
\Cmax \frac{\Nhat_k }{N\lambdabar} \frac{\epsilon}{2} .
\end{align*}
On the other hand, for any $k$ such that  $\Nhat_k = 0$, $\abs{Z^{\sf LOO}_k(\alpha,\bp_0) - Z^{\sf LOO}_k(\alpha_j,\bp_0)}  = 0$.  
In total, this implies 
\(
\| \mathbf Z^{\sf LOO}(\alpha,\bp_0) - \mathbf Z^{\sf LOO}(\alpha_j,\bp_0) \|^2_2 \leq  \frac{\epsilon^2}{4} \frac{\Cmax^2}{N^2\lambdabar^2} \| \bm{\Nhat} \|_2^2,
\)
 which implies 
$ \| \mathbf Z^{\sf LOO}(\alpha,\bp_0) - \mathbf Z^{\sf LOO}(\alpha_j,\bp_0) \| \leq \frac{\epsilon}{2}  \| \mathbf F^{\sf LOO} \|_2$, as was to be proven.    

Thus, the total size of the covering is at most 
\[
1 + M 
\ \leq \ \edit{
    2 + \frac{8 L^2 \alpha_{\max}}{\gamma(1+\Nhat_{\min}) \Cmax \epsilon}
\ =\ 
    2 + 
    \frac{\Nhat_{\max}}{1+\Nhat_{\min}}  \frac{32 L^4}{\Cmax^2 \gamma^2 \epsilon^2}.
    }\]
This completes the proof.  \endproof
}

\subsubsection{Maximal deviation bounds.}
We next use the above lemmas to bound the maximal deviations of interest via \cref{thm:pollard}.  

\begin{lemma}[Bounding the Maximal Deviations] \label{lem:UCDeviations}
Suppose $\frac{4L^2}{C\gamma} \geq 1$.  Then, under the assumptions of \cref{thm:FixedPointShrinkage}, there exists a universal constant $\const$ such that for any $0 < \delta< 1/2$, the following two statements each hold (separately) with probability at least $1-\delta$:
\begin{align*}
 \sup_{\alpha \geq 0 } \abs{\frac{1}{K }\sum_{k=1}^K \prns{Z_k(\alpha,\bp_0) - \E[Z_k(\alpha,\bp_0) ]} } 
&\leq  \edit{
\   
\const \cdot 
L \sqrt{\frac{\Cmax}{\gamma}} 
  \cdot \left(\frac{\lambdamax}{\lambdamin} \right)^{5/4}  \cdot \frac{ \log(1/\delta) \cdot \sqrt{\log(K)}}{\sqrt K } 
, }\\
 \sup_{\alpha \geq 0 } \abs{\frac{1}{K }\sum_{k=1}^K \prns{Z^{\sf LOO}_k(\alpha,\bp_0) - \E[Z^{\sf LOO}_k(\alpha,\bp_0) ] } } 
&\leq  \edit{
\  
\const \cdot 
L \sqrt{\frac{\Cmax}{\gamma}} \cdot 
\left( \frac{\lambdamax}{\lambdamin}\right)^{5/4} \cdot
\frac{ \log^2 (1/\delta)\cdot  \log^{3/2}(K)}{ \sqrt K }
.}
\end{align*}
\end{lemma}
\proof{Proof.}
{  \blockedit
To prove the first inequality, our strategy will be to apply \cref{thm:pollard} to the process $\Fperf$.  
To that end, we first bound the variable $J$ in \cref{eq:JDudley}.  Recall by 
\cref{lem:Envelopes}, the size of the envelope is at most $C\frac{\| \bm{\lambda} \|_2}{\lambdabar}$.  
Using the bound on the packing numbers from \cref{lem:PackSmooth},
\[
J \ \leq \ 9 \Cmax \frac{\| \bm{\lambda} \|_2}{\lambdabar} \ \int_{0}^1 \sqrt{  \log\left(   2 + 
    \frac{\Nhat_{\max}}{\Nhat_{\min} + 1}  \frac{32 L^4}{\Cmax^2 \gamma^2 \epsilon^2}
 \right) } d\epsilon 
\ \leq \ 
 9 \Cmax \frac{\| \bm{\lambda} \|_2}{\lambdabar} \ \int_{0}^1 \sqrt{  \log\left(\frac{t}{\epsilon^2}
 \right) } d\epsilon
 \]
where the second inequality uses $2 \leq 2/\epsilon^2$ and 
\(
t = 2 + 
    \frac{\Nhat_{\max}}{\Nhat_{\min} + 1}  \frac{32 L^4}{\Cmax^2 \gamma^2 } 
\) 
Substitute $\log(t /\epsilon^2) = 2 \log( \sqrt t/\epsilon)$ in the integral above, and then apply \cref{lemma:integral}, yielding
\begin{align*}
J  &\ \leq \ 
 9\sqrt 2 \cdot  \Cmax \frac{\| \bm{\lambda} \|_2}{\lambdabar}
\left( 
\sqrt{\pi}/2 + \sqrt{ \log \left(\sqrt {2 + 
    \frac{\Nhat_{\max}}{\Nhat_{\min} + 1}  \frac{32 L^4}{\Cmax^2 \gamma^2 } } \right) }
\right)
\\&\ \leq\ 
9\sqrt 2 \cdot  \Cmax \frac{\| \bm{\lambda} \|_2}{\lambdabar}
(\sqrt{\pi}+1)
\sqrt{ \log \left(\sqrt {2 + 
    \frac{\Nhat_{\max}}{\Nhat_{\min} + 1}  \frac{32 L^4}{\Cmax^2 \gamma^2 } } \right) }
,
\end{align*} 
where in the second inequality we have used $\sqrt{\pi  \log \left(\sqrt {2 + 
    \frac{\Nhat_{\max}}{\Nhat_{\min} + 1}  \frac{32 L^4}{\Cmax^2 \gamma^2 } } \right) }
\geq \sqrt{\pi \log(\sqrt 2 )} > \sqrt \pi /2$.  Thus, 
taking the $p$-norm of both sides and rounding up the leading constant shows that there exists a universal constant $\const_1$ such that 
\begin{equation} \label{eq:JSmoothFixed}
\|J\|_p  \ \leq \ 
\const_1 \cdot  \Cmax \frac{\| \bm{\lambda} \|_2}{\lambdabar}
 \magd{\sqrt{ \log \left(\sqrt {2 + 
    \frac{\Nhat_{\max}}{\Nhat_{\min} + 1}  \frac{32 L^4}{\Cmax^2 \gamma^2 } } \right) }}_p.
\end{equation}
We next bound the $p$-norm on the right. Invoke \cref{lem:RelatingOrlicz4} Part~\ref{PsiNormSqrtLogY} with $Y = \sqrt {2 + 
    \frac{\Nhat_{\max}}{\Nhat_{\min} + 1}  \frac{32 L^4}{\Cmax^2 \gamma^2 } } \geq \sqrt 2\geq1$. 
    Notice $\sqrt 2 \cdot \sqrt{\E[Y]} \geq 1$, which implies 
   \[
   \max(1, \sqrt{\E[Y]}/2) \ \leq \ 
   1 + \sqrt{\E[Y]}/{2}  \ \leq \ 
   (\sqrt{ 2} + 1) \sqrt{\E[Y]}/2 \leq \sqrt{\E[Y]}.
  \]
Hence, the norm on the right-hand side of \cref{eq:JSmoothFixed} is at most     
\begin{align*}
&5^{1/p} \left( \frac{p}{2e} \right)^{1/2}  \sqrt{\E\left[\sqrt {2 + 
    \frac{\Nhat_{\max}}{\Nhat_{\min} + 1}  \frac{32 L^4}{\Cmax^2 \gamma^2 } } \right]}
\\
& \quad \ \leq \ 
5^{1/p} \left( \frac{p}{2e} \right)^{1/2}  \sqrt[4]{2 + 
    \E\left[\frac{\Nhat_{\max}}{\Nhat_{\min} + 1}\right]  \cdot \frac{32 L^4}{\Cmax^2 \gamma^2 } },
&&(\text{Jensen's Inequality})
\\
& \quad \ \leq \ 
5^{1/p} \left( \frac{p}{2e} \right)^{1/2}  \sqrt[4]{2 + 
    \frac{32 L^4}{\Cmax^2 \gamma^2 }  \cdot 6 \left(\frac{12}{e}\right)^2 \frac{\lambdamax}{\lambdamin} \log^2 K  },
&&(\text{\cref{lem:PropertiesOfPoisson} Part~\ref{PNormSqrtRatio} })
\end{align*}    

We next use the assumptions on the parameters to rewrite this bound more simply.  
By the assumption that $\frac{4L^2}{C\gamma} \geq 1$, we have $\frac{32 L^4}{\Cmax^2 \gamma^2} \geq 2$.  
Moreover, since $K \geq 2$, $(\frac{12}{e} \log K)^2  \geq 1$.  Hence, 
the term under the square root is at most 
\(
\frac{64 L^4}{\Cmax^2 \gamma^2 }  \cdot 6 \left(\frac{12}{e}\right)^2 \frac{\lambdamax}{\lambdamin} \log^2 K.
\)

Substituting and simplifying thus shows there exists a universal constant $\const_2$ such that 
\begin{align*}
\magd{\sqrt{ \log \left(\sqrt {2 + 
    \frac{\Nhat_{\max}}{\Nhat_{\min} + 1}  \frac{32 L^4}{\Cmax^2 \gamma^2 } } \right) }}_p
&   \ \leq \  \const_2 \cdot 
5^{1/p} \sqrt p \cdot   \frac{L }{\sqrt{\Cmax \gamma} } \left(\frac{\lambdamax}{\lambdamin}\right)^{1/4} \sqrt{\log K}  
\end{align*}   

Hence, substituting above into \cref{eq:JSmoothFixed} shows there exists a universal constant $\const_3$ such that 
\[
\|J\|_p  \ \leq \ 
\const_3 \cdot  L \sqrt{\frac{\Cmax}{\gamma}} \cdot   \frac{\lambdamax}{\lambdamin}^{5/4} \cdot 
5^{\frac{1}{p}} \sqrt p \cdot 
\sqrt{K \log K}  
.
\]

Finally, applying 
\cref{thm:pollard}  yields
\[
 \sup_{\alpha \geq 0 } \abs{\frac{1}{K }\sum_{k=1}^K \prns{Z_k(\alpha,\bp_0) - \E[Z_k(\alpha,\bp_0) ]} } 
 \ \leq \ 
\const_3   \cdot 
\left(\frac{25}{\delta}\right)^{1/p} p  \cdot L \sqrt{\frac{C}{\gamma}} \left(\frac{\lambdamax}{\lambdamin}\right)^{5/4}
\frac{\sqrt{\log K }}{\sqrt{ K }}.
\]
This expression is minimized to first order by taking $p = 2 \log(1/\delta) \geq 1$ and observing $\left(\frac{25}{\delta}\right)^{\frac{1}{2 \log(1/\delta)}}$ is at most a constant for $0 < \delta < \frac{1}{2}$. Substituting and simplifying proves the first result.  

The proof of the second result is very similar, applying \cref{thm:pollard} to the process $\Floo$.  The only key difference is the envelope of this process is now $\frac{C}{N \lambdabar} \| \bm{\Nhat} \|_2 \leq \frac{C \sqrt K}{N \lambdabar} \Nhat_{\max}$ (cf. \cref{lem:Envelopes}).  
Thus, following the same steps that lead to \cref{eq:JSmoothFixed} but with this envelope shows that $J$ for this process satisfies
\begin{align*}
\|J\|_p  &\ \leq \ 
\const_4 \cdot   \frac{ \Cmax \sqrt K }{N \lambdabar}
 \magd{ \Nhat_{\max} \cdot \sqrt{ \log \left(\sqrt {2 + 
    \frac{\Nhat_{\max}}{\Nhat_{\min} + 1}  \frac{32 L^4}{\Cmax^2 \gamma^2 } } \right) }}_p
\\ &\ \leq \ 
\const_4 \cdot   \frac{ \Cmax \sqrt K }{N \lambdabar}
 \magd{ \Nhat_{\max}}_{2p} \cdot \magd{\sqrt{ \log \left(\sqrt {2 + 
    \frac{\Nhat_{\max}}{\Nhat_{\min} + 1}  \frac{32 L^4}{\Cmax^2 \gamma^2 } } \right) }}_{2p},
\end{align*}
for some constant $\const_4$, where the second inequality follows from H\"older's Inequality (cf. \cref{lem:NormProduct}).

Following an argument entirely analogous to the one that followed \cref{eq:JSmoothFixed} but with $p$ replaced by $2p$ shows 
\[
\magd{\sqrt{ \log \left(\sqrt {2 + 
    \frac{\Nhat_{\max}}{\Nhat_{\min} + 1}  \frac{32 L^4}{\Cmax^2 \gamma^2 } } \right) }}_{2p}
  \ \leq \ \const_5 \cdot
5^{\frac{1}{2p}}\cdot \sqrt p  \cdot 
 \frac{L }{\sqrt{\Cmax \gamma} } \cdot  \left( \frac{\lambdamax}{\lambdamin} \right)^{1/4} \cdot  \sqrt{\log K}  
\]

We bound $\| \Nhat_{\max} \|_{2p}$ using \cref{lem:PropertiesOfPoisson} Part~\ref{PNormNhatmax}.  

Then combining these bounds proves
\[ 
\| J \|_p 
\ \leq \ 
\const_6 \cdot  L   \sqrt{\frac{\Cmax }{\gamma} } \cdot \frac{\lambdamax}{\lambdamin}^{5/4} 
	6^{1/p} p^{3/2} \cdot  \sqrt{K}  \log^{3/2}(K).
\]
Applying \cref{thm:pollard}, substituting $p = 2 \log(1/\delta) > 1$ and simplifying yields the result.  \endproof
}

\subsubsection{Proof of \cref{thm:FixedPointShrinkage}}
We now can prove our main result:
\proof{Proof of \cref{thm:FixedPointShrinkage}.}
Combining \cref{lem:ConditionsForOptimality,lem:UCDeviations} shows 
if $\frac{4L^2}{\Cmax \gamma} \geq 1$, then there exists a universal constant $\const$ such that 
\begin{align*}
\edit{{\sf Sub Opt}_{\bp_0,K}(\alphaLOO_{\bp_0})}   \ \leq  \ 
\const \cdot 
L \sqrt{\frac{\Cmax}{\gamma}} \cdot 
\left( \frac{\lambdamax}{\lambdamin}\right)^{5/4} \cdot
\frac{ \log^{2} (1/\delta)\cdot  \log^{3/2}(K)}{ \sqrt K }.
\end{align*}
If $\frac{4L^2}{\Cmax \gamma} < 1$, we can always increase $L$ until $\frac{4L^2}{\Cmax \gamma} = 1$ as the larger $L$ remains a valid Lipschitz constant.  Increasing the leading constant in this case proves the theorem.  
\endproof

\subsection{Deferred Proofs from \cref{sec:SmoothCostsGeneral}: Shrunken-SAA with Data-Driven Anchors for Strongly-Convex Problems}
\label{sec:SmoothCostsGeneralAppendix}

Our strategy to proving \cref{thm:SmoothThmGeneralAnchor,thm:SmoothThmLooAnchor} is similar to proving to \cref{thm:FixedPointShrinkage} except that our process is now indexed by both $\alpha\geq0$ and $\bq\in \edit{ \mathcal P}$.

\subsubsection{Maximal deviation bounds.}
Our first step is to use \cref{lem:ContinuityAlpha}, part~\ref{part:ContinuityinP0} to reduce bounding the maximal deviations of $\Zperf(\cdot,\cdot),\,\Zloo(\cdot,\cdot)$ to bounding the maximal deviations of $\Zperf(\cdot,\bq),\,\Zloo(\cdot,\bq)$ for a finite number of fixed anchors $\bq \in \mathcal P$.
\begin{lemma}[Reduction to Maximal Deviations with Fixed Anchor] \label{lem:MaximalDeviationReduction}
Under the assumptions of \cref{thm:SmoothThmGeneralAnchor}, if $\{\bq^1,\dots,\bq^M\}$ is an $\epsilon_0$-covering of $\edit{\mathcal P}$ with respect to $\ell_1$, then
\begin{align}
\label{eq:TruePerfGeneralDecomposition}
 \sup_{\alpha \geq 0, \bq \in \edit{\imh} } \abs{ \overline Z(\alpha, \bq) - \E[\overline Z(\alpha, \bq) ] }  \ &\leq  \ 
\edit{\frac{2L^2}{\gamma}} \epsilon_0    +  \max_{j=1,\dots,M}\sup_{{\alpha \geq 0}} \abs{ \overline Z(\alpha, \bq^j) - \E[\overline Z(\alpha, \bq^j) ] },\\
\label{eq:LOOGeneralDecomposition}
\sup_{\alpha \geq 0, \bq \in \edit{\imh} } \abs{ \overline Z^{\sf LOO}(\alpha, \bq) - \E[\overline Z^{\sf LOO}(\alpha, \bq) ] }  \ &\leq  \ 
    \edit{\edit{\frac{2L^2}{\gamma} \frac{\Nhat_{\text{avg}}}{N \lambdabar}}\epsilon_0}
    \\\notag&\phantom{\leq\ }+  \max_{j=1,\dots,M}\sup_{{\alpha \geq 0}} \abs{ \overline Z^{\sf LOO}(\alpha, \bq^j) - \E[\overline Z^{\sf LOO}(\alpha, \bq^j) ] }. 
\end{align}
\end{lemma}
\proof{Proof.}%
Consider the first inequality. Fix some $\bq \in \edit{\mP}$, and suppose $\bq^j$ is the closest member of the covering.  
Then,
{\blockedit
\begin{align*}
\abs{Z_k(\alpha, \bq) - Z_k(\alpha, \bq^j)} 
& \ \leq \ 
\frac{\lambda_k}{\lambdabar} \abs{ {\bp_k}^\top \left( \bc_k( \bx_k(\alpha, \bq,  \bfmhat_k)) - \bc_k( \bx_k(\alpha, \bq^j,  \bfmhat_k)) \right) }
\\ &\  \leq \ 
L \cdot \frac{\lambda_k}{\lambdabar} \left\|  \bx_k(\alpha, \bq,  \bfmhat_k) -  \bx_k(\alpha, \bq^j,  \bfmhat_k)  \right\|_2
        && (\text{Lipschitz Continuity})
\\ & \ \leq \ 
	\frac{L^2}{\gamma} \magd{\bq-\bq^j}_1 \frac{\lambda_k}{\lambdabar}
        && (\text{\cref{lem:ContinuityAlpha}, part~\ref{part:ContinuityinP0}})
\\ & \ \leq \ 
	\frac{L^2}{\gamma} \epsilon_0 \frac{\lambda_k}{\lambdabar}
\end{align*}
Averaging over $k$ shows 
\(
\abs{\overline Z(\alpha, \bq) - \overline Z(\alpha, \bq^j)}  \ \leq \ 	\frac{L^2}{\gamma} \epsilon_0.
\)
}
By Jensen's inequality, this bound also implies that 
\(
\abs{\E[\overline Z(\alpha, \bq)] - \E[\overline Z(\alpha, \bq^j)] } 
	\leq \E\left[ \abs{\overline Z(\alpha, \bq) - \overline Z(\alpha, \bq^j)} \right] 
	\leq  \frac{L^2}{\gamma} \epsilon_0.
\)
Hence, by the triangle inequality, 
\begin{align*}
\abs{ \overline Z(\alpha, \bq) - \E\left[ \overline Z(\alpha, \bq) \right]} 
& \ \leq \ 
\abs{ \overline Z(\alpha, \bq) - \overline Z(\alpha, \bq^j) } 
+ 
\abs{ \E\left[ \overline Z(\alpha, \bq)  - \overline Z(\alpha, \bq^j) \right]} 
+ 
\abs{ \overline Z(\alpha, \bq^j) - \E\left[ \overline Z(\alpha, \bq^j) \right]} .
\\
& \  \leq \  \frac{2 L^2}{\gamma} \epsilon_0  + \abs{ \overline Z(\alpha, \bq^j) - \E\left[ \overline Z(\alpha, \bq^j) \right]} .
\end{align*}
Substituting yields the first inequality in the result.

We next prove the second inequality. Fix some $\bq \in \edit{\mP}$, and suppose $\bq^j$ is the closest member of the covering.  Then,
\begin{align*}
\Big| &\overline Z^{\sf LOO}(\alpha, \bq) - \overline Z^{\sf LOO}(\alpha, \bq^j) \Big| 
\\
& \leq 
    \frac{1}{KN \lambdabar} \sum_{k=1}^K  
    \sum_{i=1}^d  \mhat_{ki} \abs{  c_{ki}( \bx_k(\alpha, \bq,  \bfmhat_k - \be_i ) ) - c_{ki}( \bx_k(\alpha, \bq^j,  \bfmhat_k - \be_i ) )  }
\\
& \leq
    \frac{L}{K N \lambdabar } \sum_{k=1}^K     \sum_{i=1}^d  \mhat_{ki} \left\|  \bx_k(\alpha, \bq,  \bfmhat_k - \be_i ) -  \bx_k(\alpha, \bq^j,  \bfmhat_k - \be_i )  \right\|_2
        && (\text{Lipschitz Continuity})
\\
& \leq
    \edit{\frac{L^2}{N \lambdabar \gamma} \magd{\bq-\bq^j}_1\frac1K\sum_{k=1}^K     \Nhat_k}
        && (\text{\cref{lem:ContinuityAlpha}, part~\ref{part:ContinuityinP0}})
\\
& \leq
    \edit{\frac{L^2}{\gamma} \frac{\Nhat_{\text{avg}}}{N \lambdabar}}\epsilon_0
\end{align*}
By Jensen's inequality, this further implies that $\abs{\E[\overline Z^{\sf LOO}(\alpha, \bq)] - \E[\overline Z^{\sf LOO}(\alpha, \bq^j)] } \leq \E\left[ \abs{\overline Z^{\sf LOO}(\alpha, \bq) - \overline Z^{\sf LOO}(\alpha, \bq^j)} \right]  \leq          \edit{\edit{\frac{L^2}{\gamma} \frac{\Nhat_{\text{avg}}}{N \lambdabar}}\epsilon_0} $.
Using the triangle inequality as before and
applying the two bounds above yields our second inequality in the result.
\endproof

\vskip 12 pt
We next use the above lemmas to bound the maximal deviations of interest via \cref{thm:pollard}:
{\blockedit
\begin{lemma}[Bounding Maximal Deviations General Anchors] \label{lem:ULLNGeneralAnchors}
Under the assumptions of \cref{thm:SmoothThmGeneralAnchor}, there exists a universal constant $\const$ such that for any $0 < \delta < \frac{1}{2}$,  the following two statements each hold (separately) with probability at least $1-\delta$:
\begin{align*} \hspace{-20pt}
\sup_{\alpha \geq 0, \ \bq \in \edit{\mathcal P} } \abs{ \frac{1}{K} \sum_{k=1}^K Z_k(\alpha,\bq) - \E[Z_k(\alpha,\bq)] } 
&\ \leq \ 
    \const  \cdot 
    \max\left(\Cmax, \ \frac{L^2}{\gamma} + L \sqrt{\frac{\Cmax}{\gamma}} \right)  
    \left( \frac{\lambdamax}{\lambdamin} \right)^{5/4}  
    \frac{d_0 \log^{3/2}(K)  \log(1/\delta) }{\sqrt K},
\\ \hspace{-20pt}
\sup_{\alpha \geq 0, \ \bq \in \mP} \abs{ \frac{1}{K} \sum_{k=1}^K Z^{\sf LOO}_k(\alpha,\bq) - \E[Z^{\sf LOO}_k(\alpha,\bq)] } 
&\ \leq \ 
    \const  \cdot 
    \max\left( \Cmax, \ \frac{L^2}{\gamma}  \ + \  L\sqrt{\frac{\Cmax}{\gamma}} \right) \left( \frac{\lambdamax}{\lambdamin} \right)^{5/4}
    \frac{d_0^2 \log^{7/2}(K)  \log^2(1/\delta) }{\sqrt K}.
\end{align*}
\end{lemma}
\proof{Proof.}
First consider the case $\frac{4L^2}{C\gamma} \geq 1$.  Fix some $0 < \epsilon_0 < \frac{1}{2}$ and consider a minimal $\epsilon_0$-covering of $\mathcal P$ with respect to $\ell_1$.  Denote its size by $M$.  Necessarily, $M \leq D_1(\epsilon_0, \mathcal P)$ (cf. \citealp[pg. 10]{pollard1990empirical}).  Apply \cref{lem:MaximalDeviationReduction} with this covering, and then apply the first part of \cref{lem:UCDeviations} with $\delta \leftarrow \delta/M$ to bound the remaining suprema.  This shows that there exists a constant $\const_1$ such that with probability at least $1-\delta$, 
\[
\sup_{\alpha \geq 0, \ \bq \in \edit{\mathcal P} } \abs{ \frac{1}{K} \sum_{k=1}^K Z_k(\alpha,\bq) - \E[Z_k(\alpha,\bq)] } 
\ \leq \ 
    \const_1  \cdot 
    \frac{L^2}{\gamma} \epsilon_0 \ + \  \const_1 L\sqrt{\frac{\Cmax}{\gamma}} \left( \frac{\lambdamax}{\lambdamin} \right)^{5/4}
    \frac{\log^{1/2}(K) }{\sqrt K} \cdot \log\prns{\frac{D_1(\epsilon_0, \mathcal P)}{\delta}}.
\]
Directly optimizing the choice of $\epsilon_0$ appears difficult.  We instead take the (suboptimal) choice $\epsilon_0 = \frac{1}{2 \sqrt K}$ and note $\epsilon_0 < \frac{1}{2}$ since $K \geq 2$.  Furthermore, by assumptions on the parameters, $d_0 \geq 1$, $2 \log K \geq 1$ and $2 \log(1/\delta) \geq 1$.  Hence, 
\begin{align*}
\log(D_1(\epsilon_0, \mathcal P)/\delta) 
& \ \leq \ 
\log(1/\delta) + d_0 \log(1/\epsilon_0)
\\ & \ = \ 
\log(1/\delta) + d_0 \log 2 + \frac{d_0}{2} \log K 
\\ &\ \leq \ 
2 d_0 \log K \log(1/\delta) + 
2 d_0 \log K \log (1/\delta) + 
d_0 \log K \log (1/\delta) 
\\& \ =  \ 
5 d_0 \log K \log (1/\delta).
\end{align*}
Substituting above shows there exists a constant $\const_2$ such that  
\begin{align*}
\sup_{\alpha \geq 0, \ \bq \in \edit{\mathcal P} } \abs{ \frac{1}{K} \sum_{k=1}^K Z_k(\alpha,\bq) - \E[Z_k(\alpha,\bq)] } 
& \ \leq \ 
    \const_2  \cdot 
    \frac{L^2}{\gamma\sqrt K }  \ + \  \const_2 L\sqrt{\frac{\Cmax}{\gamma}} \left( \frac{\lambdamax}{\lambdamin} \right)^{5/4}
    \frac{d_0 \log^{3/2}(K)  \log(1/\delta) }{\sqrt K},
\\ 
& \ \leq \ 
    \const_3  \cdot 
    \left( \frac{L^2}{\gamma}  \ + \  L\sqrt{\frac{\Cmax}{\gamma}} \right)  \left( \frac{\lambdamax}{\lambdamin} \right)^{5/4}
    \frac{d_0 \log^{3/2}(K)  \log(1/\delta) }{\sqrt K},
\end{align*}
by collecting constants.  

In the case $\frac{4L^2}{C\gamma} < 1$, we can always increase $L$ until $\frac{4L^2}{C\gamma} = 1$ as the larger $L$ remains a valid Lipschitz constant.  Substituting this increased $L$ yields the leading term $3\Cmax/4$ and proves the first inequality.  

The proof of the second inequality is very similar.  Assume $\frac{4L^2}{C\gamma} \geq 1$.  Again, applying \cref{lem:MaximalDeviationReduction} over an $\epsilon_0$-covering and using \cref{lem:UCDeviations} with $\delta \leftarrow \frac{\delta}{2M}$ to bound the remaining suprema shows that with probability at least $1-\delta/2$, 
\begin{align*}
&\sup_{\alpha \geq 0, \ \bq \in \edit{\mathcal P} } \abs{ \frac{1}{K} \sum_{k=1}^K Z_k(\alpha,\bq) - \E[Z_k(\alpha,\bq)] } 
\\& \quad \ \leq \ 
    \const_4  \cdot 
    \frac{L^2}{\gamma} \frac{\Nhatbar}{N \lambdabar} \epsilon_0 \ + \  \const_4 L\sqrt{\frac{\Cmax}{\gamma}} \left( \frac{\lambdamax}{\lambdamin} \right)^{5/4}
    \frac{\log^{3/2}(K) }{\sqrt K} \cdot \log^2\left(\frac{2D_1(\epsilon_0, \mathcal P)}{\delta}\right).
\end{align*}
Take the (suboptimal) choice $\epsilon_0 = \frac{1}{2 \sqrt K}$.  The same simplifications from above show that 
\begin{align*}
\log(2D_1(\epsilon_0, \mathcal P)/\delta) \ \leq \ \log 2 + 5 d_0 \log K \log (1/\delta) \ \leq \ 7 d_0 \log K \log(1/\delta), 
\end{align*}
whereby with probability at least $1-\delta/2$, 
\begin{align*}
\sup_{\alpha \geq 0, \ \bq \in \edit{\mathcal P} } \abs{ \frac{1}{K} \sum_{k=1}^K Z_k(\alpha,\bq) - \E[Z_k(\alpha,\bq)] } 
\\ \quad & \ \leq \ 
    \const_5  \cdot 
    \frac{L^2}{\gamma\sqrt K } \frac{\Nhatbar}{N \lambdabar}  \ + \  \const_5 L\sqrt{\frac{\Cmax}{\gamma}} \left( \frac{\lambdamax}{\lambdamin} \right)^{5/4}
    \frac{d_0^2 \log^{7/2}(K)  \log^2(1/\delta) }{\sqrt K}.
\end{align*}

It remains to bound the fraction $\frac{\Nhatbar}{N\lambdabar} = \frac{K \Nhatbar}{KN\lambdabar}$.  Notice $K\Nhatbar \sim \text{Poisson}(KN\lambdabar)$.  From \cref{lem:PropertiesOfPoisson} Part~\ref{Psi1Poisson} 
applied to $K\Nhatbar$ and Markov's inequality, we have that with probability at least $1-\delta/2$, $\frac{\Nhatbar}{N \lambdabar} \leq \log(4/\delta)$.  

Substitute this bound above, apply the union bound and collect constants to show that with probability at least $1-\delta$
\begin{align*}
\sup_{\alpha \geq 0, \ \bq \in \edit{\mathcal P} } \abs{ \frac{1}{K} \sum_{k=1}^K Z_k(\alpha,\bq) - \E[Z_k(\alpha,\bq)] } 
& \ \leq \ 
    \const_6  \cdot 
    \left( \frac{L^2}{\gamma}  \ + \  L\sqrt{\frac{\Cmax}{\gamma}} \right) \left( \frac{\lambdamax}{\lambdamin} \right)^{5/4}
    \frac{d_0^2 \log^{7/2}(K)  \log^2(1/\delta) }{\sqrt K}.
\end{align*}

In the case $\frac{4L^2}{C\gamma} < 1$, we can again increase $L$ until $\frac{4L^2}{C\gamma} = 1$ since the larger $L$ is still a valid Lipschitz constant.  Substituting this increased $L$ yields the leading term $3\Cmax/4$ and proves the second claim.  
\endproof
}

\subsubsection{Proofs of \cref{thm:SmoothThmGeneralAnchor,thm:SmoothThmLooAnchor}.}
{\blockedit We can now prove the main results of the section via our previously outlined strategy.  
\proof{Proof of \cref{thm:SmoothThmGeneralAnchor,thm:SmoothThmLooAnchor}.}
The proofs of both theorems are identical.  
For both theorems, by \cref{lem:ConditionsForOptimality}, the quantity to be bounded is bounded by the sum of the same two maximal deviations. These are in turn bounded by \cref{lem:ULLNGeneralAnchors}. Instantiating each bound for $\delta \leftarrow \delta/2$, adding the right hand sides and applying the union bound yields a bound on the sub-optimality.  Collecting dominant terms yields the result.\endproof}

\subsection{Proof of \cref{thm:FixedPointShrinkageDiscrete}: Shrunken-SAA with Fixed Anchors for Discrete Problems}\label{sec:FixedPointShrinkageDiscreteAppendix}

We first use \cref{thm:SizeOfDiscreteSets} proven in \cref{sec:DisceteFixedAnchor} to prove the following bounds on the maximal deviations of interest via \cref{thm:pollard}.
{\blockedit
\begin{lemma}[Bounding Maximal Deviations for Discrete Problems] \label{lem:ULLNDiscrete}
Under the assumptions of \cref{thm:FixedPointShrinkageDiscrete}, there exists a constant $\const$ such that for any $0 < \delta < 1/2$, 
the following two statements hold (separately) each with probability at least $1-\delta$:
\begin{align*}
\sup_{\alpha \geq 0} \abs{ \frac{1}{K} \sum_{k=1}^K Z_k(\alpha,\bp_0) - \E[Z_k(\alpha,\bp_0)] } 
&\ \leq \ 
\const  \cdot \Cmax \frac{\lambdamax}{\lambdamin} \cdot \sqrt{ \log\left( \sum_{k=1}^K \abs{\X_k}  \right) }
                                 \cdot \frac{\sqrt{\log\left( \frac{1}{\delta} \right)}  }{\sqrt K},
\\
\sup_{\alpha \geq 0} \abs{ \frac{1}{K} \sum_{k=1}^K Z^{\sf LOO}_k(\alpha,\bp_0) - \E[Z^{\sf LOO}_k(\alpha,\bp_0)] } 
& \ \leq \ 
\const \cdot \Cmax  \frac{\lambdamax}{\lambdamin} \cdot 
\sqrt{ \log\left( N_{\max} \sum_{k=1}^K \abs{\X_k} \right)} \cdot 
	 \frac{ \log^{3/2}(K) \cdot \log^{3/2}(1/\delta)}{\sqrt K}.
\end{align*}

\end{lemma}
\proof{Proof.}%
Consider the first inequality.  We first bound the variable $J$ in \cref{eq:JDudley} corresponding to the process $\Fperf$ with the envelope given by \cref{lem:Envelopes}.
By \cref{thm:SizeOfDiscreteSets},
\[
J \ \leq  \ 9 \Cmax \cdot \frac{\lambdamax}{\lambdamin} \cdot \sqrt K \sqrt{ \log\left( 2 \sum_{k=1}^K \abs{\X_k} \right)},
\]
where we have upper bounded $\| \bm\lambda \|_2 \leq \lambdamax \sqrt K$.  
From \cref{thm:pollard}, there exists a constant $\const_1$ such that with probability at least $1-\delta$, 
\begin{align*}
\sup_{\alpha \geq 0} \abs{ \frac{1}{K} \sum_{k=1}^K Z_k(\alpha,\bp_0) - \E[Z_k(\alpha,\bp_0)] } & \leq \const_1 \cdot
    \left( \frac{5}{\delta} \right)^{1/p} p^{1/2} \cdot \Cmax \frac{\lambdamax}{\lambdamin} \cdot \sqrt{ \frac{ \log\left( 2 \sum_{k=1}^K \abs{\X_k} \right)}{K}}.
\end{align*}
Let $p = 2 \log(1/\delta) > 1$, and collect constants to complete the proof. 

The proof of the second inequality is similar but uses different envelopes (cf. \cref{lem:Envelopes}) and the larger packing numbers of \cref{thm:SizeOfDiscreteSets}.  Specifically, we note that $\min(d, \Nhat_k) \leq \Nhat_{\max}$ and  $\| \bm{\Nhat} \|_2 \leq \Nhat_{\max} \sqrt K$, and
bound $J$ as 
\[
J \ \leq \ 9 \frac{\Cmax \sqrt K}{N \lambdabar} \Nhat_{\max}  \sqrt{ \log\left( 1 + 2 \Nhat_{\max} \sum_{k=1}^K \abs{\X_k} \right) }.
\]

Recall $N_{\max} \equiv N \lambdamax \geq N \lambdamin \geq 1$.  Thus, 
we can upper bound the logarithm as 
\begin{align*}
\log\left( 1 + 2 \Nhat_{\max} \sum_{k=1}^K \abs{\X_k} \right) 
& \ \leq \ 
\log\left( 6 N_{\max} \sum_{k=1}^K \abs{\X_k} + 2 \Nhat_{\max} \sum_{k=1}^K \abs{\X_k} \right) 
\\
& \ = \ 
\underbrace{\log\left( 2 N_{\max} \sum_{k=1}^K \abs{\X_k} \right)}_{\geq \log 4}  + \underbrace{\log \left( 3 + \frac{\Nhat_{\max}}{N_{\max}} \right)}_{\geq \log 3 } 
\\ 
\ &\  \leq \ 
2 \log\left( 2 N_{\max} \sum_{k=1}^K \abs{\X_k} \right) \cdot \log \left( 3 + \frac{\Nhat_{\max}}{N_{\max}} \right), 
\end{align*}
where the last inequality follows because $a + b \leq 2 a b$ when $a, b \geq 1$ 

Substituting above and  taking the $p$-norm shows there exists a constant $\const_2$ such that 
\begin{align*}
\| J\|_p  &\ \leq \ \const_2 \cdot  \frac{\Cmax \sqrt K}{N \lambdabar}   
\sqrt{ \log\left( 2 N_{\max} \sum_{k=1}^K \abs{\X_k} \right)} \cdot  \magd{\Nhat_{\max}  \sqrt{ \log \left( 3 + \frac{\Nhat_{\max}}{N_{\max}} \right) }}_{p} 
\\ & \ \leq \ 
\const_2 \cdot \frac{\Cmax \sqrt K}{N \lambdabar}   
\sqrt{ \log\left( 2 N_{\max} \sum_{k=1}^K \abs{\X_k} \right)} \cdot 
 \magd{\Nhat_{\max}}_{2p}  \cdot \magd{\sqrt{ \log \left( 3 + \frac{\Nhat_{\max}}{N_{\max}} \right) }}_{2p} ,
\end{align*}
where the second inequality follows from H\"older's Inequality (cf. \cref{lem:NormProduct}) We next bound these two $2p$-norms.

We bound the second $2p$-norm using \cref{lem:RelatingOrlicz4} Part~\ref{PsiNormSqrtLogY} with $Y = 3 + \frac{\Nhat_{\max}}{N_{\max}} > 3$, yielding 
\begin{align*}
\magd{\sqrt{ \log \left( 3 + \frac{\Nhat_{\max}}{N_{\max}} \right) }}_{2p} 
&\ \leq \ 
5^{\frac{1}{2p}} \sqrt{\frac{p}{e}} \max\left(1, \frac{1}{2} \sqrt{ 3 + \E\left[  \frac{\Nhat_{\max}}{N_{\max}}  \right] } \right)
\\
& \ \leq \ 
5^{\frac{1}{2p}} \sqrt{\frac{p}{e}}\max\left(1, \frac{1}{2} \sqrt{ 3 + \frac{36}{e} \log K } \right) 
&&\text{(\cref{lem:PropertiesOfPoisson} Part~\ref{PNormNhatmax})}
\\ 
& \ \leq 2 \cdot \ 5^{\frac{1}{2p}} \sqrt{p}\sqrt{\log K }, 
\end{align*}
since $K \geq 2$.

Similarly, bound $ \magd{\Nhat_{\max}}_{2p}$ using \cref{lem:PropertiesOfPoisson} Part~\ref{PNormNhatmax}.  

Combining shows
\begin{align*}
\| J \|_p \ \leq \ 
\const_3 \cdot \frac{\Cmax N_{\max} \sqrt K}{N \lambdabar}   
\sqrt{ \log\left( 2 N_{\max} \sum_{k=1}^K \abs{\X_k} \right)} \cdot 
	6^{\frac{1}{p}}  p^{3/2} \cdot   \log^{3/2}(K).
\end{align*}
Applying \cref{thm:pollard} and substituting $p = 2 \log(1/\delta) > 1$ proves the second inequality.  
\endproof

}

We can now prove the main result of the section.  
\begin{proof}{Proof of \cref{thm:FixedPointShrinkageDiscrete}.}
\Cref{lem:ULLNDiscrete} bound the maximal deviations in \cref{lem:ConditionsForOptimality}.  Instantiating them for $\delta \leftarrow \delta/2$, adding their righthand sides and applying the union bound bounds the sub-optimality.  Collecting dominant terms proves the result.  
\end{proof}

\subsection{Deferred Proofs from \cref{sec:DiscreteGeneral}: Shrunken-SAA with Data-Driven Anchors for Discrete Problems.}
\label{sec:DiscreteGeneralAppendix}
As a first step towards our proof, we prove \cref{thm:SizeOfDiscreteSetsGeneral}.
{\blockedit Recall the $m \equiv \sum_{k=1}^K \binom{\abs{ \X_k}}{2}$ hyperplanes defined in \cref{sec:DiscreteGeneral}:
\begin{equation}\notag
\edit{H_{kij}=\braces{\btheta\in\R{d_0}:
\left(V\btheta +\bfmhat_{k} \right)^\top \left( \bc_k( \bx_{ki}) - \bc_k( \bx_{kj}) \right)
= \bm 0
},\quad \forall\;k=1,\dots,K,\,i\neq j=1,\dots, \abs{\X_k}.}
\end{equation}
In words, for $\btheta$ on $H_{kij}$ we are indifferent between $\bx_{ki}$ and $\bx_{kj}$ when using $\btheta$ in \cref{eq:ThetaParam}.  On either side, we strictly prefer one solution.

For any fixed $\btheta \in \R{d_0}$, we considered the polyhedron induced by the equality constraints of those hyperplanes containing $\btheta$, and the inequality constraints defined by the side on which $\btheta$ lies for the remaining hyperplanes.  
We call such polyhedra \emph{fully-specified} because they are defined by their relationship to \emph{all} $m$ hyperplanes in the arrangement.  Because this polyhedron lives in \edit{$\R{d_0}$}, it necessarily has dimension \edit{$j \leq d_0$}.  For example the shaded region in \cref{fig:DiscretePackingGeneral} is a fully-specified polyhedron with $j=2$, the bold line segment has $j =1$ and the bold point has $j = 0$.  
As argued in the main text, to bound \edit{$\abs{\FperfqmP}$} it suffices to count the number of $j$-dimensional fully-specified polyhedron in the arrangement of the above $m$ hyperplanes for all $0 \leq j \leq d_0$.}

Counting the polyhedra induced by hyperplane arrangements is a classical problem in geometry.  For example, it is well-known that the number of ${\edit{d_0}}$-dimensional, fully-specified polyhedra in a hyperplane arrangement with $m$ hyperplanes in $\R {\edit{d_0}}$ is at most $\sum_{i=0}^{\edit{d_0}} \binom{m}{i}$ \citep[Prop. 2.4]{stanley2004introduction}.  We first use this result to bound the total number of polyhedra in an arbitrary arrangement with $m$ hyperplanes in $\R {\edit{d_0}}$.

\begin{lemma}[Number of Fully-Specified Polyhedra]\label{lem:FullySpecifiedPolyhedra}
In a hyperplane arrangement with $m$ hyperplanes in $\R {\edit{d_0}}$, the number of fully-specified polyhedra is at most 
\[
\sum_{j=0}^{\edit{d_0}} \binom{m}{{\edit{d_0}}-j} \sum_{i=0}^j \binom{m-{\edit{d_0}} + j }{i}  \ \leq \ (1 + 2 m )^{\edit{d_0}}.
\]
\end{lemma}
\begin{proof}{Proof of \cref{lem:FullySpecifiedPolyhedra}}  
Each fully-specified polyhedron has some dimension, $0 \leq j \leq {\edit{d_0}}$. We will count the number of such fully-specified polyhedra by counting for each dimension $j$.

Fix some $0 \leq j \leq {\edit{d_0}}$.  Notice that each $j$-dimensional polyhedron lives in a $j$-dimensional subspace defined by ${\edit{d_0}}-j$ linearly independent hyperplanes from the arrangement.  
There are at most  $\binom{m}{{\edit{d_0}}-j}$ ways to choose these linearly independent ${\edit{d_0}}-j$ hyperplanes.  Next project the remaining hyperplanes onto this subspace which yields at most $m-{\edit{d_0}}+j$ non-trivial hyperplanes in the subspace, i.e., hyperplanes that are neither the whole subspace nor the empty set. These non-trivial hyperplanes ``cut up" the subspace into various polyhedra, including $j$-dimensional, fully-specified polyhedra. 
By \cite[Prop. 2.4]{stanley2004introduction}, the number of $j$-dimensional, fully-specified polyhedra in this hyerplane arrangement of at most $m-{\edit{d_0}}+j$ hyperplanes in $j$-dimensional space is at most $\sum_{i=0}^j \binom{m-{\edit{d_0}}+j}{i}$.  In summary, it follows that there are at most $\binom{m}{{\edit{d_0}}-j} \sum_{i=0}^j \binom{m-{\edit{d_0}}+j}{i}$ $j$-dimensional, fully-specified polyhedra in the arrangement.  

Summing over $j$ gives the lefthand side of the bound in the lemma.  

For the righthand side, recall that 
\[
\sum_{i=0}^j \binom{m -{\edit{d_0}} + j}{i}  
\ \ \leq  \ \ 
\sum_{i=0}^j (m-{\edit{d_0}} + j)^i \cdot 1 ^{m-{\edit{d_0}} + j - i } 
\ \ \leq \  \ (1 + m -{\edit{d_0}} +j )^ j 
\ \ \leq  \ \  (1 + m) ^ j, 
\]
where the penultimate inequality is the binomial expansion and the last follow because $j \leq {\edit{d_0}}$.  Next,
\begin{align*}
\sum_{j=0}^{\edit{d_0}} \binom{m}{{\edit{d_0}}-j} \sum_{i=0}^j \binom{m-{\edit{d_0}} + j }{i} 
&  \ \leq \ 
\sum_{j=0}^{\edit{d_0}} \binom{m}{{\edit{d_0}}-j} (1 + m)^j
\\
&  \ \leq  \ \sum_{j=0}^{\edit{d_0}} m^{{\edit{d_0}}-j} (1 + m)^j  && 
\\
&  \ = \  (1 + 2m )^{\edit{d_0}},
\end{align*}
where the last equality is again the binomial expansion.  
\end{proof}

\vskip 12pt
We can now bound the cardinality of the relevant solution sets.  
\begin{proof}{Proof of \cref{thm:SizeOfDiscreteSetsGeneral}.}
Recall there are $m = \sum_{k=1}^K \binom{\abs{\X_k}}{2}$ hyperplanes in the arrangement \cref{eq:totalorderinghyperplanes} \edit{in $\R{d_0}$}, and the number of fully-specified polyhedra in this arrangement upper-bounds \edit{$\abs{\FperfqmP}$}.  Noting $1 + 2m  = 1 + \sum_{k=1}^K \abs{\X_k}(\abs{\X_k} - 1) \leq \sum_{k=1}^K \abs{\X_k}^2$ yields the first bound.

A similar argument can be used to bound \edit{$\abs{ \FlooqmP }$}.  
{\blockedit 
Suppose first $\Nhat_{\max} = 0$.  Then this set has size $1$.  On the other hand, if $\Nhat_{\max} > 0$, 
let $\mathcal I_k=\{i=1,\dots,d:\mhat_{ki}>0\}$, so that
\begin{align} \notag
&\abs{ \FlooqmP}
\leq
 \abs{ \left\{ \left(\bx_k(\alpha, \bq, \bfmhat_k - \be_i ) \right)_{k = 1,\dots,K,i\in\mathcal I_k}  \ : \ \bq \in \mP, \alpha \geq 0 \right\}  } 
\\ \label{eq:ySolLOO}
&\qquad\qquad\leq
\abs{ \left\{ \left( \bx_k(\|V\btheta\|_1, V\btheta/\|V\btheta\|_1, \bfmhat_k - \be_i) \right)_{k = 1,\dots,K,i\in\mathcal I_k} \ : \ \btheta \in \R{d_0},V\btheta\in\R d_+ \right\} }.
\end{align}
We then consider the arrangement generated by}
\begin{align*}
\edit{H_{kijl}=\braces{\btheta\in\R {d_0}:
\left( \bm V\btheta +\bfmhat_{k} - \be_l \right)^\top \left( \bc_k( \bx_{ki}) - \bc_k( \bx_{kj}) \right)
=
\bm 0
}},
\end{align*}
for all $ k= 1, \ldots, K$, $i, j = 1, \ldots, \abs{\X_k}$ with $i \neq j$, and \edit{$l \in\mathcal I_k$}.  
Notice \edit{that since $\abs{ \mathcal I_k} \leq\Nhat_k$ there are at most $\Nhat_{\max} \sum_{k=1}^K \binom{\abs{\X_k}}{2}$} such hyperplanes.  Moreover, \edit{$\abs{\FlooqmP}$} is upper-bounded by the number of fully-specified polyhedra in this arrangement.  
Note that \edit{$1+ 2 \Nhat_{\max} \sum_{k=1}^K \binom{\abs{\X_k}}{2} = 1 + \Nhat_{\max} \sum_{k=1}^K \abs{\X_k} (\abs{\X_k} -1 )  \leq \Nhat_{\max} \sum_{k=1}^K \abs{\X_k}^2$}.  
\edit{Adding $1$ covers the case $\Nhat_{\max} = 0$.}  
Plugging in this value into \cref{lem:FullySpecifiedPolyhedra} yields the second bound above.
\end{proof}

\subsubsection{Maximal Deviation Bounds.}
We next use \cref{thm:SizeOfDiscreteSetsGeneral} to bound the maximal deviations of interest via \cref{thm:pollard}.
\begin{lemma}[Bounding Maximal Deviations, Discrete Case, General Anchors]  \label{lem:ULLNDiscreteGeneral}   
Under the assumptions of \cref{thm:ShrinkageDiscreteGeneralAnchor}, there exists a constant $\const$ such that 
for any $0 < \delta < \frac{1}{2}$, both of the following statements hold (separately) with probability at least $1-\delta$:
\begin{align*}
\sup_{\alpha \geq 0,\ \bq\in\edit{\mP}} \abs{ \frac{1}{K} \sum_{k=1}^K Z_k(\alpha,\bq) - \E[Z_k(\alpha,\bq)] } 
&\ \leq \ 
\const \cdot   \Cmax \frac{\lambdamax}{\lambdamin} \cdot 
                                  \sqrt{\edit{d_0}\log\left(\sum_{k=1}^K \abs{\X_k}  \right) }  \cdot
                                 \frac{\sqrt{\log\left( 1/\delta \right)}}{\sqrt K}.
\\
\sup_{\alpha \geq 0,\ \bq\in\edit{\mP}} \abs{ \frac{1}{K} \sum_{k=1}^K Z^{\sf LOO}_k(\alpha,\bq) - \E[Z^{\sf LOO}_k(\alpha,\bq)] } 
&\ \leq \ 
\const \cdot 
\Cmax  \frac{\lambdamax}{\lambdamin}  \sqrt{ \edit{d_0} \log\left( N_{\max} \sum_{k=1}^K \abs{\X_k}  \right) }  
	    \cdot  \frac{\log^{3/2}(K)\log^{2}(1/\delta) }{\sqrt K}.
\end{align*}

\end{lemma}
\proof{Proof.}
Using \cref{lem:Envelopes,thm:SizeOfDiscreteSetsGeneral} to bound the variable $J$ in \cref{eq:JDudley}  and since $\prns{\sum_{k=1}^K\abs{\X_k}^2}^{\edit{d_0}}\leq \prns{\sum_{k=1}^K\abs{\X_k}}^{2\edit{d_0}}$, proves
\[
\| J\| _p  \ \leq  \ 9 \Cmax \frac{\lambdamax}{\lambdamin} \sqrt K \sqrt{ 2\edit{d_0}\log\left(\sum_{k=1}^K \abs{\X_k} \right) }.
\]
{\blockedit  \label{ReviewerClarity}
Next apply \cref{thm:pollard} and let $p = 2 \log(1/\delta)$ to prove the first statement.   

For the second inequality, we follow a similar strategy with the appropriate envelope and packing number (cf. \cref{lem:Envelopes,thm:SizeOfDiscreteSetsGeneral}).  In this case, 
\[
J \ \leq \ \frac{9 C \Nhat_{\max}\sqrt K }{N \lambdabar} \cdot \sqrt{ \log\left( 1 + \Nhat_{\max}^{d_0} \left( \sum_{k=1}^K \abs{\X_k} \right)^{2d_0} \right) }.
\]
Consider the inner logarithm, and let $\X_{\rm{tot}} \equiv \sum_{k=1}^K \abs{\X_k} \geq K \geq 2$.  Then, 
\begin{align*}
\log\left( 1 + \Nhat_{\max}^{d_0} \X^{2d_0}_{\rm tot }\right) 
&\ \leq \ 
d_0 \log\left(1 + \Nhat_{\max} \X^{2}_{\rm tot }\right)  && (\text{since } \Nhat_{\max} \X_{\rm tot}^2 > 0)
\\
& \ \leq d_0 \log\left( 3 N_{\max} \X^{2}_{\rm tot } + \Nhat_{\max} \X^{2}_{\rm tot }\right) 
\\
& \ \leq \ d_0 \left( 
\underbrace{\log\left( N_{\max} \X^{2}_{\rm tot }\right)}_{\geq \log(4)}  +  \underbrace{\log \left( 3  + \frac{\Nhat_{\max}}{N_{\max}} \right)}_{\geq \log(3)}
\right)
\\
& \ \leq \ 
2d_0  \log\left( N_{\max} \X^{2}_{\rm tot }\right) \cdot 
 \log \left( 3  + \frac{\Nhat_{\max}}{N_{\max}} \right),
 \end{align*}
where the last inequality follows because $a + b \leq 2 a b$ for $a, b \geq 1$.  

Substituting above shows 
\begin{align*}
\| J\|_p  
&\ \leq \ \frac{9 C \sqrt K }{N \lambdabar} \cdot \sqrt{ 2 d_0  \log\left( N_{\max} \X^{2}_{\rm tot }\right) }  \magd{ \Nhat_{\max} \cdot  \sqrt{ 
 \log \left( 3  + \frac{\Nhat_{\max}}{N_{\max}} \right)
 } }_p .
\\
&\ \leq \ \frac{9 C \sqrt K }{N \lambdabar} \cdot \sqrt{ 2 d_0 \log\left( N_{\max} \X^{2}_{\rm tot }\right) } \cdot  \magd{ \Nhat_{\max} }_{2p} \cdot  \magd{\sqrt{ 
 \log \left( 3  + \frac{\Nhat_{\max}}{N_{\max}} \right)
 } }_{2p},
\end{align*}

We next bound these norms.  The first is bounded by \cref{lem:PropertiesOfPoisson} Part~\ref{PNormNhatmax}.  The second was bounded in the proof of \cref{lem:ULLNDiscrete} as
\begin{align*}
\magd{\sqrt{ \log \left( 3 + \frac{\Nhat_{\max}}{N_{\max}} \right) }}_{2p} 
& \ \leq 2 \cdot \ 5^{\frac{1}{2p}} \sqrt{p}\sqrt{\log K }. 
\end{align*}
Combining proves
\[
\| J\|_p  
\ \leq \ \const_3 \cdot 
\Cmax \cdot   \frac{\lambdamax}{\lambdamin}  \sqrt{ d_0 \log\left( N_{\max} \X^{2}_{\rm tot }\right) } \cdot  
	6^{\frac{1}{p}} p^{3/2}   \cdot  \log^{3/2}(K) \sqrt K,
\]
for some constant $\const_3$.  
Now apply \cref{thm:pollard} and substitute $p = 2 \log(1/\delta)$ to prove the second inequality.  
\endproof
}

%

\subsubsection{Proofs of \cref{thm:ShrinkageDiscreteGeneralAnchor,thm:ShrinkageDiscreteLooAnchor}.}
{\blockedit We can now prove the main results of the section via our usual strategy.  
\begin{proof}{Proof of \cref{thm:ShrinkageDiscreteGeneralAnchor,thm:ShrinkageDiscreteLooAnchor}.}  
The proofs of both theorems are identical.
For both theorems, by \cref{lem:ConditionsForOptimality}, the quantity to be bounded is bounded by the sum of the same two maximal deviations. These are in turn bounded by \cref{lem:ULLNDiscreteGeneral}. Instantiating each bound for $\delta \leftarrow \delta/2$, adding the right hand sides and applying the union bound yields a bound on the sub-optimality.  Collecting dominant terms yields the result.  
\end{proof}}

\section{Contrasting the Sub-Optimality-Stability Bias-Variance Tradeoffs}
\label{sec:BiasVariance}
We here expand on the discussion from \cref{sec:StabilityIntuition} comparing the Sub-Optimality-Stability tradeoff to the classic bias-variance tradeoff. As mentioned in \cref{sec:StabilityIntuition}, one important distinction is that the former applies to general optimization problems. In the following we will show that they are different even when we restrict to the case of MSE (cf. \cref{ex:MSE}).

To be more precise, fix the cost functions $c_k(x, \xi) = (x - \xi)^2$, let $\mu_k$ and $\sigma^2_k$ denote the mean and variance of $\xi_k \in \Rl$ and assume $\lambda_k = 1$ for all $k$ for simplicity.  There are at least two ways to interpret the classical bias-variance tradeoff in context of \cref{ass:RandomData}. First, we can decompose conditionally on $\bm{\Nhat}$, yielding
\begin{align*}
\Eb{ \Zperf(\alpha, \bp_0) \mid \bm{ \Nhat} } \ \  = \ \  \frac{1}{K} \sum_{k=1}^K \underbrace{\left(\frac{\alpha}{\Nhat_k + \alpha}\right)^2 (\mu_k - \mu_{k0})^2 }_\text{Conditional Bias Squared} \ \ + \ \  
                                    \underbrace{\left(\frac{\Nhat_k}{\Nhat_k + \alpha}\right)^2 \frac{\sigma_k^2}{\Nhat_k} }_{\text{Conditional Variance}}, 
\end{align*}
where $\mu_{k0} = \bp_0^\top \ba_k$.
Taking expectations of both sides yields the identity for $\alpha > 0$
\begin{equation} \label{eq:CondBiasVariance}
\Eb{ \Zperf(\alpha, \bp_0)} \ \ = \ \ \underbrace{\frac{1}{K} \sum_{k=1}^K \Eb{\left(\frac{\alpha}{\Nhat_k + \alpha}\right)^2} (\mu_k - \mu_{k0})^2}_{\text{Expected Conditional Bias Squared}}  \ \ +\ \  
            \underbrace{\Eb{\frac{\Nhat_k}{(\Nhat_k + \alpha)^2}} \sigma_k^2}_{\text{Expected Conditional Variance}} . 
\end{equation}
This perspective is perhaps most appropriate if view \cref{ass:RandomData} as a smoothing that randomizes over instances.  

Alternatively, we can apply the bias-variance decomposition unconditionally, yielding for $\alpha > 0$,
\begin{align} \notag
\Eb{ \Zperf(\alpha, \bp_0) } 
& = 
 \frac{1}{K} \sum_{k=1}^K \left( \Eb{x_k(\alpha, \bp_0, \muhat_k) - \mu_k} \right)^2 + \text{Var}(x_k(\alpha, \bp_0, \muhat_k) ),
\\  \label{eq:UncondBiasVariance}
& = 
 \frac{1}{K} \sum_{k=1}^K 
    \underbrace{ \left( \Eb{ \frac{\alpha}{\Nhat_k + \alpha}} \right)^2 (\mu_{0k} - \mu_k)^2}_{\text{Bias Squared}}  + \underbrace{\text{Var}(x_k(\alpha, \bp_0, \muhat_k) )}_{\text{Variance}},
\end{align}
(We can, if desired, evaluate the second term using the law of total variance after conditioning on $\Nhat_k$, but this expression will not be needed in what follows.)  This perspective is perhaps most appropriate if we view the randomization of $\Nhat_k$ as intrinsic to the data-generating process.  

Finally, from \cref{cor:steinchen} and our previous comments, we have that 
\begin{align*}
\Eb{ \Zperf(\alpha, \bp_0) } = \frac{1}{N \lambdabar} \left( \Eb{\text{SAA-SubOptimality}(\alpha) } +                                                           \Eb{\text{Instability}(\alpha) }    
 + \Eb{\text{SAA}(0)  }
\right),
\end{align*}
where, again, $\text{SAA}(0)$ does not depend on $\alpha$.  A straightforward calculation yields, 
\begin{lemma}[SAA-Sub-Optimality for MSE] \label{lem:SAASubOpt}
For $\alpha > 0$, we have
\begin{align*}
\op{SAA-SubOpt}(\alpha) &=     \frac{1}{K} \sum_{k=1}^K \Nhat_k \left( \frac{\alpha}{\Nhat_k + \alpha}\right)^2 ( \muhat_k - \mu_{k0} )^2
\\
\Eb{\op{SAA-SubOpt}(\alpha) } &=     \frac{1}{K}\sum_{k=1}^K \Eb{ \Nhat_k \left( \frac{\alpha}{\Nhat_k + \alpha}\right)^2}  ( \mu_k - \mu_{k0} )^2 + \frac{1}{K}\sum_{k=1}^K \Eb{\left( \frac{\alpha}{\Nhat_k + \alpha}\right)^2} \sigma_k^2,
\end{align*}
where $\muhat_k$ is the sample mean for the $k^\text{th}$ subproblem.  
\end{lemma}
\begin{proof}{Proof of \cref{lem:SAASubOpt}}
By definition, the $k^\text{th}$ term of $\op{SAA-SubOpt}(\alpha)$ is 
\begin{align*}
\sum_{i=1}^d \mhat_{ki}\left( c_{ki}( x_k(\alpha, \bp_0, \bfmhat_k) ) - c_{ki}(x_k(0, \bp_0, \bfmhat_k) ) \right)
& = 
\Nhat_k \sum_{i=1}^d \bphat_{ki}\left( c_{ki}( x_k(\alpha, \bp_0, \bfmhat_k) ) - c_{ki}(x_k(0, \bp_0, \bfmhat_k) ) \right)
\\
& = \Nhat_k \left(  \Eb{ (\xihat_k - \muhat_k(\alpha) )^2 \mid \bfmhat_k }  + \Eb{ ( \xihat_k - \muhat_k)^2 \mid \bfmhat_k }\right)
\end{align*}
where $\bx_k(\alpha, \bp_0, \bfmhat_k) = \muhat_k(\alpha) \equiv \frac{\alpha}{\Nhat_k + \alpha} \mu_{k0} + \frac{\Nhat_k}{\Nhat_k + \alpha} \muhat_k$, and $\xihat_k \sim \bphat_k$.  

Note 
\(
\Eb{ (\xihat_k - \muhat_k(\alpha) )^2 \mid \bfmhat_k }  = (\muhat_k - \muhat_k(\alpha))^2 + \hat{\sigma}_k^2,
\)
where $\hat{\sigma}_k^2$ is the variance of $\xihat_k \mid \bfmhat_k$.  Similarly, 
\(
\Eb{ (\xihat_k - \muhat_k )^2 \mid \bfmhat_k }  = \hat{\sigma}_k^2.
\)
Hence from above, the $k^\text{th}$ term of $\op{SAA-SubOpt}(\alpha)$ is 
\(
\Nhat_k (\muhat_k - \muhat_k(\alpha))^2.
\)
Using the definition of $\muhat_k(\alpha)$ we have
\(
(\muhat_k - \muhat_k(\alpha))^2 = \left(\frac{\alpha}{\Nhat_k + \alpha}\right)^2 (\mu_0 - \muhat_k)^2.
\)
Summing across the $k$ terms yields the expression for $\op{SAA-SubOpt}(\alpha)$ in the lemma.  

Now consider taking the conditional expectation of the $k^\text{th}$ term of $\op{SAA-SubOpt}(\alpha)$ where we condition on $\bm{\Nhat}$.  From our previous expression, this is simply
\begin{align*}
\Nhat_k \left(\frac{\alpha}{\Nhat_k + \alpha}\right)^2 \Eb{(\mu_0 - \muhat_k)^2 \mid \bm{\Nhat} } 
&= 
\Nhat_k \left(\frac{\alpha}{\Nhat_k + \alpha}\right)^2 \left( ( \mu_0 - \mu_k)^2 + \frac{\sigma_k^2 }{\Nhat_k} \right).
\\ &=
\Nhat_k \left(\frac{\alpha}{\Nhat_k + \alpha}\right)^2  ( \mu_0 - \mu_k)^2 + \left(\frac{\alpha}{\Nhat_k + \alpha}\right)^2 \sigma_k^2.
\end{align*}
Taking expectations and then averaging over $k$ yields the expression for $\Eb{\op{SAA-SubOpt}(\alpha)}$, completing the lemma.  
\end{proof}

By inspection, $\frac{1}{N \lambdabar} \E[ \op{SAA-SubOpt}(\alpha)]$ involves a non-zero term that depends on both $\sigma_k^2$ and $\alpha$.  Consequently, it must differ from the bias-squared term in \cref{eq:UncondBiasVariance} and the expected conditional bias-squared term in \cref{eq:CondBiasVariance}.  In particular, since the difference  depends on $\alpha$ and $\op{SAA}(0)$ does not depend on $\alpha$, the difference is not solely due to the treatment of this constant.   Finally, since each of the identities decomposes the same quantity $\Eb{\Zperf{\alpha, \bp_0}}$, it follows that the bias-variance tradeoff and the Sub-Optimality-Instability Tradeoff are fundamentally different for this example.


\section{Computational Details and Additional Numerical Experiments}
\label{sec:AdditionalFigures}

\subsection{Simulation Set-up for \cref{fig:IllustratingDataPooling} }
\label{sec:SimSetUp}
For $d = 10$, we generate $5,000$ distributions $\bp_k$ according to a uniform distribution on the simplex and additional $5,000$ distributions $\bp_k$ according to the Dirichlet distribution with parameter $(3, \ldots, 3)$, for a total of $K=10,000$ subproblems.  We take $\lambda_k = 1$ for all $k$.  Across all runs, these $\bp_k$ and $\lambda_k$ are fixed.  Then, for each run, for each $k$, we then generate $\Nhat_k = 20$ data points independently according to \cref{eq:MultinomialCounts}.  We train each of our policies on these data, and evaluate against the true $\bp_k$.  Results are averaged across $10,000$ runs.

\subsection{Additional Figures from \cref{ex:BabyNewsvendorTake2}.}
\Cref{fig:TradeOffLOO} shows the companion figures for \cref{ex:BabyNewsvendorTake2} from \cref{sec:StabilityIntuition}.

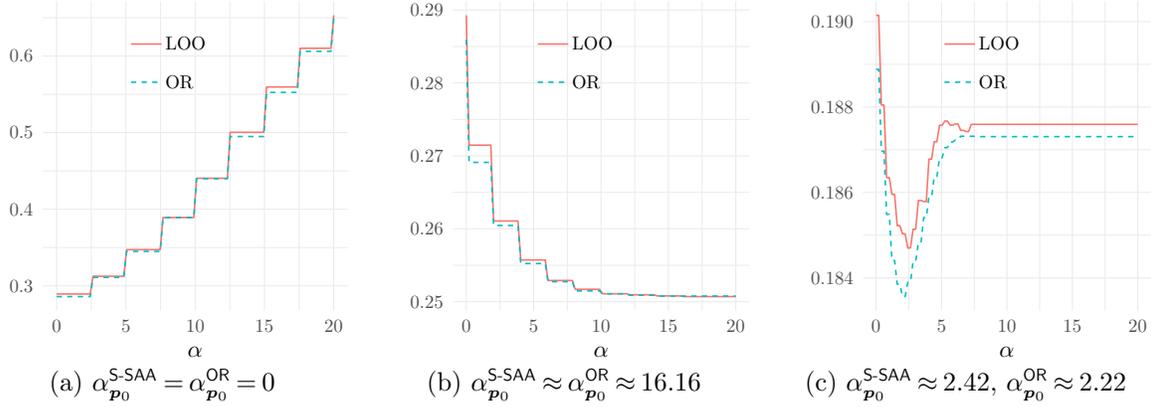
\begin{figure}
    \centering
    \begin{subfigure}[b]{0.3\textwidth}
\begin{tikzpicture}[x=1pt,y=1pt]
\definecolor{fillColor}{RGB}{255,255,255}
\path[use as bounding box,fill=fillColor,fill opacity=0.00] (0,0) rectangle (144.54,144.54);
\begin{scope}
\path[clip] ( 25.23, 23.41) rectangle (140.54,140.54);
\definecolor{drawColor}{gray}{0.92}

\path[draw=drawColor,line width= 0.2pt,line join=round] ( 25.23, 47.31) --
	(140.54, 47.31);

\path[draw=drawColor,line width= 0.2pt,line join=round] ( 25.23, 76.30) --
	(140.54, 76.30);

\path[draw=drawColor,line width= 0.2pt,line join=round] ( 25.23,105.30) --
	(140.54,105.30);

\path[draw=drawColor,line width= 0.2pt,line join=round] ( 25.23,134.30) --
	(140.54,134.30);

\path[draw=drawColor,line width= 0.2pt,line join=round] ( 43.57, 23.41) --
	( 43.57,140.54);

\path[draw=drawColor,line width= 0.2pt,line join=round] ( 69.78, 23.41) --
	( 69.78,140.54);

\path[draw=drawColor,line width= 0.2pt,line join=round] ( 95.99, 23.41) --
	( 95.99,140.54);

\path[draw=drawColor,line width= 0.2pt,line join=round] (122.20, 23.41) --
	(122.20,140.54);

\path[draw=drawColor,line width= 0.4pt,line join=round] ( 25.23, 32.81) --
	(140.54, 32.81);

\path[draw=drawColor,line width= 0.4pt,line join=round] ( 25.23, 61.81) --
	(140.54, 61.81);

\path[draw=drawColor,line width= 0.4pt,line join=round] ( 25.23, 90.80) --
	(140.54, 90.80);

\path[draw=drawColor,line width= 0.4pt,line join=round] ( 25.23,119.80) --
	(140.54,119.80);

\path[draw=drawColor,line width= 0.4pt,line join=round] ( 30.47, 23.41) --
	( 30.47,140.54);

\path[draw=drawColor,line width= 0.4pt,line join=round] ( 56.68, 23.41) --
	( 56.68,140.54);

\path[draw=drawColor,line width= 0.4pt,line join=round] ( 82.88, 23.41) --
	( 82.88,140.54);

\path[draw=drawColor,line width= 0.4pt,line join=round] (109.09, 23.41) --
	(109.09,140.54);

\path[draw=drawColor,line width= 0.4pt,line join=round] (135.30, 23.41) --
	(135.30,140.54);
\definecolor{drawColor}{RGB}{248,118,109}

\path[draw=drawColor,line width= 0.6pt,line join=round] ( 30.47, 29.71) --
	( 31.53, 29.71) --
	( 32.59, 29.71) --
	( 33.65, 29.71) --
	( 34.71, 29.71) --
	( 35.77, 29.71) --
	( 36.82, 29.71) --
	( 37.88, 29.71) --
	( 38.94, 29.71) --
	( 40.00, 29.71) --
	( 41.06, 29.71) --
	( 42.12, 29.71) --
	( 43.18, 29.71) --
	( 44.24, 36.44) --
	( 45.30, 36.44) --
	( 46.35, 36.44) --
	( 47.41, 36.44) --
	( 48.47, 36.44) --
	( 49.53, 36.44) --
	( 50.59, 36.44) --
	( 51.65, 36.44) --
	( 52.71, 36.44) --
	( 53.77, 36.44) --
	( 54.82, 36.44) --
	( 55.88, 36.44) --
	( 56.94, 46.50) --
	( 58.00, 46.50) --
	( 59.06, 46.50) --
	( 60.12, 46.50) --
	( 61.18, 46.50) --
	( 62.24, 46.50) --
	( 63.30, 46.50) --
	( 64.35, 46.50) --
	( 65.41, 46.50) --
	( 66.47, 46.50) --
	( 67.53, 46.50) --
	( 68.59, 46.50) --
	( 69.65, 46.50) --
	( 70.71, 58.65) --
	( 71.77, 58.65) --
	( 72.83, 58.65) --
	( 73.88, 58.65) --
	( 74.94, 58.65) --
	( 76.00, 58.65) --
	( 77.06, 58.65) --
	( 78.12, 58.65) --
	( 79.18, 58.65) --
	( 80.24, 58.65) --
	( 81.30, 58.65) --
	( 82.36, 58.65) --
	( 83.41, 73.53) --
	( 84.47, 73.53) --
	( 85.53, 73.53) --
	( 86.59, 73.53) --
	( 87.65, 73.53) --
	( 88.71, 73.53) --
	( 89.77, 73.53) --
	( 90.83, 73.53) --
	( 91.89, 73.53) --
	( 92.94, 73.53) --
	( 94.00, 73.53) --
	( 95.06, 73.53) --
	( 96.12, 90.92) --
	( 97.18, 90.92) --
	( 98.24, 90.92) --
	( 99.30, 90.92) --
	(100.36, 90.92) --
	(101.41, 90.92) --
	(102.47, 90.92) --
	(103.53, 90.92) --
	(104.59, 90.92) --
	(105.65, 90.92) --
	(106.71, 90.92) --
	(107.77, 90.92) --
	(108.83, 90.92) --
	(109.89,108.08) --
	(110.94,108.08) --
	(112.00,108.08) --
	(113.06,108.08) --
	(114.12,108.08) --
	(115.18,108.08) --
	(116.24,108.08) --
	(117.30,108.08) --
	(118.36,108.08) --
	(119.42,108.08) --
	(120.47,108.08) --
	(121.53,108.08) --
	(122.59,122.73) --
	(123.65,122.73) --
	(124.71,122.73) --
	(125.77,122.73) --
	(126.83,122.73) --
	(127.89,122.73) --
	(128.95,122.73) --
	(130.00,122.73) --
	(131.06,122.73) --
	(132.12,122.73) --
	(133.18,122.73) --
	(134.24,122.73) --
	(135.30,135.22);
\definecolor{drawColor}{RGB}{0,191,196}

\path[draw=drawColor,line width= 0.6pt,dash pattern=on 2pt off 2pt ,line join=round] ( 30.47, 28.73) --
	( 31.53, 28.73) --
	( 32.59, 28.73) --
	( 33.65, 28.73) --
	( 34.71, 28.73) --
	( 35.77, 28.73) --
	( 36.82, 28.73) --
	( 37.88, 28.73) --
	( 38.94, 28.73) --
	( 40.00, 28.73) --
	( 41.06, 28.73) --
	( 42.12, 28.73) --
	( 43.18, 28.73) --
	( 44.24, 35.99) --
	( 45.30, 35.99) --
	( 46.35, 35.99) --
	( 47.41, 35.99) --
	( 48.47, 35.99) --
	( 49.53, 35.99) --
	( 50.59, 35.99) --
	( 51.65, 35.99) --
	( 52.71, 35.99) --
	( 53.77, 35.99) --
	( 54.82, 35.99) --
	( 55.88, 35.99) --
	( 56.94, 45.82) --
	( 58.00, 45.82) --
	( 59.06, 45.82) --
	( 60.12, 45.82) --
	( 61.18, 45.82) --
	( 62.24, 45.82) --
	( 63.30, 45.82) --
	( 64.35, 45.82) --
	( 65.41, 45.82) --
	( 66.47, 45.82) --
	( 67.53, 45.82) --
	( 68.59, 45.82) --
	( 69.65, 45.82) --
	( 70.71, 58.64) --
	( 71.77, 58.64) --
	( 72.83, 58.64) --
	( 73.88, 58.64) --
	( 74.94, 58.64) --
	( 76.00, 58.64) --
	( 77.06, 58.64) --
	( 78.12, 58.64) --
	( 79.18, 58.64) --
	( 80.24, 58.64) --
	( 81.30, 58.64) --
	( 82.36, 58.64) --
	( 83.41, 73.37) --
	( 84.47, 73.37) --
	( 85.53, 73.37) --
	( 86.59, 73.37) --
	( 87.65, 73.37) --
	( 88.71, 73.37) --
	( 89.77, 73.37) --
	( 90.83, 73.37) --
	( 91.89, 73.37) --
	( 92.94, 73.37) --
	( 94.00, 73.37) --
	( 95.06, 73.37) --
	( 96.12, 89.34) --
	( 97.18, 89.34) --
	( 98.24, 89.34) --
	( 99.30, 89.34) --
	(100.36, 89.34) --
	(101.41, 89.34) --
	(102.47, 89.34) --
	(103.53, 89.34) --
	(104.59, 89.34) --
	(105.65, 89.34) --
	(106.71, 89.34) --
	(107.77, 89.34) --
	(108.83, 89.34) --
	(109.89,106.03) --
	(110.94,106.03) --
	(112.00,106.03) --
	(113.06,106.03) --
	(114.12,106.03) --
	(115.18,106.03) --
	(116.24,106.03) --
	(117.30,106.03) --
	(118.36,106.03) --
	(119.42,106.03) --
	(120.47,106.03) --
	(121.53,106.03) --
	(122.59,121.58) --
	(123.65,121.58) --
	(124.71,121.58) --
	(125.77,121.58) --
	(126.83,121.58) --
	(127.89,121.58) --
	(128.95,121.58) --
	(130.00,121.58) --
	(131.06,121.58) --
	(132.12,121.58) --
	(133.18,121.58) --
	(134.24,121.58) --
	(135.30,134.01);
\end{scope}
\begin{scope}
\path[clip] (  0.00,  0.00) rectangle (144.54,144.54);
\definecolor{drawColor}{gray}{0.30}

\node[text=drawColor,anchor=base east,inner sep=0pt, outer sep=0pt, scale=  0.64] at ( 21.63, 30.60) {0.3};

\node[text=drawColor,anchor=base east,inner sep=0pt, outer sep=0pt, scale=  0.64] at ( 21.63, 59.60) {0.4};

\node[text=drawColor,anchor=base east,inner sep=0pt, outer sep=0pt, scale=  0.64] at ( 21.63, 88.60) {0.5};

\node[text=drawColor,anchor=base east,inner sep=0pt, outer sep=0pt, scale=  0.64] at ( 21.63,117.59) {0.6};
\end{scope}
\begin{scope}
\path[clip] (  0.00,  0.00) rectangle (144.54,144.54);
\definecolor{drawColor}{gray}{0.30}

\node[text=drawColor,anchor=base,inner sep=0pt, outer sep=0pt, scale=  0.64] at ( 30.47, 15.40) {0};

\node[text=drawColor,anchor=base,inner sep=0pt, outer sep=0pt, scale=  0.64] at ( 56.68, 15.40) {5};

\node[text=drawColor,anchor=base,inner sep=0pt, outer sep=0pt, scale=  0.64] at ( 82.88, 15.40) {10};

\node[text=drawColor,anchor=base,inner sep=0pt, outer sep=0pt, scale=  0.64] at (109.09, 15.40) {15};

\node[text=drawColor,anchor=base,inner sep=0pt, outer sep=0pt, scale=  0.64] at (135.30, 15.40) {20};
\end{scope}
\begin{scope}
\path[clip] (  0.00,  0.00) rectangle (144.54,144.54);
\definecolor{drawColor}{RGB}{0,0,0}

\node[text=drawColor,anchor=base,inner sep=0pt, outer sep=0pt, scale=  0.80] at ( 82.88,  5.94) {$\alpha$};
\end{scope}
\begin{scope}
\path[clip] (  0.00,  0.00) rectangle (144.54,144.54);
\definecolor{drawColor}{RGB}{248,118,109}

\path[draw=drawColor,line width= 0.6pt,line join=round] ( 58.60,124.34) -- ( 70.16,124.34);
\end{scope}
\begin{scope}
\path[clip] (  0.00,  0.00) rectangle (144.54,144.54);
\definecolor{drawColor}{RGB}{0,191,196}

\path[draw=drawColor,line width= 0.6pt,dash pattern=on 2pt off 2pt ,line join=round] ( 58.60,109.89) -- ( 70.16,109.89);
\end{scope}
\begin{scope}
\path[clip] (  0.00,  0.00) rectangle (144.54,144.54);
\definecolor{drawColor}{RGB}{0,0,0}

\node[text=drawColor,anchor=base west,inner sep=0pt, outer sep=0pt, scale=  0.64] at ( 71.60,122.14) {LOO};
\end{scope}
\begin{scope}
\path[clip] (  0.00,  0.00) rectangle (144.54,144.54);
\definecolor{drawColor}{RGB}{0,0,0}

\node[text=drawColor,anchor=base west,inner sep=0pt, outer sep=0pt, scale=  0.64] at ( 71.60,107.68) {OR};
\end{scope}
\end{tikzpicture} \vspace{-20pt}
        \caption{$\alphaLOO_{\bp_0} = \alphaOR_{\bp_0} = 0$}
    \end{subfigure}
    ~ 
    \begin{subfigure}[b]{0.3\textwidth}
\begin{tikzpicture}[x=1pt,y=1pt]
\definecolor{fillColor}{RGB}{255,255,255}
\path[use as bounding box,fill=fillColor,fill opacity=0.00] (0,0) rectangle (144.54,144.54);
\begin{scope}
\path[clip] ( 28.43, 23.41) rectangle (140.54,140.54);
\definecolor{drawColor}{gray}{0.92}

\path[draw=drawColor,line width= 0.2pt,line join=round] ( 28.43, 40.60) --
	(140.54, 40.60);

\path[draw=drawColor,line width= 0.2pt,line join=round] ( 28.43, 68.20) --
	(140.54, 68.20);

\path[draw=drawColor,line width= 0.2pt,line join=round] ( 28.43, 95.80) --
	(140.54, 95.80);

\path[draw=drawColor,line width= 0.2pt,line join=round] ( 28.43,123.40) --
	(140.54,123.40);

\path[draw=drawColor,line width= 0.2pt,line join=round] ( 46.26, 23.41) --
	( 46.26,140.54);

\path[draw=drawColor,line width= 0.2pt,line join=round] ( 71.74, 23.41) --
	( 71.74,140.54);

\path[draw=drawColor,line width= 0.2pt,line join=round] ( 97.22, 23.41) --
	( 97.22,140.54);

\path[draw=drawColor,line width= 0.2pt,line join=round] (122.70, 23.41) --
	(122.70,140.54);

\path[draw=drawColor,line width= 0.4pt,line join=round] ( 28.43, 26.80) --
	(140.54, 26.80);

\path[draw=drawColor,line width= 0.4pt,line join=round] ( 28.43, 54.40) --
	(140.54, 54.40);

\path[draw=drawColor,line width= 0.4pt,line join=round] ( 28.43, 82.00) --
	(140.54, 82.00);

\path[draw=drawColor,line width= 0.4pt,line join=round] ( 28.43,109.60) --
	(140.54,109.60);

\path[draw=drawColor,line width= 0.4pt,line join=round] ( 28.43,137.20) --
	(140.54,137.20);

\path[draw=drawColor,line width= 0.4pt,line join=round] ( 33.52, 23.41) --
	( 33.52,140.54);

\path[draw=drawColor,line width= 0.4pt,line join=round] ( 59.00, 23.41) --
	( 59.00,140.54);

\path[draw=drawColor,line width= 0.4pt,line join=round] ( 84.48, 23.41) --
	( 84.48,140.54);

\path[draw=drawColor,line width= 0.4pt,line join=round] (109.96, 23.41) --
	(109.96,140.54);

\path[draw=drawColor,line width= 0.4pt,line join=round] (135.44, 23.41) --
	(135.44,140.54);
\definecolor{drawColor}{RGB}{248,118,109}

\path[draw=drawColor,line width= 0.6pt,line join=round] ( 33.52,135.22) --
	( 34.55, 86.06) --
	( 35.58, 86.06) --
	( 36.61, 86.06) --
	( 37.64, 86.06) --
	( 38.67, 86.06) --
	( 39.70, 86.06) --
	( 40.73, 86.06) --
	( 41.76, 86.06) --
	( 42.79, 86.06) --
	( 43.82, 57.35) --
	( 44.85, 57.35) --
	( 45.88, 57.35) --
	( 46.91, 57.35) --
	( 47.94, 57.35) --
	( 48.97, 57.35) --
	( 50.00, 57.35) --
	( 51.03, 57.35) --
	( 52.06, 57.35) --
	( 53.08, 57.35) --
	( 54.11, 42.59) --
	( 55.14, 42.59) --
	( 56.17, 42.59) --
	( 57.20, 42.59) --
	( 58.23, 42.59) --
	( 59.26, 42.59) --
	( 60.29, 42.59) --
	( 61.32, 42.59) --
	( 62.35, 42.59) --
	( 63.38, 42.59) --
	( 64.41, 34.83) --
	( 65.44, 34.83) --
	( 66.47, 34.83) --
	( 67.50, 34.83) --
	( 68.53, 34.83) --
	( 69.56, 34.83) --
	( 70.59, 34.83) --
	( 71.62, 34.83) --
	( 72.65, 34.83) --
	( 73.67, 34.83) --
	( 74.70, 31.49) --
	( 75.73, 31.49) --
	( 76.76, 31.49) --
	( 77.79, 31.49) --
	( 78.82, 31.49) --
	( 79.85, 31.49) --
	( 80.88, 31.49) --
	( 81.91, 31.49) --
	( 82.94, 31.49) --
	( 83.97, 31.49) --
	( 85.00, 29.81) --
	( 86.03, 29.81) --
	( 87.06, 29.81) --
	( 88.09, 29.81) --
	( 89.12, 29.81) --
	( 90.15, 29.81) --
	( 91.18, 29.81) --
	( 92.21, 29.81) --
	( 93.24, 29.81) --
	( 94.26, 29.81) --
	( 95.29, 29.42) --
	( 96.32, 29.42) --
	( 97.35, 29.42) --
	( 98.38, 29.42) --
	( 99.41, 29.42) --
	(100.44, 29.42) --
	(101.47, 29.42) --
	(102.50, 29.42) --
	(103.53, 29.42) --
	(104.56, 29.42) --
	(105.59, 28.98) --
	(106.62, 28.98) --
	(107.65, 28.98) --
	(108.68, 28.98) --
	(109.71, 28.98) --
	(110.74, 28.98) --
	(111.77, 28.98) --
	(112.80, 28.98) --
	(113.82, 28.98) --
	(114.85, 28.98) --
	(115.88, 28.73) --
	(116.91, 28.73) --
	(117.94, 28.73) --
	(118.97, 28.73) --
	(120.00, 28.73) --
	(121.03, 28.73) --
	(122.06, 28.73) --
	(123.09, 28.73) --
	(124.12, 28.73) --
	(125.15, 28.73) --
	(126.18, 28.73) --
	(127.21, 28.73) --
	(128.24, 28.73) --
	(129.27, 28.73) --
	(130.30, 28.73) --
	(131.33, 28.73) --
	(132.36, 28.73) --
	(133.39, 28.73) --
	(134.41, 28.73) --
	(135.44, 28.73);
\definecolor{drawColor}{RGB}{0,191,196}

\path[draw=drawColor,line width= 0.6pt,dash pattern=on 2pt off 2pt ,line join=round] ( 33.52,125.98) --
	( 34.55, 79.53) --
	( 35.58, 79.53) --
	( 36.61, 79.53) --
	( 37.64, 79.53) --
	( 38.67, 79.53) --
	( 39.70, 79.53) --
	( 40.73, 79.53) --
	( 41.76, 79.53) --
	( 42.79, 79.53) --
	( 43.82, 55.66) --
	( 44.85, 55.66) --
	( 45.88, 55.66) --
	( 46.91, 55.66) --
	( 47.94, 55.66) --
	( 48.97, 55.66) --
	( 50.00, 55.66) --
	( 51.03, 55.66) --
	( 52.06, 55.66) --
	( 53.08, 55.66) --
	( 54.11, 41.26) --
	( 55.14, 41.26) --
	( 56.17, 41.26) --
	( 57.20, 41.26) --
	( 58.23, 41.26) --
	( 59.26, 41.26) --
	( 60.29, 41.26) --
	( 61.32, 41.26) --
	( 62.35, 41.26) --
	( 63.38, 41.26) --
	( 64.41, 34.39) --
	( 65.44, 34.39) --
	( 66.47, 34.39) --
	( 67.50, 34.39) --
	( 68.53, 34.39) --
	( 69.56, 34.39) --
	( 70.59, 34.39) --
	( 71.62, 34.39) --
	( 72.65, 34.39) --
	( 73.67, 34.39) --
	( 74.70, 30.88) --
	( 75.73, 30.88) --
	( 76.76, 30.88) --
	( 77.79, 30.88) --
	( 78.82, 30.88) --
	( 79.85, 30.88) --
	( 80.88, 30.88) --
	( 81.91, 30.88) --
	( 82.94, 30.88) --
	( 83.97, 30.88) --
	( 85.00, 29.79) --
	( 86.03, 29.79) --
	( 87.06, 29.79) --
	( 88.09, 29.79) --
	( 89.12, 29.79) --
	( 90.15, 29.79) --
	( 91.18, 29.79) --
	( 92.21, 29.79) --
	( 93.24, 29.79) --
	( 94.26, 29.79) --
	( 95.29, 29.23) --
	( 96.32, 29.23) --
	( 97.35, 29.23) --
	( 98.38, 29.23) --
	( 99.41, 29.23) --
	(100.44, 29.23) --
	(101.47, 29.23) --
	(102.50, 29.23) --
	(103.53, 29.23) --
	(104.56, 29.23) --
	(105.59, 28.93) --
	(106.62, 28.93) --
	(107.65, 28.93) --
	(108.68, 28.93) --
	(109.71, 28.93) --
	(110.74, 28.93) --
	(111.77, 28.93) --
	(112.80, 28.93) --
	(113.82, 28.93) --
	(114.85, 28.93) --
	(115.88, 28.93) --
	(116.91, 28.93) --
	(117.94, 28.93) --
	(118.97, 28.93) --
	(120.00, 28.93) --
	(121.03, 28.93) --
	(122.06, 28.93) --
	(123.09, 28.93) --
	(124.12, 28.93) --
	(125.15, 28.93) --
	(126.18, 28.93) --
	(127.21, 28.93) --
	(128.24, 28.93) --
	(129.27, 28.93) --
	(130.30, 28.93) --
	(131.33, 28.93) --
	(132.36, 28.93) --
	(133.39, 28.93) --
	(134.41, 28.93) --
	(135.44, 28.93);
\end{scope}
\begin{scope}
\path[clip] (  0.00,  0.00) rectangle (144.54,144.54);
\definecolor{drawColor}{gray}{0.30}

\node[text=drawColor,anchor=base east,inner sep=0pt, outer sep=0pt, scale=  0.64] at ( 24.83, 24.59) {0.25};

\node[text=drawColor,anchor=base east,inner sep=0pt, outer sep=0pt, scale=  0.64] at ( 24.83, 52.20) {0.26};

\node[text=drawColor,anchor=base east,inner sep=0pt, outer sep=0pt, scale=  0.64] at ( 24.83, 79.80) {0.27};

\node[text=drawColor,anchor=base east,inner sep=0pt, outer sep=0pt, scale=  0.64] at ( 24.83,107.40) {0.28};

\node[text=drawColor,anchor=base east,inner sep=0pt, outer sep=0pt, scale=  0.64] at ( 24.83,135.00) {0.29};
\end{scope}
\begin{scope}
\path[clip] (  0.00,  0.00) rectangle (144.54,144.54);
\definecolor{drawColor}{gray}{0.30}

\node[text=drawColor,anchor=base,inner sep=0pt, outer sep=0pt, scale=  0.64] at ( 33.52, 15.40) {0};

\node[text=drawColor,anchor=base,inner sep=0pt, outer sep=0pt, scale=  0.64] at ( 59.00, 15.40) {5};

\node[text=drawColor,anchor=base,inner sep=0pt, outer sep=0pt, scale=  0.64] at ( 84.48, 15.40) {10};

\node[text=drawColor,anchor=base,inner sep=0pt, outer sep=0pt, scale=  0.64] at (109.96, 15.40) {15};

\node[text=drawColor,anchor=base,inner sep=0pt, outer sep=0pt, scale=  0.64] at (135.44, 15.40) {20};
\end{scope}
\begin{scope}
\path[clip] (  0.00,  0.00) rectangle (144.54,144.54);
\definecolor{drawColor}{RGB}{0,0,0}

\node[text=drawColor,anchor=base,inner sep=0pt, outer sep=0pt, scale=  0.80] at ( 84.48,  5.94) {$\alpha$};
\end{scope}
\begin{scope}
\path[clip] (  0.00,  0.00) rectangle (144.54,144.54);
\definecolor{drawColor}{RGB}{248,118,109}

\path[draw=drawColor,line width= 0.6pt,line join=round] ( 60.52,124.34) -- ( 72.08,124.34);
\end{scope}
\begin{scope}
\path[clip] (  0.00,  0.00) rectangle (144.54,144.54);
\definecolor{drawColor}{RGB}{0,191,196}

\path[draw=drawColor,line width= 0.6pt,dash pattern=on 2pt off 2pt ,line join=round] ( 60.52,109.89) -- ( 72.08,109.89);
\end{scope}
\begin{scope}
\path[clip] (  0.00,  0.00) rectangle (144.54,144.54);
\definecolor{drawColor}{RGB}{0,0,0}

\node[text=drawColor,anchor=base west,inner sep=0pt, outer sep=0pt, scale=  0.64] at ( 73.52,122.14) {LOO};
\end{scope}
\begin{scope}
\path[clip] (  0.00,  0.00) rectangle (144.54,144.54);
\definecolor{drawColor}{RGB}{0,0,0}

\node[text=drawColor,anchor=base west,inner sep=0pt, outer sep=0pt, scale=  0.64] at ( 73.52,107.68) {OR};
\end{scope}
\end{tikzpicture} \vspace{-20pt}
        \caption{$\alphaLOO_{\bp_0} \approx \alphaOR_{\bp_0} \approx 16.16$}
    \end{subfigure}
    ~ 
    \begin{subfigure}[b]{0.3\textwidth}
\begin{tikzpicture}[x=1pt,y=1pt]
\definecolor{fillColor}{RGB}{255,255,255}
\path[use as bounding box,fill=fillColor,fill opacity=0.00] (0,0) rectangle (144.54,144.54);
\begin{scope}
\path[clip] ( 31.63, 23.41) rectangle (140.54,140.54);
\definecolor{drawColor}{gray}{0.92}

\path[draw=drawColor,line width= 0.2pt,line join=round] ( 31.63, 51.98) --
	(140.54, 51.98);

\path[draw=drawColor,line width= 0.2pt,line join=round] ( 31.63, 84.32) --
	(140.54, 84.32);

\path[draw=drawColor,line width= 0.2pt,line join=round] ( 31.63,116.66) --
	(140.54,116.66);

\path[draw=drawColor,line width= 0.2pt,line join=round] ( 48.95, 23.41) --
	( 48.95,140.54);

\path[draw=drawColor,line width= 0.2pt,line join=round] ( 73.71, 23.41) --
	( 73.71,140.54);

\path[draw=drawColor,line width= 0.2pt,line join=round] ( 98.46, 23.41) --
	( 98.46,140.54);

\path[draw=drawColor,line width= 0.2pt,line join=round] (123.21, 23.41) --
	(123.21,140.54);

\path[draw=drawColor,line width= 0.4pt,line join=round] ( 31.63, 35.81) --
	(140.54, 35.81);

\path[draw=drawColor,line width= 0.4pt,line join=round] ( 31.63, 68.15) --
	(140.54, 68.15);

\path[draw=drawColor,line width= 0.4pt,line join=round] ( 31.63,100.49) --
	(140.54,100.49);

\path[draw=drawColor,line width= 0.4pt,line join=round] ( 31.63,132.83) --
	(140.54,132.83);

\path[draw=drawColor,line width= 0.4pt,line join=round] ( 36.58, 23.41) --
	( 36.58,140.54);

\path[draw=drawColor,line width= 0.4pt,line join=round] ( 61.33, 23.41) --
	( 61.33,140.54);

\path[draw=drawColor,line width= 0.4pt,line join=round] ( 86.08, 23.41) --
	( 86.08,140.54);

\path[draw=drawColor,line width= 0.4pt,line join=round] (110.84, 23.41) --
	(110.84,140.54);

\path[draw=drawColor,line width= 0.4pt,line join=round] (135.59, 23.41) --
	(135.59,140.54);
\definecolor{drawColor}{RGB}{248,118,109}

\path[draw=drawColor,line width= 0.6pt,line join=round] ( 36.58,135.22) --
	( 37.58,135.22) --
	( 38.58,101.26) --
	( 39.58,101.26) --
	( 40.58, 73.77) --
	( 41.58, 73.77) --
	( 42.58, 67.46) --
	( 43.58, 67.46) --
	( 44.58, 55.57) --
	( 45.58, 55.57) --
	( 46.58, 52.50) --
	( 47.58, 52.50) --
	( 48.58, 47.13) --
	( 49.58, 47.13) --
	( 50.58, 54.24) --
	( 51.58, 54.24) --
	( 52.58, 65.12) --
	( 53.58, 65.12) --
	( 54.58, 64.71) --
	( 55.58, 64.71) --
	( 56.58, 80.76) --
	( 57.58, 80.76) --
	( 58.58, 87.39) --
	( 59.58, 87.39) --
	( 60.58, 93.58) --
	( 61.58, 93.58) --
	( 62.58, 95.19) --
	( 63.58, 95.19) --
	( 64.58, 93.62) --
	( 65.58, 93.62) --
	( 66.58, 94.10) --
	( 67.58, 94.10) --
	( 68.58, 91.68) --
	( 69.58, 91.68) --
	( 70.58, 91.19) --
	( 71.58, 91.19) --
	( 72.58, 93.94) --
	( 73.58, 93.94) --
	( 74.58, 93.94) --
	( 75.58, 93.94) --
	( 76.58, 93.94) --
	( 77.58, 93.94) --
	( 78.58, 93.94) --
	( 79.58, 93.94) --
	( 80.58, 93.94) --
	( 81.58, 93.94) --
	( 82.58, 93.94) --
	( 83.58, 93.94) --
	( 84.58, 93.94) --
	( 85.58, 93.94) --
	( 86.58, 93.94) --
	( 87.58, 93.94) --
	( 88.58, 93.94) --
	( 89.58, 93.94) --
	( 90.58, 93.94) --
	( 91.58, 93.94) --
	( 92.58, 93.94) --
	( 93.58, 93.94) --
	( 94.58, 93.94) --
	( 95.59, 93.94) --
	( 96.59, 93.94) --
	( 97.59, 93.94) --
	( 98.59, 93.94) --
	( 99.59, 93.94) --
	(100.59, 93.94) --
	(101.59, 93.94) --
	(102.59, 93.94) --
	(103.59, 93.94) --
	(104.59, 93.94) --
	(105.59, 93.94) --
	(106.59, 93.94) --
	(107.59, 93.94) --
	(108.59, 93.94) --
	(109.59, 93.94) --
	(110.59, 93.94) --
	(111.59, 93.94) --
	(112.59, 93.94) --
	(113.59, 93.94) --
	(114.59, 93.94) --
	(115.59, 93.94) --
	(116.59, 93.94) --
	(117.59, 93.94) --
	(118.59, 93.94) --
	(119.59, 93.94) --
	(120.59, 93.94) --
	(121.59, 93.94) --
	(122.59, 93.94) --
	(123.59, 93.94) --
	(124.59, 93.94) --
	(125.59, 93.94) --
	(126.59, 93.94) --
	(127.59, 93.94) --
	(128.59, 93.94) --
	(129.59, 93.94) --
	(130.59, 93.94) --
	(131.59, 93.94) --
	(132.59, 93.94) --
	(133.59, 93.94) --
	(134.59, 93.94) --
	(135.59, 93.94);
\definecolor{drawColor}{RGB}{0,191,196}

\path[draw=drawColor,line width= 0.6pt,dash pattern=on 2pt off 2pt ,line join=round] ( 36.58,114.76) --
	( 37.58,114.76) --
	( 38.58, 83.80) --
	( 39.58, 83.80) --
	( 40.58, 59.84) --
	( 41.58, 59.84) --
	( 42.58, 42.40) --
	( 43.58, 42.40) --
	( 44.58, 33.42) --
	( 45.58, 33.42) --
	( 46.58, 28.73) --
	( 47.58, 28.73) --
	( 48.58, 34.15) --
	( 49.58, 34.15) --
	( 50.58, 42.21) --
	( 51.58, 42.21) --
	( 52.58, 49.04) --
	( 53.58, 49.04) --
	( 54.58, 58.63) --
	( 55.58, 58.63) --
	( 56.58, 66.58) --
	( 57.58, 66.58) --
	( 58.58, 74.07) --
	( 59.58, 74.07) --
	( 60.58, 81.26) --
	( 61.58, 81.26) --
	( 62.58, 85.12) --
	( 63.58, 85.12) --
	( 64.58, 87.28) --
	( 65.58, 87.28) --
	( 66.58, 87.99) --
	( 67.58, 87.99) --
	( 68.58, 89.43) --
	( 69.58, 89.43) --
	( 70.58, 89.43) --
	( 71.58, 89.43) --
	( 72.58, 89.43) --
	( 73.58, 89.43) --
	( 74.58, 89.29) --
	( 75.58, 89.29) --
	( 76.58, 89.29) --
	( 77.58, 89.29) --
	( 78.58, 89.29) --
	( 79.58, 89.29) --
	( 80.58, 89.29) --
	( 81.58, 89.29) --
	( 82.58, 89.29) --
	( 83.58, 89.29) --
	( 84.58, 89.29) --
	( 85.58, 89.29) --
	( 86.58, 89.29) --
	( 87.58, 89.29) --
	( 88.58, 89.29) --
	( 89.58, 89.29) --
	( 90.58, 89.29) --
	( 91.58, 89.29) --
	( 92.58, 89.29) --
	( 93.58, 89.29) --
	( 94.58, 89.29) --
	( 95.59, 89.29) --
	( 96.59, 89.29) --
	( 97.59, 89.29) --
	( 98.59, 89.29) --
	( 99.59, 89.29) --
	(100.59, 89.29) --
	(101.59, 89.29) --
	(102.59, 89.29) --
	(103.59, 89.29) --
	(104.59, 89.29) --
	(105.59, 89.29) --
	(106.59, 89.29) --
	(107.59, 89.29) --
	(108.59, 89.29) --
	(109.59, 89.29) --
	(110.59, 89.29) --
	(111.59, 89.29) --
	(112.59, 89.29) --
	(113.59, 89.29) --
	(114.59, 89.29) --
	(115.59, 89.29) --
	(116.59, 89.29) --
	(117.59, 89.29) --
	(118.59, 89.29) --
	(119.59, 89.29) --
	(120.59, 89.29) --
	(121.59, 89.29) --
	(122.59, 89.29) --
	(123.59, 89.29) --
	(124.59, 89.29) --
	(125.59, 89.29) --
	(126.59, 89.29) --
	(127.59, 89.29) --
	(128.59, 89.29) --
	(129.59, 89.29) --
	(130.59, 89.29) --
	(131.59, 89.29) --
	(132.59, 89.29) --
	(133.59, 89.29) --
	(134.59, 89.29) --
	(135.59, 89.29);
\end{scope}
\begin{scope}
\path[clip] (  0.00,  0.00) rectangle (144.54,144.54);
\definecolor{drawColor}{gray}{0.30}

\node[text=drawColor,anchor=base east,inner sep=0pt, outer sep=0pt, scale=  0.64] at ( 28.03, 33.60) {0.184};

\node[text=drawColor,anchor=base east,inner sep=0pt, outer sep=0pt, scale=  0.64] at ( 28.03, 65.94) {0.186};

\node[text=drawColor,anchor=base east,inner sep=0pt, outer sep=0pt, scale=  0.64] at ( 28.03, 98.28) {0.188};

\node[text=drawColor,anchor=base east,inner sep=0pt, outer sep=0pt, scale=  0.64] at ( 28.03,130.63) {0.190};
\end{scope}
\begin{scope}
\path[clip] (  0.00,  0.00) rectangle (144.54,144.54);
\definecolor{drawColor}{gray}{0.30}

\node[text=drawColor,anchor=base,inner sep=0pt, outer sep=0pt, scale=  0.64] at ( 36.58, 15.40) {0};

\node[text=drawColor,anchor=base,inner sep=0pt, outer sep=0pt, scale=  0.64] at ( 61.33, 15.40) {5};

\node[text=drawColor,anchor=base,inner sep=0pt, outer sep=0pt, scale=  0.64] at ( 86.08, 15.40) {10};

\node[text=drawColor,anchor=base,inner sep=0pt, outer sep=0pt, scale=  0.64] at (110.84, 15.40) {15};

\node[text=drawColor,anchor=base,inner sep=0pt, outer sep=0pt, scale=  0.64] at (135.59, 15.40) {20};
\end{scope}
\begin{scope}
\path[clip] (  0.00,  0.00) rectangle (144.54,144.54);
\definecolor{drawColor}{RGB}{0,0,0}

\node[text=drawColor,anchor=base,inner sep=0pt, outer sep=0pt, scale=  0.80] at ( 86.08,  5.94) {$\alpha$};
\end{scope}
\begin{scope}
\path[clip] (  0.00,  0.00) rectangle (144.54,144.54);
\definecolor{drawColor}{RGB}{248,118,109}

\path[draw=drawColor,line width= 0.6pt,line join=round] ( 62.44,124.34) -- ( 74.00,124.34);
\end{scope}
\begin{scope}
\path[clip] (  0.00,  0.00) rectangle (144.54,144.54);
\definecolor{drawColor}{RGB}{0,191,196}

\path[draw=drawColor,line width= 0.6pt,dash pattern=on 2pt off 2pt ,line join=round] ( 62.44,109.89) -- ( 74.00,109.89);
\end{scope}
\begin{scope}
\path[clip] (  0.00,  0.00) rectangle (144.54,144.54);
\definecolor{drawColor}{RGB}{0,0,0}

\node[text=drawColor,anchor=base west,inner sep=0pt, outer sep=0pt, scale=  0.64] at ( 75.44,122.14) {LOO};
\end{scope}
\begin{scope}
\path[clip] (  0.00,  0.00) rectangle (144.54,144.54);
\definecolor{drawColor}{RGB}{0,0,0}

\node[text=drawColor,anchor=base west,inner sep=0pt, outer sep=0pt, scale=  0.64] at ( 75.44,107.68) {OR};
\end{scope}
\end{tikzpicture} \vspace{-20pt}
        \caption{$\alphaLOO_{\bp_0} \approx 2.42$, $\alphaOR_{\bp_0} \approx 2.22$}
    \end{subfigure}    \caption{\textbf{LOO and Oracle Curves}.  
We consider $K = 10,000$ newsvendors where $p_{k1} \sim \text{Uniform}[.6, .9]$, $\Nhat_k \sim \text{Poisson}(10)$.
    We consider a single data draw.  The values of $p_{01}$ and the critical fractile $s$ \edit{are $(p_{01}, s) = (.3, .5)$, $(p_{01}, s) = (.75, .5)$, and $(p_{01}, s) = (.3, .2)$, respectively.\label{correctPanelLabel}}  In the first panel, instability initially increases, and there is no benefit to pooling.  In the second and third,  instability is decreasing and there is a benefit to pooling.}\label{fig:TradeOffLOO}
\end{figure}

\subsection{Implementation Details for Computational Experiments from \cref{sec:Numerics}}
\label{sec:DataCleaning}
On average, less than $2.5\%$ of stores are open on weekends, and hence we drop all weekends from our dataset.  Similarly, the data exhibits a mild upward linear trend at a rate of $215$ units a year (approximately 3.7\% increase per year), with a p-value $< .001$.  This trend is likely due to inflation and growing GDP over the time frame.  We remove this trend using simple ordinary least squares.  Finally, many stores engage in promotional activities periodically throughout the month of December leading up to Christmas.  These promotions distort sales in the surrounding period.  Hence we drop data for the month of December from our dataset.  

Throughout, $\alphaOR_{\bp_0}, \alphaLOO_{\bp_0}$ are obtained by exhaustively searching a grid of length $120$ points from $0$ to $180$.  The grand-mean  and Beta variants are obtained similarly.  Notice when $\Nhat_k = 10$, a value of $\alpha = 180$ amounts to having $18$ times more weight on the anchor point than the data, itself.  Unless otherwise specified in an experiment, $d = 20$ and $\Nhat_k = 10$ (fixed, non-random for all $k$).

{\blockedit  \label{KSPolicy}
The ``KS" policy described in the main-text corresponds to solving a data-driven distributionally robust version of the newsvendor problem, namely, 
\[
\bx_k^{\sf KS}(\rho_k, \S_k) \ \in \  \min_{x} \sup_{\P \in \mathcal P^{\sf KS}(\rho_k, \S_k)} \E_{\xi \sim \P}\left[ \max\braces{\frac{s}{1-s} (\xi-x), (x-\xi)}\right],
\] 
where the ambiguity set $\mathcal P^{\sf KS}(\rho_k, \S_k)$ is the Kolmogorov-Smirnov ball around the empirical distribution, i.e., 
\[
\mathcal P^{\sf KS}(\rho_k, \S_k) \equiv \left\{ \P \ : \ \sup_{t \in \Rl} \abs{ \P(\xi \leq t ) - \frac{1}{\Nhat_k} \sum_{j=1}^{\Nhat_k} \I{\xihat_{jk} \leq t } } \leq \rho_k  \right\}.
\]
This ambiguity set enjoys strong statistical guarantees in the large-sample setting, and, for the special case of the newsvendor problem, $\bx_k^{\sf KS}(\rho_k, \S_k)$ can be evaluated in closed-form \citep{bertsimas2018robust}.  For these reasons, we employ it in our experiments as a strong, distributionally robust benchmark.  Throughout, we select the parameters $\rho_k$ in a decoupled fashion, using $5$-fold cross-validation on $\S_k$ to select $\rho_k$ for each $k$.  

As mentioned, our ``Beta" policies use data-driven anchors selected from $\mP$, the class of all (scaled) Beta-distributions.  More specifically, this class consists of all $\text{Beta}\left( \frac{\mu}{1-\mu} \theta_2, \theta_2\right)$ distributions  with mean $\mu \in \{ 1e-6, .05, .1, \ldots, 1 \}$ and shape parameter $\theta_2 \in \{0, .05, .1, .15, \ldots, 3 \}$.  (In cases where $d < \infty$, we discretize this distribution into $d$ equal sized bins on $[0, 1]$.)  This beta-distribution should be interpreted as the distribution of the \emph{normalized} demand at the $k^\text{th}$ store.  Said differently, when shrinking the $k^\text{th}$ problem, we shrink to the un-normalized demand, i.e., towards the distribution of $\xihat_{k,\min} + (\xihat_{k, \min} - \xihat_{k,\max}) \cdot \text{Beta}\left( \frac{\mu}{1-\mu} \theta_2, \theta_2\right)$.  

}

\subsection{Summary of Historical Dataset}
{\blockedit \Cref{fig:SampleDemands} illustrates typical demand distributions at our stores as described in \cref{sec:Numerics}.  The stores display significant heterogeneity.}
\begin{figure}
    \centering
    \begin{subfigure}[b]{0.49\textwidth}
	\ifdraft{Figures Goes hEre} {\input{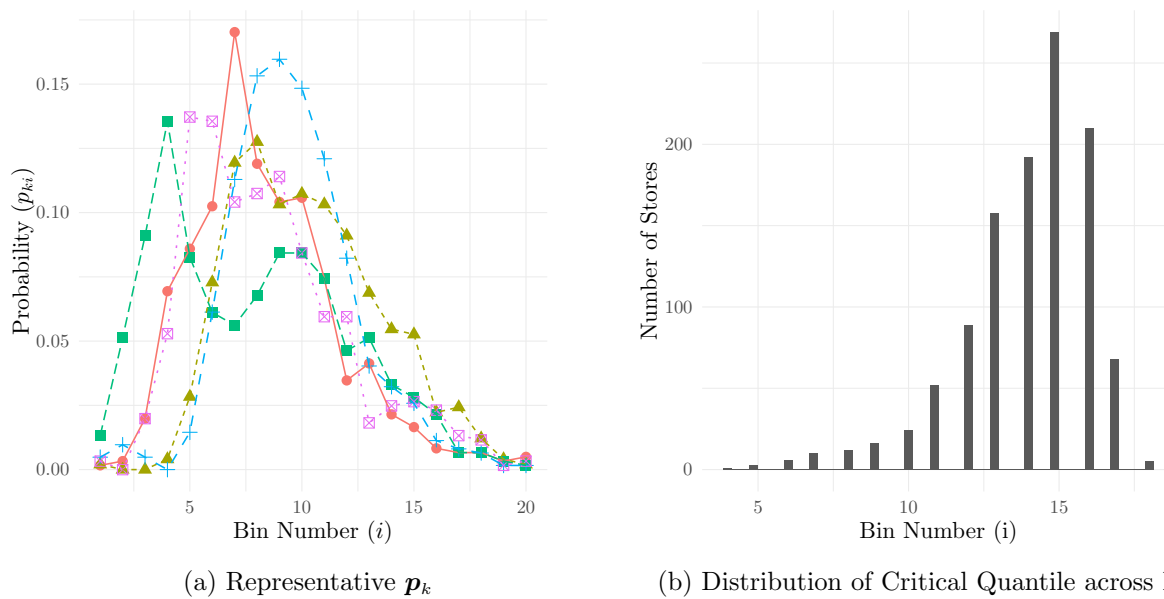}}
        \caption{Representative $\bp_k$}
    \end{subfigure}
~%
    \begin{subfigure}[b]{0.49\textwidth}

	\ifdraft{FigureGoesHere}{ 
\begin{tikzpicture}[x=1pt,y=1pt]
\definecolor{fillColor}{RGB}{255,255,255}
\path[use as bounding box,fill=fillColor,fill opacity=0.00] (0,0) rectangle (209.58,209.58);
\begin{scope}
\path[clip] ( 26.65, 23.41) rectangle (205.58,205.58);
\definecolor{drawColor}{gray}{0.92}

\path[draw=drawColor,line width= 0.2pt,line join=round] ( 26.65, 62.47) --
	(205.58, 62.47);

\path[draw=drawColor,line width= 0.2pt,line join=round] ( 26.65,124.04) --
	(205.58,124.04);

\path[draw=drawColor,line width= 0.2pt,line join=round] ( 26.65,185.60) --
	(205.58,185.60);

\path[draw=drawColor,line width= 0.2pt,line join=round] ( 76.26, 23.41) --
	( 76.26,205.58);

\path[draw=drawColor,line width= 0.2pt,line join=round] (133.20, 23.41) --
	(133.20,205.58);

\path[draw=drawColor,line width= 0.2pt,line join=round] (190.13, 23.41) --
	(190.13,205.58);

\path[draw=drawColor,line width= 0.4pt,line join=round] ( 26.65, 31.69) --
	(205.58, 31.69);

\path[draw=drawColor,line width= 0.4pt,line join=round] ( 26.65, 93.25) --
	(205.58, 93.25);

\path[draw=drawColor,line width= 0.4pt,line join=round] ( 26.65,154.82) --
	(205.58,154.82);

\path[draw=drawColor,line width= 0.4pt,line join=round] ( 47.80, 23.41) --
	( 47.80,205.58);

\path[draw=drawColor,line width= 0.4pt,line join=round] (104.73, 23.41) --
	(104.73,205.58);

\path[draw=drawColor,line width= 0.4pt,line join=round] (161.66, 23.41) --
	(161.66,205.58);
\definecolor{fillColor}{gray}{0.35}

\path[fill=fillColor] ( 34.78, 31.69) rectangle ( 38.04, 32.30);

\path[fill=fillColor] ( 38.04, 31.69) rectangle ( 41.29, 31.69);

\path[fill=fillColor] ( 41.29, 31.69) rectangle ( 44.54, 31.69);

\path[fill=fillColor] ( 44.54, 31.69) rectangle ( 47.80, 33.53);

\path[fill=fillColor] ( 47.80, 31.69) rectangle ( 51.05, 31.69);

\path[fill=fillColor] ( 51.05, 31.69) rectangle ( 54.30, 31.69);

\path[fill=fillColor] ( 54.30, 31.69) rectangle ( 57.56, 31.69);

\path[fill=fillColor] ( 57.56, 31.69) rectangle ( 60.81, 35.38);

\path[fill=fillColor] ( 60.81, 31.69) rectangle ( 64.06, 31.69);

\path[fill=fillColor] ( 64.06, 31.69) rectangle ( 67.32, 31.69);

\path[fill=fillColor] ( 67.32, 31.69) rectangle ( 70.57, 37.84);

\path[fill=fillColor] ( 70.57, 31.69) rectangle ( 73.82, 31.69);

\path[fill=fillColor] ( 73.82, 31.69) rectangle ( 77.08, 31.69);

\path[fill=fillColor] ( 77.08, 31.69) rectangle ( 80.33, 31.69);

\path[fill=fillColor] ( 80.33, 31.69) rectangle ( 83.58, 39.07);

\path[fill=fillColor] ( 83.58, 31.69) rectangle ( 86.84, 31.69);

\path[fill=fillColor] ( 86.84, 31.69) rectangle ( 90.09, 31.69);

\path[fill=fillColor] ( 90.09, 31.69) rectangle ( 93.34, 41.54);

\path[fill=fillColor] ( 93.34, 31.69) rectangle ( 96.60, 31.69);

\path[fill=fillColor] ( 96.60, 31.69) rectangle ( 99.85, 31.69);

\path[fill=fillColor] ( 99.85, 31.69) rectangle (103.10, 31.69);

\path[fill=fillColor] (103.10, 31.69) rectangle (106.36, 46.46);

\path[fill=fillColor] (106.36, 31.69) rectangle (109.61, 31.69);

\path[fill=fillColor] (109.61, 31.69) rectangle (112.86, 31.69);

\path[fill=fillColor] (112.86, 31.69) rectangle (116.12, 63.70);

\path[fill=fillColor] (116.12, 31.69) rectangle (119.37, 31.69);

\path[fill=fillColor] (119.37, 31.69) rectangle (122.62, 31.69);

\path[fill=fillColor] (122.62, 31.69) rectangle (125.88, 31.69);

\path[fill=fillColor] (125.88, 31.69) rectangle (129.13, 86.48);

\path[fill=fillColor] (129.13, 31.69) rectangle (132.38, 31.69);

\path[fill=fillColor] (132.38, 31.69) rectangle (135.64, 31.69);

\path[fill=fillColor] (135.64, 31.69) rectangle (138.89,128.96);

\path[fill=fillColor] (138.89, 31.69) rectangle (142.14, 31.69);

\path[fill=fillColor] (142.14, 31.69) rectangle (145.40, 31.69);

\path[fill=fillColor] (145.40, 31.69) rectangle (148.65, 31.69);

\path[fill=fillColor] (148.65, 31.69) rectangle (151.90,149.90);

\path[fill=fillColor] (151.90, 31.69) rectangle (155.16, 31.69);

\path[fill=fillColor] (155.16, 31.69) rectangle (158.41, 31.69);

\path[fill=fillColor] (158.41, 31.69) rectangle (161.66,197.30);

\path[fill=fillColor] (161.66, 31.69) rectangle (164.92, 31.69);

\path[fill=fillColor] (164.92, 31.69) rectangle (168.17, 31.69);

\path[fill=fillColor] (168.17, 31.69) rectangle (171.42, 31.69);

\path[fill=fillColor] (171.42, 31.69) rectangle (174.68,160.98);

\path[fill=fillColor] (174.68, 31.69) rectangle (177.93, 31.69);

\path[fill=fillColor] (177.93, 31.69) rectangle (181.18, 31.69);

\path[fill=fillColor] (181.18, 31.69) rectangle (184.44, 73.55);

\path[fill=fillColor] (184.44, 31.69) rectangle (187.69, 31.69);

\path[fill=fillColor] (187.69, 31.69) rectangle (190.94, 31.69);

\path[fill=fillColor] (190.94, 31.69) rectangle (194.20, 31.69);

\path[fill=fillColor] (194.20, 31.69) rectangle (197.45, 34.76);
\end{scope}
\begin{scope}
\path[clip] (  0.00,  0.00) rectangle (209.58,209.58);
\definecolor{drawColor}{gray}{0.30}

\node[text=drawColor,anchor=base east,inner sep=0pt, outer sep=0pt, scale=  0.64] at ( 23.05, 29.48) {0};

\node[text=drawColor,anchor=base east,inner sep=0pt, outer sep=0pt, scale=  0.64] at ( 23.05, 91.05) {100};

\node[text=drawColor,anchor=base east,inner sep=0pt, outer sep=0pt, scale=  0.64] at ( 23.05,152.62) {200};
\end{scope}
\begin{scope}
\path[clip] (  0.00,  0.00) rectangle (209.58,209.58);
\definecolor{drawColor}{gray}{0.30}

\node[text=drawColor,anchor=base,inner sep=0pt, outer sep=0pt, scale=  0.64] at ( 47.80, 15.40) {5};

\node[text=drawColor,anchor=base,inner sep=0pt, outer sep=0pt, scale=  0.64] at (104.73, 15.40) {10};

\node[text=drawColor,anchor=base,inner sep=0pt, outer sep=0pt, scale=  0.64] at (161.66, 15.40) {15};
\end{scope}
\begin{scope}
\path[clip] (  0.00,  0.00) rectangle (209.58,209.58);
\definecolor{drawColor}{RGB}{0,0,0}

\node[text=drawColor,anchor=base,inner sep=0pt, outer sep=0pt, scale=  0.80] at (116.12,  5.94) {Bin Number (i)};
\end{scope}
\begin{scope}
\path[clip] (  0.00,  0.00) rectangle (209.58,209.58);
\definecolor{drawColor}{RGB}{0,0,0}

\node[text=drawColor,rotate= 90.00,anchor=base,inner sep=0pt, outer sep=0pt, scale=  0.80] at (  9.51,114.49) {Number of Stores};
\end{scope}
\end{tikzpicture} }
        \caption{Distribution of Critical Quantile across Bins}
    \end{subfigure}
\caption{\textbf{Heterogeneity in $\bp_k$ across stores.}
The left panel shows some representative (discretized) distributions $\bp_k$ when $d = 20$ for several stores.  The right panel shows a histogram of the number of stores whose critical quantile occurs in each bin.}\label{fig:SampleDemands}
\end{figure}

The first panel of \cref{fig:AvgDailyDemands} shows the average daily demand by store for each of the $1,115$ stores in our dataset.  The second panel shows estimates of the demand distributions at a few stores.  We stress that the individual demand distributions exhibit markedly different means, variances and skewness.  

\begin{figure}
    \centering
    \begin{subfigure}[b]{0.45\textwidth}
	\ifdraft{FiguresGoesHere}{
\begin{tikzpicture}[x=1pt,y=1pt]
\definecolor{fillColor}{RGB}{255,255,255}
\path[use as bounding box,fill=fillColor,fill opacity=0.00] (0,0) rectangle (216.81,216.81);
\begin{scope}
\path[clip] ( 23.45, 23.41) rectangle (212.81,212.81);
\definecolor{drawColor}{gray}{0.92}

\path[draw=drawColor,line width= 0.2pt,line join=round] ( 23.45, 56.97) --
	(212.81, 56.97);

\path[draw=drawColor,line width= 0.2pt,line join=round] ( 23.45,106.88) --
	(212.81,106.88);

\path[draw=drawColor,line width= 0.2pt,line join=round] ( 23.45,156.79) --
	(212.81,156.79);

\path[draw=drawColor,line width= 0.2pt,line join=round] ( 23.45,206.70) --
	(212.81,206.70);

\path[draw=drawColor,line width= 0.2pt,line join=round] ( 26.67, 23.41) --
	( 26.67,212.81);

\path[draw=drawColor,line width= 0.2pt,line join=round] ( 68.67, 23.41) --
	( 68.67,212.81);

\path[draw=drawColor,line width= 0.2pt,line join=round] (110.68, 23.41) --
	(110.68,212.81);

\path[draw=drawColor,line width= 0.2pt,line join=round] (152.68, 23.41) --
	(152.68,212.81);

\path[draw=drawColor,line width= 0.2pt,line join=round] (194.69, 23.41) --
	(194.69,212.81);

\path[draw=drawColor,line width= 0.4pt,line join=round] ( 23.45, 32.01) --
	(212.81, 32.01);

\path[draw=drawColor,line width= 0.4pt,line join=round] ( 23.45, 81.92) --
	(212.81, 81.92);

\path[draw=drawColor,line width= 0.4pt,line join=round] ( 23.45,131.83) --
	(212.81,131.83);

\path[draw=drawColor,line width= 0.4pt,line join=round] ( 23.45,181.74) --
	(212.81,181.74);

\path[draw=drawColor,line width= 0.4pt,line join=round] ( 47.67, 23.41) --
	( 47.67,212.81);

\path[draw=drawColor,line width= 0.4pt,line join=round] ( 89.67, 23.41) --
	( 89.67,212.81);

\path[draw=drawColor,line width= 0.4pt,line join=round] (131.68, 23.41) --
	(131.68,212.81);

\path[draw=drawColor,line width= 0.4pt,line join=round] (173.68, 23.41) --
	(173.68,212.81);
\definecolor{fillColor}{gray}{0.35}

\path[fill=fillColor] ( 32.06, 32.01) rectangle ( 34.35, 42.00);

\path[fill=fillColor] ( 34.35, 32.01) rectangle ( 36.65, 49.48);

\path[fill=fillColor] ( 36.65, 32.01) rectangle ( 38.95, 59.46);

\path[fill=fillColor] ( 38.95, 32.01) rectangle ( 41.24, 71.94);

\path[fill=fillColor] ( 41.24, 32.01) rectangle ( 43.54, 79.43);

\path[fill=fillColor] ( 43.54, 32.01) rectangle ( 45.83, 94.40);

\path[fill=fillColor] ( 45.83, 32.01) rectangle ( 48.13,136.82);

\path[fill=fillColor] ( 48.13, 32.01) rectangle ( 50.42,171.76);

\path[fill=fillColor] ( 50.42, 32.01) rectangle ( 52.72,174.26);

\path[fill=fillColor] ( 52.72, 32.01) rectangle ( 55.01,176.75);

\path[fill=fillColor] ( 55.01, 32.01) rectangle ( 57.31,189.23);

\path[fill=fillColor] ( 57.31, 32.01) rectangle ( 59.60,199.21);

\path[fill=fillColor] ( 59.60, 32.01) rectangle ( 61.90,169.26);

\path[fill=fillColor] ( 61.90, 32.01) rectangle ( 64.19,171.76);

\path[fill=fillColor] ( 64.19, 32.01) rectangle ( 66.49,201.71);

\path[fill=fillColor] ( 66.49, 32.01) rectangle ( 68.78,174.26);

\path[fill=fillColor] ( 68.78, 32.01) rectangle ( 71.08,204.20);

\path[fill=fillColor] ( 71.08, 32.01) rectangle ( 73.37,136.82);

\path[fill=fillColor] ( 73.37, 32.01) rectangle ( 75.67,161.78);

\path[fill=fillColor] ( 75.67, 32.01) rectangle ( 77.96,124.35);

\path[fill=fillColor] ( 77.96, 32.01) rectangle ( 80.26,114.36);

\path[fill=fillColor] ( 80.26, 32.01) rectangle ( 82.55,104.38);

\path[fill=fillColor] ( 82.55, 32.01) rectangle ( 84.85, 99.39);

\path[fill=fillColor] ( 84.85, 32.01) rectangle ( 87.15, 76.93);

\path[fill=fillColor] ( 87.15, 32.01) rectangle ( 89.44, 86.91);

\path[fill=fillColor] ( 89.44, 32.01) rectangle ( 91.74, 79.43);

\path[fill=fillColor] ( 91.74, 32.01) rectangle ( 94.03, 71.94);

\path[fill=fillColor] ( 94.03, 32.01) rectangle ( 96.33, 44.49);

\path[fill=fillColor] ( 96.33, 32.01) rectangle ( 98.62, 56.97);

\path[fill=fillColor] ( 98.62, 32.01) rectangle (100.92, 56.97);

\path[fill=fillColor] (100.92, 32.01) rectangle (103.21, 66.95);

\path[fill=fillColor] (103.21, 32.01) rectangle (105.51, 44.49);

\path[fill=fillColor] (105.51, 32.01) rectangle (107.80, 51.98);

\path[fill=fillColor] (107.80, 32.01) rectangle (110.10, 42.00);

\path[fill=fillColor] (110.10, 32.01) rectangle (112.39, 34.51);

\path[fill=fillColor] (112.39, 32.01) rectangle (114.69, 34.51);

\path[fill=fillColor] (114.69, 32.01) rectangle (116.98, 34.51);

\path[fill=fillColor] (116.98, 32.01) rectangle (119.28, 37.01);

\path[fill=fillColor] (119.28, 32.01) rectangle (121.57, 46.99);

\path[fill=fillColor] (121.57, 32.01) rectangle (123.87, 34.51);

\path[fill=fillColor] (123.87, 32.01) rectangle (126.16, 34.51);

\path[fill=fillColor] (126.16, 32.01) rectangle (128.46, 34.51);

\path[fill=fillColor] (128.46, 32.01) rectangle (130.75, 39.50);

\path[fill=fillColor] (130.75, 32.01) rectangle (133.05, 37.01);

\path[fill=fillColor] (133.05, 32.01) rectangle (135.35, 39.50);

\path[fill=fillColor] (135.35, 32.01) rectangle (137.64, 34.51);

\path[fill=fillColor] (137.64, 32.01) rectangle (139.94, 32.01);

\path[fill=fillColor] (139.94, 32.01) rectangle (142.23, 37.01);

\path[fill=fillColor] (142.23, 32.01) rectangle (144.53, 32.01);

\path[fill=fillColor] (144.53, 32.01) rectangle (146.82, 32.01);

\path[fill=fillColor] (146.82, 32.01) rectangle (149.12, 32.01);

\path[fill=fillColor] (149.12, 32.01) rectangle (151.41, 32.01);

\path[fill=fillColor] (151.41, 32.01) rectangle (153.71, 34.51);

\path[fill=fillColor] (153.71, 32.01) rectangle (156.00, 32.01);

\path[fill=fillColor] (156.00, 32.01) rectangle (158.30, 32.01);

\path[fill=fillColor] (158.30, 32.01) rectangle (160.59, 34.51);

\path[fill=fillColor] (160.59, 32.01) rectangle (162.89, 37.01);

\path[fill=fillColor] (162.89, 32.01) rectangle (165.18, 34.51);

\path[fill=fillColor] (165.18, 32.01) rectangle (167.48, 34.51);

\path[fill=fillColor] (167.48, 32.01) rectangle (169.77, 32.01);

\path[fill=fillColor] (169.77, 32.01) rectangle (172.07, 32.01);

\path[fill=fillColor] (172.07, 32.01) rectangle (174.36, 34.51);

\path[fill=fillColor] (174.36, 32.01) rectangle (176.66, 32.01);

\path[fill=fillColor] (176.66, 32.01) rectangle (178.96, 34.51);

\path[fill=fillColor] (178.96, 32.01) rectangle (181.25, 32.01);

\path[fill=fillColor] (181.25, 32.01) rectangle (183.55, 32.01);

\path[fill=fillColor] (183.55, 32.01) rectangle (185.84, 32.01);

\path[fill=fillColor] (185.84, 32.01) rectangle (188.14, 32.01);

\path[fill=fillColor] (188.14, 32.01) rectangle (190.43, 32.01);

\path[fill=fillColor] (190.43, 32.01) rectangle (192.73, 32.01);

\path[fill=fillColor] (192.73, 32.01) rectangle (195.02, 32.01);

\path[fill=fillColor] (195.02, 32.01) rectangle (197.32, 32.01);

\path[fill=fillColor] (197.32, 32.01) rectangle (199.61, 32.01);

\path[fill=fillColor] (199.61, 32.01) rectangle (201.91, 32.01);

\path[fill=fillColor] (201.91, 32.01) rectangle (204.20, 34.51);
\end{scope}
\begin{scope}
\path[clip] (  0.00,  0.00) rectangle (216.81,216.81);
\definecolor{drawColor}{gray}{0.30}

\node[text=drawColor,anchor=base east,inner sep=0pt, outer sep=0pt, scale=  0.64] at ( 19.85, 29.81) {0};

\node[text=drawColor,anchor=base east,inner sep=0pt, outer sep=0pt, scale=  0.64] at ( 19.85, 79.72) {20};

\node[text=drawColor,anchor=base east,inner sep=0pt, outer sep=0pt, scale=  0.64] at ( 19.85,129.63) {40};

\node[text=drawColor,anchor=base east,inner sep=0pt, outer sep=0pt, scale=  0.64] at ( 19.85,179.54) {60};
\end{scope}
\begin{scope}
\path[clip] (  0.00,  0.00) rectangle (216.81,216.81);
\definecolor{drawColor}{gray}{0.30}

\node[text=drawColor,anchor=base,inner sep=0pt, outer sep=0pt, scale=  0.64] at ( 47.67, 15.40) {5,000};

\node[text=drawColor,anchor=base,inner sep=0pt, outer sep=0pt, scale=  0.64] at ( 89.67, 15.40) {10,000};

\node[text=drawColor,anchor=base,inner sep=0pt, outer sep=0pt, scale=  0.64] at (131.68, 15.40) {15,000};

\node[text=drawColor,anchor=base,inner sep=0pt, outer sep=0pt, scale=  0.64] at (173.68, 15.40) {20,000};
\end{scope}
\end{tikzpicture} \vspace{-20pt} }
        \caption{\textbf{Distribution of Average Daily Demand}}
    \end{subfigure}
    ~ 
    \begin{subfigure}[b]{0.45\textwidth}
	\ifdraft{FiguresGoesHere}{\input{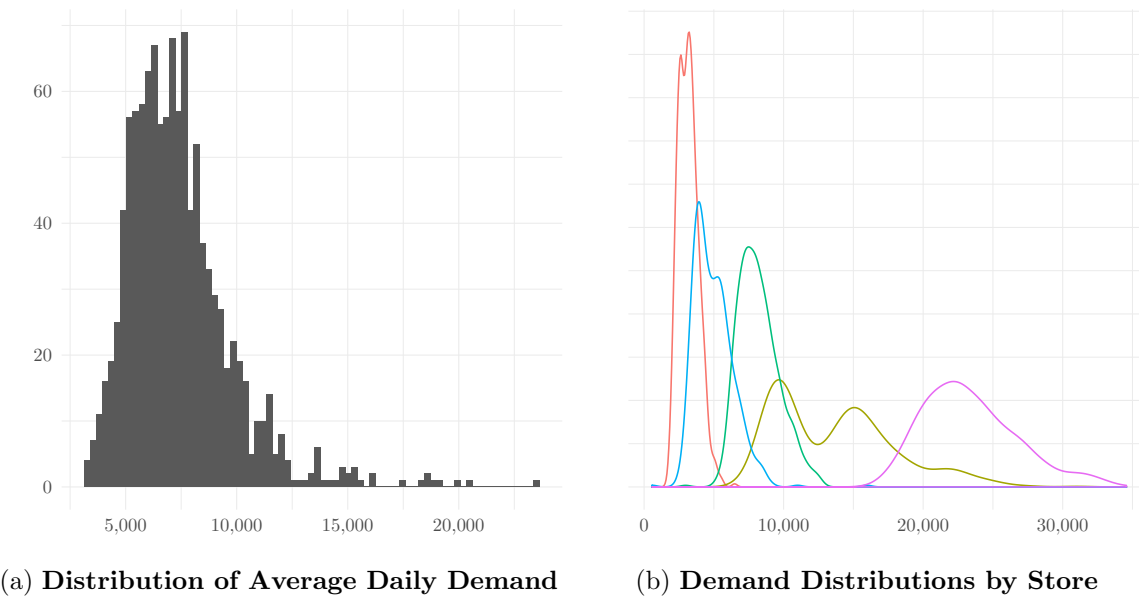}\vspace{-20pt}}
        \caption{\textbf{Demand Distributions by Store}}
    \end{subfigure}
    \caption{\textbf{Heterogeneity in Store Demand}.  
	The first panel shows a histogram of average daily demand by store  across $1,115$ stores in a European drugstore chain.
	The second panel shows estimates of the demand distribution at a few representative stores.  
	}
	\label{fig:AvgDailyDemands}
\end{figure}

\subsection{Additional Figures from \cref{sec:SyntheticData,sec:NumericallyTestingAssRandomData}.}
{\blockedit
The relative performance improvement over all SAA for all of our policies from the experiment in \cref{sec:SyntheticData} is displayed in \cref{tab:SynthDataTrueAllPolicies,tab:SynthDataFalseAllPolicies} for the case where $\Nhat_k$ is random and non-random, respectively.  To ease comparison, policies that shrink to the same type of anchor are grouped together.  Notice qualitative features are similar in both tables. 
\begin{table} \blockedit
\caption{\textbf{Relative Performance Improvement over SAA (\%), $\Nhat_k \sim \op{Poisson}(10)$}.  \\
Performance using simulated data as described in \cref{sec:SyntheticData}.
\label{tab:SynthDataTrueAllPolicies}}  
\centering
\begin{tabular}{rrrrrrrrrrr}
  \toprule
   & \multicolumn{2}{c}{Beta} & \multicolumn{3}{c}{Grand-Mean} & \multicolumn{3}{c}{Fixed (Uniform)} & \multicolumn{2}{c}{Decoupled}
\\
K & Oracle & S-SAA & Oracle & S-SAA & JS & Oracle & S-SAA & JS & SAA & KS\\ 
  \midrule
 10 & 17.20 & 10.19 & 15.28 & 12.61 & 9.83 & 12.49 & 8.43 & 5.53 &   0 & -8.71 \\ 
   32 & 11.02 & 6.42 & 9.44 & 7.07 & 4.00 & 6.09 & 3.48 & 0.56 &   0 & -12.05 \\ 
   64 & 11.34 & 8.57 & 10.17 & 8.71 & 5.20 & 7.40 & 6.65 & 1.24 &   0 & -11.57 \\ 
  128 & 13.04 & 11.75 & 12.38 & 11.68 & 5.27 & 9.38 & 9.34 & 1.37 &   0 & -11.49 \\ 
  256 & 13.10 & 12.37 & 12.66 & 12.27 & 4.94 & 9.66 & 9.66 & 0.92 &   0 & -10.71 \\ 
  362 & 13.08 & 12.57 & 12.69 & 12.43 & 5.13 & 9.71 & 9.71 & 0.36 &   0 & -10.43 \\ 
  431 & 13.26 & 12.80 & 12.91 & 12.68 & 5.13 & 9.95 & 9.95 & 0.46 &   0 & -10.25 \\ 
  512 & 12.95 & 12.48 & 12.50 & 12.29 & 5.21 & 9.67 & 9.67 & 0.27 &   0 & -10.64 \\ 
  609 & 13.12 & 12.72 & 12.69 & 12.49 & 5.32 & 9.82 & 9.82 & 0.20 &   0 & -10.57 \\ 
  724 & 13.21 & 12.85 & 12.80 & 12.63 & 5.39 & 9.97 & 9.97 & 0.17 &   0 & -10.43 \\ 
  861 & 13.35 & 13.04 & 12.95 & 12.78 & 5.40 & 10.08 & 10.08 & 0.13 &   0 & -10.46 \\ 
  1024 & 13.07 & 12.79 & 12.67 & 12.52 & 5.29 & 9.78 & 9.78 & 0.05 &   0 & -10.62 \\ 
  1115 & 13.12 & 12.86 & 12.73 & 12.58 & 5.27 & 9.82 & 9.82 & 0.05 &   0 & -10.68 \\ 
   \bottomrule
\end{tabular}

\end{table}

\begin{table} \blockedit
\caption{\textbf{Relative Performance Improvement over SAA (\%), $\Nhat_k =10$ (non-random)}.  \\
Performance using simulated data as described in \cref{sec:SyntheticData}.
\label{tab:SynthDataFalseAllPolicies}}  
\centering
\begin{tabular}{rrrrrrrrrrr}
  \toprule
   & \multicolumn{2}{c}{Beta} & \multicolumn{3}{c}{Grand-Mean} & \multicolumn{3}{c}{Fixed (Uniform)} & \multicolumn{2}{c}{Decoupled}
\\
K & Oracle & S-SAA & Oracle & S-SAA & JS & Oracle & S-SAA & JS & SAA & KS\\ 
  \midrule
 10 & 13.07 & 6.89 & 11.13 & 8.54 & 6.46 & 10.12 & 7.42 & 4.46 &   0 & -13.89 \\ 
   32 & 7.37 & 3.19 & 6.19 & 3.49 & 1.00 & 4.52 & 2.36 & 0.05 &   0 & -17.05 \\ 
   64 & 7.09 & 4.75 & 6.27 & 4.70 & 1.27 & 4.70 & 3.88 & 0.16 &   0 & -16.93 \\ 
  128 & 8.71 & 7.62 & 8.28 & 7.69 & 1.28 & 6.43 & 6.35 & 0.47 &   0 & -17.20 \\ 
  256 & 8.92 & 8.25 & 8.67 & 8.37 & 1.13 & 6.68 & 6.68 & 0.25 &   0 & -16.26 \\ 
  362 & 8.93 & 8.47 & 8.71 & 8.50 & 1.08 & 6.65 & 6.65 & 0.03 &   0 & -16.06 \\ 
  431 & 9.11 & 8.75 & 8.95 & 8.78 & 1.26 & 6.83 & 6.83 & 0.03 &   0 & -15.92 \\ 
  512 & 8.87 & 8.52 & 8.57 & 8.40 & 1.55 & 6.69 & 6.69 & 0.00 &   0 & -16.22 \\ 
  609 & 9.03 & 8.70 & 8.70 & 8.53 & 1.47 & 6.83 & 6.83 & 0.00 &   0 & -16.33 \\ 
  724 & 9.16 & 8.88 & 8.88 & 8.73 & 1.57 & 6.98 & 6.98 & 0.00 &   0 & -16.21 \\ 
  861 & 9.42 & 9.15 & 9.12 & 8.98 & 1.61 & 7.27 & 7.27 & 0.00 &   0 & -16.26 \\ 
  1024 & 9.19 & 8.96 & 8.86 & 8.74 & 1.67 & 7.02 & 7.02 & 0.00 &   0 & -16.45 \\ 
  1115 & 9.22 & 8.98 & 8.90 & 8.77 & 1.62 & 7.05 & 7.05 & 0.00 &   0 & -16.49 \\ 
   \bottomrule
\end{tabular}

\end{table}
}

\Cref{fig:StdDevPlot} shows the standard deviation of each of our methods on simulated data from \cref{sec:SyntheticData} as a function of $K$, both when \cref{ass:RandomData} holds and when it is violated and the amount of data is fixed.  Performance is again quite similar in both cases. 

\begin{figure}
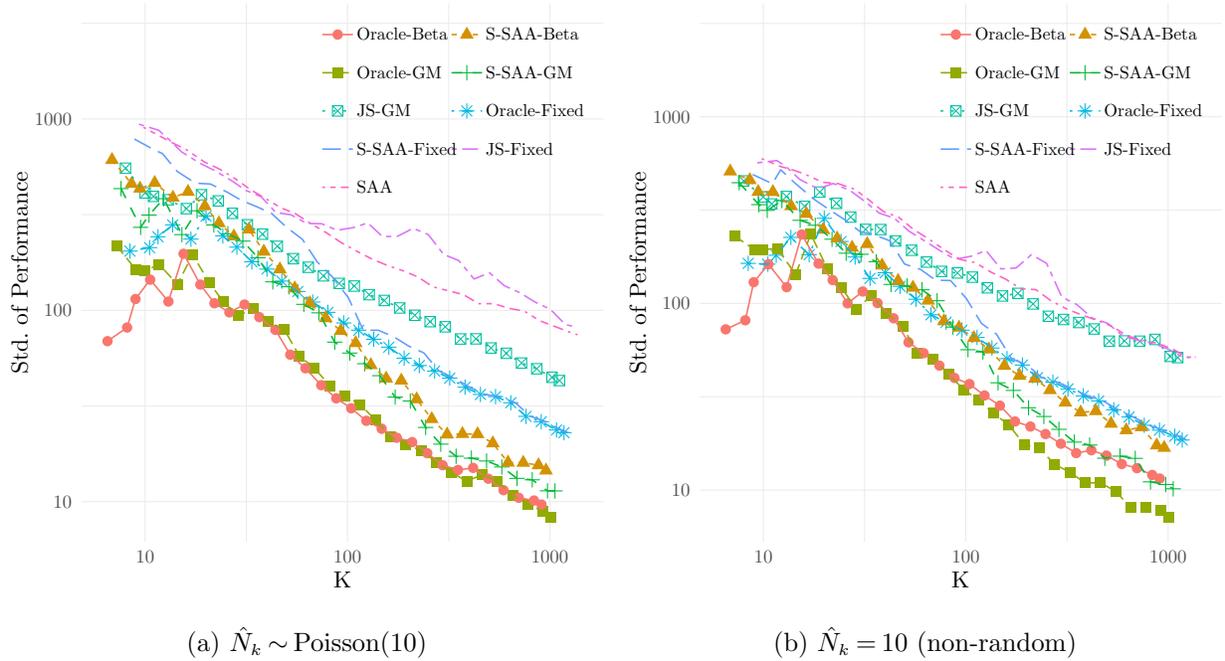

    \centering
    \begin{subfigure}[b]{0.49\textwidth}
	\ifdraft{Figure Goes Here}{ 
	\input{./Figures/synDataStdDev_True.tex} \vspace{-10pt}
	} 
        \caption{$\Nhat_k \sim \op{Poisson}(10)$}
    \end{subfigure}
    \begin{subfigure}[b]{0.49\textwidth}
	\ifdraft{Figure Goes Here}{
	\input{./Figures/synDataStdDev_False.tex} \vspace{-10pt}
	}  
        \caption{$\Nhat_k = 10$ (non-random)}
    \end{subfigure}
\caption{\textbf{Standard Deviation of Performance}
For each method, the standard deviation of converges to zero because performance concentrates at its expectation as $K \rightarrow \infty$.  Notice that our Shrunken-SAA methods are less variable than the decoupled SAA solution because pooling increases stability.} 
\label{fig:StdDevPlot}
\end{figure}

\Cref{fig:AlphaPlot} shows the average amount of pooling by method by $K$ on our simulated data set from \cref{sec:SyntheticData}, both when \cref{ass:RandomData} holds and when the amount of data is fixed.  Again, in both cases the performance is quite similar, and we see that both Shrunken-SAA and the oracle method when using $\bphat^{\rm GM}$ shrink more than the other methods.  
\begin{figure}
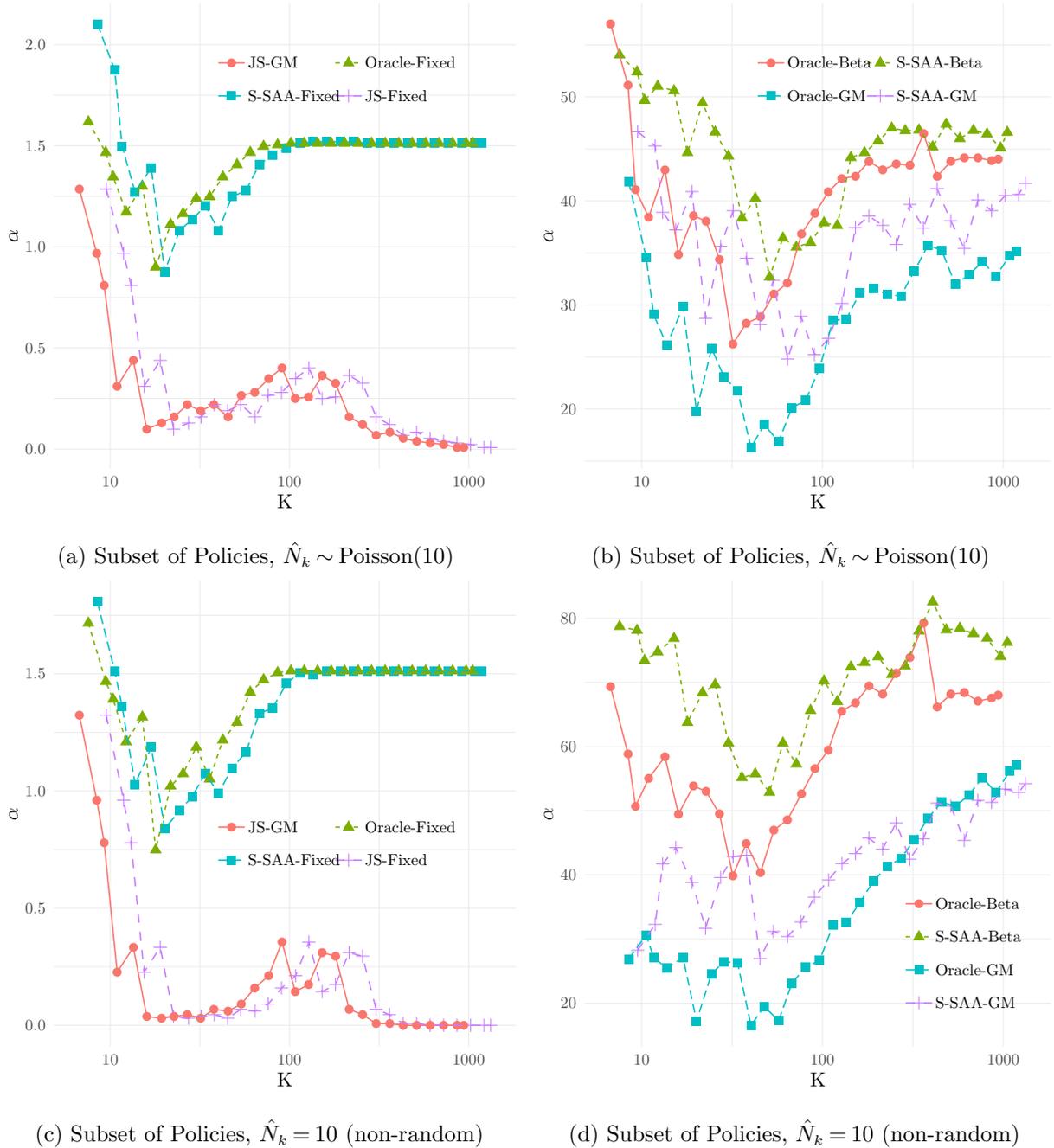

    \centering
    \begin{subfigure}[c]{.49\textwidth}
	\ifdraft{FigureGoesHere}{\input{./Figures/synDataAlphaNonGM_True.tex} \vspace{-10pt}} 
	\caption{Subset of Policies, $\Nhat_k \sim \op{Poisson}(10)$ \label{fig:Subset1Random}}
    \end{subfigure}
    \begin{subfigure}[c]{0.49\textwidth}
	\ifdraft{FiguresGoesHere}{\input{./Figures/synDataAlphaGM_True.tex} \vspace{-10pt}} 
	\caption{Subset of Policies, $\Nhat_k \sim \op{Poisson}(10)$ \label{fig:Subset2Random}}
    \end{subfigure}
    \begin{subfigure}[c]{.49\textwidth}
	\ifdraft{\vspace{10pt} FiguresGoesHere}{\input{./Figures/synDataAlphaNonGM_False.tex}\vspace{-10pt} } 
	\caption{Subset of Policies, $\Nhat_k = 10$ (non-random)}
    \end{subfigure}
    \begin{subfigure}[c]{0.49\textwidth}
	\ifdraft{\vspace{10pt} FiguresGoesHere}{\input{./Figures/synDataAlphaGM_False.tex}\vspace{-10pt} } 
	\caption{Subset of Policies, $\Nhat_k = 10$ (non-random)}
    \end{subfigure}
\caption{\textbf{Amount of Pooling by Method}
We plot the amount of data-pooling ($\alpha$) for each of the above methods (plotted separately for clarity).  In panels~a) and b), the amount of data follows \cref{ass:RandomData}.  In the remainder, it is fixed.  In general, optimization-aware methods shrink much more aggressively in both instances.}\label{fig:AlphaPlot}
\end{figure}

\subsection{Additional Figures from \cref{historical}:  Historical Backtest}
{\blockedit
\Cref{tab:HistDataTrueAllPolicies20} shows the relative performance improvement over SAA for all of our policies in the historical data experiment described in \cref{historical} with $d=20$.  For convenience, policies with the same type of anchor are grouped together for comparison.  
\begin{table}[ht] \blockedit
\caption{\textbf{Relative Performance Improvement over SAA (\%), Historical Data}.  \\
Performance using historical data as described in \cref{historical}, $d = 20$.
\label{tab:HistDataTrueAllPolicies20}}  
\centering
\begin{tabular}{rrrrrrrrrrr}
  \toprule
   & \multicolumn{2}{c}{Beta} & \multicolumn{3}{c}{Grand-Mean} & \multicolumn{3}{c}{Fixed (Uniform)} & \multicolumn{2}{c}{Decoupled}
\\
K & Oracle & S-SAA & Oracle & S-SAA & JS & Oracle & S-SAA & JS & SAA & KS\\ 
  \midrule
 10 & 18.96 & 4.72 & 13.99 & 8.98 & 8.16 & 11.82 & 5.04 & 4.13 &   0 & -12.56 \\ 
   32 & 11.34 & 4.17 & 8.65 & 3.96 & 1.63 & 5.83 & 2.32 & 0.26 &   0 & -14.62 \\ 
   64 & 10.47 & 6.25 & 8.74 & 6.22 & 2.44 & 5.99 & 4.70 & 0.17 &   0 & -14.02 \\ 
  128 & 11.88 & 9.92 & 11.10 & 9.92 & 2.55 & 8.44 & 8.44 & 0.38 &   0 & -13.25 \\ 
  256 & 11.92 & 10.89 & 11.44 & 10.98 & 2.38 & 9.06 & 9.06 & 0.59 &   0 & -12.60 \\ 
  362 & 11.49 & 10.78 & 11.16 & 10.81 & 2.08 & 8.67 & 8.67 & 0.00 &   0 & -12.44 \\ 
  431 & 11.55 & 10.89 & 11.25 & 10.95 & 2.25 & 8.72 & 8.72 & 0.00 &   0 & -12.28 \\ 
  512 & 11.12 & 10.43 & 10.73 & 10.48 & 2.49 & 8.38 & 8.38 & 0.00 &   0 & -12.50 \\ 
  609 & 11.19 & 10.57 & 10.81 & 10.58 & 2.57 & 8.48 & 8.48 & 0.00 &   0 & -12.42 \\ 
  724 & 11.25 & 10.77 & 10.94 & 10.79 & 2.65 & 8.62 & 8.62 & 0.00 &   0 & -12.31 \\ 
  861 & 11.40 & 11.01 & 11.12 & 10.96 & 2.61 & 8.75 & 8.75 & 0.00 &   0 & -12.47 \\ 
  1024 & 11.20 & 10.85 & 10.93 & 10.80 & 2.58 & 8.59 & 8.59 & 0.00 &   0 & -12.56 \\ 
  1115 & 11.30 & 10.95 & 11.05 & 10.94 & 2.55 & 8.68 & 8.68 & 0.00 &   0 & -12.61 \\ 
   \bottomrule
\end{tabular}

\end{table}
}

{\blockedit 

\subsection{Performance as $d \rightarrow \infty$}
\label{sec:Infinited}
Recall that the Shrunken-SAA algorithm, does not require that the random variables $\bxi_k$ have discrete support \edit{(cf. \cref{rem:Computational})}.
Consequently, we next study the robustness of Shrunken-SAA to $d$, the number of support points of $\bxi_k$.  

To this end, we increase $d$ from our base case.  \Cref{fig:varyingDBacktest2} below shows results for $d = 50$ and $ d = \infty$, i.e., \edit{not performing any discretization.} \label{ReviewerNoDiscretization}  
The complete set of policies can be seen in \cref{tab:HistDataTrueAllPolicies50,tab:HistDataTrueAllPoliciesInfinity} below.

\begin{figure}
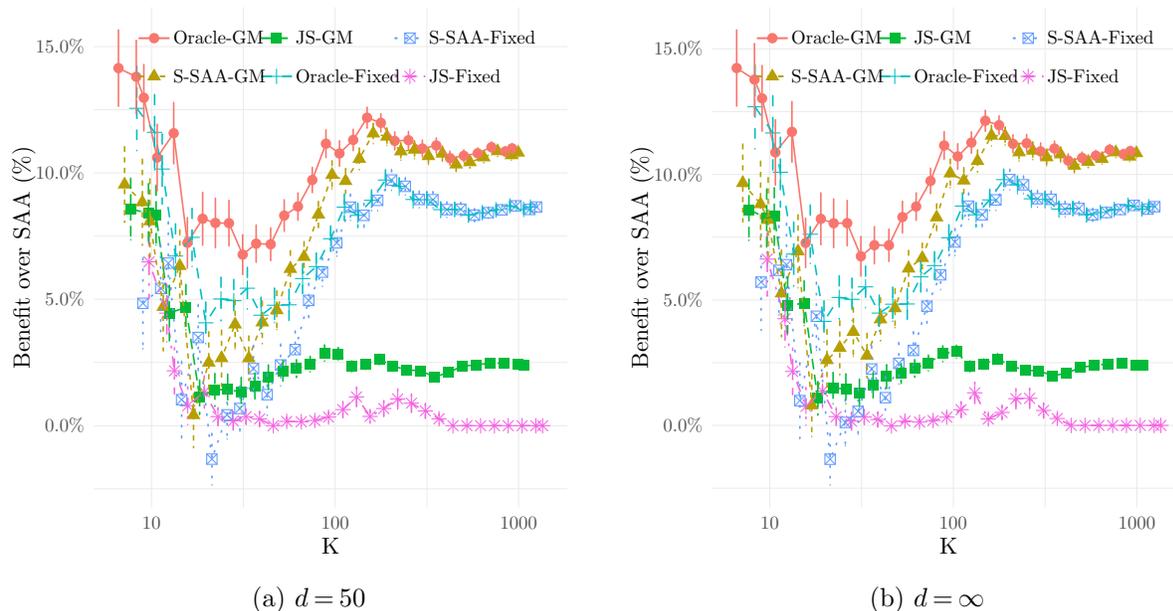

    \centering
    \begin{subfigure}{.49\textwidth}
	\ifdraft{FiguresGoesHere}{
	\input{./Figures/backtest_50.tex}
	}
        \caption{$d=50$}
    \end{subfigure}
    \begin{subfigure}{.49\textwidth}
        \ifdraft{FiguresGoesHere}{
        \input{./Figures/backtest_1000.tex}
        }
        \caption{$d=\infty$}
    \end{subfigure} 
        \caption{\textbf{Robustness to choice of $d$.}
        Performance of policies on our historical data.  In the first panel, $d=50$.  In the second panel, the distributions $\P_k$ are treated as continuous in the Shrunken-SAA algorithm, i.e., $d=\infty$.  Error bars show $\pm 1$ standard error.  The differences between the plots are essentially indiscernible.}\label{fig:varyingDBacktest2}
\end{figure} 
{\blockedit
\begin{table}[ht] \blockedit
\caption{\textbf{Relative Performance Improvement over SAA (\%)}.  \\
Performance using historical data as described in \cref{sec:Infinited}, $d = 50$.
\label{tab:HistDataTrueAllPolicies50}}  
\centering
\begin{tabular}{rrrrrrrrrrr}
  \toprule
   & \multicolumn{2}{c}{Beta} & \multicolumn{3}{c}{Grand-Mean} & \multicolumn{3}{c}{Fixed (Uniform)} & \multicolumn{2}{c}{Decoupled}
\\
K & Oracle & S-SAA & Oracle & S-SAA & JS & Oracle & S-SAA & JS & SAA & KS\\ 
  \midrule
 10 & 18.39 & 6.38 & 13.81 & 8.84 & 8.41 & 11.61 & 5.43 & 4.45 &   0 & -17.50 \\ 
   32 & 10.73 & 3.44 & 8.01 & 4.01 & 1.33 & 5.44 & 2.26 & 0.27 &   0 & -18.56 \\ 
   64 & 10.01 & 6.51 & 8.31 & 6.21 & 2.26 & 5.82 & 4.96 & 0.21 &   0 & -18.02 \\ 
  128 & 11.60 & 9.96 & 10.77 & 9.68 & 2.35 & 8.32 & 8.32 & 0.37 &   0 & -16.03 \\ 
  256 & 11.73 & 10.86 & 11.27 & 10.86 & 2.21 & 8.95 & 8.95 & 0.58 &   0 & -15.86 \\ 
  362 & 11.36 & 10.80 & 10.97 & 10.67 & 1.93 & 8.55 & 8.55 & 0.00 &   0 & -16.01 \\ 
  431 & 11.44 & 10.92 & 11.09 & 10.77 & 2.11 & 8.56 & 8.56 & 0.00 &   0 & -15.69 \\ 
  512 & 11.02 & 10.57 & 10.59 & 10.34 & 2.35 & 8.33 & 8.33 & 0.00 &   0 & -16.14 \\ 
  609 & 11.08 & 10.70 & 10.68 & 10.43 & 2.40 & 8.43 & 8.43 & 0.00 &   0 & -16.03 \\ 
  724 & 11.15 & 10.82 & 10.78 & 10.62 & 2.47 & 8.54 & 8.54 & 0.00 &   0 & -15.76 \\ 
  861 & 11.36 & 11.09 & 11.02 & 10.86 & 2.47 & 8.70 & 8.70 & 0.00 &   0 & -15.72 \\ 
  1024 & 11.20 & 10.93 & 10.85 & 10.70 & 2.44 & 8.56 & 8.56 & 0.00 &   0 & -15.98 \\ 
  1115 & 11.32 & 11.08 & 10.97 & 10.81 & 2.40 & 8.65 & 8.65 & 0.00 &   0 & -15.85 \\ 
   \bottomrule
\end{tabular}

\end{table}
The performance is nearly identical to the case of $d = 20$.  To make this clearer, in the second panel of \cref{fig:varyingDBacktest} we plot the performance of our Shrunken-SAA methods for varying $d$.  Again, the differences are quite small.  In our opinion, these results suggest that \edit{the performance of Shrunken-SAA is quite robust to size of the support of $\bxi_k$, and is still effective if $\bxi_k$ may be continuous.}
\begin{figure}[t!]%
\centering%
\begin{minipage}[m]{0.49\textwidth}
	\ifdraft{FiguresGoesHere}{
	\input{./Figures/comparing_d.tex}
	}
\end{minipage}%
\begin{minipage}[m]{0.5\textwidth}
    \captionof{figure}{\textbf{Robustness to $d$ on Historical Data.}
    We limit attention to the Shrunken-SAA policies and compare them on the same historical datasets for $d=20, 50, \infty$.  \edit{The performance of each variant is insensitive to $d$}.
	}
	\label{fig:varyingDBacktest}
\end{minipage}%
\vspace{-20pt}
\end{figure}

}

{\blockedit
\begin{table}[ht] \blockedit
\caption{\textbf{Relative Performance Improvement over SAA (\%)}.  \\
Performance using historical data as described in \cref{sec:Infinited}, $d = \infty$.
\label{tab:HistDataTrueAllPoliciesInfinity}}  
\centering
\begin{tabular}{rrrrrrrrrrr}
  \toprule
   & \multicolumn{2}{c}{Beta} & \multicolumn{3}{c}{Grand-Mean} & \multicolumn{3}{c}{Fixed (Uniform)} & \multicolumn{2}{c}{Decoupled}
\\
K & Oracle & S-SAA & Oracle & S-SAA & JS & Oracle & S-SAA & JS & SAA & KS\\ 
  \midrule
 10 & 18.13 & 6.39 & 13.78 & 8.82 & 8.26 & 11.65 & 6.18 & 4.25 &   0 & -23.75 \\ 
   32 & 10.60 & 3.71 & 8.06 & 3.73 & 1.29 & 5.52 & 2.25 & 0.24 &   0 & -26.99 \\ 
   64 & 9.94 & 6.41 & 8.30 & 6.25 & 2.27 & 5.92 & 4.75 & 0.19 &   0 & -25.18 \\ 
  128 & 11.52 & 9.98 & 10.72 & 9.76 & 2.35 & 8.39 & 8.39 & 0.26 &   0 & -23.38 \\ 
  256 & 11.66 & 10.96 & 11.22 & 10.88 & 2.19 & 9.03 & 9.03 & 0.60 &   0 & -21.83 \\ 
  362 & 11.31 & 10.83 & 10.92 & 10.68 & 1.98 & 8.62 & 8.62 & 0.00 &   0 & -21.49 \\ 
  431 & 11.38 & 10.98 & 11.03 & 10.80 & 2.08 & 8.64 & 8.64 & 0.00 &   0 & -20.74 \\ 
  512 & 10.99 & 10.65 & 10.55 & 10.35 & 2.31 & 8.39 & 8.39 & 0.00 &   0 & -21.51 \\ 
  609 & 11.06 & 10.75 & 10.66 & 10.50 & 2.39 & 8.48 & 8.48 & 0.00 &   0 & -21.43 \\ 
  724 & 11.13 & 10.92 & 10.75 & 10.62 & 2.45 & 8.60 & 8.60 & 0.00 &   0 & -21.09 \\ 
  861 & 11.34 & 11.17 & 10.99 & 10.88 & 2.48 & 8.77 & 8.77 & 0.00 &   0 & -20.76 \\ 
  1024 & 11.17 & 11.00 & 10.81 & 10.71 & 2.40 & 8.62 & 8.62 & 0.00 &   0 & -20.88 \\ 
  1115 & 11.29 & 11.16 & 10.93 & 10.85 & 2.38 & 8.71 & 8.71 & 0.00 &   0 & -20.83 \\ 
   \bottomrule
\end{tabular}

\end{table}

}
}

{\blockedit
\subsection{Performance as $N \rightarrow \infty$}
\label{sec:LargeSample}

We next study the performance of our methods as we increase $\Nhat_k$.  Recall in the experiment above, $\Nhat_k = 10$, with some instances having fewer training points due to missing values.  In \cref{fig:varyingNBacktest} we consider $\Nhat_k = 20$ days and $\Nhat_k = 40$ days for training (again with some instances having fewer data points), and let $d= \infty$.   (See also 
\cref{tab:varyingNBacktest20,tab:varyingNBacktest40} for all benchmarks.)
 As $\Nhat_k$ increases for all $k$, SAA, itself, converges in performance to the full-information optimum.  Consequently, there is ``less-room" to improve upon SAA, and we see that for $\Nhat_k = 40$, our methods still improve upon decoupling, but by a smaller amount.    We also note that the JS-GM variant performs relatively better than for small $\Nhat_k$.  We intuit this is because as $\Nhat_k \rightarrow \infty$, the empirical distribution $\phat_k$ converges in probability to the true distribution $\bp_k$, i.e., the variance of $\bphat_k$ around $\bp_k$ decreases.  For large enough $\Nhat_k$, this variance is a ``second order" concern, and hence accounting for discrepancy  in the mean (which is how $\alphaJS_{\bp_0}$ is chosen) captures most of the benefits.  This viewpoint accords more generally with intuition that estimate-then-optimize procedures work well in environments with high signal-to-noise ratios.  

In summary, we believe these preliminary studies support the idea that Shrunken-SAA retains many of SAA's strong large-sample properties, but still offers a marginal benefit for large $K$.  
}

\begin{figure}
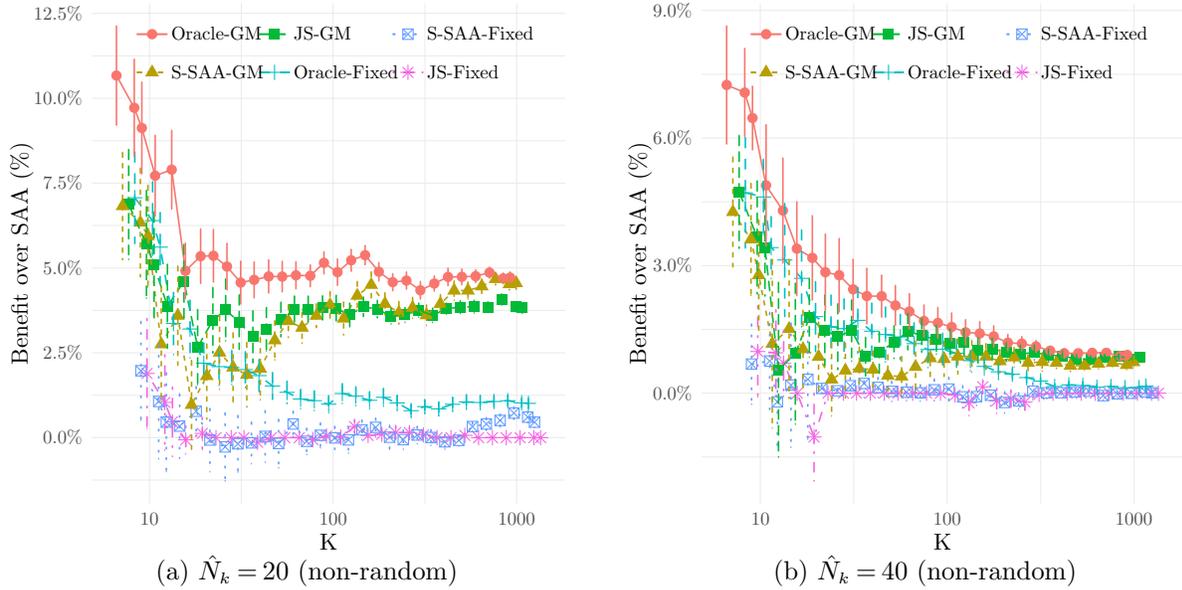

    \centering
    \begin{subfigure}{.49\textwidth}
	\ifdraft{ Figures Goes Here}{
	\input{./Figures/HistPlotN20.tex} \vspace{-10pt}
	}
    \caption{\blockedit $\Nhat_k = 20$ (non-random)}
    \end{subfigure}
    \begin{subfigure}{.49\textwidth}
	\ifdraft{Figures Goes Here}{
	\input{./Figures/HistPlotN40.tex}\vspace{-10pt}
	}
    \caption{\blockedit $\Nhat_k = 40$ (non-random)}
    \end{subfigure} 
        \caption{\textbf{Dependence on $N$.}
    Evaluated on historical data with $d  = \infty$.  Error bars show $\pm 1$ standard error.   
}\label{fig:varyingNBacktest}
\end{figure} 

{\blockedit
\begin{table}[ht] \blockedit
\caption{\textbf{Relative Performance Improvement over SAA (\%) when $N = 20$}.  \\
    Evaluated on historical data with $d  = \infty$.
\label{tab:varyingNBacktest20}}  
\centering
\begin{tabular}{rrrrrrrrrrr}
  \toprule
   & \multicolumn{2}{c}{Beta} & \multicolumn{3}{c}{Grand-Mean} & \multicolumn{3}{c}{Fixed (Uniform)} & \multicolumn{2}{c}{Decoupled}
\\
K & Oracle & S-SAA & Oracle & S-SAA & JS & Oracle & S-SAA & JS & SAA & KS\\ 
  \midrule
 10 & 15.09 & 5.58 & 10.72 & 7.57 & 7.29 & 8.26 & 2.90 & 1.03 &   0 & -7.79 \\ 
   32 & 8.74 & 2.76 & 5.85 & 2.96 & 3.81 & 2.81 & 0.29 & -0.14 &   0 & -16.33 \\ 
   64 & 6.73 & 3.62 & 5.27 & 3.87 & 3.98 & 1.74 & -0.17 & -0.06 &   0 & -16.46 \\ 
  128 & 6.32 & 4.53 & 5.21 & 3.93 & 3.89 & 1.51 & 0.18 & 0.02 &   0 & -16.94 \\ 
  256 & 5.87 & 4.32 & 4.93 & 3.92 & 3.80 & 1.23 & 0.02 & -0.13 &   0 & -15.75 \\ 
  362 & 5.68 & 4.56 & 4.75 & 4.04 & 3.89 & 1.28 & 0.05 & 0.00 &   0 & -15.72 \\ 
  431 & 5.65 & 4.59 & 4.81 & 4.23 & 3.96 & 1.32 & 0.23 & 0.07 &   0 & -15.69 \\ 
  512 & 5.72 & 4.79 & 4.98 & 4.55 & 4.00 & 1.35 & 0.50 & 0.00 &   0 & -16.09 \\ 
  609 & 5.58 & 4.88 & 4.89 & 4.56 & 3.94 & 1.28 & 0.58 & 0.00 &   0 & -16.33 \\ 
  724 & 5.54 & 4.85 & 4.87 & 4.56 & 3.97 & 1.29 & 0.62 & 0.00 &   0 & -16.28 \\ 
  861 & 5.54 & 5.09 & 5.04 & 4.82 & 4.09 & 1.32 & 0.74 & 0.00 &   0 & -16.34 \\ 
  1024 & 5.50 & 4.97 & 4.93 & 4.75 & 3.95 & 1.26 & 0.45 & 0.00 &   0 & -16.58 \\ 
  1115 & 5.44 & 4.95 & 4.90 & 4.75 & 3.90 & 1.19 & 0.30 & 0.00 &   0 & -16.50 \\ 
   \bottomrule
\end{tabular}

\end{table}

}

{\blockedit
\begin{table}[ht] \blockedit
\caption{\textbf{Relative Performance Improvement over SAA (\%) when $N = 40$}.  \\
    Evaluated on historical data with $d  = \infty$.
\label{tab:varyingNBacktest40}}  
\centering
\begin{tabular}{rrrrrrrrrrr}
  \toprule
   & \multicolumn{2}{c}{Beta} & \multicolumn{3}{c}{Grand-Mean} & \multicolumn{3}{c}{Fixed (Uniform)} & \multicolumn{2}{c}{Decoupled}
\\
K & Oracle & S-SAA & Oracle & S-SAA & JS & Oracle & S-SAA & JS & SAA & KS\\ 
  \midrule
 10 & 10.09 & 1.21 & 7.07 & 3.62 & 3.66 & 4.61 & 0.76 & 0.95 &   0 & -3.58 \\ 
   32 & 5.44 & 0.63 & 2.77 & 0.53 & 1.47 & 1.71 & 0.23 & 0.00 &   0 & -9.09 \\ 
   64 & 3.71 & 0.80 & 2.06 & 0.39 & 1.45 & 1.17 & 0.01 & 0.00 &   0 & -6.20 \\ 
  128 & 2.55 & 1.44 & 1.56 & 0.86 & 1.20 & 0.77 & -0.08 & 0.15 &   0 & -5.58 \\ 
  256 & 2.08 & 1.51 & 1.19 & 0.84 & 0.92 & 0.37 & 0.05 & -0.03 &   0 & -5.40 \\ 
  362 & 1.92 & 1.16 & 1.10 & 0.74 & 0.89 & 0.16 & 0.01 & 0.00 &   0 & -5.35 \\ 
  431 & 1.85 & 1.32 & 1.00 & 0.72 & 0.88 & 0.19 & 0.03 & 0.00 &   0 & -5.23 \\ 
  512 & 1.67 & 1.22 & 0.94 & 0.65 & 0.80 & 0.18 & 0.05 & 0.00 &   0 & -5.06 \\ 
  609 & 1.49 & 1.11 & 0.94 & 0.65 & 0.80 & 0.14 & -0.06 & 0.00 &   0 & -5.20 \\ 
  724 & 1.57 & 1.22 & 0.95 & 0.69 & 0.84 & 0.16 & -0.02 & 0.00 &   0 & -5.04 \\ 
  861 & 1.49 & 1.16 & 0.95 & 0.71 & 0.86 & 0.13 & -0.00 & 0.00 &   0 & -4.99 \\ 
  1024 & 1.41 & 1.13 & 0.87 & 0.67 & 0.81 & 0.13 & 0.04 & 0.00 &   0 & -5.05 \\ 
  1115 & 1.44 & 1.19 & 0.90 & 0.73 & 0.84 & 0.17 & 0.01 & 0.00 &   0 & -4.99 \\ 
   \bottomrule
\end{tabular}

\end{table}

}

{\blockedit
\subsection{Other Forms of Cross-Validation}
\label{sec:KFoldCrossVal}
Our theoretical development of Shrunken-SAA naturally motivated our Modified-LOO procedure in \cref{alg:ssaa}.  
When $K \Nhat_{\rm avg}$ is very large, however, LOO may be computationally demanding, and simpler $5$-fold or $10$-fold cross-validation methods might be preferred.  We next study the performance of \cref{alg:ssaa} when we replace the Modified-LOO Cross-Validation step by a simpler Modified $\kappa$-fold cross-validation step, where $\kappa \in \{2, 5, 10\}$.  Here the qualifier ``Modified" indicates that, as in \cref{alg:ssaa}, we do not update the anchor (even if it depends on the data) for each fold.  

\Cref{fig:CVTest} shows some indicative results under the synthetic data setting of \cref{sec:SyntheticData} in the case $d = \infty$ (continuous data).  We consider both $\Nhat_k \sim \op{Poisson}(10)$ (left panels) or $\Nhat_k = 10$ (right panels).  In both settings, each form of cross-validation converges to oracle performance qualitatively similarly to the LOO performance.  We have repeated this test for other values of $N$ and $d$ and with our historical data setting of \cref{historical}, and largely observe similar results.  In summary, this suggests empirically that when computational budgets require it, Shrunken-SAA can safely be implemented with other forms of cross-validation.
}
\begin{figure}
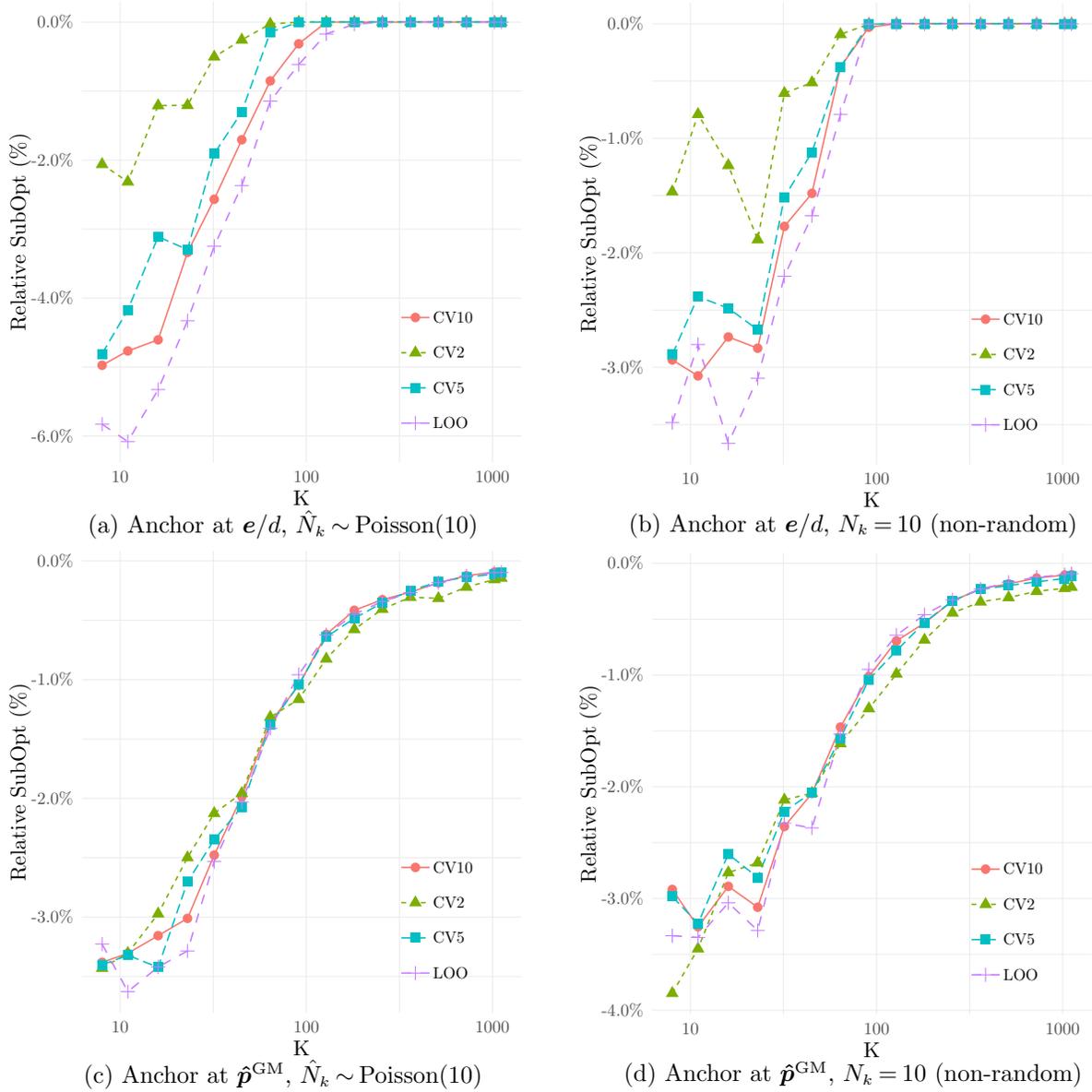

    \centering
    \begin{subfigure}[c]{.49\textwidth}
	\ifdraft{FigureGoesHere}{
	\input{./Figures/compareCV_Fixed_PoissonTrue.tex} \vspace{-10pt}
	} 
	\caption{Anchor at $\be/d$, $\Nhat_k \sim \op{Poisson}(10)$}
    \end{subfigure}
    \begin{subfigure}[c]{.49\textwidth}
	\ifdraft{\vspace{10pt} FiguresGoesHere}{
	\input{./Figures/compareCV_Fixed_PoissonFalse.tex}\vspace{-10pt} 
	} 
	\caption{Anchor at $\be/d$, $N_k = 10$ (non-random)}
    \end{subfigure}
    \begin{subfigure}[c]{0.49\textwidth}
	\ifdraft{FiguresGoesHere}{
	\input{./Figures/compareCV_GM_PoissonTrue.tex} \vspace{-10pt}
	} 
	\caption{Anchor at $\bphat^{\rm GM}$, $\Nhat_k \sim \op{Poisson}(10)$}
    \end{subfigure}
    \begin{subfigure}[c]{0.49\textwidth}
	\ifdraft{\vspace{10pt} FiguresGoesHere}{
	\input{./Figures/compareCV_GM_PoissonFalse.tex}\vspace{-10pt} 
	} 
	\caption{Anchor at $\bphat^{\rm GM}$, $N_k = 10$ (non-random)}
    \end{subfigure}
\caption{\textbf{Other Types of Cross-Validation}
We compare variants of the Shrunken-SAA procedure that leverage (modified) $2$, $5$ or $10$-fold cross-validation instead of LOO cross-validation in \cref{alg:ssaa}, both when the amount of data is random (\cref{ass:RandomData} hold) and when it is fixed.  In each case, all methods of cross-validation seem to converge to oracle optimality.  Plots show relative suboptimality to oracle performance.
}\label{fig:CVTest}
\end{figure}
%

{\blockedit

\section{Extension to Continuous Distributions}
\label{sec:ContinuousExtension}

In this section, we extend our results from \cref{sec:SmoothCostsFixedAnchor,sec:SmoothCostsGeneral} to the case where the random variables $\bxi_k$ may have continuous support and discuss the challenges of similar extensions for discrete problems.  Specifically, we no longer assume $\bxi_k \in \{\ba_{k1}, \ldots, \ba_{kd}\}$, i.e., that $\bxi_k$ is supported on a finite set. Instead we allow any compact support.

\begin{assumption}[Compact Support for $\P_k$ and $h(\S)$]\label{asn:compactsupport}
There exists a compact set $\Xi \subseteq \R{\ell}$ such that, for each $k=1,\dots,K$, $\bxi_k \sim \P_k$ is an $\ell$-dimensional real random vector whose support is contained in $\Xi$, and, with probability $1$ with respect $\S$, $h(\S) \in \mP$ and has support contained in $\Xi$. 
\end{assumption}

As mentioned in \cref{remark:continuousanalysis}, our proof technique will be to consider a discretized system whose performance is arbitrarily close to the true, continuous system and invoke our results for this discretized system.  In order to construct an arbitrarily close discretized system, we will require some additional continuity on the cost functions.

\begin{assumption}[Equicontinuity]\label{asn:equicont}
For each $k=1,\dots,K$, $\{c_k(\bx, \bxi_k ) \ : \ \bx\in\X_k\}$ is equicontinuous in $\bxi$ for all $\bxi \in \Xi$. Namely, for every $\epsilon>0,\,\bxi\in\Xi$ there exists $\delta>0$ such that $\abs{c_k(\bx, \bxi)-c_k(\bx, \bxi')}\leq \epsilon$ for all $\bx\in\X,\,\fmagd{\bxi-\bxi'}\leq\delta$.
\end{assumption}

\begin{remark}
Notice that in principle, $c_k(\bx, \bxi_k )$ need only be defined for $\bxi_k$ in the support of $\P_k$. Assuming that it is defined and equicontinuous on the larger $\Xi$ is without loss of generality via the Tietze continuous extension theorem \citep[Theorem 3.2]{munkres1974topology}. \halmos
\end{remark}

Finally, we assume the same assumptions on the cost functions as in \cref{sec:SmoothCostsFixedAnchor,sec:SmoothCostsGeneral}. We restate these below in terms of $\bxi$ that may not be finitely supported.

\begin{assumption}[\textbf{Bounded, Lipschitz, Strongly-Convex Optimization}]\label{asn:smooth_cts}
There exists $L,\gamma$ such that 
$c_{k}(\bx,\bxi)$ are $\gamma$-strongly convex and $L$-Lipschitz over $\X_k$, and, moreover, $\X_k$ is non-empty and convex, for all $k=1,\dots,K$, and $\bxi\in\Xi$.
\end{assumption}

For clarity, we repeat the definitions of some of our primitives, but now in terms of general distributions and data sets $\S_k$ and $\S$:
\begin{align*}
\bx_k(\alpha,\pQ,\S_k)&\in
\arg\min_{ \bx_k \in \X_k} \  
    \sum_{j=1}^{\Nhat_k} c_k(\bx_k, \bxihat_{kj} ) + \alpha \E_{\bxi_k \sim \pQ}\left[ c_k(\bx_k, \bxi_k) \right],
    \\
\overline Z_K(\alpha,\pQ)&=\frac1K\sum_{k=1}^KZ_k(\alpha,\pQ),
\ \text{ where } \ 
Z_k(\alpha,\pQ) =\frac{\lambda_k}{\lambdabar}\E_{\bxi_k\sim\P_k}[c_k(\bx_k(\alpha,\pQ,\S_k),\bxi_k)],
\\
\overline Z^{\LOO}_K(\alpha,\pQ)&=\frac1K\sum_{k=1}^KZ^\LOO_k(\alpha,\pQ),
\ \text{ where } \ 
Z^{\LOO}_k(\alpha,\pQ) =\frac{1}{N\lambdabar} \sum_{j=1}^{\Nhat_k}  c_{k}( \bx_k(\alpha, \pQ,  \S_k\backslash\{\bxihat_{kj}\} ), \bxihat_{kj} ).
\end{align*}
Notice $\bx_k(\alpha,\pQ,\S_k)$ is precisely as in \cref{alg:ssaa}. 

The oracle pooling amount for a specified $h(\cdot)$ is given by
\begin{align*}
\alphaOR_h&\in\argmin_{\alpha\geq0}\Zperf(\alpha,h(\S)),
\end{align*}
and the simultaneous oracle pooling amount and oracle anchor within a class $\mP$ is given by
\begin{align*}
(\alphaOR_\mP,\,\pQ^\OR_\mP) &\in \argmin_{\alpha \geq 0,\,\pQ\in\mP} \Zperf(\alpha,\pQ).
\end{align*}
Again, we will measure performance of a policy relative to these oracle benchmarks:  
\begin{align*}
{\sf Sub Opt}_{h,K}(\alpha)
	&= \Zperf(\alpha,h(\S)) -  \Zperf(\alphaOR_h,h(\S)),
\\
{\sf Sub Opt}_{\mP,K}(\alpha,\pQ)&= \Zperf(\alpha,\pQ) -  \Zperf(\alphaOR_\mP,\pQ^\OR_\mP).
\end{align*}
For convenience, we again often refer to the constant function $\S\mapsto\pQ$ as just $\pQ$.
Notice that in the special case that $\bxi_k$ has finite, discrete support, each of these above definitions is equivalent to our original definitions in \cref{sec:Model}.

We can now prove an extension of \cref{thm:FixedPointShrinkage} to the case of continuous random variables.
\begin{theorem}\label{thm:ctsdist_fixed_anchor}
{\rm \textbf{(Shrunken-SAA with Fixed Anchors for Strongly-Convex Problems and Continuous Distributions)}}
Fix any $\P_0$.  Suppose \cref{ass:RandomData,asn:compactsupport,asn:equicont,asn:smooth_cts} hold, $K \geq 2$ and $N\lambdamin \geq 1$.
Then, there exists a universal constant $\const$ such that for any $0 < \delta < 1/2$, with probability at least $1-\delta$, we have that
 \[
    {\sf Sub Opt}_{\P_0,K}(\alphaLOO_{\P_0})
     \ \leq  \ 
    \const \cdot
    \max\left(\Cmax, L \sqrt{\frac{\Cmax}{\gamma}} \right) \cdot 
    \left( \frac{\lambdamax}{\lambdamin}\right)^{5/4} \cdot
    \frac{ \log^{2} (1/\delta)\cdot  \log^{3/2}(K)}{ \sqrt K }.
\]
\end{theorem}

The first step in the proof of \cref{thm:ctsdist_fixed_anchor} is to construct our approximate discrete system:
\begin{lemma}[A Discrete Approximate System]   \label{lem:Discretization}
Suppose \cref{asn:compactsupport,asn:equicont,asn:smooth_cts} hold.  Then, for any $\epsilon > 0$, there exists a finite partition $B^\prime_1, \ldots, B^\prime_d$ of $\Xi$ and, for each $k = 1, \dots, K$, random variables $\bxi_k^\disc$ supported on $\{\ba_{k1},\dots,\ba_{kd}\}$ such that 
\begin{enumerate} [label=\roman*)]
\item $\ba_{ki} \in B^\prime_i$,
\item $\P( \bxi_{k}^\disc = \ba_{ki} ) = \P(\bxi_k \in B^\prime_i)$, and
\item $ \abs{ c_k(\bx_k, \bxi_k) - c_k(\bx, \ba_{ki}) } \leq \epsilon$ for all $k = 1, \ldots, K$, $i = 1, \ldots, d$,
$\bxi_k \in B^\prime_i$
 and $\bx \in \X_k$.  
\end{enumerate}
\end{lemma}
\proof{Proof.}
    Since $K$ is finite,  \cref{asn:equicont} implies that the larger set $\{c_k(\bx;\bxi):\bx\in\X_k,\,k=1,\dots,K\}$ is equicontinuous in $\bxi$ for all $\bxi \in\Xi$.  
    In other words,  for every $\bxi\in\Xi$, there exists $\delta(\bxi)>0$ such that $\abs{c_k(\bx;\bxi)-c_k(\bx;\bxi')}\leq\epsilon$ for all $k=1,\dots,K$ and $\bx\in\X_k$ whenever $\fmagd{\bxi-\bxi'}\leq\delta(\bxi)$  and $\bxi'\in\Xi$. Let $B(\bxi)=\{\bxi'\in\Xi:\fmagd{\bxi-\bxi'}\leq\delta(\bxi)\}$.  
    Then $\bigcup_{\bxi \in \Xi} B(\bxi)$ necessarily covers $\Xi$.  Since $\Xi$ is compact, there exists a finite subcover, $B(\bxi_{1}),\dots,B(\bxi_{d})$.  
    We construct a partition from this finite subcover, namely, 
   \[
    B^\prime_i = B(\bxi_i) \setminus \bigcup_{1 \leq j \leq i-1} B(\bxi_j) \cap B(\bxi_i).
    \]
    In words, $B^\prime_i$ is the same as $B(\bxi_i)$ but omits any point that was already covered by a previous set.  
    Let $\chi:\Xi\to\{1,\dots,d\}$ be the indicator of this partition, i.e., $\bxi \in B^\prime_{\chi(\bxi)}$ for all $\bxi \in \Xi$. 
    
    Now let
    \begin{align*}
    \ba_{ki} & \ \equiv \ \bxi_i, \quad \text{for $i=1, \ldots, d$, and $k = 1, \ldots, K$}
	\end{align*}
    and define the discrete random variable $\bxi_k^\disc$ such that $\P( \bxi_{k}^\disc = \ba_{ki} ) = \P(\bxi_k \in B^\prime_i)$.  
    
    Then the first two claims in the lemma are immediate.  For the last, notice by construction of the partition, $\| \bxi_k - \ba_{ki} \| \leq \delta(\ba_{ki})$ so that the third claim holds by equicontinuity.  
\endproof

\vskip 8pt
We will now apply our existing analysis to the discretized system.  For clarity, given $B_i',\bxihat_k^\disc$ as in \cref{lem:Discretization}, we define
    \begin{align*}
    c_{ki}(\bx) &\ \equiv \ c_k(\bx, \ba_{ki}), \quad &&\forall i =1, \ldots, d, \  k = 1, \ldots, K,
    \\
    p_{ki} &\ \equiv \ \P_k(B^\prime_i ),\quad  &&\forall i =1, \ldots, d, \  k = 1, \ldots, K,
    \\
    \mhat_{ki} &\ \equiv \ \sum_{j=1}^{\Nhat_k}\indic{\bxihat_{kj}\in B^\prime_i}, \quad &&\forall i =1, \ldots, d, \  k = 1, \ldots, K
	\\
	\bx^\disc_k(\alpha,\bq,\bfmhat_k) &\ \in \ \arg\min_{ \bx_k \in \X_k}(\bfmhat_k+\alpha \bq)^T\bc_k(\bx_k), \quad &&\forall k = 1, \ldots, K.
    \end{align*}
One can confirm directly that, under \cref{asn:smooth_cts}, $c_{ki}(\bx)$ are each $\Cmax$-bounded, $L$-Lipschitz, and $\gamma$-strongly convex for every $k,i$.
Finally, for any distribution $\pQ$ on $\R{\ell}$, define its discretization $ \dze(\pQ) = \left(\pQ(B^\prime_1),\dots,\pQ(B^\prime_d ) \right)\in\Delta_d$.  

The next step of the proof establishes that the policies $\bx^\disc_k(\cdot, \cdot, \cdot)$ of the discretized system are suitably close to the policies $\bx_k(\cdot \cdot, \cdot)$ of the original, continuous system. 
\begin{lemma}[Bounding Differences in Policies]  \label{lem:DiscretizedPolicyClose}
Suppose \cref{asn:compactsupport,asn:equicont,asn:smooth_cts} hold.   For given $\epsilon > 0$, consider the discretization given by \cref{lem:Discretization}.  Then for any $\pQ \in \mP$ and data set $\S_k$, 
\[
\magd{\bx_k(\alpha,\pQ,\S_k)-\bx^\disc_k(\alpha,\dze(\pQ),\bfmhat_k)}_2\leq \sqrt{\frac{2\epsilon}\gamma}.
\]
\end{lemma}
\proof{Proof.}
Define  
\begin{align*}
f^\disc_k(\bx_k)& \ \equiv \ 
\prns{\frac{\bfmhat_k+\alpha \dze(\pQ)}{{\Nhat_k+\alpha}}}^T\bc_k(\bx_k)
 \ =\ 
\frac1{\Nhat_k+\alpha}\sum_{j=1}^{\Nhat_k} c_k(\bx_k, \xi_{\chi(\bxihat_{kj})} ) +\frac\alpha{\Nhat_k+\alpha} \E_{\bxi_k \sim \pQ}\left[ c_k(\bx_k, \ba_{k,\chi(\bxi_k)}) \right], 
\\
f^\cts_k(\bx_k) &\ \equiv \ \frac1{\Nhat_k+\alpha}\sum_{j=1}^{\Nhat_k} c_k(\bx_k, \bxihat_{kj} ) +\frac\alpha{\Nhat_k+\alpha} \E_{\bxi_k \sim \pQ}\left[ c_k(\bx_k, \bxi_k) \right].
\end{align*}
Using \cref{lem:Discretization}  part iii) and the triangle inequality, we have that 
\(
\abs{ f^\disc_k(\bx_k) - f^\disc_k(\bx_k) } \leq \epsilon
\)
for all $\bx_k \in \X_k$, and all $k$.  

By construction $f^\disc_k$ and $f^\cts_k$ are both $\gamma$-strongly convex, and   $\bx^\disc_k(\alpha,\bq,\bfmhat_k)$ and $\bx_k(\alpha,\pQ,\S_k)$ are their respective optimizers.  Hence, we can use an argument similar to \cref{lem:ContinuityAlpha} to show that 
$\bx^\disc_k(\alpha,\bq,\bfmhat_k)$ and $\bx_k(\alpha,\pQ,\S_k)$ are close.  More specifically, by strong-convexity
\begin{align*}
f^\cts_k(\bx^\disc_k(\alpha,\dze(\pQ),\bfmhat_k))-f^\cts_k(\bx_k(\alpha,\pQ,\S_k))&\geq\frac\gamma2\magd{\bx_k(\alpha,\pQ,\S_k)-\bx^\disc_k(\alpha,\dze(\pQ),\bfmhat_k)}_2^2\\
f^\disc_k(\bx_k(\alpha,\pQ,\S_k))-f^\disc_k(\bx^\disc_k(\alpha,\dze(\pQ),\bfmhat_k))&\geq\frac\gamma2\magd{\bx_k(\alpha,\pQ,\S_k)-\bx^\disc_k(\alpha,\dze(\pQ),\bfmhat_k)}_2^2.
\end{align*}
Combining, we obtain
\begin{align*}
\gamma \magd{\bx_k(\alpha,\pQ,\S_k)-\bx^\disc_k(\alpha,\dze(\pQ),\bfmhat_k)}_2^2  & \ \leq \ 
\abs{ f^\cts_k(\bx^\disc_k(\alpha,\dze(\pQ),\bfmhat_k))- f^\disc_k(\bx^\disc_k(\alpha,\dze(\pQ),\bfmhat_k))}  
\\ & \quad + 
	\abs{f^\disc_k(\bx_k(\alpha,\pQ,\S_k))- f^\cts_k(\bx_k(\alpha,\pQ,\S_k))}
\\
& \ \leq \ 
2 \epsilon.
\end{align*}
Rearranging proves the result.  
\endproof

Finally, we introduce discrete analogues of our usual stochastic processes
\begin{align*}
\overline Z^\disc_K(\alpha,\bq)& \ = \ \frac1K\sum_{k=1}^KZ^\disc_k(\alpha,\bq) \ \text{ where } \ Z^\disc_k(\alpha,\bq)\ = \ \frac{\lambda_k}{\lambdabar}\bp_k^T\bc_k(\bx^\disc_k(\alpha,\bq)),
\\
\overline Z^{\LOO,\disc}_K(\alpha,\bq)& \ = \ \frac1K\sum_{k=1}^KZ^\disc_k(\alpha,\bq) \ \text{ where } \ Z^{\LOO,\disc}_k(\alpha,\bq) \ = \ \frac{1}{N\lambdabar} \sum_{i=1}^d  \mhat_{ki} c_{ki}( \bx^\disc_k(\alpha, \bq,  \bfmhat_k - \be_i ) ).
\end{align*}

We can now prove our first main result.  

\proof{Proof of \cref{thm:ctsdist_fixed_anchor}.}
Fix any $k$, and $\bx_k, \by_k \in \X_k$.  Then, 
\begin{align*}
\abs{ \E_{\bxi_k \sim \P_k} \left[  c_k(\bx_k, \bxi_k) \right] - \bp_k^\top \bc_k(\by_k) } 
&\ = \ 
\abs{ 
\E_{\bxi_k \sim \P_k} \left[  c_k(\bx_k, \bxi_k) \right] - \E_{\bxi_k \sim \P_k} \left[ \sum_{i=1}^d \I{\bxi_k \in B^\prime_i }c_k(\by_k, \ba_{ki}) \right]
}
\\
& \ = \ 
\abs{ 
 \E_{\bxi_k \sim \P_k} \left[  \sum_{i=1}^d c_k(\bx_k, \bxi_k) \I{\bxi_k \in B^\prime_i}\right] 
 - \E_{\bxi_k \sim \P_k} \left[ \sum_{i=1}^d \I{\bxi_k \in B^\prime_i }c_k(\by_k, \ba_{ki}) \right]
}
\\ & \ \leq \ 
\E_{\bxi_k \sim \P_k} \left[  \sum_{i=1}^d   \I{\bxi_k \in B^\prime_i} \abs{  c_k(\bx_k, \bxi_k)   -  c_k(\by_k, \ba_{ki}) } \right], 
\end{align*}
where the first equality uses the definition of $p_{ki}$, the second equality uses that $B^\prime_i$ form a partition, the last inequality uses the triangle inequality.
Now, whenever $\bxi_k \in B^\prime_i$, 
\[
\abs{  c_k(\bx_k, \bxi_k)   -  c_k(\by_k, \ba_{ki}) } 
\ \leq \ 
\abs{  c_k(\bx_k, \bxi_k)   -  c_k(\bx_k, \ba_{ki}) }
+ \abs{  c_k(\bx_k, \ba_{ki})   -  c_k(\by_k, \ba_{ki}) }
\ \leq \ \epsilon + L \| \bx_k - \by_k\|_2.
\]
Substituting above shows 
\[
\abs{ \E_{\bxi_k \sim \P_k} \left[  c_k(\bx_k, \bxi_k) \right] - \bp_k^\top \bc_k(\by_k) } 
\ \leq \ 
\epsilon + L \| \bx_k - \by_k\|_2,
\]
by construction of $B^\prime_i$ and the \cref{asn:smooth_cts}.

Now for any $\alpha$, $\pQ$, we can instantiate this inequality with $\bx_k \leftarrow \bx(\alpha, \pQ, \S_k)$ an $\by_k \leftarrow \bx^\disc_k(\alpha, \dze(\pQ), \bfmhat_k)$ to see that 
\[\hspace{-20pt}
\abs{Z_k(\alpha,\pQ)-\overline Z^\disc_k(\alpha,\dze(\pQ))} 
 \ \leq \ 
 \frac{\lambda_k}{\lambdabar} \left( \epsilon + L \|  \bx(\alpha, \pQ, \S_k) - \bx^\disc_k(\alpha, \dze(\pQ), \bfmhat_k) \|  \right) 
 \ \leq \ 
 \frac{\lambda_k}{\lambdabar} \left( \epsilon + L \sqrt{\frac{2\epsilon}{\gamma}}  \right),
\]
by \cref{lem:DiscretizedPolicyClose}.  Averaging over $k$ proves 
\begin{align*}
\abs{\overline Z_K(\alpha,\pQ)-\overline Z^\disc_K(\alpha,\dze(\pQ))}&\leq L\sqrt{\frac{2\epsilon}\gamma} + \epsilon.
\end{align*}

An entirely analogous argument yields
\begin{align*}
\abs{\overline Z^{\LOO}_K(\alpha,\pQ)-\overline Z^{\LOO,\disc}_K(\alpha,\dze(\pQ))}&\leq \left( L\sqrt{\frac{2\epsilon}\gamma} + {\epsilon} \right)  \frac{\Nhatbar}{N\lambdabar}.
\end{align*}
Notice $K\Nhatbar \sim \text{Poisson}(KN\lambdabar)$.  From \cref{lem:PropertiesOfPoisson} Part~\ref{Psi1Poisson} 
applied to $K\Nhatbar$ and Markov's inequality, we have that with probability at least $1-\delta/2$, $\frac{\Nhatbar}{N \lambdabar} \leq \log(4/\delta)$. 

Now suppose $\frac{4L^2}{C\gamma} \geq 1$.  Then, applying \cref{lem:UCDeviations} and \cref{cor:steinchen}  to $\overline Z^{\disc}_K$  and $\overline Z^{\LOO,\disc}_K$ with $\bp_0\leftarrow \dze(\P_0)$ and $\delta\leftarrow\delta/4$ shows that 
there exists a universal constant $\const_1$ such that 
with probability at least $1-\delta/2$,
\[
\sup_{\alpha\geq0}\abs{\overline Z^\disc_K(\alpha,\dze(\P_0))-\overline Z^{\LOO,\disc}_K(\alpha,\dze(\P_0))}\leq
\const_1 \cdot
    L \sqrt{\frac{\Cmax}{\gamma}} \cdot 
    \left( \frac{\lambdamax}{\lambdamin}\right)^{5/4} \cdot
    \frac{ \log^{2} (1/\delta)\cdot  \log^{3/2}(K)}{ \sqrt K },
\]
Therefore, with probability at least $1-\delta$,
\begin{align} \notag
&\sup_{\alpha\geq0}\abs{\overline Z_K(\alpha,\P_0)-\overline Z^{\LOO}_K(\alpha,\P_0)}
\\ \label{eq:MaxDiscrepancy} 
& \  \leq  \ 
\const_1 \cdot
    L \sqrt{\frac{\Cmax}{\gamma}} \cdot 
    \left( \frac{\lambdamax}{\lambdamin}\right)^{5/4} \cdot
    \frac{ \log^{2} (1/\delta)\cdot  \log^{3/2}(K)}{ \sqrt K }+\left( L\sqrt{\frac{2\epsilon}\gamma}  + \epsilon \right) \prns{1+\log(4/\delta)},
\end{align}
Similar to \cref{lem:ConditionsForOptimality}, $\overline Z^\LOO_K(\alphaLOO_{\P_0}, \P_0) \leq \overline Z^\LOO_K(\alphaOR_{\P_0}, \P_0)$ implies that 
\begin{align*}
{\sf Sub Opt}_{\P_0,K}(\alphaLOO_{\P_0}) &\ \leq\ 
\overline Z_k(\alphaLOO_{\P_0}, \P_0) - \overline Z^\LOO_k(\alphaLOO_{\P_0}, \P_0) 
\\ & \quad \ +\  
\overline Z^\LOO_k(\alphaOR{\P_0}, \P_0) -   \overline Z_k(\alphaOR_{\P_0}, \P_0) 
\\
& \ \leq \  2 \sup_{\alpha\geq0}\abs{\overline Z_K(\alpha,\P_0)-\overline Z^{\LOO}_K(\alpha,\P_0)}, 
\end{align*}
and, hence, ${\sf Sub Opt}_{\P_0,K}(\alphaLOO_{\P_0})$ is at most twice \cref{eq:MaxDiscrepancy}.  

Finally, recall the choice of $\epsilon>0$ was arbitrary.  Thus, taking a limit $\epsilon \rightarrow 0$, shows that there exists a constant $\const_2$ such that 
\[
{\sf Sub Opt}_{\P_0,K}(\alphaLOO_{\P_0}) \leq
\const_2 \cdot
    L \sqrt{\frac{\Cmax}{\gamma}} \cdot 
    \left( \frac{\lambdamax}{\lambdamin}\right)^{5/4} \cdot
    \frac{ \log^{2} (1/\delta)\cdot  \log^{3/2}(K)}{ \sqrt K }.,
\]

In the case that that $\frac{4L^2}{C\gamma}< 1$, we can always increase $L$ until $\frac{4L^2}{C\gamma} = 1$, since the larger $L$ is still a valid Lipschitz constant.  Substituting this larger $L$ above and collecting constants  proves the theorem. 
\endproof

The same key idea can also be used to prove analogues of \cref{thm:SmoothThmGeneralAnchor,thm:SmoothThmLooAnchor}.  
In the case of continuous distributions, we measure the complexity of $\mP$ by its packing number with respect to total variation distance.  Specifically, let $D_{\sf TV}(\epsilon, \mP)$ be the largest number of elements of $\mP$ that are each at least $\epsilon$ separated in total-variation distance.  

\begin{theorem}\label{thm:ctsdist_general_anchor}
{\rm \textbf{(Shrunken-SAA with Data-Driven Anchors for Strongly-Convex Problems and Continuous Distributions)}}
Fix any $h(\cdot)$.  Suppose \cref{ass:RandomData,asn:compactsupport,asn:equicont,asn:smooth_cts} hold, $K \geq 2$ and $N\lambdamin \geq 1$.
Suppose moreover that there exists $d_0$ such that for any $0<\epsilon<1/2$, $\log D_{\sf TV}(\epsilon,\mP)\leq d_0\log(1/\epsilon)$.
Then, there exists a universal constant $\const$ such that for any $0 < \delta < 1/2$, with probability at least $1-\delta$, we have that
 \[
    {\sf Sub Opt}_{h,K}(\alphaLOO_{h})
     \ \leq  \ 
    \const \cdot
    \max\left( \Cmax, \ \frac{L^2}{\gamma}  \ + \  L\sqrt{\frac{\Cmax}{\gamma}} \right) \left( \frac{\lambdamax}{\lambdamin} \right)^{5/4}
        \frac{d_0^2 \log^{7/2}(K)  \log^2(1/\delta) }{\sqrt K}.
\]
\end{theorem}
\proof{Proof.}
Fix an $\epsilon >0$.  We apply the  same discretization as in the proof of \cref{thm:ctsdist_fixed_anchor}. Let $\mP^\disc = \{ \dze(\pQ) \ : \ \pQ \in \mP \}$.  Since 
$\magd{\dze(\pQ)-\dze(\pQ')}_1\leq 2\magd{\pQ-\pQ'}_{\sf TV}$, we have that $\log D_1(\epsilon , \mP^\disc) \leq 2 d_0 \log(1/\epsilon)$.
Thus, the assumptions of \cref{thm:SmoothThmGeneralAnchor}  hold for $\mP \leftarrow \mP^\disc$ and $d_0 \leftarrow 2d_0$, 
and we can apply \cref{lem:ULLNGeneralAnchors} to bound the maximal deviations in the discrete system.  

The remainder of the proof follows the proof of \cref{thm:ctsdist_fixed_anchor} closely.  Specifically, we bound the difference between the discrete system and the original continuous system, and then bound $    {\sf Sub Opt}_{h,K}(\alphaLOO_{h})$  by twice the maximal deviations and take a limit as $\epsilon \rightarrow 0$ to yield the result. \endproof
\begin{theorem}\label{thm:ctsdist_loo_anchor}
{\rm \textbf{(Shrunken-SAA with $\hloo$ for Strongly-Convex Problems and Continuous Distributions)}}
Under the assumptions of \cref{thm:ctsdist_general_anchor},
there exists a universal constant $\const$ such that for any $0 < \delta < 1/2$, with probability at least $1-\delta$, we have that
 \[
    {\sf Sub Opt}_{\mP,K}(\alphaLOO_{\hloo},\hloo(\bfmhat))
     \ \leq  \ 
    \const \cdot
    \max\left( \Cmax, \ \frac{L^2}{\gamma}  \ + \  L\sqrt{\frac{\Cmax}{\gamma}} \right) \left( \frac{\lambdamax}{\lambdamin} \right)^{5/4}
        \frac{d_0^2 \log^{7/2}(K)  \log^2(1/\delta) }{\sqrt K}.
\]
\end{theorem}
\proof{Proof.}
The proof is the same as \cref{thm:ctsdist_general_anchor}.
\endproof

\begin{remark}[Challenges with Discrete Problems]\label{remark:discretecontinuoushard}
Proving similar extensions for continuous distributions and discrete problems poses some technical challenges .  
The key issue appears to be establishing an analogue of \cref{lem:DiscretizedPolicyClose}.  Indeed, without further assumptions, it is not clear that the set of policies $\left\{ \left( \bx_k(\alpha,\pQ,\S_k) \right)_{k=1}^K : \alpha \geq 0, \ \pQ \in \mP \right\}$ as indexed by $\alpha$ and $\pQ$ will be identical to
$\left\{ \left(   \bx^\disc_k(\alpha,\dze(\pQ),\bfmhat_k) \right)_{k=1}^K : \alpha \geq 0, \ \pQ \in \mP \right\}$, and, it also not clear in what sense these two sets might be ``approximately equal" and under what conditions.  Thus, we leave establishing the suitable additional assumptions to analyze this case to future work.  
\end{remark}
}

\end{APPENDICES}
\end{document}
